\documentclass[12pt]{article} 
\usepackage[margin=1.2in]{geometry}

\usepackage[OT1]{fontenc}

\usepackage[american]{babel}

\usepackage[utf8]{inputenc}

\usepackage[sc]{mathpazo}

\usepackage{amsmath,amssymb,amsfonts,mathrsfs}

\usepackage{graphicx}

\usepackage{pdfpages}

\usepackage{varioref}

\usepackage{xcolor}  
\definecolor{darkbrown}{RGB}{230, 77, 0}
\definecolor{darkpurple}{RGB}{179, 0, 89}
\usepackage{hyperref}  
\hypersetup{
	colorlinks=true,    
	linkcolor=darkbrown,  
	citecolor=darkpurple, 
	filecolor=black,      
	urlcolor=darkbrown    
}

\usepackage{booktabs}

\usepackage{tikzit}

\tikzstyle{red dot}=[fill={rgb,255: red,255; green,5; blue,5}, draw=black, shape=circle]
\tikzstyle{green line}=[fill={rgb,255: red,85; green,255; blue,73}, draw=black, shape=circle, tikzit fill={rgb,255: red,85; green,255; blue,73}]
\tikzstyle{textnode}=[align=center, fill=white, draw=none]
\tikzstyle{small text}=[fill=none, draw=none, shape=circle, font={\footnotesize}, align=center]
\tikzstyle{tinytext}=[fill=none, draw=none, shape=circle, font={\tiny}]

\tikzstyle{edge redish}=[-, fill=none, draw=black, tikzit draw=black]
\tikzstyle{arrow}=[->]
\tikzstyle{soft edge}=[-, fill={rgb,255: red,168; green,168; blue,168}, draw={rgb,255: red,176; green,176; blue,176}]
\tikzstyle{fill edges}=[-, draw=black, tikzit draw=black, fill={rgb,255: red,161; green,225; blue,237}, tikzit fill={rgb,255: red,161; green,225; blue,237}]

\usepackage{appendix}

\title{Semialgebraic groups and generalized affine buildings}
\author{
	Raphael Appenzeller
}
\date{\today}

\usepackage{fancyhdr}
\usepackage{tikz-cd}
\usepackage{graphicx}
\usepackage{subcaption}
\usepackage{amsthm}
\usepackage{mathtools}
\usepackage{xfrac}
\usepackage{tikz-cd}
\usepackage{stmaryrd} 
\usepackage{wrapfig} 
\usepackage{epigraph} 
\usepackage{ wasysym } 

\newcommand{\F}{\mathbb{F} }
\newcommand{\R}{\mathbb{R} }
\newcommand{\Q}{\mathbb{Q} }
\newcommand{\N}{\mathbb{N} }
\newcommand{\Z}{\mathbb{Z} }
\newcommand{\bbH}{\mathbb{H} }
\newcommand{\C}{\mathbb{C} }
\newcommand{\G}{\mathbb{G} }
\newcommand{\D}{\mathbb{D} }
\newcommand{\X}{\mathcal{X}}
\newcommand{\Qrc}{\overline{\Q}^{\operatorname{rc}}}
\newcommand{\K}{\mathbb{K}}

\newcommand{\frakg}{\mathfrak{g}}
\newcommand{\frakp}{\mathfrak{p}}
\newcommand{\frakk}{\mathfrak{k}}
\newcommand{\fraka}{\mathfrak{a}}
\newcommand{\frakap}{\overline{\mathfrak{a}}^+}
\newcommand{\Ap}{A^+}
\newcommand{\ApF}{A_\F^+}
\newcommand{\ApR}{A_\R^+}
\newcommand{\frakn}{\mathfrak{n}}
\newcommand{\frakh}{\mathfrak{h}}
\newcommand{\frakt}{\mathfrak{t}}

\newcommand{\frakl}{\mathfrak{l}}
\newcommand{\eL}{\mathcal{L}}
\newcommand{\A}{\mathbb{A}}
\newcommand{\Fun}{\mathscr{F}} 
\newcommand{\B}{\mathcal{B}}

\newcommand{\SLnR}{\operatorname{SL}_n(\R) }

\newcommand{\GLnB}{\operatorname{GL}_n(\mathbb{B}) }
\newcommand{\GLnD}{\operatorname{GL}_n(\D) }

\newcommand{\SL}{\operatorname{SL} }

\newcommand{\Id}{\operatorname{Id} }
\newcommand{\diag}[1]{\operatorname{Diag}\left( #1_1, \ldots , #1_n\right) }

\newcommand{\tran}{{^{\mathsf{T}}\!}}
\newcommand{\KPhi}{\vphantom{\Phi}_{\K}\Phi}
\newcommand{\RPhi}{\vphantom{\Phi}_{\R}\Phi}
\newcommand{\FPhi}{\vphantom{\Phi}_{\F}\Phi}
\newcommand{\KW}{\vphantom{W}_{\K}W}
\newcommand{\RW}{\vphantom{W}_{\R}W}
\newcommand{\FW}{\vphantom{W}_{\F}W}
\newcommand{\KDelta}{\vphantom{\Delta}_{\K}\Delta}

\newcommand{\FDelta}{\vphantom{\Delta}_{\F}\Delta}

\DeclareMathOperator{\leanop}{\mathsf{Lean}}
\def\lean{\ensuremath{\leanop}} 

\numberwithin{equation}{section}
\newtheorem{theorem}{Theorem}[section]
\newtheorem{lemma}[theorem]{Lemma}
\newtheorem{proposition}[theorem]{Proposition}
\newtheorem{corollary}[theorem]{Corollary}

\newtheorem{question}[theorem]{Question}

\makeatletter
\newtheorem*{rep@theorem}{\rep@title}
\newcommand{\newreptheorem}[2]{%
	\newenvironment{rep#1}[1]{%
		\def\rep@title{#2 \ref{##1}}%
		\begin{rep@theorem}}%
		{\end{rep@theorem}}}
\makeatother

\newreptheorem{theorem}{Theorem}
\newreptheorem{proposition}{Proposition}
\newreptheorem{lemma}{Lemma}
\newreptheorem{corollary}{Corollary}

\theoremstyle{remark} 
\newtheorem{remark}{Remark}[section] 

\usepackage{framed}

\begin{document}
		\begin{titlepage}
		\begin{center}
			
			\large Diss. ETH No. 30394
			\vspace*{1.5cm}
			
			\LARGE{
				\textsc{Semialgebraic groups and \\ generalized affine buildings }}\\
			\vspace{1cm}
			\vspace{1cm}
			\large
			
			A thesis submitted to attain the degree\\
			of Doctor of Sciences\\
			(Dr. sc. ETH Z\"urich)

			\vspace{1cm}
			
			presented by\\
			
			\vspace{0.5cm}
			\textsc{Raphael Appenzeller}\\
			
			\vspace{0.5cm}
			MSc Mathematics at ETH Z\"urich\\
			born on 04.02.94\\
			
			\vspace{0.5cm}
			accepted on the recommendation of\\
			\vspace{0.5cm}
			Prof. Dr. Marc Burger, examiner\\
			Prof. Dr. Petra Schwer, co-examiner\\
			
			\vspace{2cm}
			2024
		\end{center}
		
	\end{titlepage}
	  
	  \newpage
	  
	  \pagenumbering{roman}
	  \renewcommand{\thepage}{(\roman{page})}
	  
	  \section*{Abstract}
\addcontentsline{toc}{section}{Abstract}

We develop the theory of algebraic groups over real closed fields and apply the results to construct a geometric object $\B$ and to prove that $\B$ is an affine $\Lambda$-building. Real closed fields are ordered fields that have particularly nice properties in the context of \emph{semialgebraic geometry}, where objects are defined by polynomial equalities and inequalities. Over the real numbers, algebraic groups are Lie groups and have been studied extensively in the last century. Many results can be generalized to algebraic groups over other real closed fields. We use a model theoretic transfer principle to prove generalizations of statements about semisimple Lie groups. In this direction we give proofs for the Iwasawa-decomposition $KAU$, the Cartan-decomposition $KAK$ and the Bruhat-decomposition $BWB$. For unipotent subgroups we prove the Baker-Campbell-Hausdorff formula and use it to analyse root groups. We give a proof of the Jacobson-Morozov Lemma about subgroups whose Lie algebra is isomorphic to $\mathfrak{sl}_2$ and we describe other rank 1 subgroups which are the semisimple parts of Levi-subgroups. We prove a semialgebraic version of Kostant's convexity.

Over the reals, semisimple Lie groups are closely related to the symmetry groups of symmetric spaces of non-compact type. Following an idea from \cite{KrTe}, these symmetric spaces can be described semialgebraically, which allows us to consider their semialgebraic extension over any real closed field. Starting from these non-standard symmetric spaces we use a valuation (with image some non-discrete ordered abelian group $\Lambda$) on the fields to define a $\Lambda$-pseudometric. Identifying points of distance zero results in a $\Lambda$-metric space $\B$. Assuming that the root system of the associated Lie group is reduced, we prove that $\B$ is an affine $\Lambda$-building. The proof relies on a thorough analysis of stabilizers.

We further consider the special case of the group $\operatorname{SL}_2$, where $\B$ is a $\Lambda$-tree. When $\Lambda < \R$ is a proper subgroup, $\B$ is an incomplete metric space. We prove that even after completing all the apartments, the metric space may still not be complete. We also discuss the independence of axioms of $\Lambda$-trees and $\Lambda$-buildings. We formalize the notions of $\Lambda$-metric spaces and $\Lambda$-trees and produce computer verified proofs for some of our results in the proof assistant \lean .  

\newpage

\section*{Zusammenfassung}

Wir entwickeln die Theorie der algebraischen Gruppen über reell abgeschlossenen Körpern und verwenden die Resultate, um ein geometrisches Objekt $\B$ zu konstruieren, und um zu zeigen, dass $\B$ ein affines $\Lambda$-Gebäude ist. Reell abgeschlossene Körper sind geordnete Körper, die besonders gute Eigenschaften in der \emph{semialgebraischen Geometrie} haben, wo Objekte mit polynomiellen Gleichungen und Ungleichungen definiert werden. Algebraische Gruppen über den reellen Zahlen sind Lie-Gruppen und wurden im letzten Jahrhundert ausgiebig studiert. Viele Resultate können für algebraische Gruppen über reell abgeschlossenen Körpern verallgemeinert werden. Wir verwenden ein Transfertheorem aus der Modelltheorie, um Verallgemeinerungen von Aussagen über halbeinfache Lie-Gruppen zu zeigen. In dieser Richtung geben wir Beweise für die Iwasawa-Zerlegung $KAU$, die Cartan-Zerlegung $KAK$ und die Bruhat-Zerlegung $BWB$. Für unipotente Untergruppen zeigen wir die Baker-Campbell-Hausdorff-Formel und verwenden sie, um Wurzelgruppen zu untersuchen. Wir geben einen Beweis für das Jacobson-Morozov-Lemma über Untergruppen mit Lie-Algebren, die isomorph zu $\mathfrak{sl}_2$ sind, und wir beschreiben andere Untergruppen mit Rang 1, die den halbeinfachen Teilen von Levi-Untergruppen entsprechen. Wir beweisen eine semialgebraische Version von Kostant's Konvexität. 

Über den reellen Zahlen sind halbeinfache Lie-Gruppen eng verwandt mit den Symmetriegruppen von symmetrischen Räumen von nicht-kompaktem Typ. Nach einer Idee von \cite{KrTe}, können diese symmetrischen Räume semialgebraisch beschrieben werden, was ihre semialgebraische Erweiterung zu einem beliebigen reell abgeschlossenen Körper erlaubt. Mit diesen nicht-standard symmetrischen Räumen können Bewertungen (mit Bild in einer nicht-diskreten geordneten kommutativen Gruppe $\Lambda$) der Körper verwendet werden, um eine $\Lambda$-Pseudometrik zu definieren. Werden dann Punkte der Distanz Null identifiziert, entsteht ein $\Lambda$-metrischer Raum $\B$. Unter der Annahme, dass das Wurzelsystem der entsprechenden Lie-Gruppe reduziert ist, zeigen wir, dass $\B$ ein affines $\Lambda$-Gebäude ist. Der Beweis hängt von einer gründlichen Untersuchung der Stabilisatoren ab.

Weiter untersuchen wir den Spezialfall der Gruppe $\operatorname{SL}_2$, wo $\B$ ein $\Lambda$-Baum ist. Wenn $\Lambda <\R$ eine echte Untergruppe ist, dann ist $\B$ ein unvollständiger metrsicher Raum. Wir zeigen, dass $\B$ auch nachdem alle Segmente vervollständigt wurden nicht zwingend vollständig sein muss. Weiter diskutieren wir die Unabhängigkeit der Axiome von $\Lambda$-Bäumen und $\Lambda$-Gebäuden. Wir formalisieren die Konzepte von $\Lambda$-metrischen Räumen und $\Lambda$-Bäumen und produzieren Computer-überprüfte Beweise für einige der Resultate mittels Beweisassistent \lean .

	\section*{Acknowledgements} 
	\addcontentsline{toc}{section}{Acknowledgements}
	
	I am very grateful for the support of many people during the last five years. I have learned so much from my supervisor, Marc Burger, starting all the way back in my first semester at university, when he introduced me to the hyperbolic plane. I appreciate your continued support and also the freedom you gave me to work on projects that interested me. I am thankful to Petra Schwer for showing interest in my work, inviting me to Magdeburg and Heidelberg, and for many long discussions about how to prove axiom (A2). I am looking forward to working together with you in the future. For ideas and help with the mathematics in this thesis, I would like to thank the experienced mathematicians Linus Kramer, Anne Parreau, Gabriel Pallier and Auguste Hébert.
	
	I spent a lot of time with my academic siblings, Luca De Rosa, Xenia Flamm and Victor Jaeck. I learned so much from our weekly meetings and joint projects. Luca's relaxed approach to life, Xenia's energetic ideas and the humor of Victor are the perfect ingredients for a sustainably motivating team. Numerous people enriched the time at ETH, among them fellow PhD students Jannick Krifka, Merlin Incerti-Medici, Davide Spriano, Tommaso Goldhirsch, Lisa Ricci, Paula Truöl, Francesco Fournier-Facio, Lauro Silini, Martina J\o rgensen, Hjalti Isleifsson, Alessio Cela, Aitor Iribaz-Lopez, Mireille Soergel, Fernando Camacho Cadena, Merik Niemeyer, Kevin Klinge, Konstantin Andritsch and Segev Gonen Cohen and also more senior people such as Peter Feller, Tom Ilmanen, Urs Lang, Emmanuel Kowalski, Sebastian Baader, Matthew Cordes, Benjamin Brück, Johannes Schmitt and Samir Canning. It was a pleasure to set up reading courses with you, have apr\`es-lunch-meetings and arxiv-seminars, discuss mathematics, go 3D-printing and do outreach projects, but also to spend some quality free time together. Special thanks to Victor Jaeck and Segev Gonen Cohen for reading parts of this thesis and providing valuable feedback before I handed it in.
	
	I would also like to thank my family and friends, foremost Ramona Ackermann, for their support and for keeping up with my constant excitement about mathematics.

	%

	  \newpage
	
	\hypersetup{
		colorlinks=false,    
		linkcolor=black,  
		citecolor=darkpurple, 
		filecolor=black,      
		urlcolor=darkbrown    
	}
	
	\tableofcontents
	
	\hypersetup{
		colorlinks=true,    
		linkcolor=darkbrown,  
		citecolor=darkpurple, 
		filecolor=black,      
		urlcolor=darkbrown    
	}
	
\newpage

    \pagenumbering{arabic}
    \setcounter{page}{1}

	\epigraph{\emph{C’est avec la logique que nous prouvons et avec l’intuition que nous trouvons.}}{Henry Poincaré}

	\section{Introduction and main results}
	
	The rich interplay of group theory and geometry first popularized by Klein's Erlangen program enabled enormous progress in both geometry as well as group theory in the last century. A connection between Riemannian symmetric spaces and Lie groups was developed by Cartan  in 1926, leading to the classification of symmetric spaces \cite{Car27}. While Lie groups rely on the real numbers, developments in algebraic geometry allowed the study of algebraic groups over arbitrary fields. In 1962 Jacques Tits introduced axioms for buildings as a means to connect algebraic groups with incidence geometries \cite{Tit62}.
	Buildings had and continue to have a wide impact from the classification of finite simple groups to rigidity theorems by Margulis and Mostow and to the study of character varieties.
	
	Symmetric spaces are Riemannian manifolds that at every point admit an isometry given by reversing the direction of geodesics through the point \cite{Hel2, Ebe}. Examples include the model spaces of constant curvature (Euclidean space $\mathbb{E}^n$, spheres $\mathbb{S}^n$ and hyperbolic space $\mathbb{H}^n$) but also manifolds whose curvature properties may depend on the direction. Every symmetric space decomposes as a direct product into symmetric spaces of compact type, Euclidean type and non-compact type. The isometry group of a symmetric space is a Lie group and when the space is of non-compact type, the Lie group is semisimple.

  \begin{figure}[htb]  
  	\centering  
\begin{equation*}
	\begin{tikzpicture}
	\begin{pgfonlayer}{nodelayer}
		\node [style=none] (0) at (3.5, -0.25) {};
		\node [style=none] (1) at (2.5, -1.25) {};
		\node [style=none] (2) at (2.5, -13) {};
		\node [style=none] (3) at (3.5, -14) {};
		\node [style=none] (4) at (11.75, -14) {};
		\node [style=none] (5) at (12.75, -13) {};
		\node [style=none] (6) at (12.75, -1.25) {};
		\node [style=none] (7) at (11.75, -0.25) {};
		\node [style=small text] (9) at (7.25, -0.75) {simplicial buildings};
		\node [style=small text] (11) at (6, -3) {spherical \\ buildings\\};
		\node [style=none] (13) at (4.5, -1.25) {};
		\node [style=none] (14) at (3.75, -2) {};
		\node [style=none] (15) at (7.5, -1.25) {};
		\node [style=none] (16) at (8.25, -2) {};
		\node [style=none] (17) at (8.25, -13) {};
		\node [style=none] (18) at (7.5, -13.75) {};
		\node [style=none] (19) at (4.5, -13.75) {};
		\node [style=none] (20) at (3.75, -13) {};
		\node [style=none] (21) at (1.5, -14) {};
		\node [style=none] (22) at (1.5, -1) {};
		\node [style=none] (23) at (1.5, -0.75) {};
		\node [style=tinytext] (24) at (1.5, -0.5) {$\operatorname{rank} /\dim$};
		\node [style=none] (26) at (1, -13) {};
		\node [style=none] (27) at (1, -11.5) {};
		\node [style=none] (28) at (27.75, -11.5) {};
		\node [style=none] (29) at (1, -8.5) {};
		\node [style=none] (30) at (1, -10) {};
		\node [style=none] (31) at (1, -7) {};
		\node [style=none] (40) at (27.75, -8.5) {};
		\node [style=none] (41) at (27.75, -10) {};
		\node [style=none] (42) at (27.75, -7) {};
		\node [style=small text] (43) at (10.5, -3) {discrete \\ affine \\ buildings };
		\node [style=none] (44) at (9.25, -1.25) {};
		\node [style=none] (45) at (8.5, -2) {};
		\node [style=none] (46) at (27, -1.25) {};
		\node [style=none] (47) at (27.75, -2) {};
		\node [style=none] (48) at (27.75, -12.25) {};
		\node [style=none] (49) at (27, -13) {};
		\node [style=none] (50) at (9.25, -13) {};
		\node [style=none] (51) at (8.5, -12.25) {};
		\node [style=none] (71) at (17, -1.25) {};
		\node [style=none] (72) at (17.75, -2) {};
		\node [style=none] (73) at (17.75, -12.25) {};
		\node [style=none] (74) at (17, -13) {};
		\node [style=none] (75) at (22, -1.25) {};
		\node [style=none] (76) at (22.75, -2) {};
		\node [style=none] (77) at (22.75, -12.25) {};
		\node [style=none] (78) at (22, -13) {};
		\node [style=small text] (79) at (15, -3) {Euclidean \\ buildings \\ \cite{KlLe}};
		\node [style=small text] (80) at (20.25, -3) {$\R$-buildings \\ \cite{BrTi,Par}\\};
		\node [style=small text] (81) at (25.5, -3) {$\Lambda$-buildings \\ \cite{Ben1,HIL24} \\ \cite{KrTe}};
		\node [style=small text] (83) at (10.5, -12.25) {trees};
		\node [style=small text] (84) at (20, -12.25) {$\R$-trees};
		\node [style=small text] (85) at (25.25, -12.25) {$\Lambda$-trees};
		\node [style=tinytext] (86) at (1, -13.5) {$1$};
		\node [style=tinytext] (87) at (2, -13.5) {$0$};
		\node [style=small text] (88) at (19, -0.75) {affine buildings};
		\node [style=tinytext] (89) at (2, -12.25) {$1$};
		\node [style=tinytext] (90) at (1, -12.25) {$2$};
		\node [style=tinytext] (91) at (2, -10.75) {$2$};
		\node [style=tinytext] (92) at (1, -10.75) {$3$};
		\node [style=tinytext] (93) at (2, -9.25) {$3$};
		\node [style=tinytext] (94) at (1, -9.25) {$4$};
		\node [style=tinytext] (95) at (2, -7.75) {$\vdots$};
		\node [style=tinytext] (96) at (1, -7.75) {$\vdots$};
		\node [style=none] (97) at (3.75, -11.5) {};
		\node [style=none] (98) at (8.25, -11.5) {};
		\node [style=none] (99) at (8.5, -10) {};
		\node [style=none] (100) at (12.75, -10) {};
	\end{pgfonlayer}
	\begin{pgfonlayer}{edgelayer}
		\draw (0.center) to (7.center);
		\draw [bend left=45, looseness=1.75] (7.center) to (6.center);
		\draw (6.center) to (5.center);
		\draw [bend left=45, looseness=1.75] (5.center) to (4.center);
		\draw (4.center) to (3.center);
		\draw [bend left=45, looseness=1.75] (3.center) to (2.center);
		\draw (2.center) to (1.center);
		\draw [bend left=45, looseness=1.75] (1.center) to (0.center);
		\draw (14.center) to (20.center);
		\draw [bend right=45, looseness=1.75] (20.center) to (19.center);
		\draw (19.center) to (18.center);
		\draw [bend right=45, looseness=1.75] (18.center) to (17.center);
		\draw (17.center) to (16.center);
		\draw [style=fill edges] (14.center)
			 to [bend left=45, looseness=1.75] (13.center)
			 to (15.center)
			 to [bend left=45, looseness=1.75] (16.center)
			 to (98.center)
			 to (97.center)
			 to cycle;
		\draw [style=arrow] (21.center) to (22.center);
		\draw [style=soft edge] (26.center) to (5.center);
		\draw [style=soft edge, in=360, out=180] (28.center) to (27.center);
		\draw (45.center) to (51.center);
		\draw [bend right=45, looseness=1.75] (51.center) to (50.center);
		\draw (50.center) to (49.center);
		\draw [bend right=45, looseness=1.75] (49.center) to (48.center);
		\draw (48.center) to (47.center);
		\draw [bend right=45, looseness=1.75] (47.center) to (46.center);
		\draw (46.center) to (44.center);
		\draw [style=fill edges] (99.center)
			 to (100.center)
			 to (6.center)
			 to (44.center)
			 to [bend right=45, looseness=1.75] (45.center)
			 to cycle;
		\draw [bend right=45, looseness=1.75] (74.center) to (73.center);
		\draw (73.center) to (72.center);
		\draw [bend right=45, looseness=1.75] (72.center) to (71.center);
		\draw [bend right=45, looseness=1.75] (78.center) to (77.center);
		\draw (77.center) to (76.center);
		\draw [bend right=45, looseness=1.75] (76.center) to (75.center);
		\draw [style=soft edge] (29.center) to (40.center);
		\draw [style=soft edge] (31.center) to (42.center);
		\draw [style=soft edge] (30.center) to (41.center);
	\end{pgfonlayer}
\end{tikzpicture}
\end{equation*}
  	\caption{Buildings were first considered as simplicial objects, but for affine buildings the successively more general notions of Euclidean buildings, $\R$-buildings and $\Lambda$-buildings were introduced. All spherical buildings of rank at least 3 and all discrete affine buildings of rank at least 4 are associated to an algebraic group and some field. In this thesis we consider the most general notion, that of affine $\Lambda$-buildings. }  
  	\label{fig:buildings_overview}  
  \end{figure}
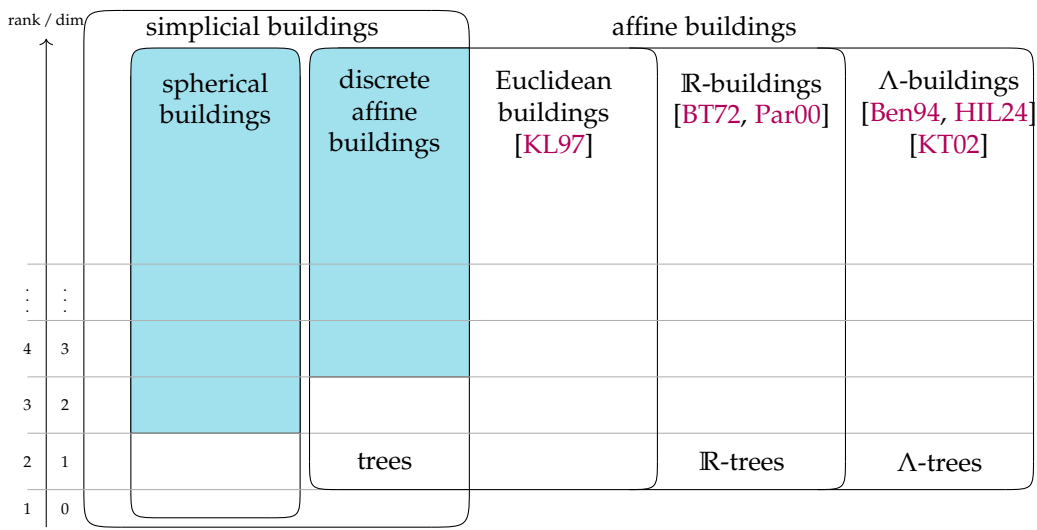
  
	\begin{figure}[t]
	\centering
	\includegraphics[scale=0.45]{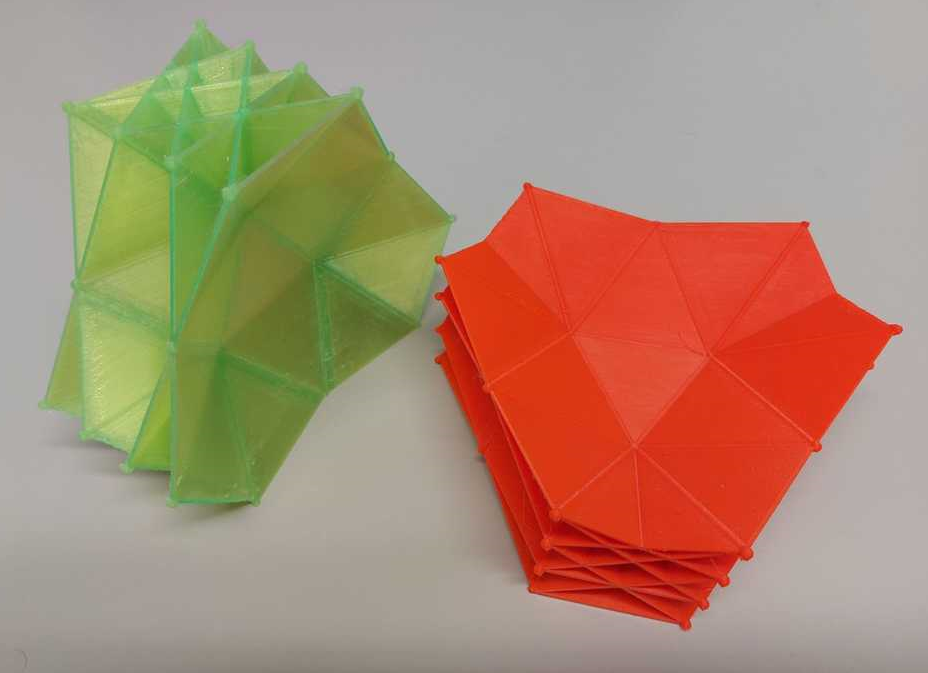}
	\caption{ Two 3D printed models of parts of a discrete affine building. The 3D models are available on \href{https://www.thingiverse.com/thing:6301564}{\texttt{www.thingiverse.com/thing:6301564}}.}
	\label{fig:affine_buildings_3D}
\end{figure}

	Buildings were first considered as simplicial complexes \cite{Ron,Bro}.  They consist of apartments glued together according to some rules of symmetry. The topological type of the apartments leads to a rough division of buildings, see Figure \ref{fig:buildings_overview}. Buildings with simplicial divisions of spheres as apartments are called spherical buildings and when the apartment is homeomorphic to $\mathbb{R}^n$, the building is called affine, for some 3D printed models of affine buildings, see Figure \ref{fig:affine_buildings_3D}. There is a large collection of buildings whose apartments may have various shapes. Spherical and affine buildings with apartments of large enough dimension are classified \cite{Wei03,Wei09}:
	they are all associated to semisimple algebraic groups over certain fields. In the case of affine buildings, the fields come with a discrete valuation.
	
	Starting with \cite{BrTi}, algebraic groups over fields with not necessarily discrete valuations were considered, which lead to the introduction of non-discrete affine buildings, or $\R$-buildings. These are defined as a set together with an atlas of apartments satisfying five axioms, though the exact choice of axioms varied by time and author. In this setting, apartments are Euclidean affine spaces equipped with the action of an affine Weyl group. Examples of $\R$-buildings were constructed and used in \cite{BrTi, KlLe, Par}.
	
	In dimension 1, discrete affine buildings correspond exactly to simplicial trees without leaves and $\R$-buildings correspond to $\R$-trees, also called real trees, without endpoints. Real trees appear in various mathematical areas, such as geometric group theory \cite{Sha87} and probability \cite{Ald91}. Real trees can also be characterized as those geodesically complete metric spaces that are $\delta$-hyperbolic for $\delta=0$, in the sense of Gromov \cite{Gro87}.
	In 1984, Morgan and Shalen \cite{MoSh} introduced $\Lambda$-trees for any linearly ordered abelian group $\Lambda$, see \cite{Chi01}. This concept generalizes simplicial trees as $\mathbb{Z}$-trees and real trees as $\mathbb{R}$-trees. When $\Lambda<\R$, $\Lambda$-trees offer a new way of describing $\R$-trees with special properties, for example a $\Q$-tree can be viewed as an $\R$-tree that only branches at rational distances from some base point. But there are also new phenomena when $\Lambda \not< \mathbb{R}$, for example when $\Lambda= \mathbb{Z}\times \mathbb{Z}$ with the lexicographical ordering, $\Lambda$-trees can be interpreted as trees of trees.
	
	The concepts of non-discrete buildings and $\Lambda$-trees were unified by Bennett \cite{Ben1}, leading to the concept of affine $\Lambda$-buildings or generalized affine buildings. Similar to how the step from $\mathbb{R}$-trees to $\Lambda$-trees required the addition of an axiom, Bennett increased the number of axioms by one to six axioms (A1) - (A6) for affine $\Lambda$-buildings. Various alternative sets of axioms were given in \cite{BeScSt}. In Chapter \ref{sec:lambda_trees} and \ref{sec:A6_independence} we will investigate the independence of some of the axioms for $\Lambda$-trees and affine $\Lambda$-buildings further.
	Bennett used homothety classes of lattices to construct an affine $\Lambda$-building associated to $\operatorname{SL}(n,\F)$ where $\F$ is any field with a valuation taking values in $\Lambda$. More recently \cite{HIL24} follow the approach in \cite{BrTi} to give examples in reductive groups. In \cite{KrTe}, Kramer and Tent suggest a construction of an affine $\Lambda$-building $\mathcal{B}$ based on a semisimple linear algebraic group defined over a non-Archimedean real closed field $\F$. A large part of this thesis is concerned with the construction and analysis of $\B$. We will construct $\B$ in Section \ref{sec:building_def} and prove that $\B$ satisfies the axioms (A1) - (A6) for affine $\Lambda$-buildings in Section \ref{sec:Bisbuilding}.
	
	There are multiple links between symmetric spaces and buildings. In 1973, Mostow \cite{Mos73} proved a rigidity result on lattices in higher rank Lie groups, now called Mostow-rigidity. A crucial step is to associate to a symmetric space of non-compact type its building at infinity, a spherical building. In 1997, Kleiner and Leeb \cite{KlLe} showed that the asymptotic cone of a symmetric space of non-compact type is a non-discrete building. They used this result to show rigidity of quasi-isometries. The asymptotic cone is a special case of the building $\B$, when the field is a Robinson field \cite{KrTe,BIPP21,BIPP23arxiv}. 
	
    A recent direction of applications lies in the theory of character varieties. Character varieties are spaces consisting of representations up to postconjugation of a finitely generated group $\Gamma$ into a Lie group $G$. A prototypical example of a component of the character variety of a surface group $\pi_1(S)$ in $\operatorname{SL}(2,\mathbb{R})$ is given by the space of marked hyperbolic structures $\mathcal{T}(S)$ on a surface $S$\footnote{The space of marked hyperbolic structures $\mathcal{T}(S)$ is often called Teichmüller space.}. Thurston \cite{Thu88}
    introduced a compactification of $\mathcal{T}(S)$ which lead to the classification of the elements of the mapping class group. The boundary points in Thurston's compactification correspond to actions of $\Gamma$ on certain $\R$-trees. Brumfiel \cite{Bru1} introduced the real spectrum compactification of $\mathcal{T}(S)$ and showed that it surjects onto Thurston's compactification. For general character varieties, Thurston's compactification was generalized to the Weyl chamber length spectrum compactification by Parreau \cite{Par12}.
    Generalizing Brumfiel's real spectrum compactification, Burger, Iozzi, Parreau and Pozetti \cite{BIPP21,BIPP23arxiv}
    introduced the real spectrum compactification of character varieties, and showed that it projects to the Weyl chamber length spectrum compactification. 
    Moreover the real spectrum compactification preserves components, 
    and therefore promises to be a useful tool in studying components. The points in the boundary of the real spectrum compactification can be identified with equivalence classes of representations $\Gamma \to G(\F)$ where $G(\F)$ is the group of $\F$-points of a semisimple algebraic group for a non-Archimedean real closed field $\F$. The group $G(\F)$ acts on the building $\B$ studied in this thesis, and hence every point in the boundary of the real spectrum compactification defines an action of $\Gamma$ on $\B$.    	
    
    Apart from applications to the study of symmetric spaces and character varieties, $\B$ is also an interesting affine $\Lambda$-building on its own. In \cite{Ben1}, an affine $\Lambda$-building of type $\tilde{\operatorname{A}}_n$ was constructed for every $\Lambda$ arising as the valuation group of a valued field. In \cite{HIL24} a building for quasi-split reductive groups was constructed. Apart from abstractly turning any affine $\R$-building into an affine $\Lambda$-building using the functoriality result of \cite{ScSt12}, the building $\B$ is the first example for general semisimple groups. Moreover in contrast to \cite{BrTi} and \cite{HIL24}, the construction of $\B$ is relatively straightforward as it does not rely on a large number of group-theoretic results. However the proof that $\B$ is an affine $\Lambda$-building does rely on such results.
	
	Central tools used in this thesis are real closed fields, semialgebraic geometry and the transfer principle, for an introduction see \cite{BCR}. The field of real numbers $\mathbb{R}$ is an example of a real closed field and all real closed fields are similar to $\mathbb{R}$ in a strong model-theoretic sense, namely the following transfer principle holds \cite{Tar51}: a first-order logic sentence holds over $\mathbb{R}$ if and only if it holds over any real closed field, or expressed in terms of model theory, the theory of real closed fields is complete. The transfer principle was used by Artin \cite{Art27} 
	to solve in the affirmative Hilbert's 17th problem on whether every positive real polynomial can be represented as a square of rational functions.

		\subsection{Semialgebraic groups and decompositions}
		
		For a detailed introduction to affine $\Lambda$-buildings, real closed fields and algebraic groups, see Sections \ref{sec:affine_L_buildings}, \ref{sec:real_closed} and \ref{sec:algebraic_groups}.
		Let $\F$ and $\K$ be real closed fields such that $\K \subseteq \F \cap \R$. Let $G$ be a Zariski-connected semisimple, self-adjoint (if $g\in G$ then $g\tran \in G$) linear algebraic $\K$-group\footnote{The whole theory could also be formulated for a group $\tilde{G}$ that is defined as a semialgebraic set, and that lies between the semialgebraic connected component of the identity $G_\F^\circ$ of $G_\F$ and $G_\F$ itself, so $G_\F^\circ < \tilde{G} < G_\F$, see \cite{BIPP23arxiv}.}
		that is a subgroup of $\operatorname{SL}(n)$. The $\K$-points of $G$ form a semialgebraic subset $G_\K \subseteq \mathbb{K}^{n\times n}$. The $\R$-extension $G_\R$ is a semisimple Lie group. We denote by $G_\F$ the $\F$-extension. We will call the study of $G_\F$ the \emph{semialgebraic setting}, as opposed to the algebraic setting (studying $G$ as an algebraic group) or the Lie group setting (studying $G_\R$). In this terminology a large part of the thesis is dedicated to the development of the semialgebraic setting.

	 Many decompositions of $G_\F$ rely on the choice of a torus. The following first result shows that all real closed fields define the same split tori.
	\begin{reptheorem}{thm:split_tori}
		A torus $S < G$ is maximal $\K$-split if and only if it is maximal $\F$-split. Moreover there is a self-adjoint maximal $\F$-split torus.
	\end{reptheorem} 
	To state decomposition theorems, we define certain subgroups of $G_\F$, for details see Section \ref{sec:decompositions}. Let $S_\K$ be the $\K$-points of a self-adjoint maximal $\K$-split torus. Let $A_\K$ be the semialgebraic connected component of $S_\K$ and $A_\F$ its semialgebraic extension. Let $K = G \cap \operatorname{SO}_n$, $K_\K$ the $\K$-points and $K_\F$ the semialgebraic extension, which coincides with the $\F$-points as $K$ is an algebraic group. Over the reals, $K_\R$ is a maximal compact subgroup. We also extend the semialgebraic groups $N_\K = \operatorname{Nor}_{K_\K}(A_\K)$ and $M_{\K} = \operatorname{Cen}_{K_{\K}}(A_\K)$ to $N_\F$ and $M_\F$. An order on the root system $\Sigma$ associated to $A_\F$ allows us to define $U_\K$ and $U_\F$ as the exponential of the sum of root spaces corresponding to positive roots.
	There are various definitions of root systems and Weyl groups in the literature. In Section \ref{sec:compatibility} we use Theorem \ref{thm:split_tori} to verify how these objects defined via the theory of algebraic groups, real Lie groups, Lie algebras and semialgebraic setting all coincide. This allows us to later use these settings interchangeably.
	\begin{repproposition}{prop:weylgroups}
		The algebraic root system $\KPhi$ is isomorphic to the root system $\Sigma$ from the real setting. The spherical Weyl groups $\KW, \FW, N_\R/M_\R, N_\F/M_\F$ and the group generated by reflections in roots of the root system $\Sigma$ are all isomorphic. 
	\end{repproposition}
	Decompositions play a central role in the theory of Lie groups, the theory of algebraic groups and in this thesis. We prove the following versions of the Iwasawa (KAU), Cartan (KAK) and Bruhat (BWB) decomposition for $G_\F$.
	\begin{reptheorem}{thm:KAU}[$G=KAU$]
		For every $g \in G_\F$, there are $k \in K_\F , \, a \in A_\F,\, u \in U_\F$ such that $g = kau$. This decomposition is unique.
	\end{reptheorem} 
	\begin{reptheorem}{thm:KAK}[$G=KAK$]
		For every $g \in G_\F$, there are $k_1,k_2 \in K_\F, \, a \in A_\F$ such that $g = k_1 a k_2$.
		In this decomposition $a$ is uniquely determined up to a conjugation by an element of $N_{\F}/M_{\F}$.
	\end{reptheorem}
	\begin{reptheorem}{thm:BWB}[$G=BWB$]
		For every $g \in G_\F$, there are $b_1,b_2 \in B_\F := M_\F A_\F U_\F$ and $n\in N_\F$ such that $g = b_1 n b_2$. In this decomposition $n$ is unique up to multiplying by an element in $M_{\F}$. For the spherical Weyl group $W_s := N_\F/M_\F$, we have a disjoint union of double cosets
		$$
		G_{\F} = \prod_{[n] \in W_s} B_{\F}nB_{\F}.
		$$ 
	\end{reptheorem}
	
	In Sections \ref{sec:BCH}, \ref{sec:U} we note that while in contrast to the setting of Lie groups the exponential map is not defined for $G_\F$, it is still defined for $U_\F$ as the elements of $U_\F$ are nilpotent. As a consequence, the Baker-Campbell-Hausdorff formula holds for elements in $U_\F$ which is useful in proofs using induction over the root system. 
	
	In Section \ref{sec:Jacobson_Morozov} we give a proof of the folklore result that the Jacobson-Morozov Lemma holds for algebraic groups.
	\begin{repproposition}{prop:Jacobson_Morozov_group}
		Let $g\in G$ be any unipotent element in a semisimple linear algebraic group $G$ over an algebraically closed field $\D$ of characteristic $0$. Then there is an algebraic subgroup $\operatorname{SL}_g < G$ with Lie algebra $\operatorname{Lie}(\operatorname{SL}_g) \cong \mathfrak{sl}_2$ and $g \in \operatorname{SL}_g$. The element $\log(g) \in \operatorname{Lie}(G)$ corresponds to 
		$$
		\begin{pmatrix}
			0 & 1 \\ 0 & 0
		\end{pmatrix}\in \mathfrak{sl}_2.
		$$
		Moreover, if $g\in G_\F$ for a field $\F \subseteq \D$, then $\operatorname{SL}_g$ can be assumed to be defined over $\F$. 
	\end{repproposition} 
	For us, the following more precise version of the Jacobson-Morozov Lemma, now in the semialgebraic setting, will be useful. Note that the assumption on $\alpha$ in Proposition \ref{prop:Jacobson_Morozov_real_closed} means that it can only be applied to reduced root systems (if $\alpha \in \Sigma$, then $2\alpha \notin \Sigma$).
	\begin{repproposition}{prop:Jacobson_Morozov_real_closed}
		Let $\F$ be a real closed field with an order compatible valuation $v \colon \F \to \Lambda \cup \{ \infty\}$, $\alpha \in \Sigma$ and assume $\frakg_{2\alpha} = 0$. Let $u \in (U_{\alpha})_\F$. Then there are $ X \in (\frakg_\alpha)_\F$ and $t \in \F$ such that $u=\exp(tX)$ and $\max_{ij}\left\{ (-v)(X_{ij}) \right\} = 0$. Let $(X,Y,H)$ be the $\mathfrak{sl}_2$-triplet of the Jacobson-Morozov Lemma in the Lie algebra setting, see Lemma \ref{lem:JM_basic}. Then there is a morphism of algebraic groups
		$
		\varphi \colon \operatorname{SL}(2,\D)  \to G 
		$
		defined over $\F$ such that $\varphi$ has finite kernel and
		\begin{align*}
			\varphi \begin{pmatrix}
				1 & t \\ 0& 1
			\end{pmatrix} =  u = \exp(tX)  \quad \text{and} \quad  \varphi \begin{pmatrix}
				1 & 0 \\ t & 1
			\end{pmatrix} = \exp(tY).
		\end{align*}
		If $\varphi$ is not injective, then $\ker (\varphi) \cong \Z/2\Z$ and $\varphi$ factors through the isomorphism
		$$
		\operatorname{PGL}(2,\D) := \operatorname{SL(2,\D)}/\ker(\varphi)  \xrightarrow{\sim} \varphi(\operatorname{SL}(2,\D))
		$$
		which is also defined over $\F$. Moreover $\varphi(g\tran) = \varphi(g)\tran$, for any $g \in \operatorname{SL}(2,\F)$.
	\end{repproposition}
	
	The Jacobson-Morozov Lemma produces subgroups with Lie algebras isomorphic to $\mathfrak{sl}_2$. The following theorem produces potentially larger rank one subgroups $L_{\pm \alpha}$ associated to a root $\alpha \in \Sigma$ in the algebraic setting. These subgroups are the semisimple parts of Levi subgroups. 
	\begin{reptheorem}{thm:levi_group}
		Let $\alpha \in \Sigma$. Then there is a connected semisimple self-adjoint linear algebraic group $L_{\pm \alpha}$ defined over $\K$ such that
		\begin{enumerate}
			\item [(i)] $\operatorname{Lie}(L_{\pm \alpha}) = ( \frakg_\alpha \oplus \frakg_{2\alpha} ) \oplus (\frakg_{-\alpha} \oplus \frakg_{-2\alpha} ) \oplus ( [\frakg_\alpha , \frakg_{-\alpha}] + [\frakg_{2\alpha}, \frakg_{-2\alpha}] )$.
			\item [(ii)] $\operatorname{Rank}_\R(L_{\pm \alpha}) = \operatorname{Rank}_\F(L_{\pm \alpha}) = 1$.
		\end{enumerate}
	\end{reptheorem}
	
	\subsection{The nonstandard symmetric space $\X_\F$}
	Now $\F$ is a non-Archimedean real closed field with order compatible valuation $v \colon \F \to \Lambda \cup \{ \infty\}$ and the groups are defined as in the previous subsection. The semialgebraic set
	$$
	P_\R := \left\{ x \in \R^{n\times n} \colon x=x\tran, \det(x)=1, x \text{ is positive definite }  \right\} 
	$$
    is a model for the symmetric space $\operatorname{SL}(n,\R)/\operatorname{SO}(n,\R)$. For the Lie group $G_\R < \operatorname{SL}(n,\R)$, the symmetric space $G_\R/K_\R$ can be realized as a totally geodesic submanifold of $P_\R$ given by the semialgebraic orbit $\X_\R := G_\R.\operatorname{Id} \subseteq P_\R$. The semialgebraic extension $\X_\F$ is called the \emph{non-standard symmetric space}. 
    
    Over the reals, the symmetric space $\X_\R$ comes equipped with a Riemannian distance, such that $G_\R$ acts by isometries: the maximal flats are isometric to a Euclidean vector space and for any two points there is a maximal flat containing the two, so the distance is given as the norm of the difference of the two. We mimic this construction in the semialgebraic setting. First we project two arbitrary points to the \emph{non-standard maximal flat} $A_\F.\operatorname{Id} \subseteq \X_\F$, which contains the fundamental Weyl chamber $A_\F^+.\operatorname{Id}$ for the semialgebraic subset 
    $$
    A_\F^+ := \left\{ a \in A_\F \colon \chi_\alpha(a) \geq 1 \text{ for all } \alpha \in \Sigma_{>0}  \right\}.
    $$
    The Cartan decomposition KAK can be used to define a Weyl chamber valued distance on $\X_\R$, the semialgebraic version of which is the following Lemma.  
	
	\begin{replemma}{lem:Cartan_projection}
		The Cartan decomposition defines a Cartan-projection 
		\begin{align*}
			\delta_\F \colon \X_\F \times \X_\F & \to \ApF\\
			(x,y) &\mapsto a
		\end{align*}
		which is invariant under the action of $G_\F$. 
	\end{replemma} 
	We then use the semialgebraic norm 
		\begin{align*}
		N_{\F} \colon A_{\F} &\to \F_{\geq 1} \\
		a & \mapsto  \prod_{\alpha \in \Sigma} \max \left\{\chi_\alpha(a), \chi_\alpha (a)^{-1}\right\},
	\end{align*}
	Over the reals, $d = \log \circ N_\R \circ \delta_\R$ defines a metric on $\X_\R$. The logarithm is not a semialgebraic function, but we can replace $\log$ by the valuation $-v$, to obtain a symmetric non-negative function 
	\begin{align*}
			d \colon \X_\F \times \X_\F &\to \Lambda \\
			(x,y) & \mapsto (-v)(N_\F (\delta_{\F}(x,y))).
	\end{align*}
	The function $d$ is not quite a $\Lambda$-metric, since there are points $x\neq y \in \X_\F$ with $d(x,y) = 0$, but we have the following.
	\begin{reptheorem}{thm:pseudodistance}
		The function $d \colon \X_\F \times \X_\F \to \Lambda$ is a pseudo-distance on $\X_\F$.
	\end{reptheorem}
    To prove the triangle inequality in Theorem \ref{thm:pseudodistance}, we use the following two results of independent interest. First, the Iwasawa decomposition KAU can be used to define the semialgebraic retraction $a_\F \colon G_\F \to A_\F$ by $g = uak \mapsto a =: a_\F(g)$. For some $b\in A_\R$, Kostant's convexity theorem then describes the set $a_\R(K_\R b)$ as a polytope in $A_\R$ given by the convex hull of the Weyl group orbit of $b$. In the semialgebraic setting, we prove the following variation for the elements $b \in A_\F^+$ for certain characters $\chi_i  \colon A_\F \to \F^\times$, $i \in \{1,2, \ldots, \dim(A_\R) \}$.
	
	\begin{reptheorem}{thm:kostant_F}
		For all $b \in \ApF$, we have
		$$
		\left\{ a \in \ApF \colon \exists k \in K_\F , a = a_\F(kb)  \right\}
		= \left\{ a\in \ApF \colon \chi_i (a) \leq \chi_i (b) \text{ for all } i \right\}.
		$$
	\end{reptheorem}
	Theorem \ref{thm:kostant_F} can then be used to show that the Iwasawa retraction $\rho \colon \X_\F \to A_\F.\operatorname{Id}$ defined by $g.\operatorname{Id} \mapsto a_\F(g).\operatorname{Id}$ decreases distances.
	\begin{reptheorem}{thm:UAKret}
		The map $\rho \colon \X_\F \to A_\F.\operatorname{Id}$ is a $d$-diminishing
		retraction to $A_\F.\operatorname{Id}$.
	\end{reptheorem}

	\subsection{The building $\B$}
	In view of Theorem \ref{thm:pseudodistance}, we can define the $\Lambda$-metric space $\B := \X_\F/\!\!\sim$ where $x \sim y$ whenever $d(x,y) = 0 \in \Lambda$. We note that $G_\F$ acts transitively by isometries on $\B$. We denote the basepoint by $o := [\operatorname{Id}] \in \B$. In Section \ref{sec:Xapartment}, we verify that $\A := A_\F.\operatorname{o} \subseteq \B$ is an apartment in the sense of affine $\Lambda$-buildings, and that $\operatorname{Nor}_{G_\F}(A_\F)/ \operatorname{Cen}_{G_\F}(A_\F)$ acts as the affine Weyl group $\A \rtimes W_s$, where $W_s$ is the spherical Weyl group of the root system $\Sigma$. The inclusion $f_0 \colon \A \to \B$ can then be used to define an atlas
	$
	\Fun = \left\{  g.f_0 \colon \A \to \B  \colon  g\in G_\F  \right\}.
	$
    The main result of this thesis is that the $\Lambda$-metric space $\B$ together with the atlas $\Fun$ is a generalized affine building. Theorem \ref{thm:B_is_building} was announced without proof as Theorem 5.7 in \cite{KrTe} and as Theorem 4.3 in \cite{KrTe04}. 
    \begin{reptheorem}{thm:B_is_building}
    	If the root system $\Sigma$ is reduced, then the pair $(\B,\Fun)$ is an affine $\Lambda$-building.
    \end{reptheorem}
    
    \begin{remark}
    	 The technical assumption that the root system $\Sigma$ is reduced is an artifact of our use of the Jacobson-Morozov Lemma, Proposition \ref{prop:Jacobson_Morozov_real_closed}. As a consequence of the classification of real Lie groups, see for example the appendix of \cite{OV90}, only Hermitian groups of non-tube type have non-reduced root systems; and they are of type $\operatorname{BC}_n$ \cite{GuWi24arxiv}.
    	 In the study of character varieties, components with discrete and faithful representations are only known to exist when the group admits a $\Theta$-positive structure. All groups admitting a $\Theta$-positive structure have reduced root systems and our theorem can be applied to points in the boundary of the real spectrum compactification in this case. 
    \end{remark}
    
	To be more precise, we verify that $(\B,\Fun)$ satisfies the following six axioms, which are equivalent to the original axioms for affine $\Lambda$-buildings by the work of \cite{BeScSt}.
	
	\begin{enumerate}
		\item [(A1)] For all $f\in \Fun$, $w\in W_a$, $f \circ w \in \Fun$.
		\item [(A2)] For all $f,f'\in \Fun$, the set $B= f^{-1}\left( f(\A)\cap f'(\A)  \right)$ is $W_a$-convex and there is a $w\in W_a$ such that
		$
		f|_B = f'\circ w |_B.
		$
		\item [(A3)] For all $x,y\in X$, there is a $f\in \Fun$ such that $x,y\in f(\A)$.
		\item [(A4)] For any sectors $s_1,s_2 \subseteq X$ there are subsectors $s_1' \subseteq s_1,\,  s_2' \subseteq s_2$ such that there is an $f\in \Fun$ with $s_1', s_2' \subseteq f(\A)$.
\item [(TI)] The function $d \colon \B \times \B \to \Lambda$ induced from the distance on apartments (via axioms (A1), (A2) and (A3)) satisfies the triangle inequality.
\item [(EC)] For $f_1,f_2 \in \Fun$, if $f_1(\A)\cap f_2(\A)$ is a half-apartment, then there exists $f_3 \in \Fun$ such that $f_i(\A)\cap f_3(\A)$ are half-apartments for $i\in \{1,2\}$. Moreover $f_3(\A)$ is the symmetric difference of $f_1(\A)$ and $f_2(\A)$ together with the boundary wall of $f_1(\A) \cap f_2(\A)$.
	\end{enumerate}
	
	The proof of Theorem \ref{thm:B_is_building} uses all the theory we have developed for $G_\F$ up to this point and relies on further analysis of $\B$, developed in Section \ref{sec:root_groups_and_BT}. The proof of (A2) is the most involved among the axioms. We adapt some of the methods in \cite{BrTi} to our setting to be able to prove (A2). For the special case of $G_\F = \operatorname{SL}(n,\F)$, a direct proof of axiom (A2) is given in Appendix \ref{sec:appendixSLn}. We will point out some results that may be of independent interest. Many of the following results are about stabilizers of certain subsets. 
	
	The stabilizer of $o \in \B$ was calculated by \cite{Tho} for the special case when $\F$ is a Robinson field. For general fields, it has been suggested by \cite{KrTe} and proven by \cite{BIPP23arxiv}. Let $O = \{ a \in \F \colon (-v)(a) \leq 0 \}$ be the valuation ring associated to the valuation $v$ and let $E_\F(O) := E_\F \cap \operatorname{SL}(n,O)$ for any semialgebraic subset $E_\F \subseteq G_\F$.
	\begin{reptheorem}{thm:stab} 
		The stabilizer of $o\in\B$ in $G_\F$ is $G_\F(O) $.
	\end{reptheorem}
	As a consequence of the Iwasawa retraction Theorem \ref{thm:UAKret}, we give an Iwasawa decomposition of the stabilizer of $o$.
	\begin{repcorollary}{cor:UAKO}
		There is an Iwasawa decomposition $G_\F(O) = U_\F(O) A_\F(O) K_\F$, meaning that for every $g \in G_\F(O)$ there are unique $u \in U_\F(O),\, a\in A_\F(O), \, k\in K_\F=K_\F(O)$ with $g=uak$.
	\end{repcorollary}
    
    In Section \ref{sec:root_groups_and_BT} we implement some ideas from \cite{BrTi} in the semialgebraic setting. For the root groups $(U_{\alpha})_\F = \exp((\frakg_\alpha \oplus \frakg_{2\alpha})_\F)$, we introduce root group valuations $\varphi_\alpha \colon U_\alpha \to \Lambda \cup \{-\infty \}$ in Subsection \ref{sec:root_group_valuation}. 
    Taking a closer look at the unipotent group, we notice that if $u \in U_\F$ sends some point of the apartment to some other point of the apartment, then the two points have to be the same.
	\begin{repproposition}{prop:UonA}
		For all $u\in U_\F$ and $a \in A_\F$,
		$
		ua.o \in \mathbb{A} \iff ua.o = a.o. 
		$
	\end{repproposition}
	In fact, the fixed point set of elements in $U_\F$ is a half-apartment. 
	\begin{repproposition}{prop:Ualphaconv}
		Let $\alpha \in \Sigma$. For $u \in (U_\alpha)_\F$ we have
		$$
		\left\{ p \in \mathbb{A} \colon u.p  \in \mathbb{A}  \right\} = \left\{ a.o \in \A \colon \varphi_\alpha(u) \leq (-v)\left(\chi_\alpha\left(a\right)\right)  \right\}
		$$
		and therefore this set is a half-apartment, when $u \neq \operatorname{Id}$.
	\end{repproposition}
	This can be upgraded to a finite intersection of half-apartments for elements in $U_\F$.
	\begin{repproposition}{prop:Uconv}
		For $u\in U_\F$ there are $k_\alpha \in \Lambda \cup \{-\infty \}$ for $\alpha \in \Sigma_{>0}$ such that
		$$
		\{p \in \A \colon u.p \in \A\} = \left\{ a.o \in \A \colon  k_\alpha \leq (-v)\left(\chi_\alpha \left(a\right)\right) \text{ for all } \alpha \in \Sigma_{>0}   \right\}
		$$
		and therefore the set of fixed points is a finite intersection of half-apartments. If $u$ fixes all of $\A$, then $u = \operatorname{Id}$. 
	\end{repproposition}
	
	We proceed to describe the pointwise stabilizers of the whole apartment $\A$, the fundamental Weyl chamber
	$$
	C_0 = \{ a.o \in \A \colon (-v)(\chi_\alpha(a)) \geq 0\},
	$$
	and the half-apartments
	$$
	H_\alpha^+ = \{ a.o \in \A \colon (-v)(\chi_\alpha(a)) \geq 0\}.
	$$
	\begin{reptheorem}{thm:NO_fixes_s0}
		The pointwise stabilizer of $C_0$ in $G_\F$ is $B_\F(O) = U_\F(O)A_\F(O)M_\F$.
	\end{reptheorem}
	
	\begin{reptheorem}{thm:stab_A}
		The pointwise stabilizer of $\A$ in $G_\F$ is $A_\F(O)M_\F$.
	\end{reptheorem}
	
	\begin{reptheorem}{thm:NalphaO_fixes_H}
		Let $\alpha \in \Sigma$. The pointwise stabilizer of $H_\alpha^+$ in $G_\F$ is $(U_\alpha)_\F(O) A_\F(O) M_\F$.
	\end{reptheorem}
	
	In Subsection \ref{sec:BT_rank_1}, we take a closer look at the rank 1 subgroup $(L_{\pm \alpha})_\F < G_\F$. We now restrict to groups $G_\F$ with reduced root systems, to be able to use the Jacobson-Morozov Lemma to find elements $m(u)\in \operatorname{Nor}_{G_\F}(A_\F)$ representing a reflection along some hyperplane
	$
	M_{\alpha,\ell} = \{ a.o \in \A \colon (-v)(\chi_\alpha(a)) = \ell \}
	$
	for $\alpha\in \Sigma$ and $\ell \in \Lambda$. A careful analysis of the rank one subgroup $L = \langle (U_\alpha)_\F, (U_{-\alpha})_\F \rangle < (L_{\pm \alpha})_\F$ results in the following decomposition of its stabilizer.
\begin{repproposition}{prop:stab_L_Omega}
	Assume $\Sigma$ is reduced. Let $\alpha \in \Sigma$, $\Omega \subseteq \A$ a non-empty finite subset. The pointwise stabilizers
	\begin{align*}
		L_\Omega &:= L\cap \operatorname{Stab}_{G_\F}(\Omega),\\
		 U_{\alpha,\Omega} &:= (U_\alpha)_\F \cap \operatorname{Stab}_{G_\F}(\Omega), \\
		  U_{-\alpha,\Omega} &:= (U_{-\alpha})_\F \cap \operatorname{Stab}_{G_\F}(\Omega), \\
		 N_\Omega &:= \operatorname{Nor}_{G_\F}(A_\F) \cap \operatorname{Stab}_{G_\F}(\Omega)
	\end{align*}
	satisfy 
	$
	L_\Omega = \langle U_{\alpha,\Omega} , U_{-\alpha,\Omega} , (N_\Omega \cap L) \rangle =  U_{\alpha, \Omega} U_{-\alpha,\Omega} (N_\Omega \cap L).
	$
\end{repproposition}
	
	In Section \ref{sec:BT_higher_rank}, we upgrade the rank 1 result to $G_\F$. For a subset $\Omega \subseteq \A$, let
	\begin{align*}
	\hat{P}_\Omega &:= \langle N_{\Omega}, U_{\alpha,\Omega}\colon \alpha \in \Sigma \rangle \\
	U_\Omega^+ &:= \langle U_{\alpha,\Omega} \colon \alpha \in \Sigma_{>0} \rangle \\ 
		U_\Omega^- &:= \langle U_{\alpha,\Omega} \colon \alpha \in \Sigma_{<0} \rangle.
	\end{align*}
	\begin{reptheorem}{thm:BTstab}
		Assume $\Sigma$ is reduced. The pointwise stabilizer of an arbitrary subset $\Omega \subseteq \A$ satisfies
		$
		\operatorname{Stab}_{G_\F}(\Omega) = \hat{P}_\Omega = U_\Omega^+ U_{\Omega}^- N_{\Omega} .
		$
	\end{reptheorem}
	
	The proof of Theorem \ref{thm:BTstab} is tightly intertwined with the proof of axiom (A2). It first relies on the following \emph{mixed Iwasawa} decomposition for $G_\F$. For $o\in \B$, let $\hat{P}_o := \hat{P}_{\{o\}}$. 
	\begin{reptheorem}{thm:BT_mixed_Iwasawa}
		Assume $\Sigma$ is reduced. Then
		$G_\F = U_\F \cdot \operatorname{Nor}_{G_\F}(A_\F) \cdot \hat{P}_o$.
	\end{reptheorem}
	With the help of the mixed Iwasawa decomposition we obtain Theorem \ref{thm:BTstab} in the case of finite subsets $\Omega \subseteq \A$. This statement for finite subsets is then enough to prove the second part of axiom (A2). The full Theorem \ref{thm:BTstab} then follows from this part of (A2) and the rest of axiom (A2) follows from the full Theorem \ref{thm:BTstab}.

	 \subsection{Further results on $\Lambda$-trees}
	 
	 To build up intuition in the study of affine $\Lambda$-buildings it often makes sense to first consider the rank one case of $\Lambda$-trees. We sketch some results about $\Lambda$-trees obtained as part of the thesis, details can be found in Section \ref{sec:lambda_trees} and \cite{App24, ADFJ24}.
	 In the axiomatic characterization, $\Lambda$-trees are defined as $\Lambda$-metric spaces satisfying three axioms (1), (2) and (3). The following theorem gives a purely algebraic condition on $\Lambda$ for when the geometric property of axiom (3) follows from (1) and (2).
	 
	 	\begin{reptheorem}{thm:main_Lambda_trees}(\cite[Theorem 1]{App24})
	 	Let $\Lambda$ be an ordered abelian group. The following are equivalent:
	 	\begin{itemize}
	 		\item [(a)] For every positive $\lambda_0 \in \Lambda$, the set
	 		$
	 		\{t \in \Lambda \colon 0 \leq 2t\leq \lambda_0 \}
	 		$
	 		has a maximum.
	 		\item [(b)] Every $\Lambda$-metric space that satisfies axioms (1) and (2) also satisfies (3).
	 	\end{itemize}
	 \end{reptheorem} 
	 
	 For the case of $G = \operatorname{SL}(2)$, Brumfiel \cite{Bru2} constructs a non-standard hyperbolic plane with a pseudometric and proves that the quotient $\B$ is a $\Lambda$-tree. Starting for simplicity with the field $\F$ of Puiseux-series, see definition in Section \ref{sec:real_closed_fields_examples}, $\B$ is a $\Q$-tree $\mathcal{T}$. Completing all the segments of $\mathcal{T}$ results in a $\R$-tree $\overline{\mathcal{T}}^{\operatorname{sc}}$. The segment completion is a special case of a general base change functor for affine $\Lambda$-buildings investigated in \cite{ScSt12}. The following theorem answers the question of whether the segment-completion $\overline{\mathcal{T}}^{\operatorname{sc}}$ is complete as a metric space. 
	 
	 \begin{reptheorem}{thm:ARFJ}(\cite[Theorem 1]{ ADFJ24})
	 	Let $\mathcal{T}$ be the $\Q$-tree constructed from the field of Puiseux-series in \cite{Bru2}. There is a Cauchy sequence in the segment completion $\overline{\mathcal{T}}^{\operatorname{sc}}$ that does not converge and hence $\overline{\mathcal{T}}^{\operatorname{sc}}$ is not complete.
	 \end{reptheorem}

	 \subsection{Formalization in the proof assistant \lean}
	 
	 Mathematics thrives on the intuition surrounding its objects but with its rigor also promises to be a pathway to absolute truth. In practice, mathematicians are human and mistakes happen. With the advent of computers, it has become possible to verify proofs by machine and to considerably improve the reliability of mathematical results. The \lean\ theorem prover is a proof assistant developed mainly by Leonardo de Moura at Microsoft Research \cite{MKADR15}. There is an extensive community-built mathematical library \texttt{mathlib} \cite{mC19}, which by now contains a large part of undergraduate mathematics. Recently some more advanced projects such as a formalization of perfectoid spaces \cite{BCM20} have been completed and a project to formalize Andrew Wiles' proof of Fermat's last theorem has been announced.
	 	
	 	As a proof of concept, some parts of this thesis were formalized in \lean. We built on the definition of ordered abelian groups already present in \texttt{mathlib}. We first formalized the notions of $\Lambda$-metric spaces and $\Lambda$-trees. We then showed some basic properties of $\Lambda$-trees before successfully formalizing that (a) implies (b) in Theorem \ref{thm:main_Lambda_trees}. This project showed that it is possible to verify new research results with a computer, however the time investment is substantial and we would therefore currently recommend not to formalize most new results. The \lean-files can be found in \cite{App23}. 
	 
	 \newpage

	\section{Real closed fields}\label{sec:real_closed}
	An \emph{ordered field} is a field together with a total order such that the sum and the product of positive elements are positive. 
	A field $\F$ is called \emph{real closed} if it satisfies one of the following equivalent conditions.
	\begin{enumerate}
		\item [(1)] There is a total order on $\F$ turning $\F$ into an ordered field such that every positive element has a square root and every polynomial of odd degree has a solution.
		\item [(2)] There is an order on $\F$ that does not extend to any proper algebraic field extension of $\F$.
		\item [(3)] $\F$ is not algebraically closed but every finite extension is algebraically closed.
		\item [(4)] $\F$ is not algebraically closed but $\F[\sqrt{-1}]$ is algebraically closed.
	\end{enumerate}

	An ordered field is called \emph{Archimedean} if every element is bounded by a natural number. The real numbers $\R$ and the subset of real algebraic numbers are examples of Archimedean real closed fields. A major tool when working with real closed fields is the following transfer principle from model theory.

	\subsection{The transfer principle}
	
	Recall that a \emph{first-order formula of ordered fields with parameters in $\F$} is a formula that contains a finite number of conjunctions $\wedge$, disjunctions $\lor$, negations $\lnot$, and universal $\forall$ or existential $\exists$ quantifiers on variables, starting from atomic formulas which are formulas of the kind $f(x_1, \ldots, x_n) = 0$ or $g(x_1, \ldots , x_n) \leq 0$, where $f$ and $g$ are polynomials with coefficients in $\F$. A first-order formula without free variables is called a \emph{sentence}. 
	\begin{theorem}\label{thm:logic} (Transfer principle, \cite{BCR})
		Let $\F$ and $\F'$ be real closed fields. Let $\varphi$ be a sentence with parameters in $\F \cap \F'$. Then $\varphi$ is true for $\F$ if and only if $\varphi$ is true for $\F'$, formally $\F \models \varphi \iff \F' \models \varphi$.
	\end{theorem}

	Let $\varphi$ be a first-order formula with parameters in some field $\K \subseteq \F \cap \R$ with $n$ free variables. Let $X$ be a subset of $\K^n$ which can be described as $X = \{x \in \K^n \colon \F \models \varphi(x)\}$. It follows from the transfer principle, that the \emph{semialgebraic extension $X_\F = \{x \in \F^n \colon \F \models \varphi(x)\}$ of $X$} depends only on $X$ and not on $\varphi$ and is thus well defined. Sets of the form $\{x \in k \colon k \models \varphi(x)\}$ for any ordered field $k \supseteq \K$ are called \emph{semialgebraic} sets.

	\subsection{Valuations}
	
	For an overview of valued fields, we recommend \cite{BaSh18}. A subring $O \subseteq \F$ of an ordered field $\F$ is an \emph{order convex subring}, if for all $a,b \in O,\  c \in \F, \ a \leq c \leq b $ implies $c \in O$. Note that every order convex subring is in particular a valuation ring: for all $a \in \F$, we have $a \in O$ or $a^{-1} \in O $. Let $(\Lambda,+)$ be an ordered abelian group. A \emph{valuation} on an ordered field $\F$ is a map $v \colon \F \to \Lambda \cup \{\infty\}$ which satisfies for all $a,b \in \F$
	\begin{itemize}
		\item [(1)] $v(a) = \infty$ if and only if $a = 0$.
		\item [(2)] $v(ab) = v(a) + v(b)$.
		\item [(3)] $v(a+b) \geq \min \{ v(a), v(b)\}$. 
	\end{itemize}
	We say that the valuation is \emph{order compatible}, if $(-v)(a) \geq (-v)(b)$ whenever $a \geq b \geq 0$. We will often use the same letter for $v$ and $v|_{\F_{>0}} \colon \F_{>0} \to \Lambda$. We will often be more interested in $(-v)$ than in $v$, as $(-v)$ is order preserving. 
	
	There is a correspondence between order convex valuation rings and order compatible valuations, which follows from a theorem by Krull \cite{Kru}.
	
	\begin{theorem}
		Every order convex valuation ring $O\subseteq \F$ gives rise to an order compatible valuation
		\begin{align*}
			v \colon \F & \to \Lambda \cup \{\infty\},
		\end{align*}
		where $\Lambda \cong \F^\times/O^\times$. On the other hand, every order compatible valuation gives rise to a order convex valuation ring
		$$
		O = \{ a \in \F \colon v(a)>0\}.
		$$
	\end{theorem}
	
	Real closed fields may admit order compatible valuation rings. The ordered abelian group $\Lambda$ of a valuation of a real closed field is always a $\Q$-vector space.

	\subsection{Examples}\label{sec:real_closed_fields_examples}
	The field of \emph{Puiseux series} over the real algebraic numbers $\Qrc$
	$$
	\F  := \left\{  \sum_{k=-\infty}^{k_0} c_k X^{\frac{k}{m}} \, \colon \, k_0, m \in \mathbb{Z}, \, m > 0 , \, c_k \in \Qrc, \, c_{k_0}\neq 0\right\},
	$$ 
	is a non-Archimedean real closed field, where the usual order on $\Qrc$ is extended by $X>r$ for all $r \in \Qrc$ \cite{BCR}. An order compatible valuation $v \colon \F \to \Lambda\cup \{\infty\} =  \Q\cup \{\infty\}$ is given by the degree
	\begin{align*}
		v\left(  \sum_{k=-\infty}^{k_0} c_k X^{k/m}\right)  =-\frac{k_0}{m}.
	\end{align*} 
	
	A \emph{non-principal ultrafilter on $\Z$} is a function $\omega \colon \mathcal{P}(\Z) \to \{0,1\}$ that satisfies
	\begin{enumerate}
		\item [(1)] $\omega(\emptyset) = 0$, $\omega(\Z) = 1$
		\item [(2)] If $A,B \subseteq \Z$ satisfy $A\cap B = \emptyset$, then $\omega(A \cup B) = \omega(A) + \omega(B)$.
		\item [(3)] All finite subsets $A \subseteq \Z$ satisfy $\omega(A) = 0$.
	\end{enumerate}
	Ultrafilters can be thought of as finitely-additive probability measures that only take values in $0$ and $1$. The existence of non-principal ultrafilters is equivalent to the axiom of choice, \cite{Hal}. For a given ultrafilter $\omega$, we define the \emph{hyperreal numbers} $\R_{\omega}$ to be the equivalence classes of infinite sequences $\R_{\omega} = \R^{\N}/ \!\! \sim$, where $x = (x_i)_{i \in \N} \sim y = (y_i)_{i \in \N}$ if $\omega(\{ i \in \N \colon x_i \neq y_i\}) = 0$ or $\omega(\{ i \in \N \colon x_i = y_i\}) = 1$. We define addition and multiplication componentwise, the multiplicative inverse is obtained by taking the inverses of all non-zero entries, turning $\F_\omega$ into a field. Considering constant sequences, the real numbers are a subfield of $\R_\omega$. The hyperreals are an ordered field with respect to the order defined by $[(x_i)_{i\in \N}] \leq [(y_i)_{i \in \N} ]$ if and only if $\omega(\{i \in \N \colon x_i \leq y_i\}) = 1$. The hyperreals are real closed, since $\R$ is. One can check that the hyperreals do not admit a valuation to a subgroup of $\R$, but they do admit a valuation to $(\Lambda, +) := ((\R_{\omega})_{\neq 0}, \cdot )$, namely the identity. The hyperreals are non-Archimedean, since the equivalence class containing $(1, 2, 3, \ldots)$ is an \emph{infinite} element, meaning it is larger than any natural number. 
	
	Let $b \in \R_{\omega}$ be an infinite element. Then
	$$
	O_b := \{ x \in \R_{\omega} \colon |x| < b^m \text{ for some } m \in \Z \}
	$$
	is an order convex subring of $\R_\omega$ with maximal ideal
	$$
	J_b := \{ x \in \R_\omega \colon |x| <  b^m \text{ for all } m \in \Z \}.
	$$
	The \emph{Robinson field} associated to the non-principal ultrafilter $\omega$ and the infinite element $b$ is the quotient $\R_{\omega,b} := O_b/J_b$ \cite{Rob96}. The Robinson field is a non-Archimedean real closed field. Note that $[b] \in \R_{\omega,b} $ is a \emph{big} element, meaning that for all $a \in \R_{\omega,b}$ there is an $n \in \N$ such that $a < b^n$. 
	
	Non-Archimedean ordered fields $\F$ with big elements admit an order compatible valuation $v \colon \R \to \R \cup \{\infty\}$ by letting $v(a)$ be the real number defined by the Dedekind cut
	\begin{align*}
		A_a &:= \left\{ \frac{p}{q} \in \Q \colon b^p  \leq a^q , q \in \Z_{>0}, p \in \Z\right\} \\
		B_a &:= \left\{ \frac{p}{q} \in \Q \colon b^p  \geq a^q , q \in \Z_{>0}, p \in \Z\right\} 
	\end{align*} 
	for $a \in \F$. Note that when $\F$ is Archimedean, every element $b>1$ is big and we can still define as above $v(a) = -\log_b|a|$, which then is the usual logarithm with base $b$. However $v$ is not a valuation in our sense, since it does not satisfy the strong triangle inequality, condition (3) in the definition.

	\section{Affine $\Lambda$-buildings}\label{sec:affine_L_buildings} 
	
	\subsection{$\Lambda$-metric spaces}
	
	An abelian group $(\Lambda, +)$ with a linear order such that $x,y \geq 0$ implies $x+y \geq 0$ is called an \emph{ordered abelian group}. If an ordered abelian group $\Lambda$ is isomorphic (as an ordered group) to a subgroup of $(\R,+)$, it is called \emph{Archimedean}. Hahn's embedding theorem classifies all ordered abelian groups.
	
	\begin{theorem}
		(\cite{Hah07}) For every ordered abelian group $(\Lambda,+)$ there is an ordered set $\Omega$ such that $\Lambda<\R^{\Omega}$ as an ordered subgroup, where 
		$$
		\R^\Omega = \{ f\colon \Omega \to \R \colon \operatorname{supp}(f) \text{ is contained in a well ordered set } \}
		$$
		is equipped with the lexicographical ordering.  
	\end{theorem} 
	Let $(\Lambda,+)$ be a non-trivial ordered abelian group,  In particular, $\Lambda$ has no torsion.
	
	We will use the following generalization of metric spaces.
	If $X$ is a set and $d \colon X \times X \to \Lambda$ a function, we call $(X,d)$ a \emph{$\Lambda$-pseudometric space} if for all $x,y,z \in X$
	\begin{enumerate}
		\item [(1)] $d(x,x) = 0$, $d(x,y) \geq 0$
		\item [(2)] $d(x,y) = d(y,x)$
		\item [(3)] $d(x,y) \leq d(x,z) + d(z,y)$. 
	\end{enumerate}
	If in addition, $d(x,y) = 0$ implies $x=y$, then $(X,d)$ is called a \emph{$\Lambda$-metric space}. The axioms are direct generalizations of the notions of (pseudo)metric spaces when $\Lambda = \R$. 
	
	\subsection{$\Lambda$-trees}\label{sec:lambda_trees}
	
	We define the notion of a $\Lambda$-tree following Chiswell's book \cite{Chi01}, where more details can be found. Let $(\Lambda,+)$ be an ordered abelian group and $(X,d)$ a $\Lambda$-metric space. A \emph{closed $\Lambda$-interval} is a set of the form
	$$
	[a,b]:=  \{\lambda \in \Lambda \colon a \leq \lambda \leq b\}
	$$
	for $a\leq b\in \Lambda$. Closed $\Lambda$-intervals with $d_\Lambda(t,t')=|t'-t|$ for $t,t'\in [a,b]$ are $\Lambda$-metric spaces. Let $X$ be any $\Lambda$-metric space. An isometric embedding $\varphi \colon [a,b] \to X$ is called a \emph{parametrization} of its image $s=\varphi([a,b])\subseteq X$, which is called a \emph{segment}. Note that the set of \emph{endpoints} $\{\varphi(a), \varphi(b)\}\subseteq X$ is independent of the parametrization $\varphi$ of the segment $s$. A $\Lambda$-metric space is \emph{geodesic} if for any two points $p,q\in X$ there exists a segment that has $p$ and $q$ as endpoints. If there is only one such segment, then $X$ is called \emph{uniquely geodesic} (or \emph{geodesically linear} in \cite{Chi01}). In a uniquely geodesic $\Lambda$-metric space, we denote a segment $s$ with endpoints $x,y \in X$ by $s=[x,y]$. 
	
\begin{figure}[t]
		\centering
		\includegraphics[scale=0.8]{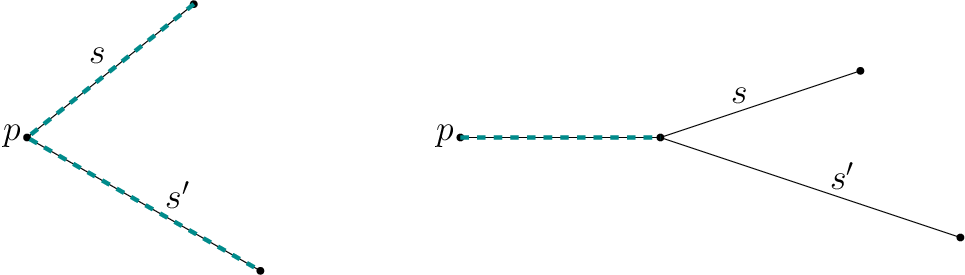}
		\caption{Illustration of axiom (2) (left) and axiom (3) (right) of $\Lambda$-trees. In axiom (2) the segments $s, s'$ are only allowed to intersect in one point.}
		\label{fig:axioms_2_3}
	\end{figure}

	A $\Lambda$-metric space $(X,d)$ is a \emph{$\Lambda$-tree} if it satisfies the following three axioms, two of which are illustrated in Figure \ref{fig:axioms_2_3}.
	\begin{itemize}
		\item[(1)] $(X,d)$ is geodesic.
		\item[(2)] If two segments $s,s'\subseteq X$ intersect in a single point $s\cap s' = \{p\}$, which is an endpoint of $s$ and $s'$, then their union $s\cup s'$ is a segment.
		\item[(3)] If two segments $s,s'\subseteq X$ have an endpoint $p$ in common, then their intersection $s\cap s'$ is a segment. 
	\end{itemize}
	
	The notion of a $\Lambda$-tree generalizes trees from graph theory as $\Z$-trees and real trees as $\R$-trees.	It is known that for $\Lambda = \Z$ and $\Lambda=\R$, axiom (3) follows from (1) and (2), see 
	\cite[Lemma I.2.3, Lemma I.3.6]{Chi01}. In \cite{App24}, we gave the following complete characterization of the groups $\Lambda$ for which axiom (3) follows from (1) and (2), answering for which $\Lambda$ condition (3) has to be included.
	
	\begin{theorem}\label{thm:main_Lambda_trees}
		Let $\Lambda$ be an ordered abelian group. The following are equivalent:
		\begin{itemize}
			\item [(a)] For every positive $\lambda_0 \in \Lambda$, the set
			$
			\{t \in \Lambda \colon 0 \leq 2t\leq \lambda_0 \}
			$
			has a maximum.
			\item [(b)] Every $\Lambda$-metric space that satisfies axioms (1) and (2) also satisfies (3).
		\end{itemize}
	\end{theorem} 

    The ordered abelian group $\Lambda = \Z[1/3] < \R$ is an example of a group for which axiom (3) is independent of (1) and (2), so axiom (3) can not be omitted from the definition for general $\Lambda$. In \cite{App24} we leave open two questions.
    
    \begin{question}\label{q:2_independent}
    	Let $\Lambda$ be a non-trivial ordered abelian group. Is axiom (2) independent of the axioms (1) and (3)?
    \end{question}
    
    \begin{question}\label{q:unique_independent}
    	Let $\Lambda$ be a non-trivial ordered abelian group. Is there a uniquely geodesic $\Lambda$-metric space that does not satisfy axiom (2)?
    \end{question}

	In applications, the algebraic condition $\Lambda = 2 \Lambda$ is often satisfied. Elsewhere in the thesis, $\Lambda$ will be the image of a valuation $v \colon \F_{>0} \to \Lambda$ of a real closed field $\F$. Since square roots exist for positive elements in $\F$, $\Lambda$ is divisible by 2.

		\subsubsection{Example, Brumfiel's tree}
	
	Let $\F$ be a non-Archimedean real closed field with valuation $v \colon \F \to \Lambda \cup \{\infty\}$ where $\Lambda \subseteq \R$. We follow \cite{Bru2} to construct a $\Lambda$-tree associated to the hyperbolic plane and $\F$. We consider the upper half-plane model $\bbH_\R := \{ (x,y) \in \R^2 \colon y >0 \}$ of the hyperbolic plane and define the \emph{non-standard hyperbolic plane} as its semialgebraic extension $\bbH_\F$. The cross ratio can be used to define the distance on $\bbH_\R$: for $p,q\in \bbH_\R$ there is a unique hyperbolic line through $p$ and $q$ given by a half circle or a ray perpendicular to the $x$-axis. Let $a,b \in \R\cup \{\infty\}$ be the two points at infinity of the line, $a$ closer to $p$ and $b$ closer to $q$. The cross-ratio
	$$
	\operatorname{cr}_\R(p,q) := \operatorname{CR}(p,q;a,b) := \frac{|a-q|\cdot |b-p|}{|a-p| \cdot |b-q|} \in \R_{\geq 1}
	$$ 
	can then be used to describe the hyperbolic distance $d_\R(p,q) = \log \operatorname{cr}_\R(p,q)$. We note that $\operatorname{cr}_\R$ is a semialgebraic function, and $d_\F(p,q) := (-v) \operatorname{cr}_\F(p,q)$ defines a $\Lambda$-pseudo-distance on $\bbH_\F$. Then $\mathcal{T}:= \bbH_\F/ \!\! \sim$, where $p\sim q$ if $d_\F(p,q)=0$ defines a $\Lambda$-metric space. 
	\begin{theorem}\label{thm:Brumfiel}
		(\cite[Theorem (28)]{Bru2})
		The $\Lambda$-metric space $\mathcal{T}$ is a $\Lambda$-tree.
	\end{theorem}
	
	This construction is an easier special case of the construction of the building $\B$ in Section \ref{sec:building_def} and follows along the same ideas: we first find a semialgebraic model for the hyperbolic plane, define the non-standard hyperbolic plane and a pseudometric on it and finally show that the quotient defines a $\Lambda$-tree. Theorem \ref{thm:Brumfiel} is a special case of our main theorem, Theorem \ref{thm:B_is_building} for $\operatorname{SL}(2)$. Our theorem is more general in two ways: we allow for higher rank symmetric spaces and thus obtain affine $\Lambda$-buildings, and we allow $\Lambda$ to be any ordered abelian group, not necessarily a subgroup of $\R$.
	
	In joint work with De Rosa, Flamm and Jaeck \cite{ADFJ24} we investigated the $\Q$-tree $\mathcal{T}$ for $\F$ the Puiseux series. Our main result states that completing all the segments of $\mathcal{T}$ does not result in a complete metric space.  
	
	\begin{theorem}\label{thm:ARFJ}
		Let $\mathcal{T}$ be the $\Q$-tree constructed from the field of Puiseux-series in \cite{Bru2}. There is a Cauchy sequence in the segment completion $\overline{\mathcal{T}}^{\operatorname{sc}}$ that does not converge and hence $\overline{\mathcal{T}}^{\operatorname{sc}}$ is not complete.
	\end{theorem}

	\subsection{Root systems}\label{sec:root_system}
	
	A detailed treatment of root systems can be found in \cite{Bou08}. Let $V$ be a finite dimensional Euclidean vector space with scalar product $\langle \cdot, \cdot \rangle$. For $\alpha \in V$, the reflection along the hyperplane $M_\alpha = \{ \beta \in V \colon \langle \alpha , \beta \rangle = 0 \}$ is given by
	$$
	r_{\alpha} (\beta) = \beta - 2 \frac{\langle \alpha , \beta \rangle}{ \langle \alpha, \alpha \rangle} \alpha . 
	$$
	A pair $(\Phi,V)$ where $\Phi \subseteq V$ is called a \emph{root system} if 
	\begin{enumerate}
		\item [(R1)] $\Phi$ is finite, symmetric ($\Phi = -\Phi$), spans $V$ and does not contain $0$.
		\item [(R2)] For every $\alpha \in \Phi$ the reflection $r_\alpha \colon V \to V$ preserves $\Phi$.
	\end{enumerate} 
	When $V$ can be determined from the context, the root system may be denoted by just $\Phi$. A root system is called \emph{crystallographic} if it satisfies the integrality condition
	\begin{enumerate}
		\item [(R3)] If $\alpha, \beta \in \Phi$, then $2\langle \alpha, \beta \rangle / \langle \alpha, \alpha \rangle \in \Z$.
	\end{enumerate} 
	For crystallographic root systems, $\operatorname{Span}_{\mathbb{Z}}(\Phi)$ is a lattice in $V$. A root system is called called \emph{reduced} if
	\begin{enumerate}
		\item [(R4)] For every $\alpha \in \Phi$, $\R \cdot \alpha \cap \Phi = \{ \pm \alpha\}$.
	\end{enumerate}
	
	\begin{figure}[h]
		\centering
		\includegraphics[width=0.8\linewidth]{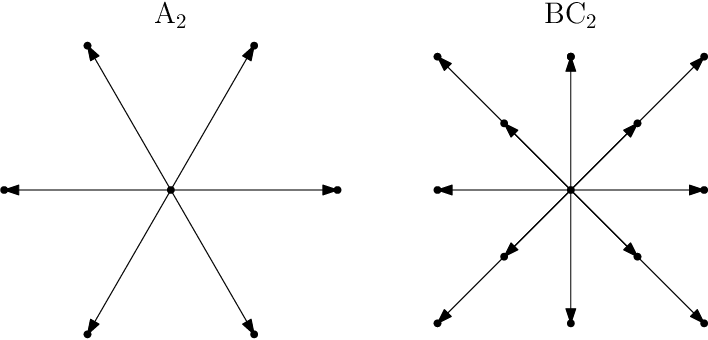}
		\caption{ Two examples of crystallographic root systems. Type $\operatorname{A}_2$ on the left, Type $\operatorname{BC}_2$ on the right. While $\operatorname{A}_2$ is reduced, $\operatorname{BC}_2$ is not. }
		\label{fig:examples_root_systems}
	\end{figure}
	
	Figure \ref{fig:examples_root_systems} gives examples of two crystallographic root systems, one of them reduced the other not.	The \emph{spherical Weyl group $W_s$} of a root system is the subgroup of isometries of $V$ generated by reflections along the hyperplanes $M_\alpha$ for $\alpha \in \Phi$. A subset $\Delta \subseteq \Phi$ is a \emph{basis of $\Phi$}, if it is a vector space basis of $V$ such that all roots $\beta \in \Phi$ can be written as $\beta = \sum_{\delta \in \Delta} \lambda_\delta \delta$ with $\lambda_\delta \in \Z$ for $\delta \in \Delta$ and such that all $\lambda_{\delta}$ have the same sign (all are non-negative or non-positive). Every root system has a basis and the spherical Weyl group acts transitively on the set of bases of $\Phi$. The cardinality of a basis $\Delta$ is called the \emph{rank} of the root system $\Phi$ and coincides with $\dim(V)$.
	 
	The connected components of $V\setminus \cup_{\alpha \in \Phi} M_\alpha$ are called \emph{(open) chambers}. The spherical Weyl group acts simply transitively on the set of chambers. Given a basis $\Delta$, the \emph{fundamental Weyl chamber of $\Delta$} is
	$$
	WC(\Delta) = \left\{ v\in V \colon \langle v, \alpha \rangle >0 \ \forall \alpha \in \Delta \right\}.
	$$
	Conversely, every chamber $C$ determines a basis $B(C)$ formed by those $\alpha \in \Phi$, that satisfy $\langle \alpha, v \rangle >0$ for all $v \in C$ and that cannot be written as a sum of other such elements of $\Phi$, see \cite[VI\S1.5]{Bou08}. A basis $\Delta$ determines the set  of \emph{positive roots} $\Sigma_{>0} \subseteq \Phi$ given by positive integer combinations of elements in the basis. Thus a basis determines a partial order on $\Sigma$ by $\alpha < \beta$ if $\beta-\alpha \in \Sigma_{>0}$.  A total order on $\Sigma$ with the same positive elements can be obtained by choosing an order on the basis and extending it lexicographically to $ \Sigma \subseteq \operatorname{Span}_{\mathbb{Z}}(\Sigma)$.

	Let now $(\Phi,V)$ be a crystallographic root system. For every $\alpha \in \Phi$ the \emph{coroots}\footnote{The coroots $\alpha^\vee$ are not to be confused with the \emph{dual roots} in $V^\star$, which are sometimes denoted by the same symbol.}  
	$$
	\alpha^{\vee} := \frac{2\alpha}{\langle \alpha, \alpha \rangle} \in V
	$$
	have the property that for all $\alpha, \beta \in \Sigma$
	$$
	\langle \alpha , \beta^\vee \rangle = \frac{2 \langle \alpha, \beta\rangle}{\langle \beta, \beta\rangle} \in \Z 
	$$
	by the crystallographic condition (R3).

	\subsection{Apartments}\label{sec:modelapartment}
	
	We introduce the model apartment $\A$. A more detailed introduction can be found in \cite{Ben1}, and including the non-crystallographic case in \cite{Sch09}. The root systems in our setting will arise from algebraic groups, and are therefore crystallographic, but possibly not reduced.
	
	Let $\Lambda$ be a non-trivial ordered abelian group. 
	 Let $(\Phi,V)$ be a crystallographic root system and $L = \operatorname{Span}_{\Z}(\Phi) \subseteq V$ its associated lattice. 
	Since both $L$ and $\Lambda$ are free $\Z$-modules, we can define the \emph{model apartment} 
	$$
	\A := L \otimes_{\Z} \Lambda .
	$$
	For a basis $\Delta = \{\alpha_1 , \ldots , \alpha_r \}$, the $\mathbb{Z}$-bilinear map
	\begin{align*}
		L \times \Lambda & \to \Lambda^{r} \\
		\left(\sum_{i=1}^r b_i \alpha_i , \lambda \right)& \mapsto (b_1 \lambda, \ldots , b_r \lambda)
	\end{align*}
    extends to a $\mathbb{Z}$-module isomorphism
	$
	\A 
	\cong \Lambda^{r}
	$, giving rise to the notation
	$$
	\A \cong \left\{ \sum_{\delta \in \Delta} \lambda_\delta \delta \colon \lambda_\delta \in \Lambda \right\}.
	$$
	Note that if $\Lambda$ is a $\Q$-vector space, then we also have
	$$
	\A = \operatorname{Span}_\Q(\Phi) \otimes_{\Q} \Lambda.
	$$
	 Since $W_s$ acts on $\Phi$, it also acts on $\A$, fixing $0 \in \A$. The apartment $\A$ itself acts by translation on $\A$. Let $T\subseteq \A$ be a subgroup of the translation group normalized by $W_s$. Then $W_a = W \ltimes T$ is called the \emph{affine Weyl group} and acts on $\A$. If $T=\A$, $W_a$ is called the \emph{full affine Weyl group}. We denote the data of a model apartment together with an affine Weyl group as $\A = \A(\Phi,\Lambda,T)$ since the root system $\Phi$, the ordered abelian group $\Lambda$ and the subgroup of translations $T$ fully determine $\A$ and the affine Weyl group.

	 For any root $\alpha \in \Phi$, the reflection $r_\alpha \colon \Phi \to \Phi$ extends to a \emph{reflection}
	 \begin{align*}
	 	r_{\alpha} \colon \A & \to \A \\
	 	\sum_{\delta \in \Delta} \lambda_\delta \delta & \mapsto \sum_{\delta \in \Delta} \lambda_\delta r_\alpha (\delta) 
	 \end{align*}
	 and for any $t\in T$, maps of the form $t \circ r_\alpha \in W_a$ are called \emph{affine reflections}. 
	 The scalar product on $\Phi$ extends naturally to the bilinear pairing
	 \begin{align*}
	 \langle \cdot , \cdot \rangle \colon 	\A \times L  & \to \Lambda \\
	 	\left( \sum_{\delta \in \Delta} \lambda_\delta \delta , \sum_{\delta' \in \Delta } \mu_{\delta'} \delta' \right) & \mapsto \sum_{\delta \in \Delta} \sum_{ \delta' \in \Delta} \mu_{\delta'} \lambda_{\delta}  \langle \delta, \delta' \rangle
	 \end{align*}
	 but we notice that the scalar product can not be defined on all of $\A \times \A$, since $\Lambda$ may not have a multiplication. For $\alpha \in \Phi$ and $\lambda \in \Lambda$ we define the \emph{affine wall}
	 $$
	 M_{\alpha,\lambda} :=
	 \left\{x \in \A \colon   \langle x, \alpha \rangle = \lambda  \right\}
	 $$
	 as well as the two \emph{affine half-spaces} 
	 \begin{align*}
	 	H_{\alpha, \lambda}^+ &= \left\{ x\in \A \colon    \langle x, \alpha \rangle \geq \lambda  \right\} \\
	 	H_{\alpha, \lambda}^- &= \left\{ x \in \A \colon  \langle x, \alpha \rangle \leq \lambda  \right\}
	 \end{align*}
	 defined by $M_{\alpha,\lambda}$. The intersection of the positive half-spaces for all $\delta \in \Delta$ is called the \emph{fundamental Weyl chamber}
	 $$
	 C_0 := \left\{ x \in \A \colon \langle x, \alpha \rangle \geq 0 \ \forall \delta \in \Delta \right\} = \bigcap_{\delta \in \Delta} H_{\delta,0}^+.
	 $$
	 The images of the fundamental Weyl chamber under the action of the affine Weyl group are called \emph{chambers}.

\begin{figure}[h]
		\centering
		\includegraphics[width=0.6\linewidth]{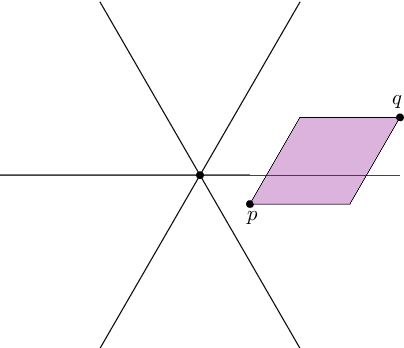}
		\caption{ Given two points $p,q \in \A$, the interval $\{ r\in \A \colon d(p,r)+d(r,q) = d(p,q)\}$ is a $W_a$-convex set defined by the intersection of four half-planes parallel to the walls in this example of type $\operatorname{A}_2$. }
		\label{fig:A2_segment}
	\end{figure}

Note that for every $\alpha \in \Phi$,
$$
\alpha^{\vee} := \frac{2\alpha}{\langle \alpha, \alpha \rangle} = \sum_{\delta \in \Delta} \langle  \alpha^\vee, \delta \rangle \delta 
= \sum_{\delta \in \Delta} 2\frac{\langle \alpha, \delta \rangle}{\langle \alpha, \alpha \rangle} \delta 
\in   L
$$
by the crystallographic condition.
     The model apartment $\A$ can be endowed with a $\Lambda$-metric. There are multiple ways to do this and we will describe one. Since $\A$ may not admit a $\Lambda$-valued scalar product, we instead use the $\Lambda$-valued norm $N \colon \A \to \Lambda$ defined by
     $$
     N( x ) = \sum_{\alpha \in \Phi_{>0}} \left| \left\langle x, \alpha^\vee \right\rangle \right|,
     $$ 
     where $|\lambda| = \max\{\lambda, -\lambda\}$ is the absolute value. Then
     $$
     d(x,y) = N(x-y)
     $$
     turns $\A$ into a $\Lambda$-metric space. Note that when $x-y \in C_0$, 
     $$
     d(x,y) =  \left\langle x-y , \sum_{\alpha \in \Phi_{>0} } \alpha^\vee \right\rangle.
     $$
     Note that the $\Lambda$-metric space $\A$ may not be uniquely geodesic: there may be two or more distinct isometric embeddings of a $\Lambda$-interval with coinciding endpoints. Since $\Lambda$ is just a group, the classical notion of convexity using linear combinations does not make sense. Instead we say that a subset $B \subseteq \A$ is \emph{$W_a$-convex}, if it is a finite intersection of affine half-spaces. Figure \ref{fig:A2_segment} illustrates that $\A$ may not be uniquely geodesic and that the union of all geodesics is a $W_a$-convex set.

	\subsection{Affine $\Lambda$-buildings}\label{sec:lambda_building}
	 
	Let $\A$ be an apartment of type $\A=\A(\Phi,\Lambda,T)$ with affine Weyl group $W_a$ as in the previous section. Let $X$ be a set and $\Fun$ a set of injective maps $\A \to X$. The set $\Fun$ is called the \emph{atlas} or \emph{apartment system} and its elements are called the \emph{charts}. The images $f(\A) $ for $f\in\Fun$ are called \emph{apartments}. The images of walls, half-spaces and chambers are called \emph{walls, half-apartments} and \emph{sectors}. 
	
\begin{figure}[h]
 		\centering
 		\includegraphics[width=0.8\linewidth]{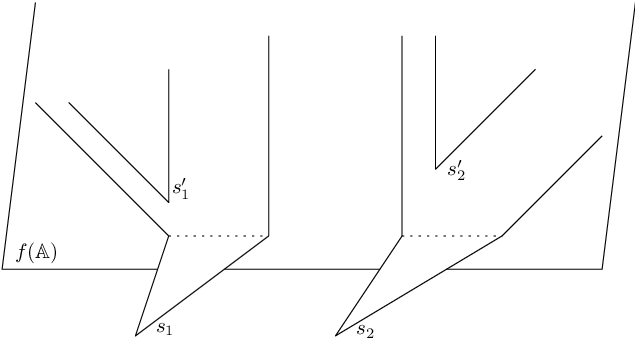}
 		\caption{Axiom (A4) states that while arbitrary sectors $s_1,s_2$ may not lie in a common flat, they contain subsectors $s_1',s_2'$ contained in a common flat $f(\A)$. }
 		\label{fig:A4}
 	\end{figure}

\begin{figure}[h]
 		\centering
 		\includegraphics[width=0.8\linewidth]{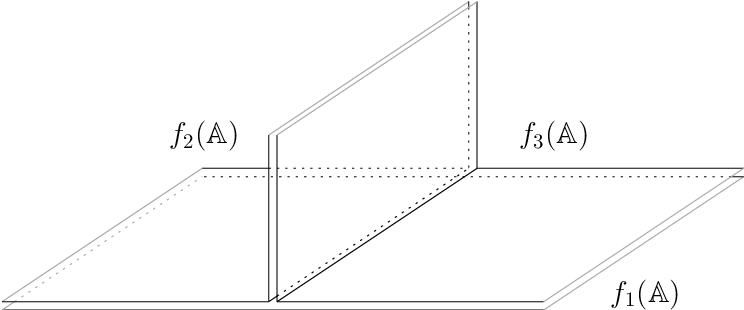}
 		\caption{Axiom (A6) states that whenever three apartments intersect pairwise in a half-apartment, then there is at least one point in the common intersection of all three.}
 		\label{fig:A6}
 	\end{figure}

	The pair $(X,\Fun)$ is called a \emph{generalized affine building} or a \emph{affine $\Lambda$-building} of type $\A = \A(\Phi,\Lambda,T)$ if the following six axioms are satisfied. Axioms (A4) and (A6) are illustrated in Figures \ref{fig:A4} and \ref{fig:A6}.
	
	\begin{enumerate}
		\item [(A1)] For all $f\in \Fun$, $w\in W_a$, $f \circ w \in \Fun$.
		\item [(A2)] For all $f,f'\in \Fun$, the set $B= f^{-1}\left( f(\A)\cap f'(\A)  \right)$ is $W_a$-convex and there is a $w\in W_a$ such that
		$
		f|_B = f'\circ w |_B.
		$
		\item [(A3)] For all $x,y\in X$, there is a $f\in \Fun$ such that $x,y\in f(\A)$.
		\item [(A4)] For any sectors $s_1,s_2 \subseteq X$ there are subsectors $s_1' \subseteq s_1, s_2' \subseteq s_2$ such that there is an $f\in \Fun$ with $s_1', s_2' \subseteq f(\A)$.
		\item [(A5)] For every $x\in X$ and $f\in \Fun$ with $x\in f(\A)$ there is a distance diminishing retraction $r_{x,f} \colon X \to f(\A)$ such that $f^{-1}(\{x\}) = \{x\}$.
		\item [(A6)] For $f_1,f_2,f_3 \in \Fun$, if $f_i(\A) \cap f_j(\A)$ are half-apartments for $i\neq j$, then $f_1(\A) \cap f_2(\A) \cap f_3(\A) \neq \emptyset$.
 	\end{enumerate}
 	
 	The distance in axiom (A5) is the function $d\colon X\times X \to \Lambda$ induced from the distance on the model apartment $\A$, whose existence follows from axioms (A1), (A2) and (A3): for any two points $x,y\in X$ we use axiom (A3) to find $f\in \Fun$ such that $x,y \in f(\A)$ and define $d(x,y) = d( f^{-1}(x) , f^{-1}(y) ) \in \Lambda $. This is well defined, by axioms (A1) and (A2). A $d$-diminishing retraction is then a function $r \colon X \to f(\A)$ satisfying $r(y)=y$ for all $y\in f(\A)$ and $$
 	d(r(x),r(y)) \leq d(x,y)
 	$$
 	for all $x,y\in X$. Note that while axioms (A1) - (A3) can be used to define a symmetric positive definite function $d$, axiom (A5) can be used to show that $d$ satisfies the triangle inequality. As shown in \cite{BeScSt}, in the presence of the other five axioms, axiom (A5) is equivalent to the triangle inequality. In fact, \cite{BeScSt} contains a number of collections of axioms that characterize affine $\Lambda$-buildings. We will use the following characterization.
 	
 	\begin{theorem}[Theorem 3.1, (5) $\implies$ (1) \cite{BeScSt}] \label{thm:equivalent_axioms}
 		Let $(X,\F)$ be a set with an atlas such that the axioms (A1), (A2), (A3) and (A4), as well as
 		\begin{enumerate}
 			\item [(TI)] The function $d$ induced from the distance in apartments satisfies the triangle inequality.
 		\end{enumerate}
 		and the exchange condition
 		\begin{enumerate}
 			\item [(EC)] For $f_1,f_2 \in \Fun$, if $f_1(\A)\cap f_2(\A)$ is a half-apartment, then there exists $f_3 \in \Fun$ such that $f_i(\A)\cap f_3(\A)$ are half-apartments for $i\in \{1,2\}$. Moreover $f_3(\A)$ is the symmetric difference of $f_1(\A)$ and $f_2(\A)$ together with the boundary wall of $f_1(\A) \cap f_2(\A)$.
 		\end{enumerate}
 		are satisfied. Then $(X,\F)$ is an affine $\Lambda$-building. 
 	\end{theorem}

	 \subsection{Independence of axiom (A6)}\label{sec:A6_independence}
	 
\begin{figure}[h]
	 	\centering
	 	\includegraphics[width=0.5\linewidth]{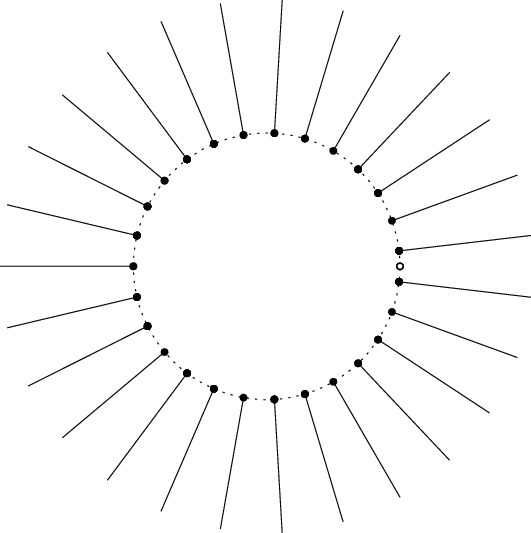}
	 	\caption{ For $\Lambda = \mathbb{Z}[1/3]$, there is a $\Lambda$-metric space $X$ that satisfies axioms (A1) - (A5), but not (A6). It is important that $1/2 \notin \Lambda$, so that there is no point opposite of $[0] \in S_\Lambda^1$, illustrated by a circle with white interior. }
	 	\label{fig:A6_counterexample}
	 \end{figure}

	 The independence of some of the axioms of affine $\Lambda$-buildings were studied in the appendix of \cite{BeScSt}. For $\Lambda = \mathbb{R}$, Parreau \cite{Par} showed that axiom (A6) follows from (A1) - (A5). When Bennett introduced affine $\Lambda$-buildings he claimed to give an example of a rank 1 $\Q$-building that satisfies axioms (A1) - (A5), but not (A6) \cite[Remark 3.3]{Ben1}. The example looked like it had the claimed properties, but did not actually admit a valid apartment system. Translated in terms of $\Lambda$-trees, the example would have been a $\Q$-tree where axioms (1) and (2) would have held, but not axiom (3). In Theorem \ref{thm:main_Lambda_trees} we showed that this is not possible. The next example shows that for $\Lambda = \Z[1/3]$, there is a $\Lambda$-metric space $X$ that admits an apartment system such that axioms (A1) - (A5) are satisfied, but (A6) is not. The example can be generalized to any ordered abelian group $\Lambda$ with $\Lambda \neq 2\Lambda$ and is illustrated in Figure \ref{fig:A6_counterexample}.
	 
	 Let $\Lambda = \Z[1/3]$. We start with a $\Lambda$-circle $S_\Lambda^1 := \Lambda/\Z$ of length 1. At every point $p \in S_\Lambda^1$ we glue a $\Lambda$-ray $r_p = [0,\infty)_\Lambda \subseteq \Lambda$ to obtain a $\Lambda$-metric space
	 $$
	 X = \left( S_\Lambda^1 \amalg \coprod_{p \in S_\Lambda^1} r_p \right) /  \! \sim
	 $$
	 where $p\in S_\Lambda^1$ is identified with $0\in r_p$. For any two points $p,q \in S_{\Lambda}^1$, a unique apartment is given by the union of $r_p$, $r_q$ and the shorter segment in $S_\Lambda^1$ connecting $p$ to $q$. Since $1/2 \notin \Lambda $, $p$ and $q$ can not lie opposite each other on $S_\Lambda^1$. Axioms (A1), (A2), (A3) and (A4) are readily checked, (A5) requires some thought and (A6) clearly does not hold. The uniqueness of the apartment is important for the validity of (A2) and is the reason why this example does not work for $\Lambda = \Q$.

	 In analogy to Theorem \ref{thm:main_Lambda_trees}, we would like to ask for which $\Lambda$ axiom (A6) is independent. Note that even for rank 1, independence of (3) does not directly translate to independence of (A6).
	 
	 \begin{question}\label{q:A6}
	 	Characterize the ordered abelian groups for which axiom (A6) is independent of the other axioms. Is there even any ordered abelian group $\Lambda$ for which axiom (A6) is independent?
	 \end{question}

\newpage

	 \section{Linear algebraic groups}\label{sec:algebraic_groups}
	 \subsection{Definitions}
	 
	 \begin{wrapfigure}{r}{0.1\textwidth}
	 	\centering
	 	\includegraphics[width=0.1\textwidth]{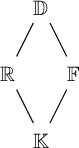}
	 \end{wrapfigure}
	 
	 Let $\K \subseteq \R$ be a subfield (usually $\K = \Q$ or $\K = \R$), $\F$ a real closed field containing $\K$ and $\D$ an algebraically closed field that contains both $\R$ and $\F$. We follow a naive approach to algebraic groups, viewing them as matrix groups. Since we are working with fields of characteristic 0, the algebraic geometry can be kept at a minimum. For an extensive introduction to algebraic groups, we refer to \cite{Bor66}, \cite{Bor}, \cite{Hum1} or \cite{Zim}. 
	 
	 The general linear group $\GLnD$ can be realized as an affine algebraic variety
	 $$
	 \GLnD \cong \left\{ \begin{pmatrix}
	 	A & \\ & t
	 \end{pmatrix} \in \D^{(n+1)\times(n+1)} \colon \det (A) \cdot t = 1 \right\}.
	 $$
	 A subgroup $G < \GLnD$ is called a \emph{linear algebraic group defined over $\K$}, if it is a set of common zeros for a set of polynomials in the coordinate ring $\K[\operatorname{GL}_n(\D)]:= \K[(x_{ij}), \det(x_{ij})^{-1}]$ of $\operatorname{GL}_n(\D)$. 
	 We will not consider more general algebraic groups and hence also call $G$ an \emph{algebraic group} or a $\K$-\emph{group}.
	 For any commutative $\K$-subalgebra $\mathbb{B} \subseteq \D$ containing $\K$, we let 
	 $$
	 \GLnB = \left\{ (a_{ij}) \in \GLnD \colon a_{ij} \in \mathbb{B} \text{ and } \det(a_{ij})^{-1} \in \mathbb{B} \right\}.
	 $$ 
	 The \emph{group of $\mathbb{B}$-points of a linear algebraic group $G$}
	 is $G_\mathbb{B} := G \cap \GLnB$. The $\R$-points of a linear algebraic group $G$ form a real Lie group \cite[Theorem 2.1]{Mil13}. 
	 
	 Viewing $G \subseteq \GLnD $ as an algebraic subset of $\D^{n+1\times n+1}$, we endow $G$ with the Zariski-topology, whose closed sets are given by common zeros of sets of polynomials in $\K[\operatorname{GL}_n(\D)]$. A connected linear algebraic group $G$ is \emph{semisimple} if every closed connected normal abelian subgroup is trivial.

 \subsubsection{Examples}
 
 The multiplicative group $\G_m = \operatorname{GL}_1(\D)$ is a linear algebraic group defined over $\Q$. Note that $\G_m$ is connected in the Zariski-topology. Its $\C$-points $(\G_m)_\C = \C \setminus \{0\}$ are connected in the Euclidean topology but its $\R$-points $(\G_m)_\R = \R \setminus \{0\}$ are not.
 
   If $V$ is a $\D$-vector space, then the group $\operatorname{GL}(V)$ of automorphisms of $V$ is an algebraic group. As $\operatorname{GL}(V)$ contains the closed connected normal abelian subgroup of scalar multiplication, $\operatorname{GL}(V)$ is not semisimple.  
  
  The real Lie groups $\SLnR,\operatorname{SO}_n(\R)$ and $\operatorname{Sp}_{2n}(\R)$ are groups of $\R$-points of the linear algebraic groups $\SL_n(\D), \operatorname{SO}_n(\D)$ and $\operatorname{Sp}_{2n}(\D)$ defined over $\Q$.

 \subsection{Morphisms and tori }

For a linear algebraic group $G$ we consider the ideal 
$$
I(G) = \left\{p \in \D\left[\left(x_{ij}\right),z\right] \colon  p\left(g,\det\left(g\right)^{-1}\right) = 0 \text{ for all } g \in G   \right\},
$$
which is finitely generated as a consequence of Hilbert's Basis theorem. The \emph{coordinate ring} of $G$ is 
$$
\D [G] = \D\left[\left(x_{ij}\right),d\right] / I(G)
$$
and its elements are called \emph{regular functions} on $G$. Given a map $\varphi\colon G \to H$ between algebraic groups $G$ and $H$, we can define the \emph{transposed map} $\varphi^\circ \colon \D[H] \to \D[G]$ by $\varphi^\circ (f) = f \circ \varphi$. A \emph{morphism} of linear algebraic groups $G$ and $H$ is a group homomorphism $\varphi \colon G \to H$ whose transposed map $\varphi^\circ$ is a ring homomorphism. If $G$ and $H$ are defined over $\K$ and $\varphi^\circ$ maps $\K[H]$ into $\K[G]$, then $\varphi$ is \emph{defined over} $\K$. An important example of a morphism $G \to G$ is the conjugation $\operatorname{Int}(g) \colon h \mapsto ghg^{-1}$ by an element $g \in G$, and if $g \in G_\K$, then $\operatorname{Int}(g)$ is defined over $\K$.

\begin{lemma}\label{lem:basic_algebraic_geometry}
	Let $\varphi \colon G \to H$ be a morphism of algebraic $\K$-groups $G < \operatorname{GL}(n,\D), H < \operatorname{GL}(m,\D)$ defined over $\K$. If we write componentwise
	\begin{align*}
		\varphi_\F \colon G_\F & \to H_\F \\
		 \begin{pmatrix}
			x_{11} & \ldots & \\ \vdots & \ddots & \\ & & x_{nn}
		\end{pmatrix} & \mapsto \begin{pmatrix}
		\varphi_{11}(x_{11}, \ldots , x_{nn}) & \ldots & \\ \vdots & \ddots & \\ & & \varphi_{mm}(x_{11}, \ldots, x_{nn})
		\end{pmatrix} 
	\end{align*}
	then all $\varphi_{ij} \colon G_\F \to \F $ are polynomials.
\end{lemma}
\begin{proof}
	Consider the polynomial representing $X_{ij} \in \K[H]$. Then 
	\begin{align*}
		\varphi^\circ (X_{ij})(x_{11}, \ldots , x_{nn}) &= X_{ij} (\varphi(x_{11}, \ldots, x_{mm})) \\
		&= X_{ij}(\varphi_{11}(x_{11}, \ldots , x_{nn}), \ldots , \varphi_{mm}(x_{11}, \ldots , x_{nn}))\\
		& = \varphi_{ij}(x_{11}, \ldots , x_{nn})
	\end{align*} 
	and thus $\varphi_{ij} \in \K[G]$ is a polynomial.
\end{proof}

\begin{proposition}\label{prop:morphisms_polynomials} 
	In characteristic $0$, a map $G\to H$ between algebraic groups $G,H$ is a morphism if and only if it is a group homomorphism whose graph is a set of common zeros of polynomials. The morphism is defined over $\K$ if and only if all the polynomials can be chosen to have coefficients in $\K$.
\end{proposition} 

 A linear algebraic group $T$ that is isomorphic to $(\G_m)^d$ is called a \emph{$d$-dimensional torus}. An element $g \in \operatorname{GL}_n(\D)$ is \emph{semisimple} if $\D^n$ is spanned by eigenvectors of $g$, or equivalently, if $g$ is diagonalizable. 
 
 \begin{theorem}[8.5 in \cite{Bor}] 
 	\label{thm:bor_torus}
 	For a Zariski-connected algebraic group $T$, the following conditions are equivalent.
 	\begin{enumerate}
 		\item $T$ is a torus.
 		\item $T$ consists only of semisimple elements.
 		\item The whole group $T$ is simultaneously diagonalizable.
 	\end{enumerate}
 \end{theorem}

Since Zariski-connectedness of algebraic groups and semisimplicity of elements is preserved under morphisms \cite[4.4(4)]{Bor}, the image of a torus under a morphism is a torus. 
 
 A morphism $\chi \colon G \mapsto \G_m$ from an algebraic group $G$ to $\G_m$ is called a \emph{character}. The set of characters of $G$ is an abelian group, denoted by $\hat{G}$. If $T \cong (\G_m)^d$ is a $d$-dimensional torus and an element $x \in T$ corresponds to $(x_1, \ldots , x_d)  \in (\G_m)^d$, then every character $\chi$ of $T$ is of the form 
 $$
 \chi(x) = x_1^{n_1} \cdots x_d^{n_d}
 $$
 for some $n_i \in \Z$ and hence $\hat{T} = \Z^d$. If there is an isomorphism $T \to (\G_m)^d$ defined over $\K$, $T$ is said to be \emph{$\K$-split}.
  
   \begin{theorem}[8.4 in \cite{Bor}.]
  	For a torus $T$ defined over $\K$, the following conditions are equivalent.
  	\begin{enumerate}
 	\item $T$ is $\K$-split.
 	\item All characters of $T$ are defined over $\K$, 
  	\end{enumerate}
  \end{theorem}

A \emph{maximal torus} $T \subseteq G$ is a subgroup that is a torus and is not properly contained in any other torus in $G$. 

\begin{theorem}[11.3 and 20.9 in \cite{Bor} ]
	In a connected algebraic group $G$, all maximal tori of $G$ are conjugate. In a connected semisimple 
	$\K$-group $G$, the maximal $\K$-split tori of $G$ are conjugate under $G_{\K}$.
\end{theorem}
The \emph{rank} of $G$ is the common dimension of the maximal tori and if $G$ is semisimple 
and defined over $\K$, the $\K$-\emph{rank} of $G$ is the common dimension of the maximal $\K$-split tori of $G$.

A morphism $\G_m \to G$ is a \emph{multiplicative one-parameter subgroup of $G$}. For a torus $T$ the one-parameter subgroups of $T$ are denoted by $X_{\star}(T)$. For $\chi \in \hat{T}$ and $\lambda \in X_{\star}(T) $, the composition $(\chi\circ\lambda) \colon \G_m \to \G_m$ is a character of the torus $\G_m$ and hence sends $x \in \G_m$ to $x^m$ for some $m \in \Z$. This defines a map $b \colon \hat{T} \times X_{\star}(T) \to \Z, (\chi,\lambda) \mapsto m$. 
\begin{proposition}{\cite[Proposition 8.7]{Bor}}\label{prop:char_cochar}
	For any torus $T$, the map 
	$$
	b\colon \hat{T} \times X_{\star}(T) \to \Z
	$$
	is a nondegenerate bilinear form.
\end{proposition}

 \subsubsection{Examples}\label{sec:ex_split}
 The $\R$-group $\operatorname{SO}(2) = \{ g \in \operatorname{SL}_2(\C) \colon gg\tran = \Id \}$ is a torus. An isomorphism with $\G_m$ is given by
 \begin{align*}
 	\G_m & \to \operatorname{SO}(2) \\
   a & \mapsto \frac{1}{2}\begin{pmatrix}
 	a+a^{-1} & i (a^{-1}-a) \\
 	i(a - a^{-1}) & a+a^{-1}
 \end{pmatrix}  \\
	 a+ib & \mapsfrom \phantom{\frac{1}{2}}  \begin{pmatrix}
	a & b \\
	c & d
\end{pmatrix}
 \end{align*}
 and the appearance of $i \in \C$ indicates that $\operatorname{SO}(2)$ is not $\R$-split. On the other hand, the $\R$-group $A = \{ g \in \operatorname{SL}_2(\C) \colon g\text{ is diagonal}\}$ is an $\R$-split torus. Since $\C$ is algebraically closed, both tori are $\D$-split. Both $\operatorname{SO}(2)$ and $A$ are maximal tori of $\operatorname{SL}_2$ and they are conjugate under $\operatorname{SL}_2$ and $\operatorname{SL}_2(\C)$, but not under $\operatorname{SL}_2(\R)$. The rank and the $\R$-rank of $\operatorname{SL}_2$ is $1$.

\subsection{The Lie algebra and the adjoint representation}

Let $e\in G$ be the neutral element of an algebraic $\K$-group $G$. There are multiple ways to define the Lie algebra of $G$, a few of which can be found for instance in Chapters 5 and 9 of \cite{Hum1}. A first explicit definition uses the description of $G$ as a matrix group
$$
G = \left\{ (a_{ij}) \in \operatorname{GL}_n(\D) \colon f(a_{ij}) = 0 \text{ for all } f \in I \right\}
$$
where $I \subset \K[(x_{ij}),\det(x_{ij})^{-1}]$ is a finite subset. For the following definition we restrict ourselves to subgroups of $\operatorname{SL}_n(\D)$, so that we can ignore the dependence on $\operatorname{det}(x_{ij})^{-1}$. For $f \in I$ we define the \emph{differential of $f$ at $e:=\operatorname{Id}$} to be the linear polynomial
$$
\operatorname{d}_{e}\!f = \sum_{ij} \frac{\partial f}{\partial x_{ij}} (e) \cdot x_{ij}
 \in \K[(x_{ij})]. 
$$
and the \emph{Zariski tangent space}  
$$
T_eG := \left\{ (a_{ij}) \in \D^{n\times n} \colon  \operatorname{d}_e\!f(a_{ij}) = 0 \text{ for all } f \in I\right\}.
$$

We now give another useful description of the tangent space in terms of derivations.

A \emph{derivation} of the coordinate ring $\D[G]$ (viewed as a $\D$-algebra) of a linear algebraic group $G$ is a linear map
$$
X \colon \D[G] \to \D[G]
$$
that satisfies the Leibnitz rule $X(ff') = X(f)f' + f X(f')$ for $f,f' \in \D[G]$. Algebraic groups act by left (and right) translation on their coordinate rings $\D[G]$. For $g \in G$, the left translation is given by $\lambda_g (f) (h) = f(g^{-1}h)$ for $f \in \D[G], h \in G$ (and the right translation by $\rho_g (f) (h) = f(hg)$). The Zariski tangent space $T_eG$ of $G$ can be identified with the set of derivations $\frakg$ of $\D[G]$ which commute with right translations
\begin{align*}
	\frakg &:= \left\{ X \colon \D[G] \to \D[G] \colon \begin{matrix}
		 X \text{ is a derivation of } \D[G] \text{ with }\\
		  X\circ \rho_g =\rho_g \circ X \text{ for all }  g \in G
	\end{matrix} \right\}.
\end{align*}
 Endowed with the Lie algebra structure defined by the bracket operation on derivations, $\frakg$ is called the \emph{Lie algebra of} $G$. For a commutative $\K$-subalgebra $\mathbb{B} \subseteq \D$ containing $\K$, we can use the identification $\frakg = T_eG$, to define the \emph{$\mathbb{B}$-points of the Lie algebra of $G$}
 $$
 \frakg_{\mathbb{B}} = T_eG \cap \mathbb{B}^{n\times n}.
 $$
 We note that the real Lie algebra $\frakg_{\R}$ is the Lie algebra of the real Lie group $G_{\R}$ \cite[Theorem 2.1]{Mil13}.
 
Let $\varphi  \colon G \to H$ be a morphism of algebraic groups $G,H$. The transposed map $\varphi^\circ \colon \D[H] \to \D[G]$ gives rise to a linear map 
\begin{align*}
	\operatorname{d}_e\!\varphi \colon  \frakg  &\to \frakh \\
	X &\mapsto X \circ \varphi^{\circ} 
\end{align*}
which is called the \emph{differential} $\operatorname{d}_e\!\varphi$ of $\varphi$ at $e\in G$. 

The group $G$ acts on itself by conjugation $\operatorname{Int}(g) \colon G \to G, h \mapsto ghg^{-1}$. The differential of $\operatorname{Int}(g)$ at $e$ is denoted $\operatorname{Ad} (g) $. The morphism $\operatorname{Ad} \colon G \to \operatorname{GL}(\frakg)$ is called the \emph{adjoint representation of} $G$ and is given by
$$
\operatorname{Ad}(g)(X) = gXg^{-1} \in T_eG
$$
when $X \in \frakg$ is viewed as an element of $T_eG$. The differential of $\operatorname{Ad}$ is called the \emph{adjoint representation of $\frakg$} and denoted by $\operatorname{ad} \colon \frakg \to \mathfrak{gl}(\frakg) $.

\subsection{Root systems and the spherical Weyl group}\label{sec:alg_root_system}

Let $G$ now be a semisimple 
algebraic group and $S$ a torus of $G$. Since $\operatorname{Ad} (S)$ is also a torus, its elements are simultaneously diagonalizable and the Lie algebra $\frakg$ decomposes into eigenspaces
$$
\frakg = \frakg_0^{(S)} \oplus\bigoplus_{\alpha \neq 0} \frakg_\alpha^{(S)}
$$
where 
$$
\frakg_\alpha^{(S)} = \{ X \in \frakg \colon \operatorname{Ad} (s) (X) = \alpha(s) \cdot X \ \text{ for all } s\in S \}
$$
for some $\alpha \in \hat{S}$. Elements $\alpha \neq 0$ with $\frakg_\alpha^{(S)} \neq 0$, are called the \emph{roots (relative to $S$)} and $\frakg_\alpha^{(S)}$ the \emph{root spaces}. The set of all roots is denoted by $\Phi(G,S)$. 
If $G$ is defined over $\K$ and $S$ is a maximal $\K$-split torus of $G$, then $\KPhi := \Phi(G,S)$ is called the set of $\K$-\emph{roots} of $G$. Since all maximal $\K$-split tori of $G$ are conjugate over $\K$ \cite[Theorem 20.9(ii)]{Bor}, $\KPhi$ only depends on $\K$ and not on the choice of maximal $\K$-split torus $S$.

\begin{theorem}[21.6 \cite{Bor}]\label{thm:alg_rootsystem}
	Let $G$ be a semisimple connected $\K$-group and $S$ a maximal $\K$-split torus of $G$. Recall that $\hat{S}\cong \Z^d$. Let $\Phi = \Phi(G,S)$. Then there is an admissible scalar product on the $\R$-vector space $V = \hat{S} \otimes_\Z \R$ such that ($V$,$\Phi$) is a crystallographic root system.
\end{theorem}

If $S$ is a maximal $\K$-split torus, the \emph{spherical Weyl group relative to $\K$} is
$$
\KW =\KW (S,G) = \operatorname{Nor}_{G}(S)/\operatorname{Cen}_{G}(S)
$$
and acts faithfully on $S$, $\hat{S}$ and $\KPhi$.

\subsection{Borel subgroups, parabolic subgroups}
Let $G$ be a connected algebraic group.
A subgroup is \emph{solvable} if its derived series consisting of iterated commutator groups terminates. A subgroup $B < G$ is a \emph{Borel subgroup} of $G$ if it is maximal among the connected solvable subgroups. 
A closed subgroup $P < G$ is \emph{parabolic} if and only if it contains a Borel subgroup. 
The minimal parabolic 
subgroups are exactly the Borel subgroups. In general, a Borel subgroup may not be defined over $\K$ and in this case the minimal parabolic $\K$-subgroups may not be Borel subgroups.

\begin{theorem}(Bruhat decomposition, \cite[Theorem 21.15]{Bor})
Let $S$ be a maximal $\K$-split torus and $P$ a minimal parabolic $\K$-subgroup containing $S$.  Denote by $\pi \colon \operatorname{Nor}_{G}(S) \to \KW$ the Weyl group projection. Then $G_\K = P_\K \operatorname{Nor}_{G_\K}(S_\K)P_\K$, in fact there is a disjoint union of double classes
$$
G_\K = \bigsqcup_{w \in  \KW } P_\K \pi^{-1}(w) P_\K.
$$
\end{theorem}

After choosing a minimal parabolic $\K$-subgroup $P$ containing a maximal $\K$-split torus, any parabolic $\K$-subgroup containing $P$ is called \emph{standard parabolic}.

\subsubsection{Examples}
The group of diagonal matrices in $\operatorname{SL}_3=\operatorname{SL}_3(\D)$ is a maximal $\Q$-split torus $T$. There are six minimal parabolic $\Q$-subgroups (which in this case are Borel subgroups) containing $T$. They are given by
\begin{align*}
	\left\{ \begin{pmatrix}
		\star & \star & \star  \\
		0 & \star & \star  \\
		0 & 0 & \star  \\
	\end{pmatrix} \in \operatorname{SL}_3 \right\} , \quad \left\{ \begin{pmatrix}
	\star & 0 & \star  \\
	\star & \star & \star  \\
	0 & 0 & \star  \\
\end{pmatrix} \in \operatorname{SL}_3 \right\} , \quad\left\{ \begin{pmatrix}
\star & \star & \star  \\
0 & \star & 0  \\
0 & \star & \star  \\
\end{pmatrix} \in \operatorname{SL}_3 \right\} , \\
\left\{ \begin{pmatrix}
	\star & 0 & 0  \\
	\star & \star & 0  \\
	\star & \star & \star  \\
\end{pmatrix} \in \operatorname{SL}_3 \right\} , \quad \left\{ \begin{pmatrix}
	\star & 0 & 0  \\
	\star & \star & \star  \\
	\star & 0 & \star  \\
\end{pmatrix} \in \operatorname{SL}_3 \right\} , \quad\left\{ \begin{pmatrix}
	\star & \star & 0  \\
	0 & \star & 0  \\
	\star & \star & \star  \\
\end{pmatrix} \in \operatorname{SL}_3 \right\} 
\end{align*}
and correspond to the Weyl chambers on which the spherical Weyl group acts transitively by conjugation. A corresponding set of representatives of $\KW$ is given by 
	\begin{align*}
		&\begin{pmatrix}
			1 & 0 & 0 \\
			0 & 1 & 0 \\
			0 & 0 & 1 \\
		\end{pmatrix}, \quad
		\begin{pmatrix}
			0 & -1 & 0 \\
			1 & \phantom{-}0 & 0 \\
			0 & \phantom{-}0 & 1 \\
		\end{pmatrix}, \quad
		\begin{pmatrix}
			1 & 0 & \phantom{-}0 \\
			0 & 0 & -1 \\
			0 & 1 & \phantom{-}0 \\
		\end{pmatrix}, \\
		&\begin{pmatrix}
			0 & 0 & -1 \\
			0 & 1 & \phantom{-}0 \\
			1 & 0 & \phantom{-}0 \\
		\end{pmatrix}, \quad
		\begin{pmatrix}
			0 & 0 & 1 \\
			1 & 0 & 0 \\
			0 & 1 & 0 \\
		\end{pmatrix}, \quad
		\begin{pmatrix}
			0 & 1 & 0 \\
			0 & 0 & 1 \\
			1 & 0 & 0 \\
		\end{pmatrix}.
	\end{align*}

\newpage

\section{Results about Lie algebras and split tori}
Let $G$ be a connected semisimple algebraic group defined over $\K \subseteq \R$, which is invariant under transposition. The Lie algebra $\frakg = T_eG$ is defined by finitely many linear equations with coefficients in $\K$ and we can therefore consider its $\K$-points $\frakg_\K \subseteq \K^{n\times n}$, which is a $\K$-vector space and an algebraic (hence semialgebraic) set. Let $\frakg_\F$ be the semialgebraic extension of $\frakg_\K$ for any real closed field $\F \supseteq \K$.

In this section, we recall facts about the real Lie group $G_\R$ and its Lie algebra $\frakg_\R$. We sometimes reference chapter 3 of S. Helgason's book \cite{Hel2} and chapter 6 of A. Knapp's book \cite{Kna}, a compact account of which can be found in \cite{Wis01}.

\subsection{Cartan involutions and Cartan decompositions}\label{sec:killing_involutions_decompositions}

Recall that the \emph{adjoint representation} of the Lie algebra $\frakg_\R$ 
$$
\operatorname{ad}_{\frakg_\R} := \operatorname{d}_e\!\operatorname{Ad}_{G_\R} \colon \frakg_\R \to \mathfrak{gl}(\frakg_\R)
$$
is given by $X \mapsto [X,\cdot\,]$. The \emph{Killing form} $B$ is the bilinear form on $\frakg_\R$ defined by
$$
B(X,Y) := \operatorname{tr}\left(\operatorname{ad}(X)\circ\operatorname{ad}(Y)\right)
$$
for $X,Y \in \frakg_\R$.
A Lie algebra is called \emph{simple} if it is non-abelian and does not contain any non-zero proper ideals. 
A Lie algebra is called \emph{semisimple} if it is a direct product of simple Lie algebras. By Cartan's criterion, being semisimple is equivalent to having a non-degenerate Killing form. Since $G_\R$ is semisimple, so is $\frakg_\R$.
A Lie algebra automorphism $\theta \colon \frakg_\R \to \frakg_\R$ is called an \emph{involution} if $\theta^2 = \operatorname{Id}$. For any involution $\theta \colon \frakg_\R \to \frakg_\R$ we can define a bilinear form $B_\theta$ by
$$
B_\theta(X,Y) := -B(X,\theta Y)
$$ 
for $X,Y \in \frakg_\R$. If $B_\theta$ is positive-definite, then $\theta$ is called a \emph{Cartan involution}. The following is a technical result on Cartan involutions which is stated as an exercise in \cite{Hel2} and proven in \cite[Theorem 6.16]{Kna}.
\begin{lemma} \label{lem:aligntheta}
	Let $\frakg_\R$ be a real semisimple Lie algebra, $\theta$ a Cartan involution and $\sigma$ any involution on $\frakg_\R$. Then there exists an automorphism
	$
	\varphi \in \operatorname{Ad}(G_\R) 
	$
	such that $\varphi \theta \varphi^{-1}$ commutes with $\sigma$.
\end{lemma}
An application of the preceding Lemma is that all Cartan involutions are conjugated, see \cite[Corollary 6.19]{Kna}.
\begin{theorem}\label{thm:cartan_inv_unique}
	Let $\frakg_\R$ be a real semisimple Lie algebra. Any two Cartan involutions of $\frakg_\R$ are conjugate via an element of $\operatorname{Ad}(G_\R)$.
\end{theorem}
\begin{proof}
	Let $\theta$ and $\theta'$ be two Cartan involutions of $\frakg_\R$. If $\theta$ and $\theta'$ commute, they have the same eigenspaces. We claim $E_{1}(\theta) = E_1(\theta')$ and $E_{-1}(\theta)=E_{-1}(\theta')$, since otherwise if for instance $\theta(X)=X$ and $\theta'(X)=-X$, then 
	$$
	B_{\theta'}(X,X) = -B(X,\theta'(X)) = - B(X,-X) = B(X,\theta(X)) = - B_{\theta}(X,X),
	$$
	but both $B_{\theta}(X,X)$ and $B_{\theta'}(X,X)$ should be positive. We conclude that if $\theta$ and $\theta'$ commute, then $\theta = \theta'$, since $\theta$ and $\theta'$ take the same values on their eigenspaces. If $\theta$ and $\theta'$ do not commute, then we can apply Lemma \ref{lem:aligntheta} to find $\varphi \in \operatorname{Ad}(G_\R)$ such that $\varphi \theta \varphi^{-1}$ commutes with $\varphi'$, and therefore $\varphi \theta \varphi^{-1} = \theta'$ by the previous argument.
\end{proof}

A decomposition of $\frakg_\R$ as a direct sum $
\frakg_\R =  \frakk \oplus \frakp.
$
is called a \emph{Cartan decomposition} if
$$
[\frakk, \frakk] \subseteq \frakk, \quad [\frakk, \frakp] \subseteq \frakp , \quad [\frakp, \frakp] \subseteq \frakk
$$
and the Killing form $B$ is negative definite on $\frakk$ and positive definite on $\frakp$. There is a correspondence between Cartan involutions and Cartan decompositions. A Cartan involution $\theta$ defines a decomposition into eigenspaces $\frakk = E_{1}(\theta)$, $\frakp = E_{-1}(\theta)$ of $\frakg_\R$. The bracket relations can be checked using that $\theta$ commutes with the bracket operation. Since then $\operatorname{ad}(\frakk)\operatorname{ad}(\frakp)(\frakk) \subseteq \frakp$ and $\operatorname{ad}(\frakk)\operatorname{ad}(\frakp)(\frakp) \subseteq \frakk$, we have 
$$
B(\frakk, \frakp) = \operatorname{Tr}(\operatorname{ad}(\frakk)\operatorname{ad}(\frakp)) = 0
$$ 
and the decomposition is orthogonal. Since $\theta$ is a Cartan involution, $B_\theta$ is positive definite, and thus $B$ is negative definite on $\frakk$ and positive definite on $\frakp$. Starting from a Cartan decomposition $\frakg_\R = \frakk \oplus \frakp$, we can define the involution
$$
\theta \colon X \mapsto \left\{\begin{matrix}
	X & \text{ if }X \in \frakk \\
	-X & \text{ if }X \in \frakp
\end{matrix} \right.
$$  
which is compatible with the bracket operation and for which $B_\theta$ is positive definite. To refine the decomposition further, we need the following Lemma.
\begin{lemma}
	Let $\frakg_\R = \frakk \oplus \frakp$ be a Cartan decomposition with Cartan involution $\theta$ and let $X \in \frakp$. The map $\operatorname{ad}(X) \colon \frakg_\R \to \frakg_\R$ is symmetric with respect to the scalar product $B_\theta$.
\end{lemma}
\begin{proof}
	We have
	\begin{align*}
		B_\theta(\operatorname{ad}(X)(Y),Z) & =-B([X,Y],\theta(Z)) = B(Y,[X,Z]) = - B(Y,[\theta(X),\theta(Z)])\\
		& = -B(Y,\theta([X,Z])) = B_\theta(Y,\operatorname{ad}(X)Z)
	\end{align*}  
for all $X \in \frakp$ and $Y,Z \in \frakg_\R$. 
\end{proof}
Let now $\frakg_\R = \frakk \oplus \frakp$ be a fixed Cartan decomposition with associated Cartan involution $\theta$. Let $\fraka \subseteq \frakp$ be a maximal abelian subspace. We can define the \emph{real rank of $G_\R$}, $\operatorname{rank}_\R(G_\R) := \operatorname{dim}(\fraka)$, which is independent of $\frakp$ as a consequence of Theorem \ref{thm:cartan_inv_unique} and independent of $\fraka$ since any two maximal abelian subspaces are conjugate to each other \cite[Theorem 6.51]{Kna}. The set $\{\operatorname{ad}(H) \colon H \in\fraka\}$ consists of symmetric, hence diagonalizable linear maps. Since they all commute, they are simultaneously diagonalizable. This results in the decomposition
$$
\frakg_\R = \frakg_0 \oplus\bigoplus_{\alpha \in \Sigma}\frakg_\alpha 
$$  
where for every $\alpha \in \fraka^\star$ in the dual space $\fraka^\star$ of $\fraka$
$$
\frakg_\alpha = \left\{ X \in \frakg_\R \colon [H,X] = \alpha(H) \cdot X \text{ for all } H \in \fraka  \right\}
$$
and where $\Sigma = \{\alpha \in \fraka^\star \colon \alpha \neq 0 \text{ and } \frakg_\alpha \neq 0\}$. The elements of $\Sigma$ are called \emph{restricted roots} and $\frakg_\alpha$ their associated \emph{restricted root spaces}. Note that in contrast to the decomposition of complex Lie algebras, the root spaces may not be one-dimensional\footnote{In the real semisimple Lie algebra $\mathfrak{sl}_2(\C)$ we have $\operatorname{dim}(\frakg_\alpha) = 2$.}. In Section \ref{sec:alg_root_system} we defined the root system of an algebraic group with a maximal split torus. In Section \ref{sec:compatibility} we will show that these two approaches give the same root spaces and root systems.
\begin{proposition}(\cite[Proposition 6.40]{Kna})\label{prop:root_decomp} 
	The \emph{restricted root space decomposition} is an orthogonal direct sum with respect to $B_\theta$ and satisfies for all $\alpha,\beta \in \Sigma$
	\begin{itemize}
		\item [(i)] $
		[\frakg_{\alpha}, \frakg_{\beta}] \subseteq \frakg_{\alpha + \beta}$,
		\item [(ii)]$\theta (\frakg_\alpha) = \frakg_{-\alpha}$ and
		\item [(iii)]$\frakg_0 = \fraka \oplus \mathfrak{z}_\frakk(\fraka)$, where $ \mathfrak{z}_\frakk(\fraka) = \left\{X\in \frakk \colon [H,X] = 0  \text{ for all } H \in \fraka \right\}$.
	\end{itemize}
\end{proposition}

The inner product $B_{\theta}$ may be restricted to $\fraka$ and used to set up an isomorphism $\fraka^{\star} \cong \fraka$ which turns $\fraka^{\star}$ into a Euclidean vector space.

\begin{theorem}(\cite[Corollary 6.53]{Kna})\label{thm:sigma_root}
	The set of roots $(\Sigma,\fraka^\star)$ is a crystallographic root system\footnote{$\Sigma$ may not be reduced, for example when $G= \operatorname{SU}(2,1)$.
		}.
\end{theorem}
By choosing an ordered basis of the root system $\Sigma$, we can define the set of positive roots $\Sigma_{>0}$. 
Then
$$
\frakn = \bigoplus_{\alpha \in \Sigma_{>0}} \frakg_\alpha
$$
is a nilpotent subalgebra of $\frakg_\R$. We have the following properties. 
\begin{lemma}\label{lem:kan_form}(\cite[Lemma 6.45]{Kna})
	Let $\frakg_\R$ be a real semisimple Lie algebra. There exists a basis of $\frakg_\R$ such that the matrices representing $\operatorname{ad}(\frakg_\R)$ have the following properties
	\begin{enumerate}
		\item The matrices of $\operatorname{ad}(\frakk)$ are skew-symmetric.
		\item The matrices of $\operatorname{ad}(\fraka)$ are diagonal.
		\item The matrices of $\operatorname{ad}(\frakn)$ are upper triangular with $0$ on the diagonal.
	\end{enumerate}
\end{lemma}
We also have the Iwasawa decomposition on the level of Lie algebras. 
\begin{theorem}(\cite[Proposition 6.43]{Kna})
	Let $\frakg_\R$ be a semisimple Lie algebra. Then $\frakg_\R = \frakk \oplus \fraka \oplus \frakn$ is a direct sum.
\end{theorem}

The following lemma about the Lie subalgebra 
$$
\frakl := ( \frakg_\alpha \oplus \frakg_{2\alpha} ) \oplus (\frakg_{-\alpha} \oplus \frakg_{-2\alpha} ) \oplus ( [\frakg_\alpha , \frakg_{-\alpha}] + [\frakg_{2\alpha}, \frakg_{-2\alpha}] )
$$
for some $\alpha \in \Sigma$ will be useful when we consider certain rank 1 subgroups in Section \ref{sec:rank1}. Note that by Cartan's criterion, $\frakl$ is semisimple, as the Killing form is the restriction of the definite Killing form of $\frakg$. Then $\frakk \cap \frakl \oplus \frakp \cap \frakl$ is a Cartan decomposition of $\frakl$.
\begin{lemma}\label{lem:levi_algebra} 
	Let $\alpha \in \Sigma$ and $X\in \frakg_\alpha  \setminus \{0\}$. Then the real rank $\operatorname{rank}_\R$ of $\frakl$ is one. A maximal abelian subspace of the symmetric part of $\frakl$ is given by $\fraka \cap \frakl = \langle [X,\theta(X)]\rangle$.
\end{lemma}
\begin{proof}
	We use the fact that $B_\theta$ gives us a scalar product on $\fraka$ allowing us to identify $\fraka \cong \fraka^\star$, sending $\alpha \in \Sigma$ to $H_\alpha$ defined by $\alpha(H)=B_\theta(H_\alpha,H)$ for all $H \in \fraka$. Let $X\in \frakg_\alpha, Y \in \frakg_{-\alpha}$ and $H \in \fraka$. Then
	\begin{align*}
		B_\theta([X,Y],H) & = -B([X,Y],\theta(H)) = B([X,Y], H) = - B(Y,[X,H])  \\
		&= B(Y,[H,X]) = B(Y,\alpha(H)X) = \alpha(H)B(Y,X) \\
		&=  B_\theta(H_\alpha, H)B(X,Y) = B_\theta(B(X,Y) \cdot H_\alpha ,H),
	\end{align*}
	where we used that $B$ is $\operatorname{ad}(X)$-invariant. 
	The element 
	$$
	W = [X,Y]- B(X,Y) \cdot H_\alpha  \in \frakg_0
	$$
	satisfies $B_\theta(W,H)=0$ for all $H \in \fraka$, hence lies in $\frakg_0$ perpendicular to $\fraka$, hence $W \in \mathfrak{z}_\frakk(\fraka) \subseteq \frakk$ by Proposition \ref{prop:root_decomp}. Similarly, one can show that for any $X' \in \frakg_{2\alpha}$ and $Y' \in \frakg_{-2\alpha}$, $[X',Y'] =W' + 2B(X',Y') \cdot H_\alpha$ for some $W' \in \mathfrak{z}_\frakk(\fraka)$. Any element $Z\in \fraka \cap \frakl$  lies in $[\frakg_\alpha , \frakg_{-\alpha} ] \oplus [ \frakg_{2\alpha} , \frakg_{-2\alpha}]$ since $\fraka \subseteq \frakg_0$ is orthogonal to $\frakg_\alpha \oplus \frakg_{2\alpha} \oplus \frakg_{-\alpha} \oplus \frakg_{-2\alpha}$. Therefore, $Z$ can be written as $Z=W+W'+ (B(X,Y)+2B(X',Y')) \cdot H_\alpha  \in \mathfrak{z}_\frakk(\fraka) \oplus \fraka$, hence $W+W'=0$ and $Z \in \langle H_\alpha \rangle$.
	
	If $Y=\theta(X)$, then 
	$$
	\theta(W)=\theta([X,\theta(X)]-B(X,Y)\cdot H_\alpha ) = [\theta(X),X]-\theta(H_\alpha) B(X,Y) = - W
	$$ 
	which implies $W \in \frakp \cap \frakk$, hence $W=0$. We have shown that $[X,\theta(X)] \in \fraka \cap \frakl \subseteq \langle H_\alpha \rangle$. Since $B$ is definite, $[X,\theta(X)] \neq 0$ and $\fraka \cap \frakl = \langle [ X, \theta(X)] \rangle$.
	
	The map $\theta|_{\frakl}$ is a Cartan-involution with Cartan decomposition $\frakl = \frakk\cap \frakl \oplus \frakp \cap \frakl$. There is a maximal abelian subspace of $\frakp \cap \frakl$ contained in $\fraka$, which is equal to $\fraka \cap \frakl = \langle [X,\theta(X)] \rangle$ by the above. The dimension of a maximal abelian subspace of $\frakp \cap \frakl$ is exactly the rank of $\frakl$ and it is one. 
\end{proof}

The following is a special case of the Jacobson-Morozov Lemma in the literature. Note that this Lemma does not work for general $X\in \frakg_\alpha \oplus \frakg_{2\alpha}$.

\begin{lemma}\label{lem:JM_basic}
	Let $\alpha \in \Sigma$ and $X \in \frakg_\alpha \setminus \{0\}$. Then there is a $Y \in \frakg_{-\alpha} \setminus \{0\}$ and $H \in \fraka$ such that $(X,Y,H)$ forms a $\mathfrak{sl}(2)$-triplet, i.e.
	$$
	[X,Y] = H, \quad [H,X] = 2X,  \quad [H,Y]=2Y.
	$$
\end{lemma}
\begin{proof}
	Let $H_\alpha \in \fraka$ be as in the proof of Lemma \ref{lem:levi_algebra} defined by $\alpha(H) = B_{\theta}(H_\alpha,H)$ for all $H \in \fraka$. Given $X \in \frakg_\alpha$, let
	$$
	Y := \frac{-2}{B_{\theta}(X,X) \cdot \alpha(H_\alpha)} \theta(X) \in \frakg_{-\alpha}
	$$
	and $H:= [X,Y]$. In the proof of Lemma \ref{lem:levi_algebra} we saw that
	$$
	[X,\theta(X)] = B(X,\theta(X))\cdot H_\alpha = -B_{\theta}(X,X) \cdot H_\alpha,
	$$ 
	whence $H \in \fraka$ and moreover
	\begin{align*}
		[H,X] &= \left[\left[X,\frac{-2}{B_{\theta}(X,X)\alpha(H_\alpha)} \theta(X)\right],X\right] \\
		&= \frac{2}{\alpha(H_\alpha)} \left[H_\alpha,X\right] = \frac{2}{\alpha(H_\alpha)} \alpha(H_\alpha) X = 2X
	\end{align*}
	and similarly $[H,Y] = -2Y$.
\end{proof}

\subsubsection{Examples}
Since $G_\R$ is invariant under transposition, $\sigma \colon g \mapsto (g^{-1})\tran$ is an involution of the Lie group whose differential $\theta := \operatorname{d}_e\! \sigma \colon X \mapsto -X\tran$ is an involution on the Lie algebra $\frakg_\R$ which decomposes $\frakg_\R= \frakp \oplus \frakk$ into a symmetric part $\frakp$ and an skew-symmetric part $\frakk$. One can check that
$$
[\frakk, \frakk] \subseteq \frakk, \quad [\frakk, \frakp] \subseteq \frakp  \quad \text{ and } \quad [\frakp, \frakp] \subseteq \frakk
$$
and that the Killing form is negative definite on $\frakk$ and positive definite on $\frakp$. Thus $\frakg_\R = \frakk \oplus \frakp$ is a Cartan decomposition and $\theta$ is a Cartan involution.

\subsection{Real Cartan subalgebras and $\R$-split subalgebras}\label{sec:real_cartan_split}

In this section we discuss Cartan subalgebras of a finite dimensional semisimple real Lie algebra $\frakg$. These results apply in particular to the Lie algebra $\frakg_\R$ of the real Lie group $G_\R$. 

An abelian subalgebra $\frakh \subseteq \frakg$ is called a \emph{Cartan subalgebra}\footnote{More generally, a Cartan subalgebra is defined to be a nilpotent self-normalizing subalgebra. In our context, we restrict to semisimple Lie algebras over fields of characteristic $0$, in which case Cartan subalgebras are abelian. } 
 if 
$$
\frakh = \operatorname{Nor}_{\frakg}(\frakh) := \left\{ X \in \frakg \colon [X,Y] \in \frakh \text{ for all } Y \in \frakh \right\}.
$$
An abelian subalgebra $\fraka \subseteq \frakg$ is called \emph{$\R$-split} if $\operatorname{ad}(X)\colon \frakg \to\frakg$ is diagonalizable over $\R$ for every $X \in \fraka$. Let $r_\R(\frakg)$ be the maximal dimension of such an $\R$-split abelian subalgebra and denote
$$
V(\frakg) := \{ \fraka \subseteq \frakg \colon \text{$\fraka$ is an $\R$-split abelian subalgebra with $\dim(\fraka) = r_\R(\frakg)$}\}.
$$ 
A \emph{maximally $\R$-split Cartan subalgebra} 
is a Cartan subalgebra containing an element of $V(\frakg_\R)$ as a subset. Let
$$
\mathcal{C} (\frakg) = \left\{  \frakh \subseteq \frakg \colon \frakh \text{ is a maximally $\R$-split Cartan subalgebra } \right\}.
$$

Knapp \cite{Kna} uses slightly different definitions and names for these objects. We will show here that they coincide. We start by showing that our notion of Cartan subalgebra coincides with the notion in \cite{Kna}.

\begin{lemma}\label{lem:defs_Cartan_subalgebra}
	Let $\mathfrak{h}$ be a subalgebra of a finite dimensional real semisimple Lie algebra $\frakg$. Then $\frakh$ is a Cartan subalgebra if and only if it satisfies one of the following conditions
	\begin{enumerate}
		\item [(i)] $\frakh$ is abelian and $\operatorname{Nor}_\frakg (\frakh) = \frakh$.
		\item [(ii)] The complexification $\frakh_\C := \frakh \oplus i \frakh$ is abelian and $\operatorname{Nor}_{\frakg_\C} (\frakh_\C) = \frakh_\C$.
		\item [(iii)] The complexification $\frakh_\C $ is nilpotent and $\operatorname{Nor}_{\frakg_\C} (\frakh_\C) = \frakh_\C$.
		\item [(iv)] The complexification $\frakh_\C $ is nilpotent and 
		$$
		\frakh_\C = \left\{ X \in \frakg_\C \colon \exists n \in \N ,\ \forall H \in \frakh_\C ,\ \operatorname{ad}(H)^n X = 0  \right\}.
		$$
		\item [(v)] The complexification $\frakh_\C$ is maximal abelian and the subset $\operatorname{ad}(\frakh_\C) \subseteq \mathfrak{gl}(\frakg_\C)$ is simultaneously diagonalizable.
	\end{enumerate}
\end{lemma}
\begin{proof}
	Notion (i) is our definition of Cartan subalgebra. For any $X = X_1 + i X_2 \in \frakg_\C$ and $Y= Y_1+iY_2 \in \frakg_\C$, we have
	$$
	\operatorname{ad}(Y)X  = [Y_1+iY_2 , X_1 + iX_2] = [Y_1,X_1] - [Y_2,X_2] +i ([Y_1, X_2] + [Y_2,X_1]).
	$$  
	From this formula we can deduce that a subalgebra of a real Lie algebra is abelian if and only if its complexification is abelian. It also follows that for any real subalgebra $\frakh \subseteq \frakg$ we have
	$$
	\operatorname{Nor}_{\frakg}(\frakh) = \frakh \iff \operatorname{Nor}_{\frakg_\C}(\frakh_\C) = \frakh_\C. 
	$$
	It follows that notions (i) and (ii) coincide. For (ii) implies (iii), it suffices to notice that every abelian subalgebra is nilpotent, the converse uses semisimplicity and is given in \cite[Proposition 2.10]{Kna}. Condition (iv) is the definition used in \cite{Kna} and the equivalence of (iii) and (iv) is given by \cite[Proposition 2.7]{Kna}. The equivalence of (iv) and (v) is given by \cite[Corollary 2.13]{Kna}.
\end{proof}

A subset $\frakh\subseteq \frakg$ is called \emph{$\theta$-stable} if there is a Cartan involution $\theta$ such that $\theta(\frakh) = \frakh$. If $\theta$ is given by the context, $\theta$-stable refers to that particular Cartan involution. 
\begin{proposition}\label{prop:Kna_real_Cartan_conj}(\cite[Prop 6.59]{Kna})
	Let $\frakg$ be the Lie algebra of a real semisimple Lie group $G$ and $\operatorname{Ad}(G)^\circ$ the connected component of $\operatorname{Ad}(G) \subseteq \operatorname{GL}(\frakg)$ in the Lie group topology. Any Cartan subalgebra $\frakh$ is conjugate via $\operatorname{Ad}(G)^\circ$ to a $\theta$-stable Cartan subalgebra.
\end{proposition}
If $\frakg = \frakp \oplus \frakk$ is the Cartan decomposition associated to a Cartan involution $\theta$, and $\frakh$ is a $\theta$-stable Cartan subalgebra, then $\fraka := \frakh \cap \frakp$ and $\frakt := \frakh \cap \frakk$ are $\theta$-stable and $\frakh = \fraka \oplus \frakt$. A Cartan subalgebra is called \emph{maximally noncompact} if $\dim(\fraka)$ is maximal. 
While all complex Cartan subalgebras are conjugated \cite[Theorem 2.15]{Kna}, for real Cartan subalgebras, the dimensions of the spaces $\fraka$ and $\frakt$ are preserved.
\begin{proposition}\label{prop:Kna_real_theta_Cartan_conj}(\cite[Prop 6.61]{Kna}) Let $\frakg$ be the Lie algebra of a semisimple real Lie group $G$ and let $K$ be a subgroup with Lie algebra $\frakk$ where $\frakg = \frakp \oplus \frakk$ is the Cartan decomposition associated to a Cartan involution $\theta$. 
	
	All maximally noncompact $\theta$-stable Cartan subalgebras are conjugate under $\operatorname{Ad}(K)^\circ$.
\end{proposition}
Together, Propositions \ref{prop:Kna_real_Cartan_conj} and \ref{prop:Kna_real_theta_Cartan_conj} allow us to extend our definition of maximally noncompact $\theta$-stable Cartan subalgebras to all Cartan subalgebras, by stipulating that a Cartan subalgebra is called \emph{maximally noncompact} if it is conjugated to a $\theta$-stable maximally noncompact Cartan subalgebra.
We will now show that the set of maximally $\R$-split Cartan subalgebras $\mathcal{C}(\frakg)$ coincides with the set of maximally noncompact Cartan subalgebras.

\begin{proposition}\label{prop:defs_max_noncompact}
	For any maximally $\R$-split subalgebra $\fraka \subseteq \frakg$ of a finite dimensional real semisimple Lie algebra $\frakg$, there is a maximally noncompact Cartan subalgebra $\frakh \supseteq \fraka$. Every maximally noncompact Cartan subalgebra contains a maximally $\R$-split subalgebra $\fraka$ and if $\frakh$ is $\theta$-stable, then $\fraka = \frakh \cap \frakp$ and $\frakh = \fraka \oplus \frakt$ in the decomposition above and $r_\R(\frakg) = \dim (\fraka)$.
\end{proposition}
\begin{proof}
	We introduce a few concepts. If $k$ is a field and $\frakg$ is a semisimple $k$-Lie algebra, then an abelian subalgebra $\frakh$ is called \emph{toral} if $\operatorname{ad}(\frakh) \subseteq \mathfrak{gl}(\frakg)$ consists of linear maps that are diagonalizable over the algebraic closure of $k$. If all the elements of $\operatorname{ad}(\frakh)$ are diagonalizable over $k$ itself, then $\frakh$ is called \emph{$k$-split toral}. If moreover $\frakh$ is maximal among all $k$-split toral subalgebras, $\frakh$ is called \emph{maximal $k$-split toral}.
	
	Let now $\frakg$ be a finite dimensional real semisimple Lie algebra. Let $\fraka \in V(\frakg)$, in our terminology $\fraka$ is maximal $\R$-split toral. Let $\frakh$ be maximal toral containing $\fraka$, which exists since $\frakg$ is finite dimensional. We want to show that $\frakh$ is a Cartan subalgebra and use characterization (v) of Lemma \ref{lem:defs_Cartan_subalgebra}. Since $\frakh$ is abelian, so is its complexification $\frakh_\C$. Since $\frakh$ is $\C$-split, all elements of $\operatorname{ad}(\frakh_\C)=\operatorname{ad}(\frakh) \oplus i \operatorname{ad}(\frakh)$ are diagonalizable over $\C$. This means that $\frakh_\C$ is a toral subalgebra of $\frakg_\C$ and by a general linear algebra fact, $\operatorname{ad}(\frakh_\C)$ is simultaneously diagonalizable. The complexification $\frakh_\C$ is also maximal abelian, since for all $X=X_1+iX_2 \in \frakg_\C$, if $[X,\frakh_\C] =0$, then we have in particular for $H \in \frakh$
	$$
	0 = [X,H] = [X_1+iX_2, H] = [X_1, H] + i [X_2, H] ,
	$$
	hence $X_1,X_2 \in \frakh$, so $X \in \frakh_\C$. This shows that $\frakh$ is a Cartan subalgebra by characterization (v).
	
	By Proposition \ref{prop:Kna_real_Cartan_conj}, there is a Cartan involution $\theta$ and a $g\in G_\R$ such that $\operatorname{Ad}(g)(\frakh)$ is a $\theta$-invariant Cartan subalgebra with decomposition $\operatorname{Ad}(g)(\frakh) = \fraka' \oplus \frakt'$ as described before. Since $\fraka' \subseteq \frakp $, it is $\R$-split. In fact, for $X=X_{\fraka'} + X_{\frakt'} \in \fraka' \oplus \frakt'$, $\operatorname{ad}(X)$ is only diagonalizable over $\R$ if $X_{\frakt'}=0$. Since diagonalizability is preserved under conjugation, $\operatorname{Ad}(\fraka) \subseteq \fraka'$ and by maximality of $\fraka$, we have $\operatorname{Ad}(\fraka) = \fraka'$. Thus $\frakh$ is maximally noncompact. 
	
	If we instead start with a maximally noncompact Cartan subalgebra $\frakh \subseteq \frakg$, we similarly obtain (possibly using a conjugation to a $\theta$-stable one) an $\R$-split subalgebra $\fraka$ with $\frakh = \fraka \oplus \frakt$. A priori $\fraka$ is maximal $\R$-split as a subalgebra of $\frakh$, but if $\fraka$ was contained in a larger subalgebra $\fraka'$ maximal $\R$-split in $\frakg$, then the above construction would result in a Cartan subalgebra $\frakh' = \fraka' \oplus \frakt'$ with $\dim(\fraka')>\dim(\fraka)$, which is impossible, since we assumed $\frakh$ to be maximally noncompact, which means $\dim(\fraka)$ is maximal. By definition, $r_\R(\frakg)$ is the maximal dimension of an $\R$-split abelian subalgebra, so $r_\R(\frakg) = \dim(\fraka)$.
\end{proof}

We now turn back to the Lie algebra $\frakg_\R$ of the $\R$-points $G_\R$ of the self-adjoint semisimple algebraic group $G$. We now collect a few properties of $V(\frakg_\R)$ and $\mathcal{C}(\frakg_\R)$. 
\begin{lemma}\label{lem:basic_facts} Recall that for the algebraic group $G$, $ \operatorname{rank}_\R(G)$ is the maximal dimension of any abelian subspace of $\frakp \subseteq \frakg_\R$, where $\frakg_\R = \frakp \oplus \frakk$ is the Cartan decomposition associated to the Cartan involution $X \mapsto -X\tran$. We have
	\begin{itemize}
		\item [(1)] $r_\R(\frakg_\R) = \operatorname{rank}_\R(G)$.
		\item [(2)] $G_\R$ acts transitively on $\mathcal{C}(\frakg_\R)$ and $ V(\frakg_\R)$.
		\item [(3)] $\mathcal{C}(\frakg_\R)$ contains a $\theta$-stable element. 
	\end{itemize}
\end{lemma}
\begin{proof}
	A maximally $\R$-split Cartan subalgebra, or by Proposition \ref{prop:defs_max_noncompact} a maximally noncompact Cartan subalgebra can be obtained by starting with a Cartan involution $\theta$, taking $\fraka \subseteq \frakp$ a maximal abelian subspace of $\frakp$ and taking $\frakt \subseteq \operatorname{Cen}_{\frakk}(\fraka)$ a maximal abelian subspace of $\operatorname{Cen}_{\frakk}(\fraka)$. Then
	$\frakh = \fraka \oplus \frakt$ is a $\theta$-stable maximally noncompact Cartan subalgebra, see \cite[Proposition 6.47]{Kna}. This shows (3).
	
	By Proposition \ref{prop:defs_max_noncompact}, $r_\R(\frakg_\R) = \dim(\fraka) $ and since $\fraka$ is a maximal abelian subspace of $\frakp$, $r_\R(\frakg_\R) = \operatorname{rank}_\R(G)$, showing (1).
	
	By Propositions \ref{prop:Kna_real_Cartan_conj} and \ref{prop:Kna_real_theta_Cartan_conj}, $\operatorname{Ad}(G_\R)$ and thus $G_\R$ act transitively on $\mathcal{C}(\frakg_\R)$. Next, let $\fraka, \fraka'\in V(\frakg_\R)$. By Proposition \ref{prop:defs_max_noncompact}, there are maximally noncompact Cartan subalgebras $\frakh, \frakh'$ such that $\fraka \subseteq \frakh$ and $\fraka' \subseteq \frakh'$. By Proposition \ref{prop:Kna_real_theta_Cartan_conj}, both $\frakh$ and $\frakh'$ are conjugated to a $\theta$-stable maximally noncompact Cartan subalgebra $\frakh''$. By Proposition \ref{prop:defs_max_noncompact}, the corresponding conjugates of $\fraka$ and $\fraka'$ coincide with the intersection $\frakh'' \cap \frakp$, hence with each other. This shows that $\operatorname{Ad}(G_\R)$ acts transitively on $V(\frakg_\R)$ and concludes the proof of (2).
\end{proof}

\subsection{Lie algebras over real closed fields}

Let $\K$ be a subfield of $\R$ and $\F$ be a real closed fields with $\K \subseteq \F$. In this section we additionally assume that $\K$ is real closed. Let $\frakg_\K  \subseteq \K^{n\times n}$ be the $\K$-points of the Lie algebra of a connected semisimple algebraic group $G$ defined over $\K$. Let $\frakg_\R$ and $\frakg_\F$ be the semialgebraic extensions. The definitions of Section \ref{sec:real_cartan_split} apply also to $\frakg_\F$:

An abelian subalgebra $\frakh \subseteq \frakg_\F$ is called \emph{Cartan subalgebra} if $\frakh = \operatorname{Nor}_{\frakg_\F}(\frakh)$. An abelian subalgebra $\fraka \subseteq \frakg_\F$ is called \emph{$\F$-split} if $\operatorname{ad}(X) \colon \frakg_\F \to \frakg_\F$ is diagonalizable over $\F$ for every $X \in \fraka$. Let $r_{\F}(\frakg_\F)$ be the maximal dimension of such an $\F$-split abelian subalgebra and denote by $V(\frakg_\F)$ the set of all $\F$-split abelian subalgebras with $\dim(\fraka) = r_\F (\frakg_\F)$. A \emph{maximally $\F$-split Cartan subalgebra} is a Cartan subalgebra containing an element of $V(\frakg_\F)$ as a subset. We denote by $\mathcal{C}(\frakg_\F)$ the set of all maximally $\F$-split Cartan subalgebras. We will now use the transfer principle to relate $\F$-split algebras to the real subalgebras studied in the previous section.

\begin{lemma} \label{lem:3}
	Whenever $\F$ and $\K$ are two real closed fields with $\K \subseteq \F$ and $\K \subseteq \R$, then
	\begin{itemize}
		\item [(1)] $r_\K(\frakg_\K) = r_\R (\frakg_\R) =  r_\F (\frakg_\F)$.
		\item [(2)] $G_\K$ acts transitively on $\mathcal{C}(\frakg_\K)$ and $V(\frakg_\K)$ and those two sets are non-empty.
		\item [(3)] \label{lem:3.3}$\mathcal{C}(\frakg_\K)$ contains a $\theta$-stable element, i.e. there is $\frakh \in \mathcal{C}(\frakg_\K)$ such that $\theta(\frakh) = \frakh$.
		\item [(4)] Let $\fraka$ and $\frakh$ be subalgebras of $\frakg_\K$. Then 
		\begin{align*}
			\fraka \in V(\frakg_\K) & \iff \fraka_\R \in V(\frakg_\R) \iff \fraka_\F \in V(\frakg_\F)
	\\
	\frakh \in \mathcal{C}(\frakg_\K) &\iff \frakh_\R \in \mathcal{C}(\frakg_\R) \iff \frakh_\F \in \mathcal{C}(\frakg_\F).	
	\end{align*}
	\end{itemize}
\end{lemma}

\begin{proof}
	Our first goal is to see $V(\frakg_\K)$ as a semi-algebraic set. Let $k$ be any field of characteristic $0$ containing $\K$. Let $\eL$ be a semisimple $k$-algebra for any field $k$, such that $\eL$ is also a semi-algebraic set defined over $\K$. Let $r \leq n := \operatorname{dim}_k(\eL)$. The set $\operatorname{Grass}_r(\eL)$ can be identified with an algebraic subset of $k^{n\times n}$, namely
	$$
	\varphi \colon \operatorname{Grass}_r(\eL) \to \left\{ A \in k^{n\times n} \colon A\tran = A, A^2=A, \operatorname{Tr}(A) = r \right\}
	$$
	sends the $k$-dimensional subspace $V$ to the orthogonal projection $\eL \to V$, see \cite[Theorem 3.4.4]{BCR}. Moreover, for any $A \in \varphi(\operatorname{Grass}_r(\eL))$ we have the description
	$$
	\varphi^{-1}(A) = \{ v \in \eL \colon Av = v \}.
	$$
	An abelian subalgebra $\fraka \subseteq \eL$ is $k$-split exactly when $\operatorname{ad}(x)$ is diagonalizable over $k$ for all $x \in \fraka$. It is enough to require that for a basis $\{v_1, \ldots , v_r\}$ of $\fraka$, the maps $\operatorname{ad}(v_i) \colon \eL \to \eL$ are diagonalizable over $k$. A linear map with matrix $M=(M_{ij})$ is diagonalizable over $k$ if and only if its characteristic polynomial decomposes into linear factors, which can be written as a first-order formula:
	$$
	f(M) = \exists x_1, \ldots , x_n \colon \forall X \colon  \det(M_{ij} - X\delta_{ij} ) = \prod_{i=1}^n (X-x_i).
	$$
	We can write the statement ``$\fraka = \varphi^{-1}(A)$ is a $k$-split abelian subalgebra'' as a first-order formula
	\begin{align*}
		& \exists v_1, \ldots , v_n \in \eL \colon \left( \forall a_1 , \ldots ,a_n \colon \sum_{i=1}^n a_i v_i = 0 \to a_1 = 0 \wedge \ldots \wedge a_n = 0\right) \wedge\\
		&\bigwedge_{i=1}^r Av_i = v_i \ \wedge  \bigwedge_{i,j=1}^r	[v_i,v_j] = 0 \ \wedge \\
		&
		\bigwedge_{\ell = 1}^r \left( \exists M_{11}, M_{12}, \ldots , M_{nn} \colon \bigwedge_{i=1}^n [v_\ell , v_i] = \sum_{j=1}^n M_{ij} e_j \, \rightarrow f(M) \right), 
	\end{align*}
	in words: there exists a basis of $\eL$ whose first $r$ vectors form a basis of $\fraka = \varphi^{-1}(A)$, such that $\fraka$ is abelian and for every $\ell \in \{1,\ldots, r\}$ $\operatorname{ad}(v_\ell)$ with matrix $M$ in this basis is diagonalizable over $k$.
	
	By quantifier elimination we get an equivalent first-order statement $g(A)$ without quantifiers and we can write 
	$$
	\varphi(V(\eL)_r) = \{A \in \varphi(\operatorname{Grass}_r(\eL)) \colon g(A) \}
	$$
	as a semi-algebraic set defined by polynomials with coefficients in $k$, where $V(\eL)_r$ denotes the set of all $k$-split abelian subalgebras of dimension $r$.

	Now for $k = \R$ we know that $V(\frakg_\R) = V(\frakg_\R)_{r_\R(\frakg_\R)}$ is non-empty, but $V(\frakg_\R)_{r}$ is empty for any $r > r_\R(\frakg_\R)$. We consider the first-order formula
	$$
	\exists A \in  \varphi( V(\frakg_k)_r ),
	$$
	which is defined over $\K$ since $\frakg_\R$ and hence the brackets are defined over $\K$. Since the formula is satisfied for $k=\R$ and $r=r_\R(\frakg_\R)$, we conclude by the transfer principle that it also holds for $k=\K$ and $k=\F$, i.e. $V(\frakg_\K) = V(\frakg)_{r_\R(\frakg_\R)} \neq \emptyset$ and $V(\frakg_\F) = V(\frakg_\F)_{r_\R(\frakg_\R)} \neq \emptyset$. For any $r > r_\R(\frakg_\R)$, we know that the formula is not satisfied for $k=\R$, therefore by the transfer principle it is also not satisfied for $k=\K$ and $k=\F$, i.e. $V(\frakg_\K)_{r} = \emptyset$ and $V(\frakg_\F)_{r} = \emptyset$ for any  $r > r_\R(\frakg_\R)$. It follows that $r_\F(\frakg_\F) = r_\R(\frakg_\R) = r_\K (\frakg_\K)$, which is statement (1) of the Lemma we are proving.
	
	We describe 
	\begin{align*}
		\mathcal{C}(\frakg_k) = \left\{ \frakh \in \operatorname{Grass}_{r_{k}(\frakg_k)}(\frakg_k) \colon 
		\begin{matrix}
			[\frakh , \frakh] = 0 , \ \operatorname{Nor}_{\frakg_k}(\frakh) = \frakh, \\
			\exists \fraka \in V(\frakg_k) \colon \fraka \subseteq \frakh 
		\end{matrix}
		\right\}
	\end{align*}
	similarly. Let $v_1, \ldots , v_n$ again be a basis of a semisimple $k$-algebra $\eL$ such that $v_1 , \ldots , v_r$ is a basis of a subalgebra $\fraka \subseteq \eL$. We have $\operatorname{Nor}_{\eL}(\fraka) = \fraka$ whenever $[v_i,v_j] \in \fraka \ \forall i \in \{1, \ldots , n\}, j \in \{1, \ldots , r\}$, i.e. given $\fraka = \varphi^{-1}(A)$ we have
	\begin{align*}
		h(A) \colon \quad &
		\bigwedge_{i=1}^n \bigwedge_{j=1}^r A[v_i , v_j] = [v_i,v_j].
	\end{align*}
	We have that $\varphi^{-1}(A) \in \mathcal{C}(\eL)$ if and only if the following first-order formula holds
	\begin{align*}
		& \exists v_1, \ldots , v_n \in \eL \colon \left( \forall a_1 , \ldots ,a_n \colon \sum_{i=1}^n a_i v_i = 0 \to \bigwedge_{i=1}^n  a_i = 0 \right) \wedge\\
		& \bigwedge_{i=1}^r Av_i = v_i \ \ \wedge\ \  
		\bigwedge_{i,j=1}^r [v_i,v_j] = 0   \ \ \wedge \ \ h(A) \ \ \wedge \\
		&  \exists  B \in V(\eL) \colon \forall v \in \eL \colon Bv=v \to Av = v. 
	\end{align*}
	Again we use quantifier elimination to to get an equivalent statement $s(A)$ without quantifiers. Thus, $\mathcal{C}(\eL)$ can be identified with the semialgebraic set
	$$
	\varphi(\mathcal{C}(\eL)) = \{ A  \in \varphi(\operatorname{Grass}_{r_k(\eL)} (\eL)) \colon s(A)\}.
	$$
	From the theory of real Lie groups we know by Lemma \ref{lem:basic_facts}(3) that $\mathcal{C}(\frakg_{\R})$ is non-empty and since $\exists A \in \varphi(\mathcal{C}(\frakg_k))$ is a first-order statement, defined over $\K$, we can use the transfer principle to conclude that $\mathcal{C}(\frakg_\K)$ and $\mathcal{C}(\frakg_{\F})$ are also non-empty. Statement (4) follows from the semi-algebraic description of the sets $V(\frakg_\K)$ and $\mathcal{C}(\frakg_\K)$ and the transfer principle.
	
	The group $G_k$ acts by conjugation on $\operatorname{Grass}_r(\frakg_k)$. The corresponding action on $\varphi(\operatorname{Grass}_r(\frakg_k))$ is given by $g.A = A \circ \operatorname{Ad}(g^{-1}) $, where $g\in G_k$, ${A \in \varphi(\operatorname{Grass}_r(\frakg_k))}$ and $\operatorname{Ad}(g) \colon \frakg_k \to \frakg_k, X \mapsto gXg^{-1}$.
	
	We know that this action is transitive on $V(\frakg_\R)$ and $\mathcal{C}(\frakg_\R)$ by Lemma \ref{lem:basic_facts}(2). As all involved sets are semi-algebraic, we can formulate transitivity as a first-order formula,
	\begin{align*}
		\forall A,B \in \varphi(V(\frakg_k)) \ \exists g \in G \colon \forall v \in \frakg_k \colon 
		A(g^{-1}vg) = Bv. 
	\end{align*}
	and conclude that $G_\K$ acts transitively on $V(\frakg_\K)$ and $\mathcal{C}(\frakg_\K)$, concluding the proof of (2).
	
	Finally, for $\theta \colon \frakg_k \to \frakg_k, X \mapsto -X\tran$ and $A \in \varphi(\mathcal{C}(\frakg_k))$, we note that $v \in \theta \left(\varphi^{-1}( A )\right)$ if and only if $A(\theta(v)) = \theta (v)$, that is $\theta A \theta v = v$. The condition $\theta(\fraka) = \fraka$ therefore corresponds to $A =\theta A\theta$. We know by Lemma \ref{lem:basic_facts}(3) that the first-order formula
	$$
	\exists A \in \varphi(\mathcal{C}(\frakg_k) ) \colon \forall v \in \frakg_k \colon Av = \theta A \theta v 
	$$
	is true for $k=\R$ and conclude that it therefore is also true for $k = \K$, proving (3). 
\end{proof}

\subsection{Split tori of algebraic groups over real closed fields}\label{sec:split_tori}

Let $G$ be a semi-simple connected algebraic group which is invariant under transposition. Let $\K$ and $\F$ be real closed fields with $\K \subseteq \R \cap \F$. Let $\frakg_\K \subseteq \K^{n \times n}$ be the $\K$-points of the Lie algebra. Let $\frakg_\R$ and $\frakg_\F$ be the semialgebraic extensions, then $\frakg_\R$ is also the Lie algebra of $G_\R$.

All subfields $\K$ of $\R$ are dense in $\R$, in the sense that $\overline{\K}=\R$ where $\overline{\K}$ is the closure of $\K$ in the usual topology of $\R$. The following generalization of this fact will be used in the proof of Theorem \ref{thm:split_tori}, the main result of this section.

\begin{lemma}\label{lem:closure_alg}
	Let $A\subseteq \R^n$ be an algebraic set defined over $\K$. Then we have $\overline{A_\K} = A_\R$ in the usual $\R^n$ topology.
\end{lemma}
\begin{proof}
	We first note that the algebraic set $A_\R$ is Zariski-closed and therefore also closed in the usual topology of $\R^n$, $\overline{A_\R} = A_\R$. We know that $A_\K \subseteq A_\R$ and therefore $\overline{A_\K} \subseteq \overline{A_\R} = A_\R$.
	
	On the other hand, let $x \in A_\R$. Since $\overline{\K} = \R$, there are $y_k \in \K^n$ and $\varepsilon_k >0$ with $|y_k - x|< \varepsilon_k$ and $\varepsilon_k \to 0$ as $k \to \infty$. Now we have the following first-order formula
	$$
	\varphi(y, \varepsilon ) \colon \quad \exists z\in A \colon \lVert z-y \rVert < \varepsilon.
	$$
	For every $k\in \N$, the formula $\varphi(y_k,\varepsilon_k)$ is true over $\R$, we can just take $z=x$ for every $k$. By the transfer principle it is also true over $\K$, which means that we have
	$
	z_k \in A_\K \colon \lVert z_k -y_k \rVert < \varepsilon_k.
	$
	We conclude
	$$
	\lVert z_k - x \rVert \leq \lVert z_k - y_k\rVert + \lVert y_k - x\rVert < 2 \varepsilon_k \to 0
	$$
	as $k \to \infty$, i.e. $x \in \overline{A_\K}$.
\end{proof}

By Subsection \ref{sec:ex_split}, a torus may be split with respect to some field, while not being split in a subfield. We now prove that this does not happen as long as all involved fields are real closed.

\begin{theorem}\label{thm:split_tori}
	Any maximal $\K$-split torus $S < G$ is maximal $\F$-split. Moreover, there is a maximal $\F$-split torus $S$ so that $g\tran = g$ for all $g \in S$.
\end{theorem}

\begin{proof}  We first find a maximal $\K$-split torus $T^s$ with $g\tran=g$ for all $g \in T^s$.
	
    Lemma \ref{lem:3}(3) states that there is a maximally $\K$-split Cartan subalgebra $\frakh \in \mathcal{C}(\frakg_\K)$ such that $\theta(\frakh) = \frakh$. We consider the complexification $\frakh_{\C} = \C \otimes_{\K} \frakh = \frakh_\R \oplus i \frakh_\R  \subseteq \frakg_{\C} = \frakg_\R \oplus i \frakg_\R$, which by Lemma \ref{lem:defs_Cartan_subalgebra}(ii) is also a Cartan subalgebra, in the sense that $\frakh_\C$ is abelian and $\operatorname{Nor}_{\frakg_\C}(\frakh_\C) = \frakh_\C$.
    Then
    $$
    T_\C := \operatorname{Nor}_{G_\C}(\frakh_\C)^{\circ} = \{ g \in G_\C \colon \operatorname{Ad}(g)(\frakh_\C) = \frakh_\C \}^{\circ}.
    $$
    forms a Zariski-connected closed subgroup of $G_\C$ with Lie algebra $\operatorname{Nor}_{\frakg_\C}(\frakh_\C) = \frakh_\C$. By \cite[Lemma 18.5]{Bor}, $T_\C$ is defined over $\K$. Since $\frakh_\C$ is finite dimensional, $T_\C$ can be written as the zero-set of finitely many polynomials and is therefore the $\C$-points of an algebraic group $T$.

	Since $\frakg_\C$ is semisimple and $\frakh_\C$ is a Cartan subalgebra, $\operatorname{ad}(H)$ is simultaneously diagonalizable for all $H \in \frakh_\C$ Lemma \ref{lem:defs_Cartan_subalgebra}(v). For $H \in \frakh_\C$, $\exp(H) \in T_\C$ and $\operatorname{Ad}(\exp(H)) = \exp(\operatorname{ad}(H))$ is diagonalizable. Possibly not all elements in $T_\C$ are of the form $\exp(H)$ with $H \in \frakh_\C$, but $\exp(\frakh_\C)$ is an open neighborhood of the identity and thus generates $T_\C$. For any $g\in T_\C$ there are $H_1, \ldots , H_n \in \frakh_\C$ such that
	\begin{align*}
		\operatorname{Ad}(g) &= \operatorname{Ad}\left(\prod_{i=1}^n \exp(H_i) \right)= \operatorname{Ad}\left(\exp \left(\sum_{i=1}^n H_i\right)\right) \\
		&= \exp\left( \operatorname{ad}\left( \sum_{i=1}^n H_i  \right) \right)
		=\exp\left( \sum_{i=1}^n \operatorname{ad}\left(  H_i  \right) \right),
	\end{align*}
	where we used that $\frakh_\C$ is abelian. Since $\operatorname{ad}(\frakh_\C)$ is simultaneously diagonalizable, $\operatorname{Ad}(g)$ is diagonalizable. Similarly we have $\operatorname{Ad}(g)\operatorname{Ad}(h) = \operatorname{Ad}(h) \operatorname{Ad}(g)$ for all $g , h\in T_\C$. Therefore $\operatorname{Ad}(T_\C)$ is simultaneously diagonalizable, meaning $T_\C$ is diagonalizable by \cite[Proposition 8.4]{Bor}. Since $T_\C$ is connected, $T$ is a torus by \cite[Proposition 8.5]{Bor}.
	
	If we restrict to symmetric matrices
	$$
	T^s = \{ g \in T \colon \sigma (g) = g^{-1} \}
	$$
	we see that $\operatorname{Lie}(T^s_\R) = \operatorname{Lie}(T_\R) \cap \frakp = \frakh_\R \cap \frakp$, where $\frakg_\R = \frakp \oplus \frakk$ is the Cartan decomposition corresponding to the standard Cartan involution $\theta$. 
	
	Let $\fraka \in V(\frakg_\K)$ be a maximal abelian $\K$-split subalgebra of $\frakg_\K$ with $\fraka \subseteq \frakh$. By Lemma \ref{lem:3}(4), we have $\frakh_\R \in \mathcal{C}(\frakg_\R)$ and $\fraka_\R \in V(\fraka_\R)$ with $\fraka_\R \subseteq \frakh_\R$. By Proposition \ref{prop:defs_max_noncompact}, $\frakh_\R$ is a maximal $\R$-split Cartan subalgebra and $\fraka_\R = \frakh_\R \cap \frakp$.

	Thus $\operatorname{Lie}(T^s_\R) = \frakh_\R\cap \frakp$ is maximal $\R$-split (or maximally noncompact in the terminology of Section \ref{sec:real_cartan_split}). We conclude that $T^s_\R$ is a maximal $\R$-split torus. 	
	Note that $T^s$ is defined over $\K$.

	We now want to show that $T^s$ is in fact $\K$-split. 
	We consider the relative root system $\RPhi := \Phi(G,T^s) $. 
	For any root $\alpha \colon T^s \to \G_m $, $\alpha \in ~\!\! \RPhi$, we know that $\alpha_\R(T^s _\K) \subseteq \K^\times$, since $\alpha$ is defined by the property that $\operatorname{Ad}(g)X =\alpha_\R(g)X$ for $g \in T^s_\K, X \in \frakg_\K$. 
	The graph of $\alpha$ is an algebraic set
	$$
	\operatorname{graph}(\alpha) = \{ (x,z) \in  T^s \times \mathbb{G}_m \colon z = \alpha(x) \}  
	$$
	defined over $\R$, our goal is to show that it is actually defined over $\K$. 
	
	Now $(x,z) \in \operatorname{graph}(\alpha)_\R$ if and only if $x \in T^s_\R$ and $z = \alpha_\R(x)$. In view Lemma \ref{lem:closure_alg}, $T^s_\R = \overline{T^s_\K}$, so $x \in T^s_\R$ is equivalent to saying that there is a sequence of $x_n \in T^s_\K$ such that $x_n \to x$ as $n \to \infty$. Since $\alpha_\R$ is continuous,
	we have $\alpha_\R(x_n) \to z$ as $n \to \infty$. We conclude that $(x,z) \in \operatorname{graph}(\alpha)_\R$ if and only if there is a sequence $(x_n, \alpha_\R(x_n)) \in \operatorname{graph}(\alpha)_\K$ with $(x_n, \alpha_\R(x_n)) \to (x,z)$ as $n \to \infty$, i.e. 
	$$
	\overline{\operatorname{graph}(\alpha)_\K} = \operatorname{graph}(\alpha)_\R,
	$$ 
	which means that $\operatorname{graph}(\alpha)_\K$ is dense in $\operatorname{graph}(\alpha)_\R$, in particular Zariski-dense (Zariski-open sets are also open in the usual topology and a set is dense if every open set intersects it). Viewing $\operatorname{graph}(\alpha)$ as an algebraic group, we can use \cite[Proposition 3.1.8]{Zim} to conclude that $\operatorname{graph}(\alpha)$ is defined over $\K$. Hence every $\alpha \in \RPhi$ is defined over $\K$, indeed every multiplicative character is defined over $\K$. By \cite[Corollary 8.2]{Bor}, this means that $T^s$ is $\K$-split. We have thus found a maximal $\K$-split torus which satisfies $g\tran = g $ for all $g \in T^s$.
	
	Next we will prove that $\operatorname{rank}_\F (G) = \operatorname{rank}_\K (G)$ as follows.
	\begin{align*}
		\operatorname{rank}_\F (G) &\leq r_\F (\frakg_\F) 
		= r_\R (\frakg_\R) = \operatorname{rank}_\K (G) \leq \operatorname{rank}_\F (G)
	\end{align*}
	We recall that $\operatorname{rank}_\F(G)$ is the dimension of a maximal $\F$-split torus and $r_\F(\frakg_\F)$ is the dimension of an element in $V(\frakg_\F)$. Indeed, the Lie algebra of any maximal $\F$-split torus is abelian and $\F$-split, i.e. contained in an element of $V(\frakg_\F)$. The first equality is Lemma \ref{lem:3}(1). The second equality is what we did in this proof: we found the maximal $\K$-split torus $T^s$ with $\operatorname{Lie}(T^s_\R) = \frakh_\R\cap\frakp \in V(\frakg_\R)$, with dimension $r_\R(\frakg_\R)$. 
	The last inequality holds because every $\K$-split torus is also $\F$-split.
	
	Let $S$ be a maximal $\K$-split torus. Then $S$ is also $\F$-split (but possibly not maximal). Let $T$ be a maximal $\F$-split torus with $S \subseteq T$. Then $\dim(S) = \operatorname{rank}_\K(G) = \operatorname{rank}_\F(G) = \dim (T)$ and since tori are connected $S=T$. Thus every maximal $\K$-split torus $S$ is also maximal $\F$-split.
\end{proof}

 \begin{corollary}
	Let $S < G$ be a maximal $\K$-split torus. Then the set $\KPhi$ of $\K$-roots of $S$ in $G$ coincides with $\FPhi$ and hence the set of standard parabolic $\K$-subgroups coincides with the set of standard parabolic $\F$-subgroups. In particular any $\F$-parabolic subgroup is $G_{\F}$-conjugate to a parabolic $\K$-subgroup.
\end{corollary}

\begin{proof}
	The set of $\K$-roots is defined as $\KPhi := \Phi(S,G)$, where $S$ is a maximal $\K$-split torus of $G$ \cite[21.1]{Bor}. Since $S$ is also $\F$-split by Theorem \ref{thm:split_tori}, we have $\KPhi := \Phi(S,G) = \FPhi $. 
	
	Following \cite[21.11]{Bor}, we choose an ordering on $\KPhi$ and fix the minimal parabolic $\K$-subgroup $P$ associated to the positive roots in $\KPhi$. Any parabolic $\K$-subgroup containing $P$ is called \emph{standard parabolic}. The standard parabolic $\K$-subgroups are in one-to-one correspondence with the subsets $I \subseteq \KDelta$ of the simple roots $\KDelta$ of $\KPhi$ \cite[Proposition 21.12]{Bor}. Since $\FDelta = \KDelta$, the standard parabolic $\K$-groups coincide with the standard parabolic $\F$-groups. By \cite[Proposition 21.12]{Bor}, any parabolic $\F$-group is conjugate to one and only one standard parabolic, by an element in $G_{\F}$. In particular every $\F$-parabolic subgroup is $G_{\F}$-conjugate to a $\K$-parabolic subgroup.
\end{proof}

\newpage

\section{Semialgebraic groups}\label{sec:decompositions}

Real Lie groups admit more decompositions than semisimple algebraic groups. In this section we extend some of the Lie group decompositions to algebraic groups over real closed fields. This includes the Iwasawa-decomposition $KAU$, Cartan-decomposition $KAK$, Bruhat-decomposition $BWB$ as well as Kostant-convexity, the Jacobson-Morozov Lemma and the study of rank 1 subgroups. We believe that these results are of independent interest. 

Our application of these results lies in the construction of an affine $\Lambda$-building in Section \ref{sec:building_def} as well as the verification of its axioms in Section \ref{sec:Bisbuilding}. The Cartan-decomposition is crucial in the definition of the building as well as in the verification of axiom (A3). The Iwasawa-decomposition together with Kostant convexity allows us to verify the triangle inequality. The Bruhat-decomposition is used in the verification of axioms (A2), (A4) and (EC). Jacobson-Morozov shows up in the verification of axioms (A2) and (EC).   

Let $\K$ and $\F$ be real closed fields such that $\K \subseteq \R \cap \F$. Let $G$ be a semisimple connected self-adjoint algebraic $\K$-group\footnote{The whole theory could also be formulated for a group $\tilde{G}$ that is defined as a semialgebraic set, and that lies between the semialgebraic connected component of the identity $G_\F^\circ$ of $G_\F$ and $G_\F$ itself, so $G_\F^\circ < \tilde{G} < G_\F$, see \cite{BIPP23arxiv}.} and $S$ a maximal $\K$-split torus, which we may assume to be self-adjoint by Theorem \ref{thm:split_tori}. We consider the $\K$-points $\frakg_\K \subseteq \K^{n\times n}$ of the Lie algebra $T_eG$. 
 Let $\frakg_{\mathbb{R}}$ be the semialgebraic extension of $\frakg_\K$. Then $\frakg_{\R}$ is the Lie algebra of the real Lie group $G_{\mathbb{R}}$. Since $G_{\R}$ is self-adjoint, so is $\frakg_{\R}$. The differential of the Lie group homomorphism $G_{\R} \to G_{\R}, g \mapsto (g^{-1})\tran$ is the standard Cartan-involution $\theta \colon \frakg_{\R} \to \frakg_{\R}, X \mapsto -X^{T}$ and leads to the Cartan decomposition into two eigenspaces $\frakg_{\R} = \frakk \oplus \frakp$ as in Section \ref{sec:killing_involutions_decompositions}. We define the algebraic $\K$-group 
$$
K = G \cap \operatorname{SO}_n
$$
whose real points $K_{\mathbb{R}}$ is a real Lie group with Lie algebra 
$$
\operatorname{Lie}(K_{\R}) = \operatorname{Lie}(G_{\R} \cap \operatorname{SO}_n(\R)) = \operatorname{Lie}(G_{\R}) \cap \operatorname{Lie}(\operatorname{SO}_n(\R)) = \frakg_{\R} \cap \mathfrak{so}_n
 = \frakk.
 $$

Let $A_{\K} = (S_{\K})^\circ$ be the semialgebraic connected component of $S_{\K}$ containing the identity. Let $A_{\R}$ be the semialgebraic extension of $A_{\K}$. We denote the Lie algebra of $A_{\R}$ by $\fraka$. 
Note that unlike for $G$ and $K$, there may not be an algebraic group $A$.

We note that $\fraka \subseteq \frakp$. By Theorem \ref{thm:split_tori}, $S$ is maximal $\K$-split as well as maximal $\R$-split and hence $\fraka$ is maximal abelian in $\frakp$. As in Section \ref{sec:killing_involutions_decompositions}, we can now define a root space decomposition 
$$
\frakg_{\R} = \frakg_0 \oplus \bigoplus_{\alpha \in \Sigma} \frakg_{\alpha}
$$
where $\Sigma = \{\alpha \in \fraka^{\star}\colon \alpha \neq 0, \frakg_{\alpha} \neq 0\}$ is a root system.

We note that
$$
N_{\R} := \operatorname{Nor}_{K_{\R}}(\fraka) = \left\{ k \in K_{\R} \colon kXk^{-1} \in \fraka \text{ for all } X \in \fraka \right\} 
$$
and
$$
M_{\R} := \operatorname{Cen}_{K_{\R}}(\fraka) = \left\{ k \in K_{\R} \colon kXk^{-1} = X \text{ for all }  X \in \fraka  \right\},
$$
are semialgebraic groups, since it suffices to verify the conditions for $X$ on a basis of $\fraka$.

A choice of ordered basis $\Delta \subseteq \Sigma$ gives a total order on $\Sigma$. We consider the Lie algebra
$$
\frakn = \bigoplus_{\alpha >0} \frakg_{\alpha}
$$
which is nilpotent by Lemma \ref{lem:kan_form}. Thus the exponential map and the logarithm are polynomials and we can define the Lie group
$$
U_{\R} = \{ g \in G_{\R} \colon \log(g) \in \frakn \}
$$
which is the group of $\R$-points of an algebraic $\K$-group $U$ defined the same way\footnote{In the theory of Lie groups, this group is often denoted by $N$, while $U$ is more common in the setting of algebraic groups.}. 

\subsection{Examples}

For the algebraic group $G = \operatorname{SL}_n$ with maximal $\K$-split torus
$$
S = \left\{ \begin{pmatrix}
	\star & & \\ & \ddots & \\ && \star
\end{pmatrix} \in \operatorname{SL}_n \right\}.
$$
we have
\begin{align*}
	K &= \operatorname{SO}_n\\
	A_{\F} &= \left\{ a = (a_{ij}) \in S_{\F} \colon a_{ii}>0   \right\} \\
	U_{\F} &= \left\{g=(g_{ij}) \in \operatorname{SL}_{n}(\F) \colon  g_{ii}= 1, \, g_{ij} = 0 \text{ for } i>j  \right\} \\
	N_{\F} &= \left\{ \text{ permutation matrices with entries in }\pm 1 \right\} \\
	M_{\F} & = \left\{ a=(a_{ij}) \in S_{\F} \colon a_{ii} \in \{ \pm 1\}    \right\} \\
	B_{\F} &= \left\{ g = (g_{ij}) \in \operatorname{SL}_n(\F) \colon g_{ij} = 0 \text{ for } i > j \right\}.
\end{align*}

\subsection{Compatibility of the root systems}\label{sec:compatibility}

In Chapter \ref{sec:alg_root_system} about algebraic groups, we defined the root system $\KPhi$ relative to a $\K$-split torus $S$ consisting of those $\alpha \in \hat{S}\setminus \{1\}$ for which 
$$
\frakg_{\alpha}^{(S)} = \{X \in \frakg \colon \operatorname{Ad}(s)X = \alpha(s)\cdot X \text{ for all } s \in S \}
$$
is non-zero. By Theorem \ref{thm:split_tori}, we have $\KPhi = \RPhi = \FPhi$. In this chapter we  defined a root system $\Sigma$ consisting of those $\alpha \in \fraka^{\star}\setminus \{0\}$ for which
$$
\frakg_\alpha = \{ X \in \frakg_{\R} \colon \operatorname{ad}(H)X = \alpha(H) \cdot X \text{ for all } H \in \fraka \}
$$
is non-empty. In this section we show that all these notions of root systems and the various notions of spherical Weyl groups coincide. Let us first see how to construct an algebraic character in $\KPhi$ from a root in $\Sigma$.

\begin{lemma}\label{lem:char}
	Let $\alpha \in \Sigma$. Then the homomorphism
	\begin{align*}
		\chi_\alpha \colon A_\R & \to \R_{>0} \\
		\exp(H) &\mapsto e^{\alpha(H)}
	\end{align*}
	is the restriction to $A_\R$ of an algebraic character in $\hat{S}$.
\end{lemma}
\begin{proof}
	We choose a basis of $\frakg$ consistent with the root decomposition
	$$
	\frakg = \frakg_0 \oplus \bigoplus_{\alpha \in \Sigma} \frakg_\alpha.
	$$
	In this basis, the matrices in $\operatorname{Ad}(S_\R)$ and $\operatorname{ad}(\fraka)$ are diagonal. In fact any diagonal entry\footnote{Note that $\frakg_\alpha$ may not be one-dimensional, but then all the diagonal entries of $\operatorname{ad}(H)$ corresponding to basis vectors in $\frakg_\alpha$ are equal, so it makes sense to talk about $(\operatorname{ad}(H))_{\alpha \alpha}$.} corresponding to the root $\alpha \in \Sigma$ satisfies
	\begin{align*}
		(	\operatorname{ad}(H))_{\alpha \alpha} & = \alpha(H)
	\end{align*}
	for all $H \in \fraka$. We define a character
	\begin{align*}
		\chi (a) &:= (\operatorname{Ad}(a) )_{\alpha \alpha} \quad \quad \text{for } a \in S
	\end{align*}
	which is algebraic by its definition. Let $\chi_\R \colon S_\R \to \R$ be the $\R$-points of $\chi$ and $\chi_\R|_{A_\R}$ its restriction to $A_\R$. We now claim that $\chi_\alpha = \chi_\R|_{A_{\R}}$: when $a = \exp(H) \in A_{\R}$ for $H \in \fraka$, we have 
	\begin{align*}
		\chi_\R (a) &= (\operatorname{Ad}(a))_{\alpha \alpha} 
		= (\operatorname{Ad}(\exp(H)))_{\alpha \alpha} \nonumber \\
		&= (\exp (\operatorname{ad}(H)))_{\alpha \alpha}
		= e^{\operatorname{ad}(H)_{\alpha \alpha}} = e^{\alpha(H)} = \chi_\alpha(a).
	\end{align*} 
\end{proof}
\begin{lemma}\label{lem:oneparam_compatibility}
	Let $\alpha \in \Sigma$. Then there is an $x_\alpha \in \fraka$ such that the homomorphism
	\begin{align*}
		t_\alpha^\R \colon \R_{>0} &\to A_\R \\
		e^s &\mapsto \exp(sx_\alpha)
	\end{align*}
	is the restriction of an algebraic one-parameter group $t_\alpha$ in $X_\star(S)$, defined over $\K$. The non-degenerate bilinear form in Proposition \ref{prop:char_cochar} is then given by
	$$
	b( \chi_\alpha, t_\beta ) = \frac{2\langle \alpha, \beta \rangle}{\langle \beta,\beta\rangle} \in \Z.
	$$
\end{lemma}
\begin{proof}
 Recall that the Killing form gives rise to a scalar product $B_\theta$ on $\fraka = \operatorname{Lie}(A_\R)$, inducing an isomorphism $\fraka^\star \cong \fraka, \gamma \mapsto H_\gamma$ with the defining property $B_\theta(H_\gamma,H) = \gamma(H)$ for all $H \in \fraka$. The coroot
 $$
 \alpha^\vee = \frac{2}{B_{\theta}(H_\alpha,H_\alpha)}\alpha
 $$
 can be used to define
 $$
 x_\alpha := H_{\alpha^\vee} = \frac{2}{B_{\theta}(H_\alpha,H_\alpha)}H_\alpha.
 $$
 We now define $t_\alpha^\R(e^s) = \exp(sx_\alpha)$ and note that
 $$
 \chi_\alpha^\R(t_\alpha^\R(e^s)) = e^{s\alpha(x_\alpha)} = (e^s)^2.
 $$
 We have to show that $t_\alpha^\R$ is algebraic. Let $\mathcal{A}(H_\alpha)$ be the smallest algebraic subgroup whose Lie algebra contains $H_\alpha$, \cite[7.1]{Bor}. Since $H_\alpha \in \K^{n\times n}$, $\mathcal{A}(H_\alpha)$ is defined over $\K$, \cite[2.1(b)]{Bor}. Note that for every $r \in \G_m$ there is a unique element $a \in \mathcal{A}(H_\alpha)$ with $\chi_\alpha(a) = r^2$. We define the map $t_\alpha \colon \G_m \to S$ by its graph
 $$
 \operatorname{graph}(t_\alpha) =  \left\{ (r,a) \in \G_m \times S \colon a \in \mathcal{A}(H_\alpha), \chi_\alpha(a) = r^2 \right\}.
 $$
 By Proposition \ref{prop:morphisms_polynomials}, $t_\alpha$ is algebraic and defined over $\K$. By the defining property of the bilinear form $b \colon \hat{S} \times X_\star(S) \to \Z$ from Proposition \ref{prop:char_cochar}, we have $\chi_\alpha (t_\beta(r)) = r^{b(\chi_\alpha, t_\beta)}$ for all $\alpha, \beta \in \Sigma$. Over $\R$ we have for $r=e^s \in \R_{>0}$
 \begin{align*}
 	\chi_\alpha^\R (t_\beta^\R(r)) = \chi_\alpha^\R (\exp(s x_\alpha)) = e^{s\alpha(x_\beta)} = r^{\frac{2\langle \alpha, \beta\rangle}{\langle \beta, \beta \rangle}},
 \end{align*}
 as claimed.
\end{proof}

\begin{proposition}\label{prop:rootsystem_compatible}
	The root systems $\KPhi$ and $\Sigma$ are isomorphic. The spherical Weyl groups $\KW$ and $N_\R / M_\R$ are isomorphic.
\end{proposition}
\begin{proof}
   We first find a compatible map $\KPhi \to \Sigma$. Let $\chi \in \KPhi$. Since $S$ is $\K$-split, $\chi$ is defined over $\K$ \cite[Corollary 8.2]{Bor}. Hence we can consider the Lie group homomorphism $\chi_{\R} \colon S_{\R} \to \R^{\times}$ and its derivative $\operatorname{d}_e\!\chi_{\R} \colon \fraka \to  \R$, which we claim to be an element in $\Sigma$. We first claim that $\operatorname{d}_e\!\chi_{\R}$ is nonzero, since otherwise $\chi_{\R}$ would be locally constant, hence only take finitely many values. But since $\chi_{R}(S_{\R})$ is Zariski dense in $\chi(S) = \G_m$ \cite[Corollary 18.3]{Bor}, this cannot be the case. Next, we claim that 
   $$
   \frakg_{\chi} := \left(\frakg_{\chi}^{(S)}\right)_{\R} \subseteq \frakg_{\operatorname{d}_e\!\chi_{\R}},
   $$
   which shows that $\KPhi \to \Sigma$ is well defined. Let $H \in \fraka$, such that $a=\exp(H) \in A_{\R} \subseteq S_{\R}$. Since $S$ is $\K$-split, we have 
   $$
   \operatorname{Ad}(\exp(H)) = \operatorname{Id} \cdot \alpha_{\R}(\exp(H))
   $$ 
   as maps $\frakg_{\chi} \to \frakg_{\chi}$. In fact, $\operatorname{Ad}\circ \exp = \operatorname{Id} \cdot (\chi_{\R}\circ\exp)$ as maps $\fraka \to \operatorname{GL}(\frakg_{\chi})$. Taking the derivative at $0 \in \fraka$, we get
   \begin{align*}
\operatorname{ad} &= \operatorname{ad} \circ \operatorname{Id} =\operatorname{d}_e\!\operatorname{Ad} \circ \operatorname{d}_0 \exp =  \operatorname{d}_0 ( \operatorname{Ad} \circ \exp ) \\
&= \operatorname{d}_0 (\operatorname{Id} \cdot (\chi_{\R}\circ\exp)) = \operatorname{Id}\cdot \operatorname{d}_e\!\chi_{\R} \circ \operatorname{d}_0\exp = \operatorname{Id}\cdot  \operatorname{d}_e\!\chi_{\R}
   \end{align*}
   as maps $\fraka \to \mathfrak{gl}(\frakg_{\chi})$. This means that for all $X \in \frakg_{\chi}$ and all $H \in \fraka$, we have $\operatorname{ad}(H)X = \operatorname{d}_e\! \chi_{\R}(H)\cdot X$ and hence $X \in \frakg_{\operatorname{d}_e\!\chi_{\R}}$.
   
   Now, we will show that the map $\alpha \mapsto \chi_\alpha$ in Lemma \ref{lem:char}, $\chi_\alpha$ now viewed as an algebraic character in $\hat{S}$, actually sends $\Sigma \to \KPhi$. Since
   $$
   (\chi_\alpha)_\R(\exp(H)) = e^{\alpha(H)},
   $$
   we can easily see that $\chi_\alpha \neq 1$, when $\alpha \neq 0$. Next we will show that
   $$
   \frakg_\alpha \subseteq \left(\frakg_{\chi_\alpha}^{(S)} \right)_{\R}.
   $$ 
   For $a=\exp(H) \in A_{\R}$, we see from the above formula that
   \begin{align*}
   	\operatorname{Ad}(a) & = \exp(\operatorname{ad}(H)) = \exp( \operatorname{Id} \cdot \alpha(H) ) 
   	= \operatorname{Id} \cdot e^{\alpha(H)} = \operatorname{Id} \cdot \chi_\alpha(a) 
   \end{align*}
   as functions on $\frakg_\alpha$. Since $A_{\R}$ is Zariski dense in $S$, $\operatorname{Ad}(s) X = \chi_\alpha (a) \cdot X$ for all $s \in S$ and $X \in \frakg_{\alpha}$.
   
   We note that the two maps between $\KPhi$ and $\Sigma$ are inverses of each other, which follows from the fact that
   $$
   \left(\frakg_{\alpha}^{(S)}\right)_{\R} \subseteq \frakg_{\operatorname{d}_e\!\alpha_{\R}} \subseteq 
   \left(\frakg_{\chi_{\operatorname{d}_e\!\alpha_\R}}^{(S)}\right)_{\R} \quad \text{ for } \alpha \in \KPhi
   $$
   and 
   $$
   \frakg_\alpha \subseteq \left(\frakg_{\chi_\alpha}^{(S)}\right)_{\R} \subseteq \frakg_{\operatorname{d}_e\!(\chi_{\alpha})_{\R}} \quad \text{ for } \alpha \in \Sigma. 
   $$
   We note that since $[\frakg_{\alpha},\frakg_\beta] \subseteq \frakg_{\alpha + \beta}$ the map $\KPhi \to \Sigma$ extends to an isomorphism on their vector spaces. For both $\KPhi$ and $\Sigma$, the Weyl groups $\KW$ and $W_s=N_{\R}/M_{\R}$ are generated by reflections along hyperplanes perpendicular to the roots \cite[Theorem 21.2]{Bor} and \cite[Proposition 7.32]{Kna}. This implies that the scalar products of two roots are preserved under the map $\KPhi \to \Sigma$ and hence the two root systems are isomorphic.
\end{proof}

We can also conclude that the various Weyl groups that can be defined coincide.

\begin{proposition} \label{prop:weylgroups}
	The following definitions of Weyl groups are isomorphic.
	\begin{itemize}
		\item [(i)] The Weyl group $W$ generated by reflections in roots of the root system $\KPhi$.
		\item [(ii)] The Weyl group $\RW = \operatorname{Nor}_K(S)/\operatorname{Cen}_K(S)$ from the theory of algebraic groups.
		\item [(iii)] The Weyl group $W_s$ generated by reflections in roots of the root system $\Sigma$.
		\item [(iv)] The Weyl group in the Lie groups setting $W(G_{\R},A_{\R}) = N_{\R}/M_{\R} $.
		\item [(v)] The Weyl group $N_{\F}/M_{\F}$ of the semialgebraic extensions.
	\end{itemize}
\end{proposition}
\begin{proof}
	By Proposition \ref{prop:rootsystem_compatible}, $\KPhi \cong \Sigma$, so (i) and (iii) coincide. 
	By \cite[Theorem 21.2]{Bor}, $W=\RW $, so the notions (i) and (ii) agree. By \cite[Proposition 7.32]{Kna}, $W_s = W(G_{\R},A_{\R})$, so the notions (iii) and (iv) agree. Let $W_s = \{w_1, \ldots , w_{|W_s|}\}$ be a list of the finitely many elements in $W_s$ considered as an abstract group. The first-order formula
	\begin{align*}
		\varphi \colon \quad & \exists n_1 , \ldots , n_{|W_s|} \in N \colon 
		\left( \bigwedge_{i,j=0 \atop i\neq j}^{|W_s|} \neg \ n_j^{-1}n_i \in  M \right)  \wedge \\
		&\left(\forall n \in N  \colon \bigvee_{i=1}^{|W_s|}  n^{-1}n_i \in M  \right) \wedge  
		\left( \bigwedge_{i,j,\ell=0 \atop w_iw_j=w_{\ell}}^{|W_s|} n_\ell^{-1}n_in_j \in M \right)
	\end{align*}
	states that there are $|W_s|$ many elements in $N$ with distinct representatives in $N/M$ and there are only $|W_s|$ many and they satisfy the same group multiplication able as $W_s$. In short it says that $N/M$ is isomorphic to $W_s$. Now, by (iv) the formula $\varphi$ holds over $\R$ and we can apply the transfer principle to get the statement over $\F$, showing that the notion in (v) gives the same Weyl group.
\end{proof}

\subsection{Iwasawa decomposition $G=KAU$}\label{sec:KAU}

We recall the classical Iwasawa-decomposition, which applies to $G_{\R}$.
\begin{theorem}\label{thm:KAU_R}(\cite[Theorem 6.46]{Kna})
	Let $G_\R$ be a semisimple Lie group. Let $\frakg_\R = \frakk \oplus \fraka \oplus \frakn$ be an Iwasawa decomposition of its Lie algebra. Let $A_\R$ and $U_\R$ be the analytic subgroups with Lie algebras $\fraka$ and $\frakn$. Then the multiplication map $K_\R \times A_\R \times U_\R \to G_\R$ is a diffeomorphism. This decomposition is unique.
\end{theorem}

We note that all the groups $G_\R, K_\R, A_\R$ and $U_\R$ are semialgebraic defined over $\K$, hence we can consider $G_\F, K_\F, A_\F$ and $U_\F \subseteq \F^{n\times n}$ for any real closed field $\F$ containing $\K$. We use the transfer principle to deduce the following semialgebraic version of the Iwasawa decomposition. We remark that by taking inverses, we also have a decomposition $G=UAK$ in addition to $G=KAU$, both of which will be called the \emph{Iwasawa decomposition}.
\begin{theorem}[$G=KAU$]\label{thm:KAU}
	For every $g \in G_{\F}$, there are $k\in K_{\F}, a\in A_{\F}, u \in U_{\F}$ such that $g=kau$. This decomposition is unique. 
\end{theorem} 
\begin{proof}
	We write the decomposition as a first-order formula
	\begin{align*}
	\varphi \colon \quad & \left(\forall g \in G \ \exists k \in K , a \in A, u\in U \colon g= kau \right)\ \wedge \\
	&\left(\forall k,k' \in K, a,a'\in A, u,u'\in U \colon kau = k'a'u' \to \left(k=k' \wedge a=a' \wedge u = u'\right) \right)	
	\end{align*}
	which holds over $\R$ by the classical Iwasawa decomposition of Theorem \ref{thm:KAU_R}.
	By the transfer principle, Theorem \ref{thm:logic}, $\varphi$ holds over all real closed fields $\F$. Note that to apply the classical Iwasawa-theorem we use that $K_\R$ is maximal compact, $\operatorname{Lie}(A_\R) = \fraka$ and $\operatorname{Lie}(U_\R)= \frakn$.
\end{proof}

\subsection{Cartan decomposition $G=KAK$}\label{sec:KAK}

We use the Cartan decomposition for real Lie groups to find an analogue statement over $\F$.

\begin{theorem}(\cite[Theorem 7.39]{Kna})\label{thm:KAK_R}
	Every element $g\in G_{\R}$ has a decomposition $g = k_1 a k_2$ with $k_1,k_2 \in K_{\R}$ and $a \in A_{\R}$. In this decomposition, $a$ is uniquely determined up to a conjugation by a member of $W(G_{\R},A_{\R})$.
\end{theorem}
\begin{proof}
	By Proposition \ref{prop:weylgroups},
	$\RW = \operatorname{Nor}_G(S)/\operatorname{Cen}_G(S)$ is isomorphic to $W(G_{\R},A_{\R}) = \operatorname{Nor}_{K_{\R}}(\fraka)/\operatorname{Cen}_{K_{\R}}(\fraka)$. Then this is the statement of \cite[Theorem 7.39]{Kna}. 
\end{proof}

\begin{theorem}[$G=KAK$]\label{thm:KAK}
	Every element $g\in G_\F$ has a decomposition $g = k_1 a k_2$ with $k_1,k_2 \in K_{\F}$ and $a \in A_{\F}$. In this decomposition, $a$ is uniquely determined up to a conjugation by a member of $N_{\F}/M_{\F}$.
\end{theorem}
\begin{proof}
	The existence part of the statement follows from the first-order formula
	$$
	\varphi \colon \quad \forall g \in G\  \exists k_1 ,k_2 \in K\ \exists a \in A \colon g =k_1ak_2,
	$$
	which holds over $\R$ by Theorem \ref{thm:KAK_R} and hence over $\F$ by the transfer principle. For uniqueness of $a$, we consider the first-order logic formula
	\begin{align*}
		\psi \colon \quad & \forall a,a' \in A, k_1,k_1',k_2,k_2' \in K \colon k_1ak_2 = k_1'a'k_2' \  \to \\
		&\left( \exists n \in N \colon a = na'n^{-1} \ \wedge \ \left(\forall n' \in N \colon a=n'a'(n')^{-1} \to n^{-1}n' \in M \right) \right)
	\end{align*}
	which states that $a$ is determined up to a conjugation by a member of $N$ and that this member is unique up to an element of $M$. Over $\R$, $\psi$ holds by Theorem \ref{thm:KAK_R}, hence $\psi$ also holds over $\F$ by the transfer principle, concluding the proof.
\end{proof}

\subsection{Bruhat decomposition $G=BWB$}\label{sec:BWB}

By \cite[page 398]{Kna}, $B_{\R} := M_{\R}A_{\R}U_{\R}$ is a closed subgroup of $G_{\R}$ and we have the following Bruhat decomposition.

\begin{theorem}(\cite[Theorem 7.40]{Kna})\label{thm:BWB_R} Every element $g \in G_{\R}$ can be written as $g=b_1nb_2$ with $b_1,b_2 \in B_{\R}$ and $n \in N_{\R}$. In this decomposition, $n$ is unique up to multiplying by an element in $M_{\R}$. Since the spherical Weyl group is $W_s=N_{\R}/M_{\R}$, we have a disjoint union of double cosets
	$$
	G_{\R} = \prod_{[n] \in W_s} B_{\R}nB_{\R}.
	$$ 
\end{theorem}

The group $B_{\R}$ is semialgebraic and the Bruhat decomposition can be extended to $G_{\F}$.

\begin{theorem}[$G=BWB$]\label{thm:BWB}
	 Every element $g \in G_{\F}$ can be written as $g=b_1nb_2$ with $b_1,b_2 \in B_{\F}$ and $n \in N_{\F}$. In this decomposition $n$ is unique up to multiplying by an element in $M_{\F}$. We have a disjoint union of double cosets
	$$
	G_{\F} = \prod_{[n] \in W_s} B_{\F}nB_{\F}.
	$$ 
\end{theorem}
\begin{proof}
	The existence of the decomposition follows directly from the transfer principle and Theorem \ref{thm:BWB_R}. For uniqueness we utilize the first-order formula
	\begin{align*}
		\varphi \colon \quad &  \forall b_1,b_2,b_1',b_2' \in B, n,n' \in N \\
		& b_1nb_2 = b_1'n'b_2' \to n^{-1}n' \in M  ,
	\end{align*}
	which holds over $\R$ by Theorem \ref{thm:BWB_R}. Since $W_s= N_\F/M_\F =N_\R/M_\R$ by Proposition \ref{prop:weylgroups}, $G_\F$ is a disjoint union as described. 
\end{proof}

Let 
$$
\ApF = \left\{ a \in A_\F \colon \chi_\alpha(a) \geq 1 \text{ for all } \alpha \in \Delta \right\},
$$
where $\chi_\alpha$ is the algebraic character associated to $\alpha \in \Delta$. Then we may choose $a \in \ApF$ in the Cartan decomposition $G_\F = K_\F \ApF K_\F$ which is then unique, since the Weyl group acts simply transitively on the set of Weyl chambers.

\subsection{Baker-Campbell-Hausdorff formula}\label{sec:BCH}

Given $X,Y \in \frakg_\R$, the Baker-Campbell-Hausdorff formula gives a formal power series description of $Z = \log(\exp(X)\exp(Y))$ in terms of $X,Y$ and iterated commutators of $X$ and $Y$, see for instance \cite[Proposition V.1]{Jac79} or \cite{Mag54}. The formal power series converges in a neighborhood of the identity \cite[Theorem 3.1 in X.3]{Hoch65}, 
but may not converge everywhere in general. If $X,Y \in \frakn_\R$, then the power series is given by a polynomial and thus converges everywhere \cite{Wei63}. Since only finitely many terms are involved, one can see directly or invoke the transfer principle to obtain the Baker-Campbell-Hausdorff formula for elements in $\frakn_\F := \bigoplus_{\alpha >0} (\frakg_\alpha)_\F$.
\begin{proposition}\label{prop:BCH}
	For every $X,Y \in \frakn_\F$, there is $Z \in \frakn_\F$ such that $\exp(X)\exp(Y) = \exp(Z)$. The element $Z$ is given by a finite sum of iterated commutators, the first terms of which are given by
	$$
	Z = X + Y + \frac{1}{2}[X,Y] + \frac{1}{12}\left(\left[X,\left[X,Y\right]\right]- \left[Y, \left[Y,X\right]\right]\right) - \frac{1}{24} \left[Y,\left[X,\left[X,Y\right]\right]\right] + \ldots
	$$ 
\end{proposition}
There are various variations in the literature. We will make use of the Zassenhaus formula
$$
\exp(X+Y) = \exp(X) \exp(Y) \exp\left(-\frac{1}{2}\left[X,Y\right]\right)\exp \left( \frac{1}{3} [Y,[X,Y]] + \frac{1}{6}[X,[X,Y]] \right) \cdots 
$$
for $X,Y \in \frakn_\F$, which can be obtained from the above by calculating $\exp(-X)\exp(X+Y) = \exp(Z) $ iteratively \cite{Mag54}.

\subsection{The unipotent group $U$}\label{sec:U}

The root spaces $\frakg_\alpha$ in the root decomposition
$$
\frakg = \frakg_0 \oplus \bigoplus_{\alpha \in \Sigma} \frakg_\alpha.
$$
corresponding to positive roots, consist of simultaneously nilpotent elements by Lemma \ref{lem:kan_form}. Therefore the group 
$$
U = \exp(\frakn) = \exp\left( \bigoplus_{\alpha \in \Sigma_{>0}} \frakg_\alpha \right)
$$
is unipotent.

Let $\alpha \in \Sigma$. When $2\alpha \notin \Sigma$, the root space $\frakg_\alpha$ is an ideal. In any case, $\frakg_\alpha \oplus \frakg_{2\alpha}$ is an ideal, where $\frakg_{2\alpha} = 0$ if $2\alpha \notin \Sigma$. Thus, for every $\alpha \in \Sigma$ there is a unipotent subgroup $U_\alpha = \exp(\frakg_\alpha \oplus \frakg_{2\alpha}) < U$, called the \emph{root group} which is also an algebraic group, since the exponential function is a polynomial on nilpotent elements. Note that $U_{2\alpha} < U_\alpha$ if both exist. 

The following Lemma shows that $A_\F$ normalizes the root groups $(U_\alpha)_\F$.

\begin{lemma} \label{lem:aexpXa}
	Let $\alpha \in \Sigma, X\in \mathfrak{g}_\alpha, X' \in \mathfrak{g}_{2\alpha}$ and $a \in A_\F $. Then 
	$$a\exp\left(X + X'\right)a^{-1} = \exp\left(\chi_\alpha (a) X + \chi_\alpha(a)^2 X'\right),$$
	where $\chi_\alpha \colon A_\F \to \mathbb{F}_{>0}$ is the algebraic character from Lemma \ref{lem:char}.
\end{lemma}
\begin{proof}
	We know 
	that $\operatorname{Ad}_a(X) = \chi_\alpha(a) X$ and $\operatorname{Ad}_a(X') = \chi_{2\alpha}(a) X' = \chi_\alpha(a)^2 X'$.
	We then use distributivity of matrix multiplication to obtain
	\begin{align*}
		a\exp(X+X')a^{-1} &= a\sum_{n=0}^\infty \frac{(X+X')^n}{n!}a^{-1} 
		=\sum_{n=0}^\infty \frac{1}{n!} a(X+X')^n a^{-1} \\
		& =\sum_{n=0}^\infty \frac{1}{n!} \left(a(X+X')a^{-1}\right)^n 
		= \exp \left(aXa^{-1} + aX'a^{-1}\right) \\
		&= \exp \left(\operatorname{Ad}_a(X) + \operatorname{Ad}_a(X')\right)
		= \exp \left(\chi_\alpha(a)X + \chi_\alpha(a)^2 X'\right)
	\end{align*}
\end{proof}
We consider $\Theta \subseteq \Sigma_{>0}$ \emph{closed under addition}, meaning that for any $\alpha, \beta \in \Theta$, if the sum $\alpha + \beta \in \Sigma$, then $\alpha + \beta \in \Theta$. For any $\Theta \subseteq \Sigma_{>0}$ closed under addition,
$$
\frakg_\Theta := \bigoplus_{\alpha \in \Theta} \frakg_\alpha 
$$
is an ideal, since $\left[ \frakg_\alpha , \frakg_\beta \right] \subseteq \frakg_{\alpha + \beta}$ for any $\alpha, \beta \in \Sigma$, see Proposition \ref{prop:root_decomp}. Hence 
$$
U_\Theta :=\exp(\frakg_\Theta)
$$
is a real algebraic group, in fact for $\Theta = \Sigma_{>0}$ we recover $U = U_{\Sigma_{>0}}$. As a consequence of the Baker-Campbell-Hausdorff-formula we obtain the following description of $(U_\Theta)_\F$.

\begin{lemma}\label{lem:BCH_consequence}
	Let $\Theta = \left\{ \alpha_1, \ldots , \alpha_k \right\} \subseteq \Sigma_{>0}$ be a subset closed under addition with $\alpha_1 > \ldots > \alpha_k$. Then
	$$	
	(U_{\Theta})_\F := \exp\left( \bigoplus_{\alpha \in \Theta} (\frakg_\alpha)_\F\right)
	= \prod_{i = 1}^k \exp\left(\left(\frakg_{\alpha_i}\right)_\F\right) = \langle u \in U_{\alpha} \colon \alpha \in \Theta \rangle.
	$$
\end{lemma}
\begin{proof}
	We start by proving
	$$
	\exp\left( \bigoplus_{\alpha \in \Theta} (\frakg_\alpha)_\F\right) \subseteq  \prod_{i = 1}^k \exp\left(\left(\frakg_{\alpha_i}\right)_\F\right) 
	$$
    using induction over $|\Theta|$. If $|\Theta| =1$, then the statement is immediate. Now assume that $\Theta = \{\alpha_1\} \cup \Theta'$ with $\alpha_1 >\beta$ for all $\beta \in \Theta'$ and such that
	$$	
	\exp\left( \bigoplus_{\alpha \in \Theta'} (\frakg_\alpha)_\F\right)
	\subseteq \prod_{i = 2}^k \exp\left(\left(\frakg_{\alpha_i}\right)_\F\right) .
	$$
	Let $X_i \in (\frakg_{\alpha_i})_\F$. 
Recall that $[\frakg_\alpha, \frakg_{\beta}] \subseteq \frakg_{\alpha + \beta}$ for all $\alpha, \beta \in \Sigma$, see Proposition \ref{prop:root_decomp}. Thus we can use the Zassenhaus-formula in Section \ref{sec:BCH} about the Baker-Campbell-Hausdorff formula to obtain the finite product
	\begin{align*}
		\exp\left( \sum_{i=1}^k X_i \right) &= \exp(X_1)\exp\left(\sum_{i=2}^k X_i\right) \exp\left( -\frac{1}{2}\left[ X_1, \sum_{i=2}^k X_i  \right] \right)  \ldots 
	\end{align*}
	and then repeatedly apply the original version in Proposition \ref{prop:BCH} to simplify the expression back to
	\begin{align*}
		\exp(X_1)\exp\left(\sum_{i=2}^{k} \tilde{X}_i\right) 
	\end{align*}
	for some new $\tilde{X}_i \in (\frakg_{\alpha_i})_\F$.  Applying the induction hypothesis, we conclude
	\begin{align*}
		\exp\left( \sum_{i=1}^k X_i \right) &\in  \exp(X_1) \prod_{i=2}^k \exp\left(\left(\frakg_{\alpha_i}\right)_\F\right) \subseteq \prod_{i = 1}^k \exp\left(\left(\frakg_{\alpha_i}\right)_\F\right).
	\end{align*}
	Next, we notice that $\prod_{i=1}^k \exp((\frakg_{\alpha_i})_\F) \subseteq \langle u \in (U_{\alpha})_\F \colon \alpha \in \Theta\rangle$. Finally, we prove the inclusion
	$$
 \langle u \in U_{\alpha} \colon \alpha \in \Theta \rangle \subseteq \exp\left( \bigoplus_{\alpha \in \Theta} (\frakg_\alpha)_\F\right)
	$$
	using induction over the word length of an element in $\langle u \in U_{\alpha} \colon \alpha \in \Theta \rangle $. If the word length is $1$, then the statement holds. Now assume that 
	$$
	v =\exp\left( \sum_{i=1}^k X_i\right) 
	$$
	for some $X_i \in (\frakg_{\alpha_i})_\F$ and consider $u=\exp(X)$ for some $X\in (\frakg_{\alpha})_\F$ and some $\alpha \in \Theta$. We apply Proposition \ref{prop:BCH} to obtain
	\begin{align*}
		uv &= \exp(X) \exp\left(\sum_{i=1}^kX_i\right) \\
		&= \exp\left( X + \sum_{i=1}^k X_i + \frac{1}{2}\left[X, \sum_{i=1}^k X_i\right] + \ldots \right) \in \exp\left( \bigoplus_{\alpha \in \Theta} (\frakg_\alpha)_\F\right)
	\end{align*}
	concluding the proof.
\end{proof}
We point out that the order of the product expression in Lemma \ref{lem:BCH_consequence} starts with $\exp ((\frakg_{\alpha})_\F)$ corresponding to the largest $\alpha = \alpha_1$ followed in decreasing order. Writing elements of $U_{\Theta}$ as the inverses of elements in $U_\Theta$, also gives an expression starting with the smallest root, followed by an increasing order. The following technical Lemma will be used in the proof of the mixed Iwasawa decomposition, Theorem \ref{thm:BT_mixed_Iwasawa}.

\begin{lemma}\label{lem:BCH_normalizer}
	Let $\Theta \subseteq \Sigma_{>0}$ be a subset closed under addition and $\alpha \in \Sigma_{>0}$ such that $\alpha > \beta$ for all $\beta \in \Theta$. Then $uU_{\Theta}u^{-1}= U_{\Theta}$ for all $u \in \exp((\frakg_{\alpha})_\F)$.
	
	For every subset $\Psi \subseteq \Sigma_{>0}$, elements $u \in U_{\Theta}$ can be expressed as $u = u' u''$ with
	$$
	u' \in \prod_{\alpha \in \Theta \cap \Psi} \exp((\frakg_{\alpha})_\F) \quad \text{ and } \quad u'' \in \prod_{\alpha \in \Theta \setminus \Psi} \exp((\frakg_{\alpha})_\F) .
	$$ 
\end{lemma}
\begin{proof}
	Let $X \in (\frakg_{\alpha})_\F$ and $\exp(Y) \in U_{\Theta}$. Then we can apply the Baker-Campbell-Hausdorff formula, Proposition \ref{prop:BCH}, twice to obtain
	\begin{align*}
		\exp(X) \exp(Y) \exp(X)^{-1} &= \exp\left(X + Y + \frac{1}{2}[X,Y] + \ldots\right) \exp(-X)  \\
		&= \exp\left( Y + \frac{1}{2} [X,Y] + \ldots\right) \in U_\Theta.
	\end{align*}
	We show the second statement using induction over the size of $\Theta$. If $|\Theta| = 1$, the statement is clear. Now consider the subset closed under addition $\Theta' := \Theta \setminus \{\alpha\} $ where $\alpha$ is the largest element of $\Theta$. For $u \in U_\Theta$, we use Lemma \ref{lem:BCH_consequence} to obtain $u = u_\alpha \bar{u}$ with $\bar{u} \in U_{\Theta'}$. If $\alpha \notin \Psi$, then the first part of this Lemma can be used to instead write $u = \bar{u} u_\alpha$. Either way, the induction assumption gives 
	$$
		\bar{u}' \in \prod_{\beta \in \Theta' \cap \Psi} \exp((\frakg_{\beta})_\F) \quad \text{ and } \quad \bar{u}'' \in \prod_{\beta \in \Theta' \setminus \Psi} \exp((\frakg_{\beta})_\F)
	$$
	such that $\bar{u} = \bar{u}' \bar{u}''$. Then $u = u_{\alpha} \bar{u}' \bar{u}''$ or $u =  \bar{u}' \bar{u}'' u_\alpha $ as required.
\end{proof}

Lemmas \ref{lem:aexpXa} and \ref{lem:BCH_consequence} can be used to prove that $A_\F$ normalizes all of $U_\F$. 

\begin{proposition} \label{prop:anainN}
	Let $\Theta \subseteq \Sigma_{>0}$ closed under addition. Then $A_\F$ normalizes $(U_\Theta)_\F$: for all $a \in A_\F, u \in (U_\Theta)_\F \colon aua^{-1} \in (U_\Theta)_\F$.
	In particular $aU_\F a^{-1} = U_\F$ for all $a \in A_\F$.
\end{proposition}
\begin{proof}
	Let $a \in A_\F$. By Lemma \ref{lem:BCH_consequence}, any $u\in (U_\Theta)_\F$ can be written as $u = u_{\alpha_1} \cdot \ldots \cdot u_{\alpha_k}$ with $u_{\alpha_i} \in (U_{\alpha_i})_\F$ where $\Theta = \{\alpha_1, \ldots , \alpha_k\}$. By Lemma \ref{lem:aexpXa}, $au_{\alpha_i}a^{-1} \in (U_{\alpha_i})_\F$, so
	$$
	aua^{-1} = au_{\alpha_1}a^{1} \cdot \ldots \cdot au_{\alpha_k}a^{-1} \in \prod_{i=1}^k (U_{\alpha_i})_\F = (U_\Theta)_\F,
	$$
	where we used Lemma \ref{lem:BCH_consequence} again.
\end{proof}

\subsection{Jacobson-Morozov Lemma for algebraic groups}\label{sec:Jacobson_Morozov}

The Jacobson-Morozov Lemma on the level of algebraic groups seems to be folklore. Since we could not find a detailed treatment in the literature we will give its statement and proof here. The classical Jacobson-Morozov Lemma applies to semisimple Lie algebras over fields of characteristic $0$. The following formulation is more general than what we proved in Lemma \ref{lem:JM_basic}.
\begin{proposition}\label{prop:Jacobson_Morozov_algebra}
	\cite[VIII \S 11.2 Prop. 2]{Bou08_789} Let $\frakg$ be a semisimple Lie algebra and $x\in \frakg$ a non-zero nilpotent element. Then there exist $h,y \in \frakg$ such that $(x,y,h)$ is an $\mathfrak{sl}_2$-triplet, meaning that 
	$$
	[h,x] = 2x, \quad [h,y] = -2y, \quad [x,y] = h
	$$
	and hence the Lie algebra generated by $x,h,y$ is isomorphic to $\mathfrak{sl}_2$.  
\end{proposition}

We will use \cite[Chapter 7]{Bor} and \cite{PR93} to show the following version of the Jacobson-Morozov Lemma for algebraic groups over an algebraically closed field $\D$ of characteristic $0$. Note that any semisimple linear algebraic group with Lie algebra $\mathfrak{sl}_2$ is isomorphic to either $\operatorname{SL}_2$ or $\operatorname{PGL}_2$ \cite[Corollary 32.2]{Hum1}. 

\begin{proposition}\label{prop:Jacobson_Morozov_group}
	Let $g\in G$ be any unipotent element in a semisimple linear algebraic group $G$ over an algebraically closed field $\D$ of characteristic $0$. Then there is an algebraic subgroup $\operatorname{SL}_g < G$ with Lie algebra $\operatorname{Lie}(\operatorname{SL}_g) \cong \mathfrak{sl}_2$ and $g \in \operatorname{SL}_g$. The element $\log(g) \in \operatorname{Lie}(G)$ corresponds to 
	$$
	\begin{pmatrix}
		0 & 1 \\ 0 & 0
	\end{pmatrix}\in \mathfrak{sl}_2.
	$$
	Moreover, if $g\in G_\F$ for a field $\F \subseteq \D$, then $\operatorname{SL}_g$ can be assumed to be defined over $\F$. 
\end{proposition} 

\begin{proof}
	Since $g \in G$ is unipotent, the nilpotent element $x=\log(g)\in \frakg $ exists. 
	By the Jacobson-Morozov Lemma, Proposition \ref{prop:Jacobson_Morozov_algebra}, there are $y,h \in \frakg$ such that $(x,y,h)$ is an $\mathfrak{sl}_2$-triplet. Let $\mathfrak{sl}_g$ denote the subalgebra of $\frakg$ generated by $x,y$ and $h$. Thus, $x \in \operatorname{sl}_g$ corresponds to
	$$
	\begin{pmatrix}
		0 & 1 \\ 0 & 0
	\end{pmatrix}\in \mathfrak{sl}_2
	$$
	under the isomorphism $\mathfrak{sl}_g \cong \mathfrak{sl}_2$.
	We follow \cite[7.1]{Bor} and define
	$$
	\mathcal{A}(\mathfrak{sl}_g) = \bigcap \left\{ H \colon H \text{ algebraic subgroup of $G$ with } \mathfrak{sl}_g \subseteq \operatorname{Lie}(H) \right\},
	$$
	which is a 
	connected, normal algebraic subgroup. Let 
	$$
	\operatorname{SL}_g = [	\mathcal{A}(\mathfrak{sl}_g) , 	\mathcal{A}(\mathfrak{sl}_g) ]
	$$ 
	be its commutator group, which is an algebraic group by \cite[2.3]{Bor} since $\mathcal{A}(\mathfrak{sl}_g)$ is connected. We then use \cite{Bor} to obtain
	\begin{align*}
		\operatorname{Lie}(\operatorname{SL}_g ) &\stackrel{7.8}{=} [ \operatorname{Lie} (\mathcal{A}(\mathfrak{sl}_g)) ,\operatorname{Lie}( \mathcal{A}(\mathfrak{sl}_g) )] \stackrel{7.9}{=} [\mathfrak{sl}_g , \mathfrak{sl}_g] = \mathfrak{sl}_g
	\end{align*}
	where the last inequality follows from $\mathfrak{sl}_g \cong \mathfrak{sl}_2$. The map
	\begin{align*}
		\alpha \colon \G_a & \to G \\
		t & \mapsto \exp(tx) = \sum_{n= 0}^\infty \frac{(tx)^n}{n!}
	\end{align*}
	is a polynomial and hence a morphism of algebraic groups. We have $\operatorname{Lie}(\alpha(\G_a)) = \langle x \rangle \subseteq \mathfrak{sl}_g$ and hence 
	$
	g \in \alpha(\G_a) \subseteq \operatorname{SL}_g
	$ 
	by \cite[7.1(2)]{Bor}.
	
	If $g \in G_\F$, $\log(g) \in \frakg_\F$, and we may assume $y,h \in \frakg_\F$ as well. Thus $\mathfrak{sl}_g \subseteq \frakg_\F$. Now $\mathcal{A}(\mathfrak{sl}_g)$ is defined over $\F$ as in \cite[2.1(b)]{Bor}. 
	Then by \cite[2.3]{Bor}, $\operatorname{SL}_g$ is defined over $\F$. 
\end{proof}

Let us return to the case of real closed fields $\K \subseteq \F \cap \R$. Recall from Section \ref{sec:U}, that $U_\F$ has subgroups $(U_\alpha)_\F$ consisting of unipotent elements for $\alpha \in \Sigma_{>0}$ defined as the semialgebraic extensions of $U_\alpha = \exp(\frakg_\alpha \oplus \frakg_{2\alpha })$. The following is another variation of the Jacobson-Morozov Lemma that only works when $\Sigma$ is reduced. 
\begin{proposition}\label{prop:Jacobson_Morozov_real_closed}
	Let $\alpha \in \Sigma$ and assume $\frakg_{2\alpha} = 0$. Let $u \in (U_{\alpha})_\F$. Then there are $ X \in (\frakg_\alpha)_\F$ and $t \in \F$ such that $u=\exp(tX)$ and $(-v)(X) = 0$. Let $(X,Y,H)$ be the $\mathfrak{sl}_2$-triplet of Lemma \ref{lem:JM_basic}. Then there is a morphism of algebraic groups
	$
	\varphi \colon \operatorname{SL}(2,\D)  \to G 
	$
	defined over $\F$ such that $\varphi$ has finite kernel and
	\begin{align*}
		\varphi \begin{pmatrix}
			1 & t \\ 0& 1
		\end{pmatrix} =  u = \exp(tX)  \quad \text{and} \quad  \varphi \begin{pmatrix}
		1 & 0 \\ t & 1
		\end{pmatrix} = \exp(tY).
	\end{align*}
	If $\varphi$ is not injective, then $\ker (\varphi) \cong \Z/2\Z$ and $\varphi$ factors through the isomorphism
	$$
	\operatorname{PGL}(2,\D) := \operatorname{SL(2,\D)}/\ker(\varphi)  \xrightarrow{\sim} \varphi(\operatorname{SL}(2,\D))
	$$
	which is also defined over $\F$. Moreover $\varphi(g\tran) = \varphi(g)\tran$, for any $g \in \operatorname{SL}(2,\F)$.
\end{proposition}
\begin{proof}
	Let $u,X,t$ as in the statement. We apply Proposition \ref{prop:Jacobson_Morozov_group} to $\exp(X) \in (U_\alpha)_\F$ to obtain an algebraic group $\operatorname{SL}_{\exp(X)}(\D) < G_\D$ defined over $\F$ with Lie algebra $\operatorname{sl}_2(\D)$. By \cite[Corollary 32.3]{Hum1}, the only algebraic groups with Lie algebra $\mathfrak{sl}_2(\D)$ are $\operatorname{SL}(2,\D)$ and $\operatorname{PGL}(2,\D) := \operatorname{SL}(2,\D)/Z(\operatorname{SL(2,\D)})$. In both cases we obtain an algebraic homomorphism $\varphi \colon \operatorname{SL}_2(\D) \to \operatorname{SL}_{\exp(X)}$ with finite kernel.
	
	We note that since $X\in (\frakg_\alpha)_\F$, the $\mathfrak{sl}_2$-triplet is described by $(X,Y,H)$ in Lemma \ref{lem:JM_basic} and we have a Lie algebra isomorphism $\mathfrak{sl}_2(\D) \cong \operatorname{Lie}(\operatorname{SL}_{\exp(X)}(\D))$. We note that the two Lie algebra isomorphisms
	 \begin{align*}
	 \varphi \colon \operatorname{SL}_2(\mathbb{D})  
	 	\xrightarrow{\sim} & \  \operatorname{SL}_{\exp(X)}(\mathbb{D}) \\
	 	\mathfrak{sl}_2(\mathbb{D}) \cong & \  \mathfrak{sl}_{\exp(X)}(\mathbb{D}) &&\hspace{-14ex} \cong \mathfrak{sl}_2(\mathbb{D}) \\
	 	& \ X &&\hspace{-14ex} \mapsfrom
	 	\begin{pmatrix}
	 		0 & 1 \\ 0 & 0
	 	\end{pmatrix}       
	 \end{align*}
	 may not coincide. Since $\operatorname{sl}_2(\D)$ does not have any outer automorphisms (\cite[Theorem 14.1]{Hum1} and \cite[Proposition D.40]{FuHa}), the two isomorphisms only differ by $\operatorname{Ad}(g)$ for some $g\in \operatorname{SL}(2,\D)$. Up to conjugation we may therefore assume that $\operatorname{d}\!\varphi \colon \mathfrak{sl}_2(\D) \to \operatorname{Lie}(\operatorname{SL}_{\exp(X)}(\D))$ maps 
	 \begin{align*}
	                  \begin{pmatrix}
	                  	0 & 1 \\ 0 & 0
	                  \end{pmatrix}        \mapsto X ,            \quad
	                  	 	\begin{pmatrix}
	                   		1 & 0 \\ 0 & -1
	                  	 	\end{pmatrix}   \mapsto   H   \quad \text{and} \quad
	                  	 	\begin{pmatrix}
	                  	 		0 & 0 \\ 1 & 0
	                  	 	\end{pmatrix}      \mapsto  Y    .
	 \end{align*}
	 Even though the conjugation may be by an element in $\operatorname{SL}_2(\D)$, the explicit description of $\operatorname{d}\!\varphi$ shows that $\varphi$ is still defined over $\F$. The description of $\exp(tX)$ and $\exp(tY)$ in the statement of the Proposition also follows. Finally, note that transposition is $\operatorname{d}\!\varphi$-equivariant and hence also $\varphi$-equivariant. 
 \end{proof} 

Since in the $\mathfrak{sl}_2$-triplet, $H \in \fraka_\F$ we obtain multiplicative one-parameter groups
\begin{align*}
	\F_{>0} & \to A_\F \\
\lambda &\mapsto \varphi\begin{pmatrix}
	\lambda & 0 \\ 0 & \lambda^{-1}
\end{pmatrix} 	
\end{align*}
which satisfy the following property that will be useful later.
\begin{lemma}\label{lem:Jacobson_Morozov_oneparam}
	Let $\alpha \in \Sigma$ such that $(\frakg_{2\alpha})_\F = 0$. Let $u\in (U_\alpha)_\F$ and $\varphi \colon \operatorname{SL}(2,\D) \to G$ as in Proposition \ref{prop:Jacobson_Morozov_real_closed}. For every $\lambda \in \F_{>0}$
	$$
	\chi_\alpha \left(\varphi \begin{pmatrix}
		\lambda & 0 \\ 0 & \lambda^{-1}
	\end{pmatrix} \right) = \lambda^2.
	$$
\end{lemma}
\begin{proof}
	We note that for all $g \in \operatorname{SL}(2,\F)$, the diagrams 
	$$
	\begin{tikzcd}
		{\mathfrak{sl}(2,\mathbb{F})} \arrow[r, "\operatorname{Ad}(g)"] \arrow[d, "\operatorname{d}\!\varphi"] & {\mathfrak{sl}(2,\mathbb{F})} \arrow[d, "\operatorname{d}\!\varphi"] &   {\operatorname{SL}(2,\mathbb{F})} \arrow[d, "\varphi"] \arrow[r, "c_g"] & {\operatorname{SL}(2,\mathbb{F})} \arrow[d, "\varphi"] \\
		\mathfrak{g}_{\mathbb{F}} \arrow[r, "\operatorname{Ad}(\varphi(g))"]                                   & \mathfrak{g}_{\mathbb{F}}                                              & G_{\mathbb{F}} \arrow[r, "c_{\varphi(g)}"]                              & G_{\mathbb{F}}                                        
	\end{tikzcd}
	$$
	commute. For
	$$
	a = \varphi \begin{pmatrix}
		\lambda & 0 \\ 0 & \lambda^{-1}
	\end{pmatrix} \quad \text{ and } \quad X = \operatorname{d}\! \varphi \begin{pmatrix}
	0 & 1 \\ 0 &0 
	\end{pmatrix} \in (\frakg_{\alpha})_\F
	$$
	we then have
	\begin{align*}
		\operatorname{Ad}\left(a \right) X &= \operatorname{d}\!\varphi \left(\operatorname{Ad}\begin{pmatrix}
		\lambda & 0 \\ 0 & \lambda^{-1}
		\end{pmatrix} \begin{pmatrix}
		0 & 1 \\ 0 & 0
		\end{pmatrix}
		\right) = \operatorname{d}\!\varphi \begin{pmatrix}
			0 & \lambda^{2} \\ 0 & 0
		\end{pmatrix} = \lambda^2 X
	\end{align*}
	and since $\chi_\alpha(a)$ is defined by $\operatorname{Ad}(a)X = \chi_{\alpha}(a)X$ we have $\chi_\alpha(a) = \lambda^{2}$.
\end{proof}

\begin{lemma}\label{lem:Jacobson_Morozov_m}
	Let $\alpha \in \Sigma$ such that $(\mathfrak{g}_{2\alpha})_\F = 0$. Let $u \in (U_\alpha)_\F $ and $X\in (\frakg_\alpha)_\F$, $t\in \F$ and $\varphi \colon \operatorname{SL}(2,\D) \to G$ as in Proposition \ref{prop:Jacobson_Morozov_real_closed}, so that $u = \exp(tX)$. The element
	$$
	m(u) := \varphi \begin{pmatrix}
		0 & t \\ -1/t & 0
	\end{pmatrix} \in G_\F
	$$ 
	is contained in $\operatorname{Nor}_{G_\F}(A_\F)$ and 
	$$
	\chi_\alpha(m(u)\cdot a\cdot m(u)^{-1}) = \chi_\alpha(a)^{-1}
	$$ 
	for any $a \in A_\F$.
	
	When $t=1$, even $m(u) \in N_\F = \operatorname{Nor}_{K_\F}(A_\F)$ and 
	$$
	m(u) \cdot (U_\alpha)_\F \cdot  m(u)^{-1} = (U_{-\alpha})_\F.
	$$
\end{lemma}
\begin{proof}
	Let $a\in A_\F$ and
	$$
	a_\perp := \varphi \begin{pmatrix}
		\sqrt{\chi_\alpha(a)} & 0 \\ 0 & \sqrt{\chi_\alpha(a)}^{-1}
	\end{pmatrix},
	$$ 
	then $a_0:= aa_\perp^{-1}$ satisfies
	$$
	\chi_\alpha (a_0) = \chi_\alpha(a) \chi_\alpha(a_\perp)^{-1} = \chi_\alpha(a) \frac{1}{\chi_\alpha(a)} = 1 
	$$
	by Lemma \ref{lem:Jacobson_Morozov_oneparam}. For any $u' \in \exp((\frakg_{\alpha})_\F)$ or $u'\in \exp((\frakg_{-\alpha})_\F)$, $a_0u'=u'a_0$ by Lemma \ref{lem:aexpXa}. We note that
	\begin{align*}
			m(u) &= \varphi \begin{pmatrix}
			0 & t \\ -1/t & 0
		\end{pmatrix} =  \varphi \left( \begin{pmatrix}
			1 & 0 \\ -1/t & 1
		\end{pmatrix}  \begin{pmatrix}
			1 & t \\ 0 & 1
		\end{pmatrix} \begin{pmatrix}
			1 & 0 \\ -1/t & 1
		\end{pmatrix}\right) \\
		& \in \exp((\frakg_{-\alpha})_\F) \cdot u \cdot \exp((\frakg_{-\alpha})_\F)
	\end{align*}
	which implies $a_0m(u) = m(u)a_0$. Now
	\begin{align*}
		m(u) \cdot a \cdot m(u)^{-1} &= m(u) \cdot a_0 \cdot a_\perp \cdot m(u)^{-1} = a_0 \cdot m(u) \cdot  a_\perp \cdot m(u)^{-1} \\
		&= a_0 \cdot \varphi \left(\begin{pmatrix}
			0 & t \\ -1/t & 0
		\end{pmatrix}\begin{pmatrix}
		\sqrt{\chi_\alpha(a)} & 0 \\ 0 & \sqrt{\chi_\alpha(a)}^{-1}
		\end{pmatrix}\begin{pmatrix}
		0 & -t \\ 1/t & 0
		\end{pmatrix}\right) \\
		&= a_0\cdot \varphi \begin{pmatrix}
		\sqrt{\chi_\alpha(a)}^{-1} & 0 \\ 0 & \sqrt{\chi_\alpha(a)}
		\end{pmatrix} = a_0 a_\perp^{-1} \in A_\F
	\end{align*}
	and thus $m(u) = \operatorname{Nor}_{G_\F}(A_\F)$. We see directly that 
	$$
	\chi_\alpha (m(u) \cdot a\cdot m(u)^{-1}) = \chi_\alpha(a_0 \cdot a_\perp^{-1}) = \chi_\alpha(a_0) \cdot \chi_\alpha(a_\perp)^{-1} = \chi_\alpha(a)^{-1}.
	$$
	Now if $t=1$ we can use that $\varphi$ preserves transposition by Proposition \ref{prop:Jacobson_Morozov_real_closed}, to show
	\begin{align*}
		m(u) \cdot m(u)\tran &= \varphi \begin{pmatrix}
			0 & 1 \\ -1 & 0
		\end{pmatrix} \left( \varphi \begin{pmatrix}
			0 & 1 \\ -1 & 0
		\end{pmatrix}\right)\tran \\
		& = \varphi \left( \begin{pmatrix}
			0 & 1 \\ -1 & 0
		\end{pmatrix}   \begin{pmatrix}
			0 & -1 \\ 1 & 0
		\end{pmatrix} \right)= \varphi\begin{pmatrix}
			1 & 0 \\ 0 & 1
		\end{pmatrix} = \operatorname{Id}
	\end{align*}
	and hence $m(u)\in K_\F$, so $m(u) \in N_{\F} = \operatorname{Nor}_{K_\F}(A_\F)$. Thus $m(u)$ is a representative of an element $w=[m(u)]$ of the spherical Weyl group $W_s = N_\F/M_\F$. While $u\in (U_{\alpha})_\F$,
	\begin{align*}
		m(u)\cdot u\cdot m(u)^{-1} &= \varphi_\F \left( \begin{pmatrix}
			0 & 1 \\ -1 & 0
		\end{pmatrix}\begin{pmatrix}
			1 & t \\ 0 & 1
		\end{pmatrix} \begin{pmatrix}
			0 & -1 \\ 1 & 0
		\end{pmatrix} \right) \\
		&= \varphi_\F \begin{pmatrix}
			1 & 0 \\ -t & 1
		\end{pmatrix} \in (U_{-\alpha})_\F,
	\end{align*}
	and thus we have $w(\alpha) = -\alpha$. Then 
$$
m(u) \cdot (U_\alpha)_\F \cdot  m(u)^{-1} = (U_{w(\alpha)})_\F =  (U_{-\alpha})_\F.
$$
\end{proof}

\subsection{Rank 1 subgroups}\label{sec:rank1}
Let $\alpha \in \Sigma$. In this section we investigate the group generated by $U_\alpha$ and $U_{-\alpha}$. When $\dim(U_\alpha) = 1$, then this group is given by the image of the Jacobson-Morozov-morphism $\varphi \colon \operatorname{SL}(2,\D) \to G$ from Proposition \ref{prop:Jacobson_Morozov_real_closed}. In general, when $\dim(U_\alpha)$ is not $1$, the group generated is larger than the image of $\varphi$, but is still rank $1$. 

\begin{theorem}\label{thm:levi_group}
	Let $\alpha \in \Sigma$. Then there is a connected semisimple self-adjoint linear algebraic group $L_{\pm \alpha}$ defined over $\K$ such that
	\begin{enumerate}
		\item [(i)] $\operatorname{Lie}(L_{\pm \alpha}) = ( \frakg_\alpha \oplus \frakg_{2\alpha} ) \oplus (\frakg_{-\alpha} \oplus \frakg_{-2\alpha} ) \oplus ( [\frakg_\alpha , \frakg_{-\alpha}] + [\frakg_{2\alpha}, \frakg_{-2\alpha}] )$.
		\item [(ii)] $\operatorname{Rank}_\R(L_{\pm \alpha}) = \operatorname{Rank}_\F(L_{\pm \alpha}) = 1$.
	\end{enumerate}
\end{theorem}
\begin{proof}
    We consider the semisimple $\D$-Lie algebra 
    $$
    \frakl :=  ( \frakg_\alpha \oplus \frakg_{2\alpha} ) \oplus (\frakg_{-\alpha} \oplus \frakg_{-2\alpha} ) \oplus ( [\frakg_\alpha , \frakg_{-\alpha}] + [\frakg_{2\alpha}, \frakg_{-2\alpha}] ).
    $$
    Similar to the proof of Proposition \ref{prop:Jacobson_Morozov_group} we follow Borel \cite[7.1]{Bor} by defining the connected normal algebraic subgroup
    $$
    \mathcal{A}(\frakl) = \bigcap \left\{ H \colon H<G \text{ is an algebraic subgroup } \frakl \subseteq \operatorname{Lie}(H) \right\}
    $$
    of $G$. We set $L_{\pm \alpha} := [\mathcal{A}(\frakl), \mathcal{A}(\frakl)]$. Then (i) follows by
    $$
    \operatorname{Lie}(L_{\pm \alpha}) = [\mathcal{A}(\frakl) , \mathcal{A}(\frakl)] = [\frakl , \frakl] = \frakl ,
    $$
    where we used that $\frakl$ is semisimple in the step $[\frakl, \frakl] = \frakl$. The algebraic group $\mathcal{A}(\frakl)$ is defined over $\K$ by \cite[2.1(b)]{Bor}. Thus $L_{\pm \alpha}$ is defined over $\K$ and connected by \cite[2.3]{Bor}. Since $\frakl$ is semisimple, so is $L_{\pm \alpha}$. Since $\theta(X) = -X\tran$ and $\theta(\frakg_\alpha) = \frakg_{-\alpha}$ and $\theta(\frakg_{2\alpha}) = \frakg_{-2\alpha}$ by Proposition \ref{prop:root_decomp}, $\theta([\frakg_\alpha, \frakg_{-\alpha}]) = [\frakg_\alpha, \frakg_{-\alpha}]$ and $\theta([\frakg_{2\alpha}, \frakg_{-2\alpha}]) = [\frakg_{2\alpha}, \frakg_{-2\alpha}]$, and hence $\frakl$ and $L_{\pm \alpha}$ are self-adjoint. 
    
    For (ii), Lemma \ref{lem:levi_algebra} tells us that the real rank of $\frakl_\R$ is $1$ and hence $\operatorname{Rank}_\R(L_{\pm \alpha}) = 1$. By the theorem on split tori, Theorem \ref{thm:split_tori}, we then have $\operatorname{Rank}_\R(L_{\pm\alpha}) = \operatorname{Rank}_\F(L_{\pm \alpha}) = 1$. 
\end{proof}

\begin{lemma}\label{lem:levi_fixes_A}
  $(L_{\pm \alpha})_\F \subseteq \operatorname{Cen}_{G_\F}(\left\{ a \in A_\F \colon \chi_\alpha (a) = 1 \right\})$. 
\end{lemma}
\begin{proof}
This is a semialgebraic statement, hence it suffices to prove it over the real numbers. Let $a \in A_\R$ with $\chi_\alpha(a) = 1$. We have 
$$
\operatorname{Lie}(L_{\pm \alpha}) = ( \frakg_\alpha \oplus \frakg_{2\alpha} ) \oplus (\frakg_{-\alpha} \oplus \frakg_{-2\alpha} ) \oplus ( [\frakg_\alpha , \frakg_{-\alpha}] \oplus [\frakg_{2\alpha}, \frakg_{-2\alpha}] ).
$$
For $X \in (\frakg_\alpha)_\R$ and $X' \in (\frakg_{2\alpha})_\R$, we have by Lemma \ref{lem:aexpXa}
\begin{align*}
	a \exp(X) a^{-1} & = \exp(\chi_\alpha(a)X ) = \exp(X) \\
	a \exp(X') a^{-1}  & = \exp(\chi_\alpha(a)^2 X') = \exp(X')
\end{align*}
and the same argument shows $a \exp(Y) a^{-1} = \exp(Y) $ for $Y \in (\frakg_{-\alpha})_\R \oplus (\frakg_{-2\alpha})_\R$. Since $\frakg_0 \cap \frakl = \mathfrak{z}_{\frakk\cap \frakl} (\fraka \cap \frakl) \oplus (\fraka \cap \frakl)$, also $a \exp(H)a^{-1} =\exp(H)$ for $H \in \frakg_0 \cap \frakl$ is clear. We have shown that 
$$
\operatorname{Ad}(a) \colon \operatorname{Lie}((L_{\pm \alpha})_\R) \to  \operatorname{Lie}((L_{\pm \alpha})_\R) 
$$
is the identity, and conjugation by $a$ is constant on connected components of $(L_{\pm})_\R$. Since $L_{\pm \alpha}$ is connected as an algebraic group and $c_a^{-1}(\operatorname{Id})$ is a closed algebraic set, 
conjugation by $a$ is the identity.
\end{proof}

We can treat $L_{\pm \alpha}$ from Theorem \ref{thm:levi_group} as an example of the theory we have developed so far, in particular we can apply group decompositions as outlined in the beginning of Section \ref{sec:decompositions}. Since $S <G$ is a self-adjoint maximal split torus, whose Lie algebra contains $[\frakg_\alpha, \frakg_{-\alpha}] \oplus [\frakg_{2\alpha}, \frakg_{-2\alpha}]$, $S_{\pm \alpha}:= S \cap L_{\pm \alpha}$ is a self-adjoint maximal split torus of $L_{\pm \alpha}$. Then $K_{\pm \alpha} = L_{\pm \alpha} \cap \operatorname{SO}_n = L_{\pm \alpha} \cap K$. Considering the $\F$-points, the semialgebraic connected component $(A_{\pm \alpha})_\K$ of $(S_{\pm \alpha})_\K$ containing the identity can be semialgebraically extended to $(A_{\pm \alpha})_\F$. The following is a version on Lemma \ref{lem:Jacobson_Morozov_oneparam}, but now for the rank 1 group $(L_{\pm \alpha})_\F$.
\begin{lemma}
	For every $t \in \F_{>0}$, there is an $a \in (A_{\pm \alpha})_\F$ such that $\chi_\alpha(a) = t$.
\end{lemma}
\begin{proof}
	This is clearly a first-order statement and it therefore suffices to show it for $\R = \F$. Let $X \in (\frakg_\alpha)_\R \setminus 0$. From Lemma \ref{lem:levi_algebra}, we know that $\exp([X, \theta(X)]) \in A_{\pm \alpha}$. Given $t \in \R_{>0}$, let
	$$
	a := \exp \left( \frac{ \log(t) }{\alpha([X,\theta(X)])}  \cdot [X,\theta(X)] \right).
	$$
	Then 
	$
	\chi_\alpha(a) = e^{\alpha(\log(a))} = e^{\log(t)} = t
	$
	by Lemma \ref{lem:char}.
\end{proof}
The root space decomposition then gives $U_{\pm \alpha}=\exp(\frakg_\alpha \oplus \frakg_{2\alpha}) = U_\alpha$, and the spherical Weyl group $W_{\pm \alpha} = (N_{\pm \alpha})_\F / (M_{\pm \alpha})_\F$ can be defined from \begin{align*}
	(N_{\pm\alpha})_\F &:= \operatorname{Nor}_{(K_{\pm \alpha})_\F}((A_{\pm \alpha})_\F)\\
	(M_{\pm\alpha})_\F &:= \operatorname{Cen}_{(K_{\pm \alpha})_\F}((A_{\pm \alpha})_\F).
\end{align*}
We note that as a consequence of Lemma \ref{lem:levi_fixes_A}, $N_{\pm \alpha} \subseteq N_\F$ and also $\operatorname{Nor}_{L_{\pm \alpha}}(A_{\pm \alpha}) \subseteq \operatorname{Nor}_{K_\F}(A_\F)$.

In the case that $\Sigma$ is reduced, there is an interpretation in terms of the Jacobson-Morozov Lemma. Given $u \in (U_{\pm \alpha})_\F$, we note that the morphism $\varphi_\F \colon \operatorname{SL}(2,\F) \to G_\F$ from Proposition \ref{prop:Jacobson_Morozov_real_closed} takes values in $(L_{\pm \alpha})_\F$. In fact, 
$$
(A_{\pm \alpha})_\F = \varphi_\F\left( \left\{ \begin{pmatrix}
	\lambda & 0 \\ 0 & \lambda^{-1}
\end{pmatrix}   \colon \lambda>0 \right\} \right)
$$
since over $\R$, both of these groups are connected and one-dimensional.
By Lemma \ref{lem:Jacobson_Morozov_m}, the element
$$
m = \varphi_\F\begin{pmatrix}
	0 & 1 \\ -1 & 0
\end{pmatrix} \in (N_{\pm \alpha})_\F
$$
is a representative of the only non-trivial element in the spherical Weyl group $W_{\pm \alpha} = \{ [\operatorname{Id}] , [m] \}$. We can now apply the Bruhat decomposition Theorem \ref{thm:BWB} to $(L_{\pm \alpha})_\F$.

\begin{corollary}\label{cor:levi_Bruhat}
	Let $(B_{\alpha})_\F := (M_{\pm \alpha})_\F (A_{\pm \alpha})_\F (U_\alpha)_\F$. Then 
	$$
	(L_{\pm \alpha})_\F = (B_{\alpha})_\F (N_{\pm \alpha})_\F (B_{\alpha})_\F ,
	$$
	and the element $m \in (N_{\pm \alpha})_\F$ is a representative of a unique element in $W_{\pm \alpha}$, so
	$$
	(L_{\pm \alpha})_\F = (B_{\alpha})_\F \ \amalg \   (B_{\alpha})_\F \cdot m \cdot (B_{ \alpha})_\F.
	$$
\end{corollary}
We will use the following variation in Lemma \ref{lem:BTuu_usu}.
\begin{corollary}\label{cor:levi_Bruhat_alternative}
	For
	\begin{align*}
		(B_{\alpha})_\F &:= (M_{\pm \alpha})_\F (A_{\pm \alpha})_\F (U_\alpha)_\F, \\
		(B_{-\alpha})_\F &:= (M_{\pm \alpha})_\F (A_{\pm \alpha})_\F (U_{-\alpha})_\F,
	\end{align*} 
	we obtain the decompositions
	\begin{align*}
		(L_{\pm \alpha})_\F &=(B_{\alpha})_\F \cdot (B_{-\alpha})_\F  \ \amalg \   m \cdot (B_{ -\alpha})_\F,  \\
		(L_{\pm \alpha})_\F &=(B_{\alpha})_\F \cdot (B_{-\alpha})_\F  \ \amalg \   (B_{\alpha})_\F \cdot m \cdot (B_{ -\alpha})_\F,  \\
		(L_{\pm \alpha})_\F &= (B_{ \alpha})_\F (N_{\pm \alpha})_\F (B_{- \alpha})_\F 
	\end{align*}
where in the last one, the element in $(N_{\pm \alpha})_\F$ is a representative of a unique element in $W_{\pm \alpha} = \left\{ \operatorname{Id}, [m] \right\}$.
\end{corollary}
\begin{proof}
	We use $m^{-1} \cdot (U_\alpha)_\F\cdot  m = (U_{-\alpha})_\F$ from Lemma \ref{lem:Jacobson_Morozov_m}. By Corollary \ref{cor:levi_Bruhat}, we choose $m^{-1}$ as the representative of $[m^{-1}]=[m] \in W_{\pm \alpha} $ and obtain
	$$
	(L_{\pm \alpha})_\F = (B_{\alpha})_\F \ \amalg \   (B_{\alpha})_\F \cdot m^{-1} \cdot (B_{ \alpha})_\F,
	$$
	so 
	$$
	(L_{\pm \alpha})_\F = m \cdot (B_{-\alpha})_\F \cdot m^{-1} \ \amalg \   (B_{\alpha})_\F  \cdot (B_{ -\alpha})_\F \cdot m^{-1}.
	$$
	Multiplying this expression on the right by $m\in (L_{\pm \alpha})_\F$ results in
	$$
	(L_{\pm \alpha})_\F =(B_{\alpha})_\F \cdot (B_{-\alpha})_\F  \ \amalg \   m \cdot (B_{ -\alpha})_\F, 
	$$ 
	and further multiplying by $(B_\alpha)_\F$ on the left results in the remaining two decompositions.
\end{proof}

\subsection{Kostant convexity}\label{sec:kostant}
Using the Iwasawa decomposition $G_\F=U_\F A_\F K_\F$, Theorem \ref{thm:KAU_R}, we associate to every $g=uak\in G_\R$ its $A$-component $a_\R(g) = a \in A_\R$. The following is Kostant's convexity theorem, which we will generalize to $G_\F$ in this chapter.
\begin{theorem}\cite[Theorem 4.1]{Kos}\label{thm:kostant_R}
	For every $b \in A_\R$, 
	$$
	\left\{ a_\R(kb) \in A_\R \colon k \in K_\R \right\} = \exp \left( \operatorname{conv}(W_s \log(b)) \right),
	$$
	where $\log \colon A_\R \to \fraka$ is the inverse of $\exp$, $W_s$ is the spherical Weyl group acting on $\fraka$ and $\operatorname{conv}(W_s \log(b))$ is the convex hull of the Weyl group orbit of $\log(b)$.
\end{theorem}
The left hand side of the equation in Theorem \ref{thm:kostant_R} is already a semialgebraic set and we will reformulate the right hand side as a semialgebraic set as well. For this, we first analyze the root system $\Sigma \subseteq \fraka^\star$.

\begin{figure}
	\centering
	\includegraphics[width=0.6\linewidth]{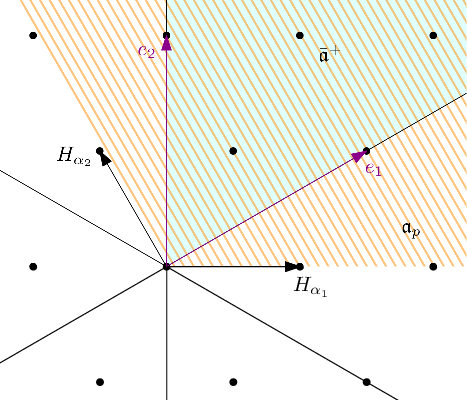}
	\caption{Root system of type $A_2$ associated to $\operatorname{SL}_3$. The convex cone $\fraka_p$ (orange stripes) can be viewed as spanned by the $H_{\alpha_i}$ or as the intersection of the half-spaces defined by the primitive vectors $e_i$.  }
	\label{fig:kostant}
\end{figure}

Recall from Section \ref{sec:killing_involutions_decompositions}, that we have a scalar product $B_\theta$ on $\fraka$, which can be used to set up an isomorphism $\fraka^\star \to \fraka, \alpha \mapsto H_\alpha$ that satisfies the defining property $\lambda(H) = B_\theta(H_\lambda, H)$ for all $H \in \fraka$. In this section we will denote $B_\theta$ as well as the corresponding scalar product on $\fraka^\star$ with brackets $\langle \cdot , \cdot \rangle$. For $\lambda \in \fraka^\star$ define\footnote{In the notation of Section \ref{sec:root_system}, $x_\alpha = H_{\alpha}^\vee$.} 
$$
x_\lambda = \frac{2}{\langle H_\lambda, H_\lambda \rangle } H_\lambda  =  \frac{2}{\langle \lambda, \lambda \rangle } H_\lambda
$$
which satisfies the property that for all $\alpha, \beta \in \Sigma$,
$$
\langle H_\alpha, x_\beta \rangle = \frac{2 \langle H_\alpha, H_\beta \rangle }{\langle H_\beta, H_\beta\rangle} \in \Z
$$
since $\Sigma$ is a crystallographic root system, Theorem \ref{thm:sigma_root}. We note that in the language of Section \ref{sec:modelapartment}, we have $x_\alpha = H_{\alpha^\vee}$ for all $\alpha \in \Sigma$. Let $\Delta = \left\{ \delta_1, \ldots , \delta_r \right\}$ be a basis of $\Sigma$ and abbreviate $x_i := x_{\delta_i}$ and $H_i := H_{\delta_i}$. Let $\frakap = \left\{ H \in \fraka \colon \delta(H) \geq 0 \text{ for all } \delta \in \Delta \right\} $. On the way to prove Theorem \ref{thm:kostant_R}, Kostant describes $\operatorname{conv}(W_s x)$ using the closed convex cone
$$
\fraka_p := \left\{ x \in \fraka \colon x = \sum_{i=1}^r \R_{\geq 0}x_i \right\}
$$
illustrated in Figure \ref{fig:kostant}.
\begin{lemma}\cite[Lemma 3.3.(2)]{Kos}\label{lem:kostant}
	Let $x,y\in \frakap$. Then 
	\begin{align*}
		y \in  \operatorname{conv}(W_s x) \quad \iff \quad  x-y \in \fraka_p.
	\end{align*}
\end{lemma}
Now we want to describe $\fraka_p$ by inequalities. The cone $\fraka_p$ is an intersection of open half-spaces defined by the half-planes
$$
E_j = \sum_{k \neq j} \R x_k .
$$  
Any vector $x\in \fraka$ orthogonal to $E_j$ satisfies $\langle x,x_i\rangle = 0$ for all $i\neq j$. Writing $x = \sum_k \lambda_{jk}H_k$, we have for all $i \neq j$
$$
\langle x, x_i \rangle = \sum_k \lambda_{jk} \langle H_k, x_i \rangle = 0.
$$
For every $i\neq j$, this is a homogeneous linear equation with variables $\lambda_{j1}, \ldots ,\lambda_{jr}$ and coefficients $\langle H_k,x_i \rangle \in \Z$. Therefore there is a rational (and hence integer) solution for the $\lambda_{jk}$.
This shows that $x\in E_j^\perp \setminus \{0\}$ may be chosen to lie in the lattice $\Gamma = \sum_{l = 1}^r \Z H_k$. There are two primitive vectors in $\Gamma \cap E_j^\perp$. Let $e_j$ be the unique one that is on the same side of $\partial E_j$ as $\fraka_p$. Thus
$$
\fraka_p = \left\{ x \in \fraka \colon \langle e_j, x \rangle \geq 0 \text{ for all } j\right\}.
$$

\begin{figure}[h]
	\centering
	\includegraphics[width=0.6\linewidth]{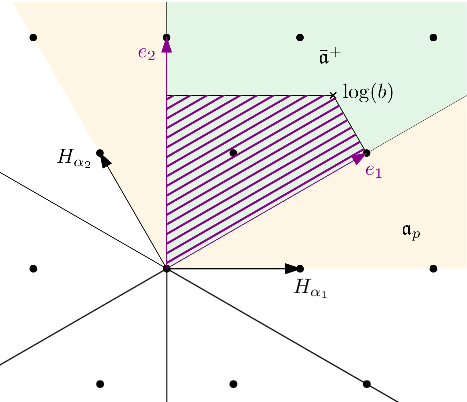}
	\caption{Root system of type $A_2$ associated to $\operatorname{SL}_3$. The convex set in Kostant's convexity Theorem \ref{thm:kostant_F} defined by inequalities is illustrated in purple. }
	\label{fig:kostant7}
\end{figure}

Under the isomorphism $\fraka^\star \cong \fraka$, the lattice $\Gamma$ corresponds to $L = \sum_{\alpha\in \Delta} \Z \alpha$ and we define $\gamma_j \in L$ to be the element corresponding to $e_j \in \Gamma$. By Lemma \ref{lem:char}, $\gamma_j$ defines an algebraic character $\chi_j  \colon A_\R \to \R $ satisfying
$$
\chi_j (\exp(H)) = e^{\gamma_j(H)}
$$
for all $H \in \fraka$. The multiplicative closed Weyl chamber is
$$
\Ap_\R : = \left\{ a \in A_{\R} \colon \chi_\delta(a) \geq 1 \text{ for all } \delta \in \Delta \right\} = \exp \frakap
$$
which is a semialgebraic set and thus has an $\F$-extension $\Ap_\F$. Using the the Iwasawa decomposition $G_\F=U_\F A_\F K_\F$, Theorem \ref{thm:KAU}, we associate to every $g=uak\in G_\F$ its $A$-component $a_\F(g) = a \in A_\F$ as before. We can now conclude the following version of Kostant's convexity theorem for $G_\F$, illustrated in Figure \ref{fig:kostant7}.

\begin{theorem}\label{thm:kostant_F}
	For all $b \in \ApF$, we have
	$$
	\left\{ a \in \ApF \colon \exists k \in K_\F , a_\F(kb) = a \right\}
	= \left\{ a\in \ApF \colon \chi_i (a) \leq \chi_i (b) \text{ for all } i \right\}.
	$$
\end{theorem}
\begin{proof}
	We first verify that the theorem holds for $\F=\R$. An element $a \in \ApR$ satisfies $a_\R(kb) = a$, for some $k \in K_\R$, if and only if $a \in \ApR \cap \exp(\operatorname{conv}(W_s \log(b)))$ by Theorem \ref{thm:kostant_R}. This means $a = \exp(H)$ for some $H \in \frakap$ and $H\in \operatorname{conv}(W_s\log(b))$. By Lemma \ref{lem:kostant}, this is equivalent to 
	$$
	\log(b) - H \in \fraka_p = \left\{ x \in \fraka \colon \langle e_j, x \rangle \geq 0 \text{ for all } j\right\},
	$$
	and taking exponents this is equivalent to
	$$
	\chi_i(ba^{-1}) \geq 1 
	$$
	for all $i$, so
	$$
	\chi_i(a) \leq \chi_i(b),
	$$
	as in the Theorem. This concludes the case $\F=\R$. Now we note that being part of either set in the statement can be formulated as a first-order formula. Since over $\R$ the two formulae imply each other, they also imply each other over $\F$ by the transfer principle. This concludes the proof. 
\end{proof}

The following Lemma will be useful later, when we prove that the Iwasawa-retraction is distance-diminishing in Theorem \ref{thm:UAKret}.

\begin{lemma}\label{lem:kostant_gammaalpha}
	For all $\eta \in L = \sum_{\delta \in \Delta} \Z \delta $, we have
	$$
	\eta = \sum_{\ell=1}^r \frac{\langle \eta, \delta_\ell\rangle}{\langle \gamma_\ell,\delta_\ell\rangle} \gamma_\ell
	$$
	and for $\eta^+ := \sum_{\alpha \in \Sigma_{>0}} \alpha$, $\eta^+$ is a positive linear combination of the $\gamma_\ell$,
	$$
	\eta^+ = \sum_{\alpha \in \Sigma_{>0}} \alpha \in \sum_{\ell=1}^r  \Q_{>0} \gamma_\ell.
	$$

\end{lemma}
\begin{proof}
	By definition $\langle x_j, e_\ell \rangle = 0 $ for all $\ell \neq j$, but $\langle x_\ell, e_\ell \rangle \neq 0$, since otherwise $e_\ell = 0$. This implies that $\left\{e_1, \ldots, e_r\right\}$ and hence $\left\{ \gamma_1, \ldots, \gamma_r \right\}$ are linearly independent: for if $x=\sum_\ell \lambda_\ell e_\ell = 0$, then $\langle x_j, x\rangle = \lambda_j \langle x_j,e_j\rangle =0$, so $\lambda_j =0$ for all $j$.
	
	Therefore we can find $n_{i\ell} \in \Q$ such that
	$$
	\delta_i = \sum_{\ell=1 }^r n_{i\ell}\gamma_\ell.
	$$
	Since $\langle \gamma_\ell,\delta_k\rangle = 0 $ when $\ell \neq k$, we have $\langle \delta_i,\delta_k \rangle = n_{ik}\langle \gamma_k, \delta_k \rangle$. Since $x_\ell \in \fraka_p$, $\langle e_\ell,x_\ell\rangle \geq 0$ and since $e_j\neq 0$ and $\langle e_\ell, x_j \rangle = 0$ when $\ell\neq j$, actually $\langle e_\ell, x_\ell\rangle >0$. Thus also $\langle \gamma_\ell,\delta_\ell \rangle >0$, so we can divide $\langle \delta_i , \delta_\ell \rangle = n_{i\ell} \langle \gamma_\ell , \delta_\ell \rangle$ by $\langle \gamma_\ell , \delta_\ell \rangle $ to get
	$$
	n_{i\ell} = \frac{\langle \delta_i ,\delta_\ell \rangle}{\langle \gamma_\ell, \delta_\ell \rangle }.
	$$
	For $\eta \in L = \sum_{\delta \in \Delta} \delta$ we have
	\begin{align*}
		\eta &= \sum_{i=1}^r \lambda_i \delta_i =  \sum_{i=1}^r \sum_{\ell=1}^r \lambda_i \frac{\langle \delta_i ,\delta_\ell \rangle}{\langle \gamma_\ell, \delta_\ell \rangle }  \gamma_\ell \\
		&= \sum_{\ell=1}^r \frac{\langle \sum_{i=1}^r \lambda_i \delta_i, \delta_\ell\rangle}{\langle \gamma_\ell, \delta_\ell \rangle} \gamma_\ell = 
		\sum_{\ell=1}^r \frac{\langle \eta, \delta_\ell \rangle }{\langle \gamma_\ell, \delta_\ell \rangle} \gamma_\ell.
	\end{align*}
	
	For $\delta \in \Delta$, recall that the reflection $\sigma_\delta$ permutes the elements of $\Sigma_{>0} \setminus \{\delta, 2\delta\}$, \cite[VI.1.6 Cor. 1]{Bou08}. Let
	$$
	\eta^+ = \sum_{\alpha \in \Sigma^{+}}\alpha.
	$$
	If $\sigma_\delta (\eta - \delta) = \eta - \delta$, for instance when $\Sigma$ is reduced, then we can use that the reflection $\sigma_\delta$ preserves the scalar product to obtain
	\begin{align*}
		\langle \eta , \delta \rangle &= \langle \sigma_\delta(\eta) , -\delta \rangle = \langle \sigma_\delta\left( \eta - \delta \right) + \sigma_\delta(\delta), -\delta \rangle \\
		&=  \langle (\eta - \delta) - \delta , -\delta \rangle 
		= -\langle \eta, \delta \rangle + 2 \langle \delta , \delta\rangle. 
	\end{align*}
	so $\langle \eta, \delta \rangle = \langle \delta, \delta \rangle >0$. If $\Sigma$ is not reduced and $\delta,2\delta \in \Sigma$, then $\sigma (\eta - 3 \delta) = \eta - 3 \delta$ and 
	\begin{align*}
		\langle \eta , \delta \rangle &= \langle \sigma_\delta(\eta) , -\delta \rangle = \langle \sigma_\delta\left( \eta - 3\delta \right) + \sigma_\delta(3\delta), -\delta \rangle \\
		&=  \langle (\eta - 3\delta) - 3\delta , -\delta \rangle
		= - \langle \eta, \delta \rangle + 6 \langle \delta , \delta \rangle , 
	\end{align*}
	so $\langle \eta, \delta \rangle = 3\langle \delta ,\delta \rangle >0$. This shows that $\eta$ is a positive linear combination of elements $\gamma_\ell$.
\end{proof}

\begin{figure}[h]
	\centering
	\includegraphics[width=0.95\linewidth]{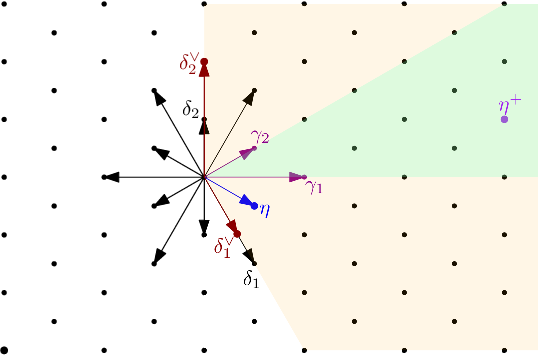}
	\caption{ The root system of type $\operatorname{G}_2$ with basis $\delta_1,\delta_2$, the coroots $\delta_1^\vee, \delta_2^\vee $ corresponding to the coroots $x_1,x_2$ and the elements $\gamma_1, \gamma_2$ spanning the Weyl chamber $\overline{(\mathfrak{a}^\star)}^+$. The element $\eta^+ := \sum_{\alpha >0} \alpha$ from Lemma \ref{lem:kostant_gammaalpha} lies in $\overline{(\mathfrak{a}^\star)}^+$, while the other candidate $\eta := \delta_1 + \delta_2$ does not.}
	\label{fig:G2_sum}
\end{figure}

Note that for the slightly simpler $\eta = \sum_{\delta \in \Delta} \delta$, $\eta$ may not be a positive linear combination of elements in $\sum \Z \gamma_\ell$, see Figure \ref{fig:G2_sum} for an example of type $\operatorname{G}_2$.
\\

 \newpage
\section{Definition of the building $\B$}\label{sec:building_def}
Let $\K$ and $\F$ be real closed fields such that $\K \subseteq \R \cap \F$. Often we assume $\F$ to be non-Archimedean with order compatible valuation $v \colon \F\to \Lambda \cup \{\infty\}$. Let $G$ be a semisimple connected self-adjoint algebraic $\K$-group and $S$ a maximal $\K$-split torus that satisfies $s=s\tran$ for all $s\in S$. In addition to the setting of Section \ref{sec:decompositions}, we assume that $G < \operatorname{SL}_n$ for some $n \in \N$. Let $A_{\F}$ be the semialgebraic extension of the semialgebraic connected component of $S_{\K}$ containing the identity and let $K=G\cap \operatorname{SO}_n$. 

For $\fraka = \operatorname{Lie}(A_\R)$, let $\Sigma \subseteq \fraka^\star$ be the root system whose elements $\alpha \in \Sigma$ correspond to $\K$-roots $\chi_\alpha \in \KPhi \subseteq \hat{S}$, see Section \ref{sec:compatibility}. Then $W_s=\KW$ is its spherical Weyl group. After choosing a basis $\Delta \subseteq \Sigma$ we let $U$ be the unipotent group associated to the positive root spaces and
$$
A_\F^+ = \left\{ a\in A_\F \colon \chi_\alpha(a) \geq 1 \text{ for all } \alpha \in \Delta \right\}.
$$

\subsection{Non-standard symmetric spaces}

In the theory of symmetric spaces, 
$$
P_\R = \left\{ x \in \R^{n\times n} \colon x=x\tran, \det(x)=1, x \text{ is positive definite }  \right\} 
$$
is a model for the symmetric space of non-compact type associated to $\operatorname{SL}(n,\R)$. The group $\operatorname{SL}(n,\R)$ acts transitively on $P_\R$ by
$$
g.x = g x g^{T}.
$$
for $g\in \operatorname{SL}(n,\R)$, $x \in P_\R$. The orbit $\X_\R = G_\R.\operatorname{Id} \subseteq P_\R$ is a closed subset and a model for the symmetric space associated to $G_\R$. We note that $P$ is a semialgebraic set defined over $\Q$ and consider its $\F$-extension $P_\F$. The action is algebraic, so the orbit can be semialgebraically extended to
$$
\X_\F = G_\F . P_\F.
$$
When $\F$ is non-Archimedean, we call $\X_\F$ the \emph{non-standard symmetric space associated to $G_\F$}. 
\begin{proposition}\label{prop:trans}
	\begin{itemize}
		\item [(a)] The group $G_\F$ acts transitively on $\X_\F$.
		\item [(b)] The stabilizer of $\operatorname{Id}\in \X_\F$ is $K_\F$.
		\item  [(c)] For any pair $x,y \in \X_\F$ there is a $g\in G_\F$ such that $g.x=\operatorname{Id}$ and $g.y$ lies in the closed Weyl chamber
		$$
		\ApF. \operatorname{Id} = \left\{ a.\operatorname{Id} \in \X_\F \colon \chi_\alpha(a) \geq 1 \text{ for all } \alpha \in \Sigma \right\}.
		$$
	\end{itemize}
\end{proposition}
\begin{proof}
	Transitivity and stabilizer of $\operatorname{Id}$ follow directly from the definitions. Use transitivity obtain $h \in G_\F$ with $h.x  =\operatorname{Id}$. Use transitivity again to obtain $h'\in G_\F$ with $h'.h.y =\operatorname{Id}$. Now decompose $h'=kak' \in G_\F = K_\F A_\F K_\F$ using the Cartan decomposition Theorem \ref{thm:KAK}, where we may assume $a^{-1}\in \ApF$ after applying an element of the spherical Weyl group $W_s$. Setting $g=k'h$ now results in the claimed
	\begin{align*}
		g.x &= k'h.x = k'.\operatorname{Id} = \operatorname{Id} \\
		g.y &= k'h.y = k'(h')^{-1}.\operatorname{Id} = a^{-1}k^{-1}.\operatorname{Id} = a^{-1}.\operatorname{Id} \in \ApF . \operatorname{Id}. 
	\end{align*}
\end{proof}

\subsection{The pseudo-distance}\label{sec:pseudodistance}

The symmetric space $\X_\R$ admits an explicit distance formula. Here we mimic this process to define a pseudo-distance on $\X_\F$. Let $x,y \in \X_\F$ be two points. We will first send $x$ and $y$ to a common flat, on which we define a multiplicative norm $N_\F$. In the real case, the logarithm would then be applied to obtain an additive distance. For the non-standard symmetric space we instead use the valuation $v \colon \F_{>0} \to \Lambda$.

By Proposition \ref{prop:trans}(a), there is for every $x,y\in \X_\F$ a $g\in G_\F$ with Cartan decomposition $g=kak' \in G_\F = K_\F A_\F K_\F$, Theorem \ref{thm:KAK}, such that $g.x = y$. The corresponding $a\in A_\F$ can be chosen to lie in $\ApF$ and is then unique. We can view this process as a map, called the \emph{Cartan projection}.

\begin{lemma}\label{lem:Cartan_projection}
	The Cartan-projection 
	\begin{align*}
		\delta_\F \colon \X_\F \times \X_\F & \to \ApF\\
		(x,y) &\mapsto a
	\end{align*}
	is well-defined and invariant under the action of $G_\F$. For all $x,y \in \X_\F$, $\delta_\F(y,x)$ is in the Weyl-group orbit of $\delta(x,y)^{-1}$.
\end{lemma} 
\begin{proof}
	We first assume $x=\operatorname{Id}$. If both $g=ka\bar{k}$ and $g' = k'a'\bar{k}'$ satisfy $y=g.x=g'.x$ with $a,a' \in \ApF$ and $k,\bar{k},k',\bar{k}' \in K_\F$, then we can apply the Cartan decomposition, Theorem \ref{thm:KAK}, to $kaa\tran k\tran = ka\bar{k}.\operatorname{Id} = k'a'\bar{k}'.\operatorname{Id} = k'a'(a')\tran(k')\tran$ to get that $aa\tran = a'(a')\tran $ as elements in $ \ApF$. Note that $aa\tran \in \ApF$, since $a=a\tran$, by our assumption on the torus $S$. Since $A_\F$ is an abelian group without torsion, $aa\tran = a'(a')\tran$ implies $a=a'$, showing that $\delta_\F(\operatorname{Id},g.\operatorname{Id}) = a$ is well defined.
	
	For general $x,y$, we can always find an $h\in G_\F$ such that $x=h.\operatorname{Id}$ and $y=g.h.\operatorname{Id}$. As the definition of $\delta_\F$ only depends on $g$ and not on $h$ in this case, $\delta_{\F}$ is invariant under $G_\F$, and thus well-defined everywhere. 
	
	If $g=kak'$ satisfies $g.x = y$, then $g^{-1} = (k')^{-1}a^{-1}k^{-1}$ satisfies $g^{-1}.y = x$. By definition, the unique element of $A_\F^+$ in the Weyl orbit of $a^{-1}$ is $\delta_\F(y,x)$.  
\end{proof}

We use a basis $\Delta$ of the root system $\Sigma$ to define a notion of positive roots $\Sigma_{>0}$. For $\alpha \in \Sigma$, let $\chi_\alpha \colon A_\R \to \R^\times$ be the corresponding ($\R$-points of the) algebraic characters. We define a continuous, semialgebraic $W_s$-invariant map
\begin{align*}
	N_{\R} \colon A_{\R} &\to \R^{\times} \\
	a & \mapsto  \prod_{\alpha \in \Sigma} \max \left\{\chi_\alpha(a), \chi_\alpha (a)^{-1}\right\}
\end{align*}
which is a multiplicative norm, meaning that for all $a,b \in A_{\R}$
\begin{itemize}
	\item [(1)] $N_{\R}(a)\geq 1$ and $N_{\R}(a) = 1$ if and only if $a=\operatorname{Id}$,
	\item [(2)] $N_{\R}(ab) \leq N_{\R}(a) N_{\R}(b)$.
\end{itemize}
We call $N_\R$ the \emph{semialgebraic norm}\footnote{There are many continuous semialgebraic multiplicative norms satisfying (1) and (2), but as norms on finite dimensional vector spaces, they are equivalent and it suffices for our purposes to fix $N_\R$.}.
Since $N_{\R}$ is semialgebraic, we can extend it to a map $N_{\F} \colon A_{\F} \to \F^{\times}$ which is still a $W_s$-invariant multiplicative norm satisfying (1) and (2) by the transfer principle and $N_\F$ is given by the same formula involving the characters.

For $\F$ non-Archimedean, we now use the Cartan projection $\delta_\F$ together with the semialgebraic norm $N_\F$ and the valuation $v \colon \F_{>0} \to \Lambda $ to define
\begin{align*}
	d \colon X_\F \times X_\F &\to \Lambda \\
	(x,y) &\mapsto (-v)(N_\F (\delta_{\F}(x,y))).
\end{align*}

We will show in Theorem \ref{thm:pseudodistance} that $d$ is a pseudo-distance on $\X_\F$. The pseudo-distance $d$ fails to be positive definite essentially due to the fact that $v$ is not injective. The proof of the triangle inequality uses Kostant's convexity theorem and the Iwasawa retraction
\begin{align*}
	\rho \colon \X_\F &\to A_\F.\operatorname{Id} \\
	g.\operatorname{Id} = uak.\operatorname{Id} & \mapsto a.\operatorname{Id}.
\end{align*}
using the Iwasawa-decomposition $G_\F = U_\F A_\F K_\F$, Theorem \ref{thm:KAU}.

\begin{lemma}\label{lem:rhoA}
	For all $a\in A_\F$, $x\in \X_\F$, $\rho(a.x) = a.\rho(x)$.
\end{lemma}
\begin{proof}
	Let $g=ua'k\in G_\F=U_\F A_\F K_\F$ such that $g.\operatorname{Id} = x$. By Proposition \ref{prop:anainN}, $aua^{-1} \in U_\F$, so
	\begin{align*}
		\rho(a.x) &= \rho(aua'k.\operatorname{Id}) = \rho((aua^{-1}aa'.\operatorname{Id}) = aa'.\operatorname{Id} = a .\rho(ua'k.\operatorname{Id}) = a .\rho(x).
	\end{align*}
\end{proof}

We use Kostant's convexity theorem to prove that $\rho$ is a $d$-diminishing retraction. 
\begin{theorem}\label{thm:UAKret}
	The map $\rho \colon \X_\F \to A_\F.\operatorname{Id}$ is a $d$-diminishing,
	$$
	\forall x,y \in \X_\F \colon d(\rho(x),\rho(y)) \leq d(x,y),
	$$ 
	retraction to $A_\F.\operatorname{Id}$.
\end{theorem} 
\begin{proof}
	It is clear that $\rho$ is a retraction, meaning that $\forall a \in A_\F, \rho(a.\operatorname{Id})=a.\operatorname{Id}$. Let $a_\F(g) = a_\F(uak)=a$ denote the $A_\F$-component of $g\in G_\F$ as in Section \ref{sec:kostant} on Kostant convexity. To show that $\rho$ is $d$-diminishing, we first claim that for all
	$
	b\in \ApF
	$
	and for all $k \in K_\F$ 
	$$
	d(\operatorname{Id}, \rho(kb.\operatorname{Id})) \leq d(\operatorname{Id}, b.\operatorname{Id}).
	$$
	For $b\in \ApF$ the set
	$$
	S_\F^b= \left\{ a\in A_\F \colon a = a_{\F}(kb) \text{ for some } k \in K_\F \right\}
	$$
	is semialgebraic and since over $\R$, $S_\R^b$ is closed under the action of the spherical Weyl group (this is a consequence of the real Kostant convexity Theorem \ref{thm:kostant_R}). The statement $W(S_\R^b) \subseteq S_\R^b$ can be formulated as a first-order formula, so $S_\F^b$ is also closed under the action of $W_s$. Note that $\rho(kb.\operatorname{Id}) = a_\F(kb).\operatorname{Id}$. While $a_\F(kb)$ may not lie in $\ApF$, there is a $w\in W$ such that $w(a_\F(kb)) \in \ApF$. We apply Theorem \ref{thm:kostant_F} to $w(a_\F(kb)) \in \ApF \cap S_\F^b$ get 
	$$
	\chi_{\gamma_i} (w(a_\F(kb))) \leq \chi_{\gamma_i}(b)
	$$
	for all $\gamma_i$ described in Section \ref{sec:kostant}. By Lemma \ref{lem:kostant_gammaalpha} for every $\alpha \in \Sigma$ there are positive rational numbers $n_{\alpha i} \in \Q_{>0}$ such that
	$$
	\alpha = \sum_{i=1}^r n_{\alpha i}\gamma_i.
	$$
	We can now prove
	\begin{align*}
		d(\operatorname{Id}, \rho(kb.\operatorname{Id})) &= \prod_{\alpha \in \Sigma_{>0}} \max \left\{ \chi_\alpha (a_\F(kb)), \chi_\alpha (a_\F(kb))^{-1}\right\} \\
		&= \prod_{\alpha \in \Sigma_{>0}}  \chi_\alpha (w(a_\F(kb))) \\
		&=  \prod_{\alpha \in \Sigma_{>0}} \prod_{i=1}^r  \chi_{ \gamma_i } (w(a_\F(kb)))^{n_{\alpha i}} \\
		&\leq \prod_{\alpha \in \Sigma_{>0}} \prod_{i=1}^r \chi_{\gamma_i} (b)^{n_{\alpha i}} \\
		&= \prod_{\alpha \in \Sigma_{>0}}  \chi_{\alpha}(b) = d(\operatorname{Id},b.\operatorname{Id}),
	\end{align*}
	where we used that $n_{\alpha i}>0$, proving the claim.
	
	Now let $x,y \in \X_\F$ arbitrary. By Proposition \ref{prop:trans}(c) and Theorem \ref{thm:KAU}, we can find $g=uak\in G_\F = U_\F A_\F K_\F$ such that $x=g.\operatorname{Id}$ and $y=g.b.\operatorname{Id}$ for some $b \in A_\F$. Now we use Lemma \ref{lem:rhoA} and the above to conclude 
	\begin{align*}
		d(\rho(x),\rho(y)) &= d(\rho(uak.\operatorname{Id}), \rho(uakb.\operatorname{Id})) 
		= d(a.\operatorname{Id}, \rho(akb.\operatorname{Id})) \\
		&= d(\operatorname{Id}, a^{-1}.\rho(akb.\operatorname{Id})) 
		= d(\operatorname{Id}, \rho(kb.\operatorname{Id})) \\
		&\leq d(\operatorname{Id}, b.\operatorname{Id}) = d(g.\operatorname{Id}, g.b.\operatorname{Id}) = d(x,y)
	\end{align*}
	concluding the proof.
\end{proof}

\begin{theorem}\label{thm:pseudodistance}
	The function $d \colon \X_\F \times \X_\F \to \Lambda$ is a pseudo-distance.
\end{theorem}
\begin{proof}
	We note that by definition, $N_\F(a) = N_\F(a^{-1})$ for all $a\in A_\F$. By Weyl group invariance and the last part of Lemma \ref{lem:Cartan_projection} we then obtain $N_\F(\delta_\F(x,y)) = N_\F(\delta_\F(y,x))$ for all $x,y\in \X_\F$, whence $d$ is symmetric.
	Since $N_\F(A_\F) \subseteq \F_{\geq 1}$, $(-v)(1)=0$ and $-v$ is order-preserving, $d$ is positive. It is also clear that $d(x,x) = 0$ for all $x\in \X_\F$. It remains to show that the triangle inequality holds. We start by analyzing the distance on the non-standard maximal flats $A_\F.\operatorname{Id}$. Let $a,b,c \in A_\F$, then we can use property (2) of the semialgebraic norm $N_\F \colon \A_\F \to \F_{>0}$ to deduce
	\begin{align*}
		d(a.\operatorname{Id},b.\operatorname{Id}) & = d(\operatorname{Id}, a^{-1}b.\operatorname{Id}) = (-v)(N_\F(a^{-1}b)) \\
		&= (-v)(N_\F(a^{-1}c c^{-1}b)) 
		\leq -(v)(N_\F(a^{-1}c) N_\F(c^{-1}b)) \\
		&= -v(N_\F(a^{-1}c)) + (-v)(N_\F(c^{-1}b)) \\
		&= d(a.\operatorname{Id},c.\operatorname{Id}) + d(c.\operatorname{Id},b.\operatorname{Id}),
	\end{align*}
	which settles the triangle inequality for points in $ A_\F.\operatorname{Id}$. For the general case we use the Iwasawa retraction from Theorem \ref{thm:UAKret}, as suggested in \cite[Lemma 1.2]{KrTe}. Let $x,y,z \in \X_\F$. By Proposition \ref{prop:trans}(c) there is a $g\in G_\F$ with $g.x ,g.y \in A_\F.\operatorname{Id}$. Then
	\begin{align*}
		d(x,y) & = d(g.x,g.y) = d(\rho(g.x),\rho(g.y)) \\
		& \leq d(\rho(g.x),\rho(g.z)) + d (\rho(g.z),\rho(g.y)) \\
		&\leq d(g.x,g.z) + d(g.z,g.y) = d(x,z) + d(z,y) 
	\end{align*}
	concludes the proof.
\end{proof}

\subsection{The apartment}\label{sec:Xapartment}

Over the reals, the orbit $A_\R.\operatorname{Id} \subseteq \X_\R$ is a maximal flat in the symmetric space $\X$. We take a closer look at the group $A_\F$ and its orbit $A_\F.\operatorname{Id} \subseteq \X_\F$ to define a space $A_\Lambda$ which will play the role of the model apartment for an affine $\Lambda$-building.
Let $O$ be an order convex valuation ring of the non-Archimedean real closed field $\F$ and $(-v) \colon \F \to \Lambda \cup \{\infty\}$ the associated order preserving valuation. 
We define the group
$$
A_{\Lambda} = A_{\F}/\{ a \in A_{\F} \colon N_{\F}(a) \in O \}.
$$
The goal of this section is to prove Theorem \ref{thm:ALisoA} which states that $A_{\Lambda}$ can be given the structure of a model apartment $\A = \A(\KPhi,\Lambda,A_{\Lambda})$, as defined in Section \ref{sec:modelapartment}.

The spherical Weyl group $W_s = N_\F/M_\F$ acts on $A_\F$ by $[k].a = kak^{-1}$ for $a\in A_\F, k\in N_\F$. Since $N_\F$ is $W_s$-invariant, the action descends to an action by automorphisms on $\A_\Lambda$.

Recall from Proposition \ref{prop:char_cochar} that there is a non-degenerate bilinear form 
\begin{align*}
	\hat{S} \times X_{\star}(S) &\to \Z \\
	(\chi, t ) &\mapsto b(\chi,t ) 
\end{align*}
with the defining property that $\chi \circ t (x) = x^{b( \chi,t )}$ for all $x \in \G_m$. 
Since $S$ is $\K$-split and hence $\F$-split by Theorem \ref{thm:split_tori}, we can take the $\F$-points and restrict to $A_\F$ and $\F_{>0}$ to obtain characters $\chi_\F \colon A_\F \to \F_{>0}$ and one-parameter subgroups $t_\F \colon \F_{>0} \to A_\F$ for $\chi \in \hat{S}$ and $t \in X_\star(S)$. Note that since $\chi_\F$ is continuous and $A_\F$ is semialgebraically connected, $\chi_\F(A_\F) \subseteq \F_{>0}$ and similarly since $t_\F$ is continuous and $\F_{>0}$ is connected, $t_\F(\F_{>0}) \subseteq A_{\F}$. We now use the valuation $(-v) \colon \F_{>0} \to \Lambda$ to set up tools to prove that $A_\Lambda$ is isomorphic to the model apartment $\A$.

\begin{proposition}\label{prop:char_cochar_Lambda}
	Let $\chi \in \operatorname{Span}_\Z(\KPhi)$ and $t \in X_\star(S)$. The $\F$-valued characters $\chi_\F$ and one-parameter subgroups $t_\F$ descend to group homomorphisms $\chi_\Lambda$ and $t_\Lambda$ such that the diagram
	$$
	\begin{tikzcd}
		\mathbb{F}_{>0} \arrow[d, "-v"] \arrow[r, "t_{\mathbb{F}}"] & A_{\mathbb{F}} \arrow[d,two heads] \arrow[r, "\chi_{\mathbb{F}}"] & \mathbb{F}_{>0} \arrow[d, "-v"]  \\
		\Lambda \arrow[r, "t_{\Lambda}"]                             & A_\Lambda \arrow[r, "\chi_\Lambda"]                     & \Lambda                                 
	\end{tikzcd}
	$$
	commutes and such that 
	$
	\chi_\Lambda \circ t_\Lambda (\lambda ) =b( \chi, t ) \cdot \lambda.
	$
\end{proposition}
\begin{proof}
	We denote by $\pi \colon A_\F \to A_\Lambda = A_\F / \{ a\in A_\F \colon N_\F(a) \in O \}$ the projection. We first show that $t_\Lambda \colon \Lambda \to A_\Lambda, \, (-v)(x) \mapsto \pi(t_\F (x))$ is well defined: we have to show that if $(-v)(x) = 0$, then $N(t_\F(x)) \in O$. Indeed, let $(-v)(x)=0$, so $x \in O^\times$. Then $ \alpha_\F (t_\F (x)) = x^{b( \alpha, t )} \in O^\times$ for all $\alpha \in \KPhi$ and thus
	\begin{align*}
		N_\F(t_\F(x)) & = \prod_{\alpha \in \KPhi} \max \left\{ \alpha_\F( t_\F(x)) , \alpha_\F(t_\F(x))^{-1} \right\} \in O.
	\end{align*}
	
	To show that $\chi_\Lambda \colon A_{\Lambda} \to \Lambda, \, \pi(a) \mapsto v(\chi_\F(a))$ is well-defined, we have to show that if $N_\F(a) \in O$, then $(-v)(\chi_\F(a)) =0$. Indeed, if $N_\F(a) \in O$, then $\alpha_\F(a) \in O$ for all $\alpha \in \KPhi$ since
	\begin{align*}
		N_\F(a) & = \prod_{\alpha \in \KPhi} \max \left\{ \chi_\alpha( a) , \chi_\alpha(a)^{-1} \right\} \in O
	\end{align*}
	is a product of elements that are larger than $1$. Since $\chi \in \operatorname{Span}_\Z(\KPhi)$, $\chi_\F(a) $ is a product of $\alpha_\F(a) \in O$ and hence in $O$.
	
	The maps $t_\Lambda$ and $\chi_\Lambda$ were defined so that the diagrams commute. For $\lambda \in \Lambda$, we can find $x \in \F_{>0}$ such that $(-v)(x)=\lambda $. Then
	\begin{align*}
		\chi_\Lambda \circ t_\Lambda (\lambda) &= \chi_\Lambda (\pi(t_\F(x ) ) )  
		= (-v)\left( \chi_\F \circ t_\F (x) \right) \\
		&= v\left(x^{b( \chi, t )} \right)
		= b( \chi , t) \cdot (-v)(x)  .
	\end{align*}
\end{proof}
Let $\Delta \subseteq \KPhi$ be a basis and define for $\delta \in \Delta$ a coroot $t_\delta \in X_\star(S)$ as in Lemma \ref{lem:oneparam_compatibility}. For general $\chi \in \operatorname{Span}_\Z (\KPhi) \subseteq \hat{S}$, we use bilinearity to define $t_\chi \in X_\star(S)$ via
\begin{align*}
	t_\chi (x) = \prod_{\delta \in \Delta} t_\delta (x)^{\lambda_\delta } \quad \text{ for } \quad
	\chi =  \prod_{\delta \in \Delta} \chi_\delta^{\lambda_\delta} \in \operatorname{Span}_{\Z}(\KPhi),\quad  \lambda_\delta \in \Z,
\end{align*}
and we note that for all $\chi, \nu \in \operatorname{Span}_{\Z}(\KPhi) $, $\chi(t_{\nu}(x)) = x^{b(\chi,t_{\nu})}$ for $x \in G_m$.

Let $r$ denote the rank of the root system $\KPhi$.
\begin{proposition}\label{prop:isoFrA}
	There is a group isomorphism 
	\begin{align*}
		f_\F \colon (\F_{>0})^r & \to A_{\F} \\
		(x_\delta)_{\delta \in \Delta} &\mapsto \prod_{\delta \in \Delta} (t_{\delta})_\F(x_\delta)
	\end{align*} 
	satisfying
	$$
	\delta (f_\F(x_1, \ldots , x_r)) = x_\delta
	$$
	for all $(x_1, \ldots , x_r) \in (\F_{>0})^r$ and $\delta \in \Delta$.
\end{proposition}
\begin{proof}
	Note that $f$ is a semialgebraic map defined over $\K$. We may write $f_\R$ or $f_\F$ for the respective semialgebraic extensions, when necessary, but forego the $\F$-subscript when there is not danger of confusion. It is clear that $f$ is a homomorphism, we will now find an inverse.
	We start by giving names to the elements of the basis $\Delta = \{\delta_1, \ldots , \delta_r \}$. We define the $r\times r$ integer matrix $M$ by $M_{ij} = b(\delta_i, t_{\delta_j} )$. Since $M$ represents a non-degenerate bilinear form, it is invertible and its inverse $M^{-1}$ has rational entries. For $j=1, \ldots , r$, we define
	$$
	\chi_j = \prod_{k=1}^r \delta_k^{(M^{-1})_{jk}} \in  \operatorname{Span}_{\Q}(\KPhi) 
	$$
	and use them to give an explicit right inverse
	\begin{align*}
		f^{-1} \colon A_{\F} & \to (\F_{>0})^r  \\
		a & \mapsto (\chi_1(a), \ldots , \chi_r(a)).
	\end{align*} 
	We calculate the $j$th entry 
	\begin{align*}
		f^{-1}(f(x_1, \ldots , x_r))_j &= f^{-1} \left( \prod_{i=1}^r t_{\delta_i} (x_i) \right)_j 
		= \chi_j\left( \prod_{i=1}^r t_{\delta_i} (x_i) \right) \\
		&=  \prod_{i=1}^r \chi_j\left( t_{\delta_i} (x_i)\right) 
		=  \prod_{i=1}^r  \prod_{k=1}^r \delta_k^{(M^{-1})_{jk}}\left(t_{\delta_i} (x_i) \right)\\
		& =  \prod_{i=1}^r  \prod_{k=1}^r (\delta_k \circ t_{\delta_i} )(x_i)^{(M^{-1})_{jk}}
		= \prod_{i=1}^r  \prod_{k=1}^r \left((x_i)^{b( \delta_k , t_{\delta_i})}\right)^{(M^{-1})_{jk}} \\
		&= \prod_{i=1}^r (x_i)^{\sum_{k=1}^r (M^{-1})_{jk} M_{ki}} = \prod_{i=1}^r (x_i)^{\operatorname{Id}_{ji} }
		= x_j.
	\end{align*}
	The existence of the right inverse proves that $f$ is surjective. We will now argue why $f$ is also injective.
	
	Note that $(\R_{>0})^r$ is a $r$-dimensional $\R$-vector space, when equipped with addition defined by multiplication $(v,w) \mapsto vw$ and scalar multiplication defined by $(\lambda,v) \mapsto e^{\lambda}v$, where $v,w \in \R_{>0}, \lambda \in \R$. Similarly, $A_\R = \exp(\fraka)$ is a $\R$-vector space since $\fraka$ is. By Proposition \ref{prop:rootsystem_compatible}, $A_{\R}$ is also $r$-dimensional. The map $f_\R$ is a surjective linear map between vector spaces of the same finite dimension, and hence is bijective.
	
	To get injectivity for $f_\F$, we write it as a first-order statement
	$$
	\varphi \colon \quad \forall x_1, \ldots , x_r \colon \bigwedge_{i=1}^r x_i > 0 \wedge \left( f(x_1, \ldots , x_r) = \operatorname{Id} \to  \bigwedge_{i=1}^r x_i = 1 \right) ,
	$$
	where we used the fact that the $t_\delta$ and thus $f$ are algebraic. By the above, $\varphi$ holds over $\R$ and thus it holds over $\F$. This concludes the proof that $f_\F$ is an isomorphism.
\end{proof}
Finally we prove that $A_{\Lambda}$ is isomorphic to $\A \cong \Lambda^r$.
\begin{theorem}\label{thm:ALisoA}
	The group isomorphism $f_\F \colon (\F_{>0})^r \cong A_\F$ descends to a group isomorphism
	\begin{align*}
		f_\Lambda \colon \Lambda^r & \to A_\Lambda \\
		(\lambda_\delta)_{\delta \in \Delta} &\mapsto \sum_{\delta \in \Delta} (t_{\delta})_{\Lambda}(\lambda_\delta).
	\end{align*}
\end{theorem}
\begin{proof}
	As a consequence of Proposition \ref{prop:char_cochar_Lambda}, the diagram
	$$
	\begin{tikzcd}
		(\mathbb{F}_{>0})^r \arrow[r, "f_{\mathbb{F}}"] \arrow[d, "v^r"'] & A_{\mathbb{F}} \arrow[d] \\
		\Lambda^r \arrow[r, "f_\Lambda"']                                 & A_{\Lambda}             
	\end{tikzcd}
	$$
	commutes. Since $v$ and the projection $\pi \colon A_\F\to A_\Lambda$ are surjective, surjectivitiy of $f_\Lambda$ follows from the surjectivity of $f_\F$, Proposition \ref{prop:isoFrA}. It remains to show that $f_\Lambda$ is injective. 
	Let $\Delta = \left\{ \delta_1, \dots , \delta_r \right\}$ and $t_i = (t_{\delta_{i}})_{\F}$. For $a \in \A_\F$ with $N_{\F}(a) \in O$, we showed in the proof of Proposition \ref{prop:char_cochar_Lambda}, that $v(\delta_j(a)) = 0$ for all $j \in \{1, \ldots , r\}$. We use the inverse function $f^{-1}_\F$ constructed in Proposition \ref{prop:isoFrA} to obtain
	\begin{align*}
		v\left(\left(f^{-1}_\F (a)\right)_{j}\right) &= v\left( \prod_{k=1}^r \delta_k^{M_{jk}^{-1}} (a) \right)
		= \sum_{k=1}^r M_{jk}^{-1} v(\delta_{k}(a)) = 0 
	\end{align*}
	where $M^{-1}_{jk}$ are some rational numbers. This shows $\ker(f_{\Lambda}) = 0 \in \Lambda^r$ and hence $f_\Lambda$ is injective.
\end{proof}

We now know that a choice of basis $\Delta \subseteq \KPhi$ leads to isomorphisms 
$$
\begin{tikzcd}
	A_{\Lambda} & \Lambda^r \arrow[l, "f_\Lambda"'] \arrow[r, "g"] & \A := \operatorname{Span}_{\mathbb{Z}}(\KPhi) \otimes_{\mathbb{Z}} \Lambda \\
\end{tikzcd}
$$
where
$$
g((\lambda_\delta)_{\delta \in \Delta}) = \sum_{\delta \in \Delta} \delta \otimes \lambda_\delta.
$$

For $A_\Lambda$ to be a model apartment of type $(\KPhi, \Lambda, \Lambda^r)$ as defined in Section \ref{sec:modelapartment}, the isomorphism $A_\Lambda \cong \A$ has to be Weyl group equivariant. The spherical Weyl group $W_s$ acts on $\KPhi$ and hence on $\A$ as described in Section \ref{sec:modelapartment}. The group $N_\F= \operatorname{Nor}_{K_\F}(A_\F)$ acts on $A_\F$ by isometries and hence on $A_\Lambda$. The spherical Weyl group $W_s$ is isomorphic to $N_\F/M_\F$ by Section \ref{sec:compatibility}. The translation group $T:=\Lambda^r$ acts on $A_\Lambda$ and $\A$ via $f_\Lambda$ and $g$ and this action coincides with the action on $\A$ defined in Section \ref{sec:modelapartment}.

\begin{proposition}\label{prop:A_equivariant_affine_Weyl_group}
	The isomorphism $g\circ f_\Lambda^{-1} \colon A_\Lambda \to \A$ is equivariant with respect to the actions of the affine Weyl group $W_s \ltimes T$.
\end{proposition}
\begin{proof}
	Let $w\in W_s$ be an element of the spherical Weyl group of $\KPhi$. For $\delta,\eta \in \Delta$, there are constants $k_{\delta,\eta} \in \Z$ such that $w(\eta) = \sum_{\delta \in \Delta} k_{\delta, \eta} \delta$. Then
	\begin{align*}
		w(g((\lambda_\delta)_\delta)) &= w \left(\sum_{\delta \in \Delta} \delta \otimes \lambda_\delta \right)
		= \sum_{\eta \in \Delta} w(\eta) \otimes \lambda_{\eta} 
		= \sum_{\eta \in \Delta} \sum_{\delta \in \Delta} k_{\delta,\eta} \delta \otimes \lambda_{\eta} \\
		&= \sum_{\delta\in \Delta} \delta \otimes \left( \sum_{\eta\in \Delta}k_{\delta,\eta}\lambda_\eta \right) 
		= g\left( \left(\sum_{\eta \in \Delta} k_{\delta,\eta} \lambda_\eta \right)_{\delta\in \Delta}\right).
	\end{align*}
	For $x\in \F_{>0}$ and $\eta \in \Delta$, we have $w((t_{\eta})_\F(x)) = (t_{w(\eta)})_\F(x)$ and hence
	\begin{align*}
		w(f_\Lambda((\lambda_\delta)_\delta))) &= w\left(\prod_{\eta\in\Delta}(t_\eta)_\Lambda(\lambda_\eta)\right)
		= \prod_{\eta\in\Delta}(t_{w(\eta)})_\Lambda(\lambda_\eta)
		= \prod_{\eta\in\Delta}(t_{\sum_{\delta} k_{\delta,\eta}\delta})_\Lambda (\lambda_\eta)
		\\	&= \prod_{\eta\in\Delta} \prod_{\delta\in\Delta} (t_\delta)_\Lambda(\lambda_\eta)^{k_{\delta,\eta}}
		= \prod_{\delta\in\Delta} (t_\delta)_\Lambda(\lambda_\eta)^{\sum_{\delta}k_{\delta,\eta}}
		= \prod_{\delta\in\Delta} (t_\delta)_\Lambda\left( \sum_{\eta}k_{\delta,\eta} \lambda_\eta\right)\\
		&= f_\Lambda \left(\left(  \sum_{\eta \in \Delta}k_{\delta,\eta} \lambda_\eta \right)_{\delta\in \Delta}\right)
	\end{align*}
	from which it follows that $g\circ f_\Lambda^{-1}$ is $W_s$-equivariant. The way we defined $T=\Lambda^r$ to act, $g\circ f_\Lambda^{-1}$ is by definition $W_s \ltimes T$-equivariant.
\end{proof}

\subsection{The building} \label{sec:building}

By Theorem \ref{thm:pseudodistance}, the non-standard symmetric space $\X_\F$ admits a $\Lambda$-pseudometric. We consider the quotient
$$
\B = \X_\F/\!\sim
$$
where $x\sim y \in \B$ when $d(x,y) = 0 \in \Lambda$. We denote the induced $\Lambda$-metric on $\B$ by the same letter $d$. We note that $G_\F$ acts by isometries on $\B$. In Section \ref{sec:Bisbuilding} we will show that the $\Lambda$-metric space $\B$ admits the structure of an affine $\Lambda$-building in certain cases, see Theorem \ref{thm:B_is_building}.

We denote the equivalence class of $\operatorname{Id}$ by $o\in \B$. The stabilizer of $o$ has been calculated by \cite{Tho} when $\F$ is a Robinson field and by \cite{KrTe} more generally. We first describe the special case of the action of $A_\F$ on $\B$. For any semi-algebraic subset $H_\F\subseteq G_\F$ we write $H_\F(O) := H_\F \cap O^{n\times n}$, where $O$ is the valuation ring.

\begin{proposition}\label{prop:stab_A}
	The following are equivalent for $a\in A_\F$.
	\begin{itemize}
		\item [(i)] $a \in \operatorname{Stab}_{G_\F}(o)$.
		\item [(ii)]  $N_\F(a) \in O$.
		\item [(iii)] $\chi_\alpha(a) \in O \ \forall \alpha \in \Sigma $.
		\item [(iv)]  $\chi_\alpha(a) \in O^\times \ \forall \alpha \in \Delta $.
		\item [(v)] $a \in A_\F(O)$.
	\end{itemize}
\end{proposition}
\begin{proof}
	If $a\in \operatorname{Stab}_{G_\F}(o)$, then $d(\operatorname{Id},a.\operatorname{Id}) = (-v)(N_\F(a)) = 0$, so $N_\F(a) \in O$. This implies by the definition of $N_\F$ that $\chi_\alpha(a) \in O \cap O^{-1} = O^\times$ for all $\alpha \in \Sigma$, in particular for all $\alpha \in \Delta$. Now use the algebraic map $f_\F \colon (\F_{>0})^r \to A_\F$ from Proposition \ref{prop:isoFrA} to write 
	$$
	a = f_\F((x_\delta)_{\delta\in \Delta}) = \prod_{\delta \in \Delta} t_\delta(x_\delta)
	$$
	for some $(x_\delta)_{\delta \in \Delta} \in (\F_{>0})^r$. Then
	$$
	\chi_\delta (a) = x_\delta \in O
	$$
	for all $\delta \in \Delta$. There is a constant $n_0 \in \N$ such that the first-order formula
	$$
	\varphi \colon \quad \forall x \ \left| t_\alpha(x)_{ij} \right| \leq n x^{n_0}, 
	$$
	holds over $\R$	(for $t_\alpha^\R(e^s) = \exp(s x_\alpha)$, take a natural number $n_0$ larger than $|(x_{\alpha})_{ij}|$ for all $i,j$.). 
	By the transfer principle, $\varphi$ holds over $\F$ as well. Since $x_\delta \in O$, $t_\delta(x_\delta) \in O^{n\times n}$ and hence $a \in A_\F(O)$.

On the other hand, if $a\in A_\F(O)$, then the linear map $\operatorname{Ad}(a) \colon \frakg \to \frakg, X \mapsto aXa^{-1}$ restricts to multiplication by $\chi_\alpha(a)$ on $\frakg_\alpha$, from which can be concluded that $\chi_\alpha(a) \in O$.	Then $N_\F(a) \in O$ and $d(o,a.o) =0$, so $a \in \operatorname{Stab}_{G_\F}(o)$.
\end{proof}

\begin{theorem}\label{thm:stab} 
	The stabilizer of $o\in\B$ in $G_\F$ is $G_\F(O) $.
\end{theorem}
\begin{proof}	
	Let $g=kak' \in G_\F = K_\F \Ap_\F K_\F$. If $g\in \operatorname{Stab}_{G_\F}(o)$, then $d(\operatorname{Id},g.\operatorname{Id}) = (-v)(N_\F(a)) = 0$. By Proposition \ref{prop:stab_A}, $a \in A_\F(O)$. Now since $K_\F = K_\F(O)$, $g = kak' \in G_\F(O)$.
	
	If on the other hand we start with a $g\in G_\F(O)$, then $a\in \A_\F(O)$, and by Proposition \ref{prop:stab_A}, $N_\F(a) \in O$ and thus also $d(\operatorname{Id},g.\operatorname{Id}) =d(\operatorname{Id},a.\operatorname{Id}) =0 $, hence $g \in \operatorname{Stab}_{G_\F}(o)$.
\end{proof}

As an application of the Iwasawa retraction, Theorem \ref{thm:UAKret}, we can give a group decomposition for the stabilizer of $o$ in $\B$.

\begin{corollary}\label{cor:UAKO}
	There is an Iwasawa decomposition $G_\F(O) = U_\F(O) A_\F(O) K_\F$, meaning that for every $g \in G_\F(O)$ there are unique $u \in U_\F(O), a\in A_\F(O), k\in K_\F=K_\F(O)$ with $g=uak$.
\end{corollary}
\begin{proof}
	Let $g=uak \in G_\F(O) \subseteq  U_\F A_\F K_\F$. We have $\rho(g.\operatorname{Id}) = a.\operatorname{Id}$. Since $g\in G_\F(O)$, we have by Theorem \ref{thm:UAKret}
	\begin{align*}
		d(\operatorname{Id},a.\operatorname{Id}) &= d(\operatorname{Id}, \rho(g.\operatorname{Id}))  \leq  d(\operatorname{Id}, g.\operatorname{Id}) = 0.
	\end{align*}
	This means that $a \in A_\F \cap G_\F(O) = A_\F(O)$. Note that since $K_\F$ stabilizes $\operatorname{Id}\in \X_\F$, $K_\F=K_\F(O)$. Since $G_\F(O)$ is a subgroup of $G_\F$, $u = gk^{-1}a^{-1} \in G_\F(O)$, so $u \in U_{\F}(O)$.
\end{proof}

\newpage

\section{Verification of the axioms for $\B$}\label{sec:Bisbuilding}

We continue in the setting of Sections \ref{sec:decompositions} and \ref{sec:building_def}: $\K$ and $\F$ are real closed fields such that $\K \subseteq \R \cap \F$ and $G$ is a semisimple connected self-adjoint algebraic $\K$-group. Let $\B = \X_\F/\!\sim$ as in Section \ref{sec:building} and $\A\cong A_\Lambda = A_\F.o  \subseteq \B$ as in Section \ref{sec:Xapartment}. Denote the inclusion by $f_0 \colon \A \to \B$. We define a set of charts
$$
\Fun = \left\{ g.f_0 \colon \A \to \B \colon g \in G_\F \right\}.
$$
The goal of this section is to show that $(\B,\Fun)$ is an affine $\Lambda$-building in the sense of Section \ref{sec:lambda_building}. 

\begin{theorem}\label{thm:B_is_building}
	If the root system $\Sigma$ of $G_\F$ is reduced, then the pair $(\B,\Fun)$ is an affine $\Lambda$-building of type $\A = \A(\KPhi,\Lambda,\Lambda^{n})$.
\end{theorem}

We will show the theorem by checking the set of axioms (A1), (A2), (A3), (A4), (TI), (EC), as described in Theorem \ref{thm:equivalent_axioms}. Axioms (A1), (A3) and (TI) are treated in Section \ref{sec:axiomsA1A3TI} and follow easily from what we have developed so far. In Section \ref{sec:root_groups_and_BT}, we will investigate the root groups and develop some theory paralleling some parts of the work of Bruhat-Tits \cite{BrTi}. We are then equipped to prove the remaining axioms in Sections \ref{sec:A2_subsubsection}, \ref{sec:A4} and \ref{sec:EC}: axiom (A2) is proven in Theorem \ref{thm:A2}, axiom (A4) in Theorem \ref{thm:A4} and axiom (EC) in Theorem \ref{thm:EC}. An alternative proof of axiom (A2) in the case of $G_\F = \operatorname{SL}(n,\F)$ is given in Appendix \ref{sec:appendixSLn}. The assumption that $\Sigma$ is reduced is directly used in the proof of axioms (A2) and (EC). Axiom (A4) uses the assumption indirectly, as it relies on the statement of (A2) in its proof. The remaining axioms (A1), (A3) and (TI) do not need the assumption.

\subsection{Axioms (A1), (A3) and (TI)}\label{sec:axiomsA1A3TI}

Three of the axioms follow from what we have done.

\begin{lemma}\label{lem:A1}
	The pair $(\B,\Fun)$ satisfies axiom
	\begin{enumerate}
		\item [(A1)] For all $f\in \Fun$, $w\in W_a$, $f \circ w \in \Fun$.
	\end{enumerate}
\end{lemma}
\begin{proof} 
	We have surjections $\operatorname{Nor}_{K_\F}(A_\F) \to W_s$ and $A_\F \to \A \cong \Lambda^n$, see Section \ref{sec:Xapartment}. 
	So for any element $(t,w) \in W_a = \Lambda^n \rtimes W_s$ we can find $a \in A_\F$ and $k\in \operatorname{Nor}_{K_\F}(A_\F)$ such that $ak.p = (t,w)(p)$ for all $p \in \A$. Let $f=g.f_0 \in \Fun$ for $g\in G_\F$. Now
	$$
	(f\circ (t,w) )(p) = g.f_0 (ak.p) = gak.p = (gak).f_0(p),
	$$
	for all $p \in \A$ and thus $f\circ (t,w) \in \Fun$, proving axiom (A1).
\end{proof}
\begin{lemma}\label{lem:A3}
	The pair $(\B,\Fun)$ satisfies axiom 
	\begin{enumerate}
		\item [(A3)] For all $p,q\in \B$, there is a $f\in \Fun$ such that $p,q\in f(\A)$.
	\end{enumerate}
\end{lemma}
\begin{proof}
	This follows from Proposition \ref{prop:trans}(c): if $[x],[y]\in \B$, for points $x,y \in \X_\F$, then there is a $g \in G_\F$ such that $g.x = \operatorname{Id} \in A_\F.\operatorname{Id} $ and $g.y \in A_\F.\operatorname{Id}$. Then $[x],[y] \in g^{-1}.f_0(\A)$ for the chart $g^{-1}.f_0 \in \Fun$.
\end{proof}

Axiom (TI) just states that $\B$ satisfies the triangle inequality. The triangle inequality was proven in Theorem \ref{thm:pseudodistance} using Kostant convexity and Iwasawa-retractions.

\subsection{Root groups and some Bruhat-Tits theory}\label{sec:root_groups_and_BT}

Before proving axiom (A2) in Section \ref{sec:A2_subsubsection}, we will develop some theory. The overarching theme of this subsection is that we want to understand pointwise stabilizers of subsets of $\A$. First we will investigate the action of the root groups $(U_\alpha)_\F$. This will give us partial results about $W_a$-convexity and allows us to describe the stabilizers of the fundamental Weyl chamber and the entire $\A$. We will then follow some ideas of \cite{BrTi, Lan96}, that allows us to describe the stabilizers of arbitrary finite subsets of $\A$.

For axiom (A4), it suffices to understand how points of the apartment $\A$ are fixed by elements of $U_\F$, which we treat in Section \ref{sec:Wconvexity_for_U}. The full development of this chapter is used in the proof of axioms (A2) and (EC), and relies on the fact that $\Sigma$ is reduced. We still develop much of the theory also in the case of non-reduced root systems and point out when we rely on the reduced case.

\subsubsection{Root group valuations} \label{sec:root_group_valuation}

Recall from Section \ref{sec:U}, that $U_\F$ has subgroups $(U_\alpha)_\F$ for $\alpha \in \Sigma_{>0}$ defined as $(U_\alpha)_\F = \exp((\frakg_\alpha)_\F \oplus (\frakg_{2\alpha})_\F)$. Since $(U_\alpha)_\F$ is unipotent, the matrix exponential $\exp \colon (\frakg_\alpha)_\F \oplus (\frakg_{2\alpha})_\F \to (U_\alpha)_\F$ and the matrix logarithm 
\begin{align*}
	\log \colon (U_\alpha)_\F &\to (\frakg_\alpha)_\F \oplus (\frakg_{2\alpha})_\F \\
	u & \mapsto \sum_{k=1}^\infty (-1)^{k+1} \frac{(u-\operatorname{Id})^k}{k}
\end{align*}
are just polynomials on $(U_\alpha)_\F$. Viewing $\frakg_\F \subseteq \mathfrak{sl}(n ,\F) \subseteq \F^{n\times n}$, we can speak about the matrix entries $Z_{ij}$ of $Z \in \frakg_\F$. Let 
$$
(-v)(Z) := \max_{ij} \{ (-v)(Z_{ij}) \}.
$$
Inspired by \cite[6.2]{BrTi} and in order to avoid talking about matrix entries too much we introduce the \emph{root group valuations} $\varphi_\alpha$
\begin{align*}
	\varphi_\alpha \colon (U_{\alpha})_\F &\to \Lambda  \cup \{ -\infty \}\\
	\exp(X+X') & \mapsto \max \left\{ \max_{i,j} \left\{(-v)(X_{ij}) \right\}, \frac{1}{2} \max_{i,j} \left\{ (-v)(X'_{ij}) \right\}\right\}
\end{align*}
where $X \in (\frakg_\alpha)_\F$ and $X' \in (\frakg_{2\alpha})_\F$. The following Lemma justifies their name.

\begin{lemma}\label{lem:BTgroup_valuation}
	For all $\alpha \in \Sigma$ and $u,v \in (U_{\alpha})_\F$
	$$
	\varphi_\alpha(uv) \leq \max \left\{ \varphi_\alpha (u) , \varphi_\alpha(v) \right\}
	$$
	and if $\varphi_\alpha(u) \neq \varphi_\alpha(v)$, then equality holds.
\end{lemma}
\begin{proof} We first claim that for all $X,Y \in \frakg_{\alpha}$,
	$$
	\frac{1}{2} (-v)\left(\left[X,Y\right]\right) \leq \max \{ (-v)(X), (-v)(Y) \}.
	$$
	Indeed in matrix entries, $[X,Y]_{ij} = \sum_{k} X_{ik}Y_{kj} - Y_{ik}X_{kj}$, and hence
	\begin{align*}
		\frac{1}{2}(-v)\left(\left[X,Y\right]\right)  &= \frac{1}{2} (-v)\left(\max_{ij} \left\{ \sum_{k} X_{ik}Y_{kj} - Y_{ik}X_{kj}  \right\} \right) \\
		&\leq \frac{1}{2} \max_{ijk}\left\{ (-v)(X_{ik} Y_{kj}) , (-v)(Y_{ik}X_{kj}) \right\} \\
		& \leq \frac{1}{2} \max_{ijk} \{ (-v)(X_{ik}) + (-v)(Y_{kj}), (-v)(Y_{ik}) + (-v)(X_{kj})  \} \\
		&\leq \frac{1}{2} \max \{ (-v)(X) + (-v)(Y) ,  
		(-v)(Y) + (-v)(X)  \} \\
		&\leq \frac{1}{2} \max \{ 2 (-v)(X), 2(-v)(Y) \} = \max \{ (-v)(X), (-v)(Y) \}.
	\end{align*}
	
	Now let $u = \exp(X+X')$ and $v=\exp(Y+Y')$ for  $X,Y \in (\frakg_\alpha)_\F$ and $X', Y' \in (\frakg_{2\alpha})_\F$. Then
	$$
	u\cdot v = \exp \left( X+Y + X'+Y'+ \frac{1}{2}[X,Y] \right)
	$$
	by the Baker-Campbell-Hausdorff-formula. By the claim
	\begin{align*}
		\varphi_\alpha(uv) &= \max \left\{ (-v)(X+Y) , \frac{1}{2} (-v)\left(X'+Y' + \frac{1}{2}[X,Y]\right) \right\} \\
		& \leq \max\left\{    (-v)(X), (-v)(Y) , \frac{1}{2}(-v)(X'), \frac{1}{2}(-v)(Y')    \right\} \\
		&= \max \{  \varphi_\alpha(u), \varphi_{\alpha}(v)  \}.
	\end{align*}
	
	If without loss of generality $\varphi_\alpha(u) < \varphi_\alpha(v)$, we distinguish two cases. If $\varphi_\alpha(v) = (-v)(Y) \geq \frac{1}{2}(-v)(Y')$, then $(-v)(X) < (-v)(Y)$ and $\frac{1}{2}(-v)(X') < (-v)(Y)$. Thus using the claim
	\begin{align*}
		\varphi_\alpha (uv) & = \max \left\{ (-v)(Y), \frac{1}{2} (-v) \left( X' + Y' + \frac{1}{2}[X,Y] \right)\right\} \\
		&= \max \left\{  (-v)(Y), \frac{1}{2}(-v)(X'), \frac{1}{2}(-v)(Y'), (-v)(X), (-v)(Y) \right\} \\
		&		= (-v)(Y) = \varphi_\alpha(v) .
	\end{align*}
	If on the other hand $\varphi_\alpha(v) = \frac{1}{2}(-v)(Y') > (-v)(Y)$, then $(-v)(X)<\frac{1}{2}(-v)(Y')$ and $(-v)(X') <  (-v)(Y')$. Thus using the claim
	\begin{align*}
	(-v)\left(X' + \frac{1}{2}[X,Y]\right) & \leq \max\left\{ (-v)(X') , 2(-v)(X), 2(-v)(Y) \right\} < (-v)(Y')
    \end{align*}
    and
	\begin{align*}
		\varphi_\alpha(uv) &=\max \left\{ (-v)(X+Y) ,\frac{1}{2}(-v)\left(Y'+ X' + \frac{1}{2}[X,Y]    \right)  \right\} \\
		&= \max \left\{ (-v)(X+Y) , \frac{1}{2} (-v)(Y')  \right\} = \frac{1}{2}(-v)(Y') = \varphi_\alpha(v).
	\end{align*}
\end{proof} 

\begin{lemma}\label{lem:orthogonal_valuation}
	Let $B_\theta$ be the scalar product defined on $\frakg_\F$ in Section \ref{sec:killing_involutions_decompositions}, then
	$$
	(-v)(X) = (-v)\left(\sqrt{B_\theta(X,X)}\right)
	$$
	for all $X \in \frakg_\F$. If $X,Y \in \frakg_\F$ are orthogonal with respect to $B_\theta$, then $(-v)(X+Y) = \max\{ (-v)(X), (-v)(Y) \}$. 
\end{lemma}
\begin{proof}
    We note that $(-v)(X) = (-v)(\lVert X\rVert_\infty)$ with the supremum norm $\lVert \cdot \rVert_\infty$. Since on $\frakg_\R$, all norms are equivalent, we can use the transfer principle to obtain a constant $k \in \N$ such that for all $X \in \frakg_\F$
	$$
	\frac{1}{k} \sqrt{B_\theta(X,X)} \leq \lVert X \rVert_\infty \leq k \sqrt{B_{\theta}(X,X)},
	$$
	from which $(-v)(X) = (-v)\left(\lVert X\rVert_\infty\right) = (-v)\left(\sqrt{B_\theta(X,X)}\right)$ follows.
	
	Now if $X,Y \in \frakg_\F$ are orthogonal, $B_\theta(X,Y) = 0$, then
	\begin{align*}
			(-v)(X+Y) &= (-v)\left(\sqrt{B_\theta (X+Y, X+Y)} \right)= \frac{1}{2}(-v)\left(B_\theta (X,X)+ B_{\theta}(Y,Y)\right) \\ 
			&= \frac{1}{2}\max\left\{ (-v)(B_{\theta}(X,X)), (-v)(B_{\theta}(Y,Y)) \right\} \\
			& = \max\{(-v)(X), (-v)(Y) \}
	\end{align*}
	where we used positive definiteness.
\end{proof}

The root group valuation allows us to describe when an element of $(U_\alpha)_\F$ fixes the base point $o$.

\begin{lemma}\label{lem:stab_Ualpha}
	Let $\alpha \in \Sigma$ and $u \in (U_\alpha)_\F$. The following are equivalent
	\begin{enumerate}
		\item [(i)] $u.o = o$
		\item [(ii)] $\varphi_\alpha(u) \leq 0$
		\item [(iii)] $u = \exp(X + X')$ for some $X\in (\frakg_\alpha)_\F(O), \ X'\in (\frakg_{2\alpha})_\F(O)$.
		\item [(iv)] $\log(u) \in \frakg_\F(O)$.
	\end{enumerate}
\end{lemma}
\begin{proof}
	Let $u = \exp(X+X')$ with $X\in (\frakg_\alpha)_\F, X'\in (\frakg_{2\alpha})_\F$. By Theorem \ref{thm:stab}, (i) $u.o = o$ is equivalent to $u \in G_\F(O)$ which by applying the logarithm, which is a polynomial, is equivalent to $X+X' \in \frakg_\F (O)$, (iv). By Proposition \ref{prop:root_decomp} the root space decomposition is orthogonal with respect to the Killing form. Since $X$ and $X'$ lie orthogonal to each other, $X,X' \in \frakg_\F(O)$ individually (iii), by Lemma \ref{lem:orthogonal_valuation}. It is then clear that all the matrix entries of $X$ and $X'$ lie in $O$, hence (ii) $\varphi_\alpha(u) \leq 0$. All these implications are equivalences.
\end{proof}

\subsubsection{$W_a$-convexity for $U_\F$}\label{sec:Wconvexity_for_U}

We first use the Iwasawa retraction to show that if $p \in \A$ is sent to $q \in \A$ by an element $u \in U_\F$, then $p=q$.

\begin{proposition} \label{prop:UonA}
	For all $u\in U_\F$ and $a \in A_\F$
	$$
	ua.o \in \mathbb{A} \iff ua.o = a.o. 
	$$
\end{proposition}
\begin{proof}
	For all $u \in U_\F$ and $a \in A_\F$. If $ua.o = a.o$, then clearly $ua.o \in \mathbb{A}$. We now have to show the converse. So let $ua.o\in \mathbb{A}$, meaning that there exists $b\in  A_\F$ such that $d(ua.o , b.o) =0$. Since the Iwasawa retraction $\rho$ is distance diminishing by Theorem \ref{thm:UAKret}, we have
	\begin{align*}
		d(b.o, a.o) &= d(b.\operatorname{Id},a.\operatorname{Id}) = d(\rho(b.\operatorname{Id}),\rho(u.a.\operatorname{Id})) \\
		&\leq d(b.\operatorname{Id},ua.\operatorname{Id}) = d(b.o, ua.o) = 0,
	\end{align*}
	and thus $a.o = b.o = ua.o \in \A$, as claimed.
\end{proof}
We now show that the set of points in $\A$ that are fixed by an element $u\in (U_\alpha)_\F$ is a half-apartment and in particular $W_a$-convex.

\begin{proposition}\label{prop:Ualphaconv}
	Let $\alpha \in \Sigma$. For $u \in (U_\alpha)_\F$ we have
	$$
	\left\{ p \in \mathbb{A} \colon u.p  \in \mathbb{A}  \right\} = \left\{ a.o \in \A \colon \varphi_\alpha(u) \leq (-v)\left(\chi_\alpha\left(a\right)\right)  \right\}
	$$
	and therefore this set is a half-apartment (with a wall parallel to the wall defined by $(-v)(\chi_\alpha) = 0$) when $u \neq \operatorname{Id}$. 
\end{proposition}
\begin{proof}
	Let $p=a.o \in \A$ and $u = \exp(X+X')$ with $X\in (\frakg_\alpha)_\F, X'\in (\frakg_{2\alpha})_\F$. We have
	\begin{align*}
		ua.o \in \A \stackrel{\text{Prop. }\ref{prop:UonA}}{\iff} & ua.o = a.o 
		\stackrel{\text{Thm. }\ref{thm:stab}}{\iff} a^{-1}ua \in G_\F(O) \\
		 \iff &a^{-1} \exp\left(X + X'\right) a \in G_\F(O) \\
		 \stackrel{\text{Lem. } \ref{lem:aexpXa}}{\iff} &\exp\left( \chi_\alpha\left(a\right)^{-1} X + \chi_\alpha(a)^{-2} X' \right) \in G_\F(O).
	\end{align*}
	Denoting $u' := \exp\left( \chi_\alpha\left(a\right)^{-1} X + \chi_\alpha(a)^{-2} X' \right)  $, the above are equivalent to $\varphi_\alpha(u') \leq 0$ by Lemma \ref{lem:stab_Ualpha}. Using the abbreviation $(-v)(Z) := \max_{ij}\{ (-v)(Z_{ij}) \}$ for $Z \in \frakg_\F$, we have
	\begin{align*}
		\varphi_\alpha(u') &= \max \left\{ (-v)\left( \chi_\alpha(a)^{-1}X  \right), \frac{1}{2} (-v)\left( \chi_\alpha(a)^{-2} X'  \right) \right\} \\
		&= \max \left\{  (-v)\left( X  \right), \frac{1}{2} (-v)\left(  X'  \right)  \right\} - (-v)(\chi_\alpha(a)) \\
		&= \varphi_\alpha(u) - (-v)(\chi_\alpha(a)) \leq 0 
	\end{align*}
	Thus we see that $ua.o \in \A$ is equivalent to $\varphi_\alpha(u) \leq (-v)(\chi_\alpha(a))$.
\end{proof}
Before upgrading the previous result to all of $U_\F$, we consider what happens to products in $U_\F$.

\begin{lemma} \label{lem:nn'}
	Let $\eta \in \Sigma_{> 0}, u\in \exp((\frakg_\eta)_\F)$ and $u' \in \exp(\bigoplus_{\alpha > \eta}(\frakg_\alpha)_\F)$. If $uu' \in G_\F(O)$, then $u\in G_\F(O)$ and $u' \in G_\F(O)$.
\end{lemma}
\begin{proof}
	For $X= \log(u)$ and $Y= \log(u')$, we consider the BCH-formula from Proposition \ref{prop:BCH}
	$$
	\exp(X)\exp(Y)= \exp\left(X+Y+\frac{1}{2}[X,Y] + \frac{1}{12}\left(\left[X,\left[X,Y\right]\right]+\left[Y,\left[Y,X\right]\right]\right) + \ldots \right)
	$$
	Now if $uu' = \exp(X)\exp(Y) \in G_\F(O)$, then
	$$
	\log(uu') = X+Y+\frac{1}{2}[X,Y] + \frac{1}{12}([X,[X,Y]]+[Y,[Y,X]]) + \ldots \in \frakn_\F(O),
	$$
	where $\frakn_\F = \bigoplus_{\lambda \in \Sigma_{ > 0}} (\frakg_\lambda)_\F$ is an orthogonal direct sum. We note that $X\in (\frakg_{\eta})_\F$ is orthogonal to the remaining terms of $\log(uu')$, hence $X \in (\frakg_\eta)_\F(O)$, see Lemma \ref{lem:orthogonal_valuation}. Thus $u= \exp(X) \in G_\F(O)$ and hence also $u' = u^{-1}uu' \in G_\F(O)$ since $G_\F(O)$ is a group.
\end{proof}

\begin{proposition}\label{prop:Uconv}
	For $u\in U_\F$ there are $k_\alpha \in \Lambda \cup \{-\infty \}$ for $\alpha \in \Sigma_{>0}$ such that
	$$
	\{p \in \A \colon u.p \in \A\} = \left\{ a.o \in \A \colon  k_\alpha \leq (-v)\left(\chi_\alpha \left(a\right)\right) \text{ for all } \alpha \in \Sigma_{>0}   \right\}
	$$
	and therefore the set of fixed points is a finite intersection of half-apartments. If $u=u_1 \cdot \ldots \cdot u_k$ for $u_i \in (U_{\alpha_i})_\F$ with $\Sigma_{>0} = \{ \alpha_1, \ldots , \alpha_k \}$ such that $\alpha_1 > \ldots > \alpha_k$, then $k_{\alpha_i} = \varphi_{\alpha_i}(u_i)$. If $u$ fixes all of $\A$, then $u = \operatorname{Id}$. 
\end{proposition}
\begin{proof}
	We use Lemma \ref{lem:BCH_consequence} to write
	$$
	u = u_1 \cdot \ldots \cdot u_k
	$$
	for some $u_i \in (U_{\alpha_i})_\F$ where $\Sigma_{>0}= \{\alpha_1, \ldots ,\alpha_k\}$ with $\alpha_1 > \ldots > \alpha_k$. By Proposition \ref{prop:UonA}, $u.p \in \A$ for $p=a.o\in \A$ if and only if $u.p=p$ and by \ref{thm:stab_A} $a^{-1}ua \in U_\F(O)$. Then
	$$
	a^{-1}ua = a^{-1}u_1a \cdot \ldots \cdot a^{-1}u_ka \in U_\F(O)
	$$ 
	and we can apply Lemma \ref{lem:nn'} repeatedly to obtain $a^{-1}u_1a \in U_\F(O)$, \ldots, $a^{-1}u_ka \in U_\F(O)$. By Proposition \ref{prop:Ualphaconv}, this implies
	$$
	k_{\alpha_i}:=\varphi_{\alpha_i}(u_i) \leq (-v)(\chi_{\alpha_i}(a))
	$$
	for all $\alpha_i \in \Sigma_{>0}$. All the previous implications are equivalences, concluding the description of the fixed point set of $u$. If $u$ fixes all of $\A$, then $a^{-1}ua \in U_\F(O)$ and $a^{-1}u_ia \in (U_{\alpha_i})_\F(O)$ for all $a \in A_\F$. This is only possible if $\varphi_{\alpha_i}(u_i) \leq (-v)(\chi_{\alpha_i}(a))$, so $u_i = \operatorname{Id}$ for all $i$ and thus $u = \operatorname{Id}$.
\end{proof}
As an application of the above, we can conclude that elements of $U_\F(O)$ fix the fundamental Weyl chamber, which is defined as 
$$
C_0 := \{a.o \in \A \colon 0 \leq (-v)(\chi_\alpha (a)) \text{ for all }\alpha \in \Sigma_{>0} \}.
$$

\begin{corollary}\label{cor:NO_fixes_s0}
	Let $u\in U_\F(O)$. Then $u.p = p$ for all $p \in C_0$.
\end{corollary}
\begin{proof}
	Elements $u\in U_\F(O)$ fix $o \in \B$, so by Proposition \ref{prop:Uconv},
	$$
	o \in \{p \in \A \colon u.p \in \A\} = \{ a.o \in \A \colon   k_\alpha \leq (-v)\left(\chi_\alpha \left(a\right)\right) \text{ for all } \alpha \in \Sigma_{>0}  \}
	$$
	from which we conclude that $k_\alpha \leq 0$ for all $\alpha \in \Sigma_{>0}$. If now $p=a.o \in C_0$, then $(-v)(\chi_\alpha(a)) \geq 0 \geq k_\alpha$, so applying Proposition \ref{prop:Uconv} again results in $u.p = p$.
\end{proof}

We also obtain that any $u\in U_\F$ can be conjugated by $A_\F$ into $G_\F(O)$. This statement will be useful in the proof of (A4).

\begin{proposition}\label{prop:exists_a_NO}
	For every $u \in U_\F$ there is an $a\in A_\F$ such that $a^{-1}ua \in U_\F(O)$.
\end{proposition}
\begin{proof}
	Use Proposition \ref{prop:Uconv} to obtain $k_\alpha \in \Lambda \cup \{-\infty\}$ for $\alpha \in \Sigma_{>0}$ such that $u$ fixes all the points in
	$$
	C:= \bigcap_{\alpha \in \Sigma_{ > 0}} H_{\alpha, k_\alpha}^+.
	$$ 
	We claim that $C$ is non-empty. Indeed, let $c\in \F_{>0}$ with $(-v)(c) > k_\alpha$ for all $\alpha >0$. Proposition \ref{prop:isoFrA} gives a group isomorphism $f_\F \colon (\F_{>0})^{|\Delta|} \cong A_\F$ with $\chi_\delta(a) = c$ for $a=f_\F(c,\ldots, c)$ and $\delta \in \Delta$. Since every $\alpha>0$ is a positive integer linear combination of elements in $\Delta$, $(-v)(\chi_\alpha(a))\geq(-v)(c)>k_\alpha$, so $a.o \in C$.
	 
By Proposition \ref{prop:anainN}, $a^{-1}ua \in U_\F$ and by Theorem \ref{thm:stab}, $a^{-1}ua.o = a^{-1}a.o=o$ implies $a^{-1}ua \in U_\F(O)$.
\end{proof}

\subsubsection{Stabilizers of apartment, half-apartments and Weyl-chambers}

In Corollary \ref{cor:NO_fixes_s0} we proved that $U_\F(O)$ fixes the fundamental Weyl chamber $C_0 = \{ a.o \in \A  \colon \chi_\alpha(a) \geq 1 \text{ for all } \alpha \in \Delta \}$, where $\Delta$ is a basis of $\Sigma$. The goal of this subsection is to describe the whole stabilizer of $C_0$ using the group $B_\F=U_\F A_\F M_\F$ from Section \ref{sec:BWB}. So far, we considered the action of $U_\F$ and $A_\F$, now we continue by investigating the action of $K_\F$. 

To determine which elements of $K_\F$ fix the standard apartment $\A \subseteq \B$, we will use some CAT(0) geometry on the symmetric space $\X_\R$, equipped with the right metric. 

\begin{theorem}\label{thm:CAT0_convex}
	(\cite[Proposition 2.2]{BH99})
	
	If $\gamma, \gamma'$ are two unit-speed geodesics in a CAT(0)-space, then the function 
	$$t \mapsto d(\gamma(t),\gamma'(t))$$ is convex.
\end{theorem}

In the non-standard symmetric space $\X_\F$, the elements of $K_\F$ that fix all points of the standard maximal flat $A_\F.\operatorname{Id}$ are exactly $\operatorname{Cen}_{K_\F}(A_\F) =: M_\F$. In the next Proposition we show that these are also exactly the elements in $K_\F$ that fix $\A \subseteq \B$ pointwise.

\begin{proposition}\label{prop:K_fixing_A_is_M}
	If $k\in K_\F$ fixes all points of $\A \subseteq \B$, then $k \in M_\F$. 
\end{proposition}
\begin{proof}
	Recall that the distance on the non-standard symmetric space $\X_\F$ is given by $d  = (-v) \circ N \circ \delta \colon \X_\F \times \X_\F \to \Lambda$. We claim that the first-order formula
		\begin{align*}
		\varphi \colon \quad \forall k \in K \colon & (\forall a \in A \colon a.\operatorname{Id}=k.a.\operatorname{Id}) 
		\vee ( \forall c >0 \ \exists a \in A \colon N(\delta(a'.\operatorname{Id},k.a.\operatorname{Id}))>c )
	\end{align*}
	holds over $\F$. The situation is illustrated in Figure \ref{fig:rotation_statement}. 
	Informally, $\varphi$ states that $k\in K_\F$ either fixes all points in the maximal flat $A_\F.\operatorname{Id}$, or there are points $a.\operatorname{Id}$ that are sent arbitrarily far away by $k$. To prove that $\varphi$ holds over $\F$, it suffices to show $\varphi$ over $\R$ by the transfer principle.
	
		\begin{figure}[h]
		\centering
		\includegraphics[width=0.6\linewidth]{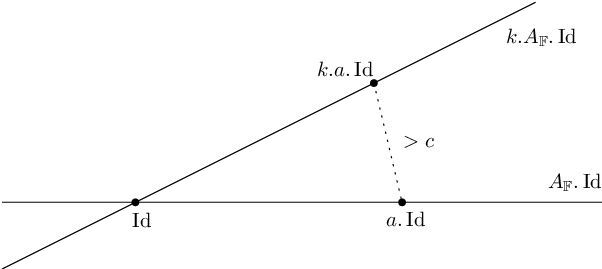}
		\caption{The first-order formula $\varphi$ states that either $k$ fixes all points in $A_\F$ or there are points $a.\operatorname{Id}, k.a.\operatorname{Id}$ whose distance is arbitrarily large. }
		\label{fig:rotation_statement}
	\end{figure}

	We note that changing the $W_s$-invariant multiplicative norm $N$ only changes $\X_\R$ up to quasi-isometry, so we may choose $N_\R$ coming from a scalar product, even if $N_\R$ is then not semialgebraic, as the truth of $\varphi$ only depends on $N_\R$ up to equivalency. So consider 
	\begin{align*}
		N_\R \colon A_\R.\operatorname{Id} &\to \R_{\geq 1} \\
		a.\operatorname{Id} &\mapsto \exp\left( \sqrt{B_\theta(\log(a),\log(a))}\right)
	\end{align*}
	for the scalar product $B_\theta$ on $\fraka$. Then by the general theory of symmetric spaces of non-compact type, $\X_\R$ with the distance $d = \log \circ N_\R \circ \delta_\R$ is a complete CAT(0)-space.
    Consider a unit-speed geodesic $\gamma \colon \R \to A_\F.\operatorname{Id} \subseteq \X_\R$ passing through $\gamma(0)=\operatorname{Id}$. Then $k.\gamma$ is also a unit-speed geodesic passing through $\operatorname{Id}$. By Theorem \ref{thm:CAT0_convex}, the function $f \colon t \mapsto d(\gamma(t),k.\gamma(t))$ is convex. Since $f$ is non-negative and $f(0)=0$, $f$ then has to be constant (hence $a.\operatorname{Id}=k.a.\operatorname{Id}$ for all $a\in A_\F$) or eventually be larger than $\log(c)$ for every $c\in \R_{>0}$ (hence there is some $a\in A_\F$ such that $N_\R(\delta_\R(a.\operatorname{Id},k.a.\operatorname{Id}))>c$).
	
	Now that $\varphi$ is established over $\F$, we consider some $k \in K_\F$ that fixes all points of $\A \subseteq \B$. Choosing $c \in \F_{>0}$ with $c\notin O$, we see that the second option in $\varphi$ cannot be true, whence $k.a.\operatorname{Id} = a.\operatorname{Id}$ for all $a\in A_\F$. This means $k \in \operatorname{Cen}_{K_\F}(A_\F)$.
\end{proof}

In Proposition \ref{prop:K_fixes_C0_in_M}, we strengthen the previous result by only requiring $k$ to fix a chamber of $\A$. We first need a preliminary result.

\begin{lemma}\label{lem:K_fixing_pm_mean}
	Let $k\in K_\F$ and $a,b\in A_\F$. If $k.a.o = a.o$, then $k.a^{-1}.o = a^{-1}.o$. If moreover $k.b.o = b.o$, then $k.\sqrt{ab}.o = \sqrt{ab}.o$. 
\end{lemma}
\begin{proof}
	We first assume that $A_\F$ consists of diagonal matrices, so we may write $a = \operatorname{Diag}(a_1, \ldots, a_n)$ and $b = \operatorname{Diag}(b_1,\ldots, b_n)$. Then
	\begin{align*}
		k.a.o = a.o & \iff a^{-1}ka \in G_\F(O) \iff \forall i,j\colon  k_{ij}\frac{a_j}{a_i} \in O \\
		&\iff \forall j,i \colon k\tran_{ij}\frac{a_i}{a_j} \in O \iff ak\tran a^{-1 } \in G_\F(O) \\
		&\iff k\tran.a^{-1}.o = a^{-1}.o \iff k.a^{-1}.o = a^{-1 }.o.
	\end{align*}
	Moreover, if $k.a.o = a.o$ and $k.b.o=b.o$, then $k_{ij}a_j/a_i \cdot k_{ij} b_j/b_i \in O$ for all $i,j$. Then also
	$$
	k_{ij} \frac{\sqrt{a_jb_j}}{\sqrt{a_ib_i}} \in O
	$$
	for all $i,j$, which translates to $k.\sqrt{ab}.o = \sqrt{ab}.o$.
	
	In general, the matrices in $A_\F$ may not be diagonal, but they are symmetric. By the spectral theorem for symmetric matrices, which holds over $\F$ by the transfer principle, $A_\F$ is orthogonally diagonalizable, meaning that there is some $Q \in \operatorname{SO}(n)$ such that $QA_\F Q\tran$ is diagonal. We can then apply the above arguments to the group $QG_\F Q\tran < \operatorname{SL}_n(\F)$. 
	For $k.a.o = a.o$ we obtain
	\begin{align*}
		a^{-1}ka \in G_\F(O) & \iff (QaQ\tran)^{-1} QkQ\tran QaQ\tran \in (QG_\F Q\tran)(O) \\
		&\iff QaQ\tran Qk\tran Q\tran (QaQ\tran)^{-1} \in (QG_\F Q\tran)(O) \\
		&\iff ak\tran a^{-1} \in G_\F(O)
	\end{align*}
	and complete the argument as above. When additionally $b^{-1}kb \in G_\F(O)$ we have
	\begin{align*}
		(QbQ\tran)^{-1} QkQ\tran QbQ\tran \in (QG_\F Q\tran)(O) 
	\end{align*}
	which by the above implies
	\begin{align*}
		(Q\sqrt{ab}Q\tran)^{-1} QkQ\tran Q\sqrt{ab}Q\tran \in (QG_\F Q\tran)(O)
	\end{align*}
	and thus$
	\sqrt{ab}^{-1}k\sqrt{ab} \in G_\F(O).
	$
\end{proof}

Let $C_0 = \{ a.o \in \A  \colon \chi_\alpha(a) \geq 1 \text{ for all } \alpha \in \Delta\}$ be the fundamental Weyl chamber associated to a basis $\Delta$ of $\Sigma$. 

\begin{proposition}\label{prop:K_fixes_C0_in_M}
	Let $k \in K_\F$ such that $k.p = p$ for all $p \in C_0$. Then $k \in M_\F$ and hence $k$ fixes all points in $\A$. In fact, if $k$ fixes all the points in $a.C_0$ for any $a \in A_\F$, then $k \in M_\F$.
\end{proposition}
\begin{proof}
	We first claim that every element $a\in A_\F$ is of the form $a=a_1a_2^{-1}$ for $a_1.o,a_2.o \in C_0$. To see this, we show that the first-order formula
	$$
	\varphi\colon \quad  \quad \forall a \in A \colon \exists a_1, a_2 \in A \colon a = a_1\cdot a_2^{-1} \wedge  \bigwedge_{\alpha \in \Sigma_{>0}} \chi_\alpha(a_1) \geq 1 \wedge \chi_\alpha(a_2) \geq 1 
	$$
	holds over $\R$ and then apply the transfer principle. Over $\R$ we can transfer the problem to the Lie algebra $\fraka_\R$ using the logarithm. We equip $\fraka$ with the distance defined by the scalar product $B_\theta$.
	 Let $H:=\log(a)$ and $R := \sqrt{B_\theta(H,H)}$. 
	Since $\mathfrak{c}_0 := \{ H \in \fraka_\R \colon \alpha(H)>0 \text{ for all }\alpha \in \Sigma_{>0} \}$ contains an open cone, it contains a ball $B_{r}(H')$ for some $r>0$ and $H' \in \mathfrak{c}_0$. Scaling the ball by the factor $R/r$, we obtain that $B_{R}(R/r\cdot H') \subseteq \mathfrak{c}_0$. As in Figure \ref{fig:scaling_statement}, we define $H_1 = R/r\cdot H'$ and $H_2=H_1-H$ which lies on the boundary of $B_{R}(R/r\cdot H')$ and hence also in $\mathfrak{c}_0$. Then $a = \exp(H_1)\exp(H_2)^{-1}$, concluding the proof of $\varphi$ over $\R$ and hence over $\F$.
	
	\begin{figure}[h]
		\centering
		\includegraphics[width=0.45\linewidth]{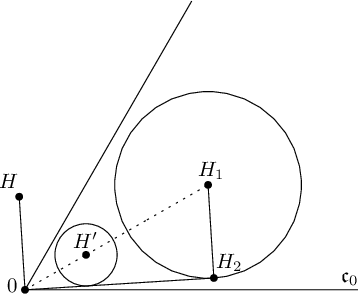}
		\caption{We use that the cone $\mathfrak{c}_0\subseteq \fraka_\mathbb{R}$ contains an open ball, which we can scale to obtain $H_1,H_2 \in \mathfrak{c}_0$ with $H=H_1-H_2$. }
		\label{fig:scaling_statement}
	\end{figure}
	
	If now $k\in K_\F$ fixes $C_0$ pointwise and $p=a.o \in \A$. Then $a=a_1a_2^{-1}$ as above with $a_1.o, a_2.o \in C_0$. Since $C_0$ is a cone, also $a_1^2.o, a_2^2.o \in C_0$, so $k.a_1^2.o = a_1^2.o$ and $k.a_2^2.o=a_2^2.o$. By Lemma \ref{lem:K_fixing_pm_mean}, then $k.a_2^{-2}.o = a_2^{-2}.o$ and 
	$$
	k.a.o = k.\sqrt{a_1^2 a_2^{-2}}.o = \sqrt{a_1^2 a_2^{-2}}.o = a.o,
	$$
	completing the first statement of the proof. If $k$ fixes $a.C_0$ for some $a\in A_\F$, a modification of the above argument similarly implies that $k$ fixes all of $\A$ pointwise.
\end{proof}

\begin{theorem}\label{thm:NO_fixes_s0}
	The pointwise stabilizer of $C_0$ in $G_\F$ is $B_\F(O) = U_\F(O)A_\F(O)M_\F$.
\end{theorem}
\begin{proof}
	If $g=uak \in  B_\F(O) = U_\F(O)A_\F(O)M_\F$ and $p \in C_0$, then $g.p=uak.p=ua.p=u.p=p$, where the last equality follows from Corollary \ref{cor:NO_fixes_s0}.
	
	If $g$ fixes $C_0$ pointwise, in particular it fixes $o \in C_0$, hence $g\in G_\F(O)$, which can be decomposed to $G_\F(O) = U_\F(O)A_\F(O)K_\F$ by Corollary \ref{cor:UAKO}. Therefore $g=uak$ with $u\in U_\F(O)$ and $a \in  A_\F(O)$. Therefore we have that $k$ fixes $C_0$ pointwise. Proposition \ref{prop:K_fixes_C0_in_M} now concludes the proof by showing $k \in M_\F$.
\end{proof}

Recall that $S <G$ is a maximal $\K$-split torus and $A_\F < S_\F$ is the semialgebraically connected component of the identity. In the following we consider the groups $T_\F := \operatorname{Cen}_{G_\F}(A_\F)$ and $T_\F(O) := T_\F \cap G_\F(O) $. Actually, $T = \operatorname{Cen}_{G}(S)$ is a maximal algebraic torus, but we do not rely on this fact in what follows. 
\begin{lemma}\label{lem:BT_T_is_MA}
	We have $T_\F = \operatorname{Cen}_{G_\F}(A_\F) = M_\F \cdot A_\F$.
\end{lemma}
\begin{proof}
	The inclusion $\supseteq$ is clear. For the other direction, let $g \in T_\F$ and choose an Iwasawa decomposition $g = nak$ with $n \in U_\F , a\in A_\F, k \in K_\F$, see Theorem \ref{thm:KAU}. For all $b \in A_\F$ we have $nak.b.o = b.nak.o = bna.o$, and thus $k.b.o = a^{-1}n^{-1}bna.o$. Denoting $\tilde{n} := (a^{-1}n^{-1}a)(a^{-1}bnb^{-1}a) \in U_\F$, where we made use of Theorem \ref{prop:anainN}. Then $k.b.o = \tilde{n}.b.o$ for all $b \in A_\F$. By Proposition \ref{prop:exists_a_NO} there is a $c\in A_\F$ such that $\tilde{\tilde{n}} := c\tilde{n}c^{-1} \in U_\F(O)$ and hence $\tilde{\tilde{n}}$ fixes the elements of $C_0$ by Theorem \ref{thm:NO_fixes_s0}. In particular, $\tilde{n}$ fixes $c.C_0$ and thus $k$ fixes $c.C_0$, since $k.p = \tilde{n}.p$ for all $p \in \A$. By Lemma \ref{prop:K_fixes_C0_in_M}, we then have $k \in M_\F$. This implies also $\tilde{n}.p=p$ for all $p\in \A$, hence $\tilde{n} = \operatorname{Id}$. Then $bnb^{-1} =n$ for all $b\in A_\F$, not only the ones with $(-v)(\chi_\alpha(b)) = 0$, hence by Lemma \ref{lem:aexpXa}, $n= \operatorname{Id}$. Thus $g = ak$ with $k \in M_\F$ and $a \in A_\F$.
\end{proof}

\begin{theorem}\label{thm:stab_A}
	The pointwise stabilizer of $\A$ in $G_\F$ is $T_\F(O) = A_\F(O)M_\F$.
\end{theorem}
\begin{proof}
	Elements of $A_\F(O)M_\F$ fix all points in $\A$. If $g\in G_\F$ fixes all points of $\A$, it fixes in particular the points in $C_0$, so $g = uak$ with $u \in U_\F(O), a \in A_\F(O), k \in M_\F$ by Theorem \ref{thm:NO_fixes_s0}. By the description of the stabilizer of $u$ in Proposition \ref{prop:Uconv}, $u$ can only fix all of $\A$ if $u=\Id$, so $g \in A_\F(O)M_\F$. By Lemma \ref{lem:BT_T_is_MA}, $T_\F(O) = A_\F(O)M_\F$.
\end{proof} 

\begin{theorem}\label{thm:NalphaO_fixes_H}
	Let $\alpha \in \Sigma$. The pointwise stabilizer of $H_\alpha^+ = \{ a.o \in \A \colon (-v)(\chi_\alpha(a)) \geq 0\}$ in $G_\F$ is $(U_\alpha)_\F(O) A_\F(O) M_\F$.
\end{theorem}
\begin{proof}
Without loss of generality, we may assume that $\Sigma$ is equipped with an order in which $\alpha>0$. Then $(U_\alpha)_\F(O) A_\F(O) M_\F \subseteq \operatorname{Stab}_{G_\F}(H_\alpha^+)$ by Proposition \ref{prop:Ualphaconv}. If $g \in \operatorname{Stab}_{G_\F}(H_\alpha^+)$, then $g = uak \in U_\alpha(O) A_\F(O) M_\F$ by Theorem \ref{thm:NO_fixes_s0} since $C_0 \subseteq H_\alpha$, so in fact $u$ fixes $H_\alpha^+$ pointwise. By Lemma \ref{lem:BCH_consequence}, $u = u_1\cdots u_k$ where $u_i \in U_{\alpha_i}$ such that $\alpha_1 > \ldots > \alpha_k > 0$. Then we can apply the refined version of Proposition \ref{prop:Uconv}, to obtain 
	$$
	\operatorname{Fix}_\A(u) = \left\{ a.o \in \A \colon k_{\alpha_i} \leq (-v)(\chi_{\alpha_i}(a)) \text{ for all } \alpha_i \in \Sigma_{>0} \right\} 
	$$
	where $k_{\alpha_i} = \varphi_{\alpha_i}(u_i)$. In our case, $H_\alpha^+ \subseteq \operatorname{Fix}_\A(u)$ and since when $k_{\alpha_i} \neq -\infty$
	$$
	H_{\alpha}^+ \subseteq \left\{ a.o \in \A \colon k_{\alpha_i} \leq (-v)(\chi_{\alpha_i}(a)) \right\}
	$$ 
	implies $\alpha = \alpha_i$, 
	$k_{\alpha_i} = -\infty$ for all $\alpha_i \neq \alpha$ and $k_\alpha = 0$, so $\varphi_{\alpha_i}(u_i) = -\infty$ implies $u_i = \operatorname{Id}$ whenever $\alpha_i \neq \alpha$, leading to $u \in (U_\alpha)_\F(O)$ concluding the proof.
\end{proof}

\begin{corollary}\label{cor:NalphaO_fixes_H_affine}
	Let $\alpha \in \Sigma, \ell \in \Lambda$. The pointwise stabilizer of the affine half-apartment
	$$
	H_{\alpha,\ell}^+ = \left\{ a.o \in \A \colon (-v)(\chi_\alpha(a) ) \geq \ell  \right\}
	$$
	in $G_\F$ is $U_{\alpha, \ell} A_\F(O) M_\F$, where $U_{\alpha, \ell} =  \{ u \in (U_\alpha)_\F \colon u \text{ stabilizes } H_{\alpha, \ell} \text{ pointwise} \}$.
\end{corollary}
\begin{proof}
	If $g\in G_\F$ stabilizes $H_{\alpha,\ell}^+$ pointwise and $a\in A_\F$ satisfies $(-v)(\chi_\alpha(a)) = \ell$, then $a^{-1} g a \in G_\F$ stabilizes $H_{\alpha, \ell}^+$ pointwise. By Theorem \ref{thm:NalphaO_fixes_H} and since $A_\F(O) M_\F \subseteq \operatorname{Cen}_{G_\F}(A_\F)$, we have $g \in a(U_\alpha)_\F(O) a^{-1} A_\F(O) M_\F$. While elements $u \in (U_\alpha)_\F(O)$ stabilize elements $b.o$ with $(-v)(\chi_\alpha(b)) = 0$ pointwise, $a ua^{-1}$ stabilize elements $ab.o$ with $(-v)(\chi_\alpha(ab)) = \ell + 0 $ pointwise, concluding the proof.
\end{proof}

\subsubsection{Bruhat-Tits theory for root groups} \label{sec:BT_root_groups}
 The goal of Subsections \ref{sec:BT_root_groups}, \ref{sec:BT_rank_1} and \ref{sec:BT_higher_rank} is to describe the pointwise stabilizer of a finite subset $\Omega \subseteq \A$ in Theorem \ref{thm:BTstab_fin}. We to obtain this, we first consider points fixed by $(U_\alpha)_\F$ in Section \ref{sec:BT_root_groups}, then points fixed by the rank one subgroups generated by $(U_\alpha)_\F$ and $(U_{-\alpha})_\F$ in Section \ref{sec:BT_rank_1}, before taking on the whole group $G_\F$ in \ref{sec:BT_higher_rank}. These subsections are inspired by the study of stabilizers in \cite[Sections 6 and 7]{BrTi}, for a good English reference see \cite{Lan96}. 
 
Let $\Omega \subseteq \A$ be any subset of the apartment $\A$. Let
\begin{align*}
	U_{\alpha, \Omega} &:= \left\{ u \in (U_\alpha)_{\F} \colon  u.p = p \text{ for all } p \in \Omega  \right\}
\end{align*}
denote the pointwise stabilizer of $\Omega$ in the group $(U_{\alpha})_\F$. The subscript $\F$ is no longer needed, since there is no corresponding group of $\R$-points. In view of Proposition \ref{prop:UonA}, Proposition \ref{prop:Ualphaconv} can be reformulated.
\begin{lemma}\label{lem:BTUalphaOmega}
	For roots $\alpha \in \Sigma$ and any subset $\Omega \subseteq \A$, we have
	\begin{align*}
		U_{\alpha,\Omega} 
		&= \left\{  u \in (U_\alpha)_\F \colon  \varphi_\alpha(u) \leq (-v)\left(\chi_\alpha\left(a\right)\right) \text{ for all } a.o \in \Omega  \right\} .
	\end{align*}
\end{lemma}
For any $\ell \in \Lambda$, denote 
$$
U_{\alpha, \ell} := \left\{ u \in (U_\alpha)_\F \colon \varphi_\alpha(u) \leq \ell \right\}.
$$
Note that if $\ell = \min_{a.o \in \Omega} \{ (-v)(\chi_\alpha(a)) \}$ and $\alpha \in \Sigma$, then Lemma \ref{lem:BTUalphaOmega} can be reformulated as
$$
U_{\alpha,\Omega} = U_{\alpha, \ell}.
$$
In the setting of \cite{BrTi}, $\ell = \inf_{a.o \in \Omega} \{ (-v)(\chi_\alpha(a)) \}$ with $U_{\alpha,\Omega} = U_{\alpha,\ell}$ always exists, since they work with $\Lambda=\R$. In our case we have to be more careful. If $\Omega$ is a finite set, this $\ell$ exists. In particular, when $|\Omega| = 1$, we have for any $a\in A_\F$
$$
U_{\alpha,\left\{ a.o \right\}} = U_{\alpha, (-v)(\chi_\alpha(a))}.
$$
\begin{lemma}\label{lem:BTaUa}
	Let $\alpha \in \Sigma, \ell \in \Lambda$ and $a \in A_\F$. Then
	$$
	a U_{\alpha , \ell} a^{-1} = U_{\alpha,\ell + (-v)(\chi_\alpha(a))}.
	$$
\end{lemma}
\begin{proof}
	Let $b \in A_\F$ with $(-v)(\chi_\alpha(b)) = \ell $, the existence of which can be concluded from Lemma \ref{lem:Jacobson_Morozov_oneparam}. We consider the single element set $\Omega = \{b.o\} $. Then $U_{\alpha, \ell} = U_{\alpha, \Omega}$. Let $u \in U_{\alpha, \ell}$. Then $aua^{-1} \in U_{\alpha, a.\Omega}$ since
	$$
	aua^{-1}.(a.b.o) = au.b.o = a.b.o  
	$$
	and since $aua^{-1} \in (U_\alpha)_\F$, see Proposition \ref{prop:anainN}. Now since 
	$$
	(-v)(\chi_{\alpha}(ab)) = (-v)(\chi_\alpha(a)) + (-v)(\chi_\alpha(b)) = \ell + (-v)(\chi_\alpha(a)) 
	$$
	we have $aU_{\alpha, \ell}a^{-1} = aU_{\alpha, \Omega}a^{-1} = U_{\alpha, a.\Omega} = U_{\alpha, \ell + (-v)(\chi_\alpha(a))}$.
\end{proof}
\begin{lemma}\label{lem:BTkUk}
	Let $\alpha \in \Sigma, \ell \in \Lambda$ and $k \in M_\F = \operatorname{Cen}_{K_\F}(A_\F)$. Then
	$$
	kU_{\alpha,\ell}k^{-1}= U_{\alpha,\ell}.
	$$
\end{lemma}
\begin{proof}
	Since $k \in M_\F$, it represents the trivial element 
	of the spherical Weyl group acting on the root system. In particular $kU_{\alpha}k^{-1} = 
	U_{\alpha}$. 
	
	Let $b\in A_\F$ with $(-v)(\chi_\alpha(b))= \ell$ and $\Omega = \{b.o\}$. Then $U_{\alpha,\ell} = U_{\alpha, \Omega}$. For any $u\in U_{\alpha,\ell}$ we have
	$$
	kuk^{-1}.b.o = ku.b.o = k.b.o = b.o,
	$$
	hence $kU_{\alpha,\Omega}k^{-1} = U_{\alpha, \Omega}$ concluding the proof.
\end{proof}

Putting the previous two results together shows that $U_{\alpha,\ell}$ is invariant under conjugation by elements in the pointwise stabilizer $T_\F(O)$.
\begin{lemma}\label{lem:BTtUt}
	Let $\alpha \in \Sigma, \ell \in \Lambda$ and $t \in T_\F(O)=M_\F A_\F(O)$. Then
	$$
	t U_{\alpha,\ell} t^{-1} = U_{\alpha,\ell}.
	$$
\end{lemma}

\subsubsection{Bruhat-Tits theory in rank 1} \label{sec:BT_rank_1}
Our goal in this section is to study the group generated by $U_{\alpha, \Omega}$ and $U_{-\alpha, \Omega}$. For this we will use Jacobson-Morozov in the form of Proposition \ref{prop:Jacobson_Morozov_real_closed}. Therefore we have to restrict ourselves to reduced root systems from now on. In this subsection we fix $\alpha \in \Sigma$ such that $(\frakg_{2\alpha})_\F = 0$ and $u \in (U_\alpha)_\F$.

When $u\neq \Id $, there is a $t \in \F$ and an $\mathfrak{sl}_2$-triplet $(X,Y,H)$ as in Proposition \ref{prop:Jacobson_Morozov_real_closed} with $(-v)(X) = (-v)(Y) = (-v)(H) = 0$ such that $u = \exp(tX)$. Moreover, there is an algebraic group homomorphism $\varphi_\F \colon \operatorname{SL}(2,\F) \to G_\F$ with finite kernel such that
\begin{align*}
	u = \exp(tX) = \varphi_\F \begin{pmatrix}
		1 & t \\ 0& 1
	\end{pmatrix}   \quad \text{and} \quad \exp(tY) = \varphi_\F \begin{pmatrix}
		1 & 0 \\ t & 1
	\end{pmatrix}.
\end{align*}
We note that $t \in O$ if and only if $u \in G_\F(O)$, using Lemma \ref{lem:stab_Ualpha}.
\begin{lemma}\label{lem:BTphiu_is_vt}
	For every $t\in \F$ we have
	$$
	\varphi_\alpha \left( \varphi_\F \begin{pmatrix}
		1 & t \\ 0 & 1
	\end{pmatrix} \right) = (-v)(t) = \varphi_{-\alpha} \left( \varphi_\F \begin{pmatrix}
		1 & 0 \\ t & 1
	\end{pmatrix} \right) .
	$$
\end{lemma}
\begin{proof}
	Since $(-v)(X) := \max_{ij}\{(-v)(X_{ij})\} = 0$, $\varphi_\alpha(u) = (-v)(tX) = (-v)(t)$. Similarly $\varphi_{-\alpha}(\exp(tY)) = (-v)(tY) = (-v)(t)$.
\end{proof}

\begin{lemma}\label{lem:BTm_uuu}
	Let $\ell := (-v)(t)$. Then
	$$
	m(u) := \varphi_\F\begin{pmatrix}
		0 & t \\ -1/t & 0
	\end{pmatrix} = \varphi_\F \left(
	\begin{pmatrix}
		1 & 0 \\ -1/t & 1
	\end{pmatrix}\begin{pmatrix}
		1 & t \\ 0 & 1
	\end{pmatrix}\begin{pmatrix}
		1 & 0 \\ -1/t & 1
	\end{pmatrix},
	\right)
	$$
	in particular $m(u) \in U_{-\alpha,-\ell} U_{\alpha, \ell} U_{-\alpha,-\ell}$.
\end{lemma}
\begin{proof}
	Let 
	$$
	u' := \varphi_\F\begin{pmatrix}
		1 & 0 \\ -1/t & 1
	\end{pmatrix}  \in (U_{-\alpha})_\F.
	$$
	The matrix expression $m(u) = u' u u'$ is a direct calculation, showing $m(u) \in (U_{-\alpha})_\F(U_\alpha )_\F (U_{-\alpha})_\F$. By Lemma \ref{lem:BTphiu_is_vt}, $\varphi_\alpha(u) = \ell$ and $\varphi_{-\alpha}(u') = (-v)(-1/t) = -\ell$, concluding the proof.
\end{proof}
The element 
$$
m(u) := \varphi_\F \begin{pmatrix}
	0 & t \\ -1/t & 0
\end{pmatrix} \in G_\F
$$
is contained in $ \operatorname{Nor}_{G_\F}(A_\F)$ by Lemma \ref{lem:Jacobson_Morozov_m} and thus a representative of an element of the affine Weyl group $W_a = \operatorname{Nor}_{G_\F}(A_\F)/\operatorname{Cen}_{G_\F}(A_\F) $. Recall that the affine Weyl group $W_a$ can be identified with $W_a = \A \rtimes W_s $, where $W_s = N_\F/M_\F= \operatorname{Nor}_{K_\F}(A_\F)/\operatorname{Cen}_{K_\F}(A_\F)$ is the spherical Weyl group.

\begin{proposition}\label{prop:BTmu_in_Wa} The action of $m(u)$ decomposes as
	$$
	m(u) = \varphi_\F \begin{pmatrix}
		t & 0 \\ 0 & t^{-1}
	\end{pmatrix} \cdot \varphi_\F \begin{pmatrix}
		0 & 1 \\ -1 & 0 
	\end{pmatrix} =: a_t\cdot m \in A_\F \cdot N_\F 
	$$
	into an affine part represented by $a_t$ and a spherical part represented by $m$.
	The element $m$ represents the reflection $r_\alpha \in W_s$ and for $a_t$ we have $(-v)(\chi_\alpha(a_t)) = 2(-v)(t) = 2\varphi_\alpha(u)$. Thus $m(u)$ represents the affine reflection along the hyperplane
	$$
	\{ a.o \in \A \colon (-v)(\chi_\alpha(a)) = \varphi_\alpha(u) \}.
	$$
\end{proposition}
\begin{proof}
	The decomposition of $m(u)$ is a direct calculation. We investigate the action of $m$. We may decompose $A_\F = (A_{\pm \alpha})_\F \cdot A_{\perp}$ as a direct product, where
	$$
	(A_{\pm \alpha})_\F = \varphi_\F \left(\left\{ \begin{pmatrix}
		\lambda & 0 \\ 0 & \lambda^{-1}
	\end{pmatrix} \colon \lambda >0 \right\}\right)
	$$
	and
	$$
	A_\perp = \{ a\in A_\F \colon \chi_\alpha(a) = 1 \}.
    $$
    The reflection $r_\alpha \colon \A \to \A$ is defined by $r_\alpha(a_\alpha.o) = a_{\alpha}^{-1}.o$ for all $a_\alpha \in (A_{\pm \alpha})_\F$ and $r_\alpha(a_\perp.o) = a_\perp.o $ for all $a_\perp \in A_\perp$. For $a_\alpha \in (A_{\pm \alpha})_\F$ there is some $\lambda >0$ such that
    \begin{align*}
    m.a_\alpha.o &= m.a_\alpha.m^{-1} . o= \varphi_\F \left(\begin{pmatrix}
    	0 & 1 \\ -1 & 0
    \end{pmatrix}  \begin{pmatrix}
    	\lambda & 0 \\ 0 & \lambda^{-1}
    \end{pmatrix}\begin{pmatrix}
    	0 & -1 \\ 1 & 0
    \end{pmatrix} \right). o \\
    &=  \varphi_\F\begin{pmatrix}
    	\lambda^{-1} & 0 \\ 0 & \lambda
    \end{pmatrix}  .o = a_\alpha^{-1}.o.
    \end{align*}
    For $a_\perp \in A_{\perp}$ we use 
    \begin{align*}
    	m &= \varphi_\F \begin{pmatrix}
    		0 & -1 \\ 1 & 0
    	\end{pmatrix} = \varphi_\F \left( \begin{pmatrix}
    	1 & 0 \\ -1 & 1
    	\end{pmatrix} \begin{pmatrix}
    	1 & 1 \\ 0 & 1
    	\end{pmatrix} \begin{pmatrix}
    	1 & 0 \\ -1 & 1
    	\end{pmatrix} \right) \\
    	& \in \exp((\frakg_{-\alpha})_\F) \cdot \exp((\frakg_{\alpha})_\F) \cdot \exp((\frakg_{-\alpha})_\F)
    \end{align*}
    and Lemma \ref{lem:aexpXa} to obtain $m. a_\perp.o = a_\perp. m.o = a_\perp.o$. Lemmas \ref{lem:Jacobson_Morozov_oneparam} and \ref{lem:BTphiu_is_vt} give $(-v)(\chi_\alpha(a_t)) = 2(-v)(t) = 2\varphi_\alpha(u)$. A point $a.o = a_\perp a_\alpha.o \in \A$ is fixed by $m(u)$ if and only if
    $$
    a.o = m(u).a.o = a_t m a_\perp a_\alpha m^{-1} .o = a_t a_\perp a_\alpha^{-1}.o 
    $$
    which is the case exactly when $(-v)(\chi_\alpha(a_t a_\alpha^{-2})) =0$, i.e. when $(-v)(\chi_\alpha(a)) = (-v)(t) = \varphi_\alpha(u)$.
\end{proof}
We obtain three corollaries from the geometric description above. Recall from Theorem \ref{thm:stab_A} that $\operatorname{Stab}_{G_\F}(\A) = M_\F A_\F(O) = T_\F(O)$.

\begin{lemma}\label{lem:BTphi_stab}
	For any $u_1,u_2 \in (U_\alpha)_\F$, 
	$$
	\varphi_\alpha(u_1) = \varphi_\alpha(u_2) \quad \iff \quad m(u_2)^{-1}m(u_1) \in T_\F(O).
	$$
	For any $u \in (U_\alpha)_\F, u' \in (U_{-\alpha})_\F$, 
	$$
	\varphi_\alpha(u) = -\varphi_{-\alpha}(u') \quad \iff \quad m(u)^{-1}m(u') \in T_\F(O).
	$$
\end{lemma}
\begin{proof}
	 We use the description of the action in Proposition \ref{prop:BTmu_in_Wa}. Both $m(u_1)$ and $m(u_2)$ act on $\A$ by an affine reflection along a hyperplane. They reflect along the same hyperplane if and only if $\varphi_\alpha(u_1) = \varphi_\alpha(u_2)$, which is the case exactly when $m(u_2)^{-1}m(u_1) \in \operatorname{Stab}_{G_\F}(\A)$. 
	
	The fixed hyperplanes of $m(u)$ are given by the conditions $(-v)(\chi_\alpha(a)) = \varphi_\alpha(u)$ and $(-v)(\chi_{-\alpha}(a)) = \varphi_{-\alpha}(u')$ respectively. The second condition can also be written as $(-v)(\chi_\alpha(a)) = -\varphi_{-\alpha}(u') = \varphi_\alpha(u)$ and hence agrees with the first.
\end{proof}

\begin{lemma}\label{lem:BTmMm_M}
	We have $m(u) \cdot M_\F \cdot m(u)^{-1} = M_\F$.
\end{lemma}
\begin{proof}
	Let $m(u)=a_t\cdot m \in A_\F \cdot N_\F$ as in Proposition \ref{prop:BTmu_in_Wa}. In the spherical Weyl group $W_s = N_\F/M_\F$ we have $mM_\F \cdot m^{-1}M_\F = M_\F$. Then 
	$$
	m(u) \cdot M_\F \cdot  m(u)^{-1} =  a_tm M_\F m^{-1}a_t^{-1} = a_t M_\F a_t^{-1} = M_\F.
	$$
\end{proof}

\begin{lemma}\label{lem:BTmUm}
	Let
$a,b \in A_\F$. If $a.o = m(u).b.o$, then $(-v)(\chi_{\alpha}(a)) = (-v)( t^2 \chi_{\alpha}(b)^{-1} )$.
Moreover
$
m(u)U_{\alpha,\ell}m(u)^{-1} = U_{-\alpha, \ell-2(-v)(t)} = U_{-\alpha, \ell-2\varphi_\alpha(u)}
$ 
for any $\ell \in \Lambda$.
\end{lemma}
\begin{proof}
	 Decomposing $m(u) =a_t \cdot m \in A_\F \cdot N_\F$ as in Proposition \ref{prop:BTmu_in_Wa}, we see $m(u)$ as consisting of a translational element $a_t$ and a representative of an element $w=[m]$ of the spherical Weyl group $W_s$. 
	
	If $a,b\in A_\F$ satisfy $a.o = m(u).b.o = a_t m b m^{-1}.o$, then
	\begin{align*}
		(-v)(\chi_\alpha(a)) &= (-v)(\chi_\alpha(a_t mbm^{-1})) 
		= (-v)(\chi_\alpha(a_t) \chi_\alpha(mbm^{-1})) \\
		&= (-v)(t^{2} \chi_\alpha(b)^{-1}).
	\end{align*}
	where the last equality comes from Lemma \ref{lem:Jacobson_Morozov_oneparam} and Lemma \ref{lem:Jacobson_Morozov_m}.
	
	As in the proof of Lemma \ref{lem:BTaUa}, we now consider $\Omega = \{b.o\}$ with $(-v)(\chi_\alpha(b)) = \ell$. Then $U_{\alpha,\ell} = U_{\alpha,\Omega}$. By Lemma \ref{lem:Jacobson_Morozov_m} and Lemma \ref{lem:aexpXa}, we have
	\begin{align*}
		m(u)\cdot (U_\alpha)_\F \cdot m(u)^{-1} & = a_t  m  (U_\alpha)_\F  m^{-1}  a_t^{-1} = a_t (U_{-\alpha})_\F a_t^{-1} = (U_{-\alpha})_\F.
	\end{align*}
	Since 
	$$	
	(m(u)\cdot u \cdot m(u)^{-1}).m(u).b.o = m(u).b.o,
	$$
	we have thus $m(u) \cdot U_{\alpha, \Omega} \cdot m(u)^{-1} = U_{-\alpha, m(u).\Omega}$. 
	Now by the comment after Lemma \ref{lem:BTUalphaOmega}, $U_{-\alpha, m(u).\Omega} = U_{-\alpha, \ell'}$ for $\ell' = (-v)(\chi_{-\alpha}(a))$ where $a.o = m(u).b.o$. We have
	$$
	\ell' = (-v)(\chi_{-\alpha}(a)) = -(-v)(\chi_\alpha(a)) = (-v)(\chi_\alpha(b)) - 2(-v)(t) = \ell -2(-v)(t)
	$$
	whence $m(u)U_{\alpha,\ell}m(u)^{-1} = U_{-\alpha,\ell-2(-v)(t)}$.
\end{proof}

We know from Lemma \ref{lem:BTm_uuu}, that $m(u)\in (U_{-\alpha})_\F u (U_{-\alpha})_\F$. Next, we will show that $m(u)$ is the only element in $\operatorname{Nor}_{G_\F}(A_\F) \cap  (U_{-\alpha})_\F u (U_{-\alpha})_\F$.

\begin{lemma}\label{lem:BTUUUm}
	Let 
	$u',u'' \in (U_{-\alpha})_\F$ such that $u'uu'' \in \operatorname{Nor}_{G_\F}(A_\F)$. Then $u'uu'' = m(u)$ and
	$$
	\varphi_{-\alpha}(u') = - \varphi_{\alpha}(u) = \varphi_{-\alpha}(u'').
	$$
\end{lemma}
\begin{proof}
	Recall that if
	$$
	u = \varphi_\F\begin{pmatrix}
		1 & t \\ 0 & 1
	\end{pmatrix}, \quad \text{and} \quad \overline{u} := \varphi_\F \begin{pmatrix}
	1 & 0 \\ 1/t & 1
	\end{pmatrix},
	$$
	then
	\begin{align*}
	m(u)& = \varphi_\F \left(\begin{pmatrix}
		1 & 0 \\ -1/t & 1
	\end{pmatrix}\begin{pmatrix}
		1 & t \\ 0 & 1
	\end{pmatrix}\begin{pmatrix}
		1 & 0 \\ -1/t & 1
	\end{pmatrix} \right) \\
	&= \overline{u}^{-1}u\overline{u}^{-1} \in (U_{-\alpha})_\F(U_{\alpha})_\F(U_{-\alpha})_\F.
	\end{align*}
	Then $u'uu'' = u'\overline{u} \cdot m(u) \cdot \overline{u} u'' \in \operatorname{Nor}_{G_\F}(A_\F)$ and $\overline{u}u'' \in (U_{-\alpha})_\F$. Now for any $\ell \leq -\varphi_{-\alpha}(\overline{u}u'')$, use Lemma \ref{lem:Jacobson_Morozov_oneparam} to obtain $a\in A_\F$ such that $(-v)(\chi_\alpha(a)) = \ell $ or equivalently $(-v)(\chi_{-\alpha}(a)) = -\ell$. By Lemma \ref{lem:BTUalphaOmega} this means that $\overline{u}u'' \in U_{-\alpha, -\ell} = U_{-\alpha, \{ a.o\}}$. Let $b \in A_\F$ be such that
	$$
	b.o = u'uu'' . a.o = u'\overline{u} \cdot m(u) \cdot \overline{u}u''.a.o = u'\overline{u} \cdot m(u) \cdot a.o
	$$ 
	and applying Proposition \ref{prop:UonA}, $b.o = m(u).a.o$. Since
	$$
	(-v)(\chi_{-\alpha}(b)) = (-v)(\chi_{\alpha}(a)) -2(-v)(t) = \ell - 2(-v)(t)  
	$$
	by Lemma \ref{lem:BTmUm}, we can use Lemma \ref{lem:BTUalphaOmega} to obtain $u'\overline{u} \in U_{-\alpha, \{ b.o\}} = U_{-\alpha, \ell - 2(-v)(t)}$ for all $\ell \leq -\varphi_{-\alpha}(\overline{u}u'')$. Therefore $\varphi_{-\alpha}(u'\overline{u}) \leq \lambda$ for all $\lambda \in \Lambda$. The only element of $(U_{-\alpha})_\F$ with this property is $u'\overline{u} = \operatorname{Id}$. 
	
	Now for all $a\in A_\F$,
	\begin{align*}
		u'uu''.a.o &= u'\overline{u} m(u) \overline{u}u''.a.o = m(u). \overline{u}u''.a.o \in \A 
	\end{align*}
	and thus $\overline{u}u''.a.o \in \A$, hence $\overline{u}u''.a.o = a.o$ by Proposition \ref{prop:UonA}. The only element of $(U_{-\alpha})_\F$ acting trivially on all of $\A$ is the identity, so $\overline{u}u'' = \operatorname{Id}$. Thus $u'uu'' = \overline{u}^{-1}u\overline{u}^{-1} = m(u)$. Moreover
	\begin{align*}
		\varphi_{-\alpha}(u')=\varphi_{-\alpha}(u'') = \varphi_{-\alpha}(\overline{u}^{-1}) = (-v)(-1/t) = -(-v)(t) = -\varphi_{\alpha}(u). 
	\end{align*}	
\end{proof}

\begin{figure}[h]
	\centering
	\includegraphics[width=1\linewidth]{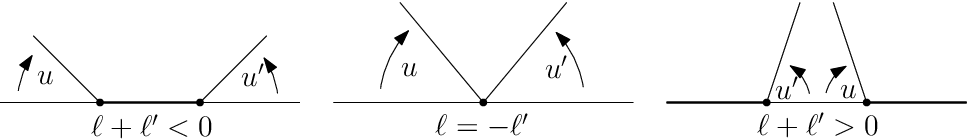}
	\caption{ Illustration of the action of elements $u\in (U_\alpha)_\F, u' \in (U_{-\alpha})_\F$ with $\ell := \varphi_{\alpha}(u) $ and $\ell' := \varphi_{-\alpha}(u')$ in the case where $\ell + \ell' < 0$ (left), $\ell = - \ell'$ (middle) and $\ell + \ell' > 0$. An element $m(u)$ of the affine Weyl group can only be generated by $u,u'$ in the case $\ell =-\ell'$. }
	\label{fig:valuationuu}
\end{figure}

The actions of $ u\in (U_\alpha)_\F, u' \in (U_{-\alpha})_\F$ depending on their root group valuations $\ell, \ell'$ are illustrated in Figure \ref{fig:valuationuu}. Only when $\ell + \ell' \leq 0$ do $u,u'$ fix a common point. To study the stabilizer we consider the cases $\ell + \ell' < 0$ and $\ell = -\ell$ separately in the next two Lemmas.
\begin{lemma}\label{lem:BTuu_usu}
	Let 
	$u'\in (U_{-\alpha})_\F$ such that $\varphi_\alpha(u) + \varphi_{-\alpha}(u') < 0$. Then there is an element $(u_1, t , u_1') \in (U_{\alpha})_\F \times T_\F 
	\times (U_{-\alpha})_\F $ such that $u'u = u_1 t u_1'$. 
	
	Moreover $\varphi_\alpha(u) = \varphi_{\alpha}(u_1)$, $\varphi_{-\alpha}(u') = \varphi_{-\alpha}(u_1')$ and $t \in T_\F(O)$.
\end{lemma}
\begin{proof}
	If $u = \operatorname{Id}$ we can choose $u_1 = \operatorname{Id}, s=\operatorname{Id}$ and $u_1' = u'$. We may thus assume $u \neq \operatorname{Id}$. The element $u'u$ lies in the group $(L_{\pm \alpha})_\F$ defined in Section \ref{sec:rank1}. By Corollary \ref{cor:levi_Bruhat_alternative},
	$$
	(L_{\pm \alpha})_\F = (B_{\alpha})_\F \cdot (B_{-\alpha})_\F  \ \amalg \   m \cdot (B_{ -\alpha})_\F,
	$$
	where $(B_{\alpha})_\F := (M_{\pm \alpha})_\F (A_{\pm \alpha})_\F (U_\alpha)_\F$ and $(B_{-\alpha})_\F := (M_{\pm \alpha})_\F (A_{\pm \alpha})_\F (U_{-\alpha})_\F$. We note that by Lemma \ref{lem:levi_fixes_A} and then Lemma \ref{lem:BT_T_is_MA} that $(M_{\pm \alpha})_\F(A_{\pm \alpha})_\F \subseteq M_\F \cdot A_\F = T_\F 
	$. 
	
	We have that $u'u \notin m \cdot (B_{-\alpha})_\F$, since otherwise there exists $u'' \in (U_{-\alpha})_\F$ such that
	$$
	u'uu'' \in m \cdot (M_{\pm \alpha})_\F (A_{\pm \alpha})_\F  \subseteq \operatorname{Nor}_{G_\F}(A_\F)
	$$
	and hence by Lemma \ref{lem:BTUUUm}, $-\varphi_\alpha(u) = \varphi_{-\alpha}(u')$ contradicting our assumption. Thus
	\begin{align*}
		u'u &\in (B_\alpha)_\F \cdot (B_{-\alpha})_\F = (U_{\alpha})_\F (A_{\pm \alpha})_\F (M_{\pm \alpha})_\F (A_{\pm \alpha}) (U_{-\alpha})_\F \\
		& \subseteq (U_{\alpha})_\F \cdot T_\F 
		\cdot (U_{-\alpha})_\F. 
	\end{align*}
	
	Let now $u_1 \in (U_\alpha)_\F,\ u_1' \in (U_{-\alpha})_\F$ and $ t \in (A_{\pm \alpha})_\F(M_{\pm \alpha})_\F (A_{\pm \alpha})_\F$ with $u'u = u_1 t u_1'$. To show $\varphi_\alpha(u)=\varphi_{\alpha}(u_1)$, we may assume as before that $u_1 \neq \operatorname{Id}$ and may thus consider the element $m(u_1) \in \operatorname{Nor}_{G_\F}(A_\F)$. Then
	$$
	m(u_1) \in (U_{-\alpha})_\F u_1 (U_{-\alpha})_\F,
	$$
	so let $u_2', u_3' \in (U_{-\alpha})_\F$ such that $u_1 = u_2' m(u_1) u_3'$ and by Lemma \ref{lem:BTUUUm} $\varphi_{-\alpha}(u_2') = - \varphi_\alpha(u_1) = \varphi_{-\alpha}(u_3')$. Then
	\begin{align*}
		u&= (u')^{-1}u_1 t u_1' = (u')^{-1} u_2' m(u_1) u_3' t u_1' \\
		&= ((u')^{-1}u_2') m(u_1)t (t^{-1} u_3' t u_1') \in (U_{-\alpha})_\F m(u_1) t (U_{-\alpha})_\F 
	\end{align*}
	where we used that $t^{-1}u_3' t \in (U_{-\alpha})_\F$ by Proposition \ref{prop:anainN} and the fact that $(M_{\pm \alpha})_\F$ acts trivially on the rank 1 root system $\{\alpha, -\alpha\}$. Applying Lemma \ref{lem:BTUUUm} again gives $\varphi_{-\alpha}((u')^{-1}u_2') = - \varphi_\alpha(u) $. Now since $\varphi_\alpha(u) + \varphi_{-\alpha}(u') < 0$, we have 
	$$
	 \varphi_{-\alpha}((u')^{-1}) = \varphi_{-\alpha}(u') < -\varphi_\alpha(u) = \varphi_{-\alpha}((u')^{-1}u_2') \leq \max\{ \varphi_{-\alpha}((u')^{-1}), \varphi_{-\alpha}(u_2) \}
	$$
	from which we can infer that $\varphi_{-\alpha}((u')^{-1}) \neq \varphi_{-\alpha}(u_2')$ and thus by the second part of Lemma \ref{lem:BTgroup_valuation}, $\varphi_{-\alpha}((u')^{-1}u_2') = \varphi_{-\alpha}(u_2')$, whence 
	$$
	\varphi_\alpha(u) = - \varphi_{-\alpha}((u')^{-1}u_2') =  -\varphi_{-\alpha}(u_2') = \varphi_\alpha (u_1).
	$$
	Proving $\varphi_{-\alpha}(u') = \varphi_{-\alpha}(u_1')$ works the same way, after replacing $u$ by $u'$, $u_1$ by $u_1'$ and $\alpha$ by $-\alpha$.
	
	From the above discussion, we obtain
	$$
	m(u_1)t = ((u_2')^{-1}u') u ((u_1')^{-1}t^{-1}(u_3')^{-1}t) \in (U_{-\alpha})_\F u (U_{-\alpha})_\F
	$$
	which by the uniqueness statement of Lemma \ref{lem:BTUUUm} implies $m(u_1)t = m(u)$. Since $\varphi_\alpha(u) = \varphi_\alpha(u_1)$, we know by Proposition \ref{prop:BTmu_in_Wa} that both $m(u_1)$ and $m(u)$ represent the same element of the affine Weyl group: the reflection $r_\alpha$ followed by a translation of length $\varphi_\alpha(u) = \varphi_\alpha(u_1)$. This means that $[t] = [m(u_1)^{-1}][m(u)] = [\operatorname{Id}] \in W_a$, in particular $t$ fixes $o$, so $t \in T_\F(O)$.
\end{proof}

\begin{proposition}\label{prop:BTL_l+l'<0} 
	Let 
	$\ell, \ell' \in \Lambda$ such that $\ell + \ell' <0$. Let 
	$$
	L_{\ell, \ell'} := \langle u \in U_{\alpha, \ell} \cup U_{-\alpha, \ell'} \rangle 
	$$
	and $H_{\ell, \ell'} := L_{\ell,\ell'} \cap T_\F(O) $. Then
	$$
	L_{\ell,\ell'} = U_{\alpha, \ell} \cdot U_{-\alpha, \ell'} \cdot H_{\ell,\ell'} .
	$$
\end{proposition}
\begin{proof}
	The inclusion $\supseteq$ is clear. We show the other direction by induction on the word length. If the length of a word in $L_{\ell, \ell'}$ is $1$, then it lies in $U_{\alpha, \ell} \cdot U_{-\alpha, \ell'} \cdot H_{\ell,\ell'}$. It remains to show that for every $u_{\pm \alpha} \in U_{\alpha, \ell} \cup U_{-\alpha,\ell'}$, and $v = u \cdot u' \cdot h \in U_{\alpha, \ell} \cdot U_{-\alpha, \ell'} \cdot H_{\ell,\ell'} $, we have $u_{\pm \alpha} v \in U_{\alpha, \ell} \cdot U_{-\alpha, \ell'} \cdot H_{\ell,\ell'}$. If $u_\alpha \in U_{\alpha,\ell}$, this is straightforward, so let $u_{-\alpha} \in U_{-\alpha,\ell'}$. 
	By Lemma \ref{lem:BTuu_usu}, we can find $(u_1,t_1,u_1') \in U_{\alpha,\ell}\times T_\F(O) \times U_{-\alpha,\ell'}$ such that $u_{-\alpha}u = u_1t_1u_1'$. Then $u_{-\alpha}v = u_1 t_1 u_1' u' h$. By Lemma \ref{lem:BTtUt}, $t_1U_{-\alpha,\ell'} t_1^{-1} = U_{-\alpha,\ell'}$, so $u_{-\alpha}v \in U_{\alpha,\ell}\cdot U_{-\alpha,\ell'} \cdot H_{\ell,\ell'}$ as required.
\end{proof}
Next up, we consider the situation $\ell + \ell' = 0$. For this we introduce
$$
U_{-\alpha, -\ell+} := \bigcup_{t < -\ell} U_{-\alpha, t} \subseteq U_{-\alpha, -\ell}
$$
which is a group by Lemma \ref{lem:BTgroup_valuation}.

\begin{proposition}\label{prop:BTL_l+l'=0} 
	Let 
	$\ell := \varphi_\alpha(u)$ and let $m(u)$ as in Lemma \ref{lem:Jacobson_Morozov_m}. Let 
	$$
	L_{\ell,- \ell} := \langle \tilde{u} \in U_{\alpha, \ell} \cup U_{-\alpha, -\ell} \rangle 
	$$
	and $H_{\ell, -\ell} := L_{\ell,-\ell} \cap T_\F(O) $. Then
	$$
	L_{\ell,-\ell} = (U_{\alpha,\ell}\cdot H_{\ell, -\ell} \cdot U_{-\alpha,-\ell+}) \cup (U_{\alpha,\ell} \cdot m(u) \cdot H_{\ell,-\ell} \cdot U_{\alpha,\ell}).
	$$
\end{proposition}
\begin{proof}
	By Lemma \ref{lem:BTm_uuu}, $m(u) \in U_{-\alpha,-\ell} U_{\alpha, \ell} U_{-\alpha,-\ell}$ and thus
	$$
	B:= (U_{\alpha,\ell}\cdot H_{\ell, -\ell} \cdot U_{-\alpha,-\ell+}) \cup (U_{\alpha,\ell} \cdot m(u) \cdot H_{\ell,-\ell} \cdot U_{\alpha,\ell}) \subseteq L_{\ell,-\ell}.
	$$
	By Lemma \ref{lem:BTmUm}, $m(u) \cdot U_{\alpha, \ell}\cdot m(u)^{-1} = U_{-\alpha, -\ell}$, which shows that $L_{\ell,-\ell}$ is generated by $U_{\alpha,\ell}$ and $m(u)$. Thus to see $L_{\ell,-\ell} \subseteq B$, it suffices to show that $B$ is a group.
By Proposition \ref{prop:BTL_l+l'<0}, $B\cdot U_{\alpha, \ell} = B$ and $B\cdot U_{-\alpha,-\ell+} = B$. By Lemma \ref{lem:BTtUt}, $B\cdot H_{\ell,-\ell} = B$. It remains to show that $B \cdot m(u) = B$. 
	
	Since $m(u)\cdot U_{\alpha,\ell}\cdot m(u)^{-1} = U_{-\alpha,-\ell}$ and by Lemma \ref{lem:BTmMm_M}, $m(u)\cdot T_\F(O) \cdot m(u)^{-1} = T_\F(O)$, we have
	$$
	(U_{\alpha,\ell} \cdot H_{\ell,-\ell}\cdot U_{-\alpha,-\ell})m(u) = U_{\alpha, \ell} \cdot m(u) \cdot  H_{\ell, -\ell} \cdot U_{\alpha,\ell} \subseteq B
	$$
	and 
	\begin{align*}
     \left( U_{\alpha,\ell}\cdot m(u) \cdot H_{\ell,-\ell} \cdot U_{\alpha,\ell} \right)  m(u)
      = U_{\alpha, \ell} \cdot  H_{\ell, -\ell} \cdot U_{-\alpha,-\ell}\cdot m(u)^2.
	\end{align*}
	By Proposition \ref{prop:BTmu_in_Wa}, $m(u)^{2}$ acts trivially on $\A$, 
	hence it lies in $ \operatorname{Cen}_{G_\F}(A_\F) \cap G_\F(O) = T_\F(O)$, so
	$$
	\left( U_{\alpha,\ell}\cdot m(u) \cdot H_{\ell,-\ell} \cdot U_{\alpha,\ell} \right)  m(u)
    = U_{\alpha, \ell} \cdot  H_{\ell, -\ell} \cdot U_{-\alpha,-\ell}.
	$$
    If $u'\in (U_{-\alpha,-\ell})_\F$ satisfies $\varphi_{-\alpha}(u') = -\ell$, then by Lemma \ref{lem:BTphi_stab}, $m(u)^{-1}m(u') \in T_\F(O) 
    $, so $m(u') \in m(u)\cdot H_{\ell,-\ell}$ since $m(u), m(u') \in L_{\ell,-\ell}$. Thus
    $$
    u' \in U_{\alpha,\ell}\cdot m(u')\cdot U_{\alpha, \ell} \subseteq U_{\alpha,\ell}\cdot  m(u) \cdot H_{\ell,-\ell} \cdot U_{\alpha,\ell}.
    $$
	and
	$$
	U_{-\alpha, -\ell} \subseteq U_{-\alpha,-\ell+} \cup U_{\alpha,\ell}\cdot m(u) \cdot H_{\ell,-\ell} \cdot U_{\alpha,\ell}.
	$$
	This shows 
	\begin{align*}
	&\left( U_{\alpha,\ell}\cdot m(u) \cdot H_{\ell,-\ell} \cdot U_{\alpha,\ell} \right)  m(u)\\
	&\subseteq \left( U_{\alpha, \ell} \cdot  H_{\ell, -\ell} \cdot U_{-\alpha,-\ell + } \right) \cup \left(U_{\alpha, \ell} \cdot  H_{\ell, -\ell} \cdot U_{\alpha,\ell } \cdot m(u)\cdot H_{\ell ,-\ell} \cdot U_{\alpha,\ell}\right) \\
	&= B.
\end{align*}
Similarly, one can show $B\cdot m(u)^{-1} \subseteq B$, which concludes the proof of $B= L_{\ell,-\ell}$.
\end{proof}

\begin{proposition}\label{prop:BTL_l+l'}  
	Let $\Omega \subseteq \A$ be a non-empty finite subset and 
	$$\ell = \min_{a.o \in \Omega} \{ (-v)(\chi_\alpha(a)) \}, \quad \ell' = \min_{a.o \in \Omega} \{ (-v)(\chi_{-\alpha}(a)) \}. 
	$$
	Let 
	$
	L_{\ell,\ell'} := \langle u \in U_{\alpha, \ell} \cup U_{-\alpha, \ell'} \rangle 
	$
	and $N_{\ell, \ell'} := L_{\ell,\ell'} \cap \operatorname{Nor}_{G_\F}(A_\F) $. Then
	$$
	L_{\ell,\ell'} = U_{\alpha, \ell} \cdot U_{-\alpha, \ell'} \cdot N_{\ell,\ell'} .
	$$
\end{proposition}
\begin{proof}
	Since $\Omega \neq \emptyset$, we have $\ell + \ell' \leq 0$. If $\ell+\ell' <0$ we can apply Proposition \ref{prop:BTL_l+l'<0} and are done. Otherwise, we have 
	$$
	L_{\ell,-\ell} = (U_{\alpha,\ell}\cdot H_{\ell, -\ell} \cdot U_{-\alpha,-\ell+}) \cup (U_{\alpha,\ell} \cdot m(u) \cdot H_{\ell,-\ell} \cdot U_{\alpha,\ell})
	$$
	by Proposition \ref{prop:BTL_l+l'=0}. As $H_{\ell,\ell'} \subseteq T_\F(O)$, we can apply Lemma \ref{lem:BTtUt} to obtain $H_{\ell,\ell'}U_{\alpha,\ell} = U_{\alpha,\ell}H_{\ell,\ell'}$ and $H_{\ell,\ell'}U_{-\alpha,\ell'} = U_{-\alpha,\ell'} H_{\ell,\ell'}$. Finally, Lemma \ref{lem:BTmUm} gives $m(u) U_{\alpha, \ell} m(u)^{-1} = U_{-\alpha, \ell- 2 \varphi_\alpha(u)} = U_{-\alpha , -\ell} = U_{-\alpha, \ell'}$.
\end{proof}

We will prove a 'mixed Iwasawa'-decomposition for the rank one subgroup $(L_{\pm \alpha})_\F$ defined in Section \ref{sec:rank1}. The compact group $K_\F$ in the Iwasawa decomposition is replaced by a subgroup of $\operatorname{Stab}_{(L_{\pm \alpha})_\F}(o)$ and an element $m$ representing the non-trivial element of the spherical Weyl group $W_{\pm \alpha} = \{[\operatorname{Id}],[m]\}$ of $(L_{\pm \alpha})_\F$.
\begin{proposition}\label{prop:BT_mixed_Iwasawa_rank_1}
	Let $L_o = \langle U_{\alpha,o} , U_{-\alpha, o} \rangle$. Then 
	\begin{align*}
		L_{\pm \alpha} &= (U_\alpha)_\F \cdot (\operatorname{Nor}_{G_\F}(A_\F) \cap (L_{\pm \alpha})_\F)\cdot L_o \\
		&= (U_\alpha)_\F \cdot (\operatorname{Cen}_{G_\F}(A_\F) \cap (L_{\pm \alpha})_\F)\cdot L_o \ \amalg \ (U_\alpha)_\F \cdot (\operatorname{Cen}_{G_\F}(A_\F) \cap (L_{\pm \alpha})_\F)\cdot m \cdot L_o.
	\end{align*}
\end{proposition}
\begin{proof}
	It is clear that the inclusion $\supseteq$ holds. For the reverse, we use Corollary \ref{cor:levi_Bruhat} to write
	$$
	(L_{\pm \alpha})_\F = (B_\alpha)_\F  \ \amalg \ (B_\alpha)_\F m (B_\alpha)_\F
	$$
	where $m$ is a representative of the non-trivial element in the spherical Weyl group $W_{\pm\alpha}$ associated with $(L_{\pm \alpha})_\F$ and $(B_\alpha)_\F = (M_{\pm \alpha})_\F (A_{\pm \alpha})_\F (U_\alpha)_\F$ as in Section \ref{sec:rank1}. If $g \in (B_\alpha)_\F$, then $g \in U_\alpha (M_{\pm \alpha})_\F (A_{\pm \alpha})_\F \subseteq (U_\alpha)_\F \cdot (\operatorname{Nor}_{G_\F}\cap (L_{\pm\alpha})_\F) $. We claim that
	$$
	m(B_\alpha)_\F = m (M_{\pm\alpha})_\F (A_{\pm \alpha})_\F (U_{\alpha})_\F \subseteq (U_{\alpha})_\F \cdot (\operatorname{Nor}_{G_\F}(A_\F) \cap (L_{\pm \alpha})_\F)\cdot L_o.
	$$
	Assuming the claim we obtain
	\begin{align*}
		(B_\alpha)_\F m (B_\alpha)_\F &\subseteq (B_\alpha)_\F (U_\alpha)_\F \cdot (\operatorname{Nor}_{G_\F}(A_\F) \cap (L_{\pm \alpha})_\F)\cdot L_o \\
		&\subseteq (U_\alpha)_\F \cdot (\operatorname{Nor}_{G_\F}(A_\F) \cap (L_{\pm \alpha})_\F)\cdot L_o
	\end{align*}
	since $(M_{\pm\alpha})_\F$ and $(A_{\pm\alpha})$ normalize $(U_\alpha)_\F$. To prove the claim, we consider an element in $m(M_{\pm \alpha})_\F (A_{\pm \alpha})_\F u$ for $u \in (U_\alpha)_\F$. If $\varphi_\alpha(u)\leq 0$, then $m(M_{\pm \alpha})_\F (A_{\pm \alpha})_\F u \subseteq (U_\alpha)_\F (\operatorname{Nor}_{G_\F}(A_\F) \cap (L_{\pm \alpha})_\F) L_o$ and we are done. If on the other hand $\varphi_\alpha(u) > 0$, then there are $u',u'' \in (U_{-\alpha})_\F$ with $u = u' \cdot m(u) \cdot u''$ and  $\varphi_{-\alpha}(u') = \varphi_{-\alpha}(u'') = -\varphi_\alpha(u) <0$ by Lemma \ref{lem:BTm_uuu}. In particular $u'' \in L_o$ and
	\begin{align*}
		m(M_{\pm \alpha})_\F (A_{\pm \alpha})_\F u &= m(M_{\pm \alpha})_\F (A_{\pm \alpha})_\F u' \cdot m(u) \cdot u'' \\
		&= mu'm^{-1} \cdot m(M_{\pm \alpha})_\F (A_{\pm \alpha})_\F m(u) \cdot u'' \\
		&\subseteq (U_\alpha)_\F \cdot (\operatorname{Nor}_{G_\F}(A_\F) \cap (L_{\pm \alpha})_\F ) \cdot P_o
	\end{align*}
	as claimed.
\end{proof}
The result of Proposition \ref{prop:BTL_l+l'} still holds, when $\Omega =\emptyset$.
\begin{lemma}\label{lem:BTL_l+l'=-infty}
	Let $\alpha \in \Sigma$, $L = \langle (U_{\alpha})_\F , (U_{-\alpha})_\F \rangle$ and 
	$N = \operatorname{Nor}_{G_\F}(A_\F) \cap L$. Then
	$$
	L = (U_{\alpha})_\F(U_{-\alpha})_\F N.
	$$
	\end{lemma}
	\begin{proof}
		Clearly $\supseteq$. We consider the word length of an element $g \in L$. If $g$ consists of at most two letters, it is clearly part of $(U_{\alpha})_\F(U_{-\alpha})_\F N$, unless $g = u'u$ for $u'\in (U_{-\alpha})_\F$ and $u\in (U_\alpha)_\F$. In this case, we write $u = u'' m(u) u'''$ for $u'',u''' \in (U_{-\alpha})_\F$. Then $\bar{u} := m(u) u''' m(u)^{-1} \in (U_{\alpha})_\F$ by Lemma \ref{lem:BTmUm} and we have $g = (u'u'')\bar{u} m(u) \in (U_{\alpha})_\F(U_{-\alpha})_\F N$.
		
		If $g$ contains at least three letters, write $g = u_n u_{n-1} \ldots u_{2} u_1$, where $u_i$ are elements in $(U_{\alpha})_\F, (U_{-\alpha})_\F$ alternating. Write $u_2 = u_2' m(u_2) u_2''$ with $u_2',u_2'' (U_{-\beta})_\F$, where $u_2\in (U_{\beta})_\F$ for $\beta \in \{\pm \alpha\}$. Using Lemma \ref{lem:BTmUm} again we obtain $\bar{u} := m(u_2) u_2'' u_1 m(u_2)^{-1} \in (U_{\beta})_\F $, and
		$g = u_n \ldots u_4 (u_3 u_2') \bar{u} m(u)$, so $g \in h N $ for some $h \in L$ with a smaller word length. Applying induction gives the result.
	\end{proof}

The following description of the stabilizer of a subset $\Omega \subseteq \A$ under the action of the rank 1 subgroup $L$ will not be used later, but is included for completeness. 

\begin{proposition}\label{prop:stab_L_Omega}
	Let $\alpha \in \Sigma$, $\Omega \subseteq \A$ a non-empty finite subset, $L = \langle (U_{\alpha})_\F , (U_{-\alpha})_\F \rangle$, $N  = \operatorname{Nor}_{G_\F}(A_\F) \cap L$ and $N_\Omega = \{g \in N \colon n.p = p \text{ for all } p \in \Omega\}$. Then
	$$
	\operatorname{Stab}_{L} (\Omega) = \langle U_{\alpha,\Omega} , U_{-\alpha,\Omega} , N_\Omega \rangle =  U_{\alpha, \Omega} U_{-\alpha,\Omega} N_\Omega.
	$$
\end{proposition}
\begin{proof}
	The inclusions $\supseteq$ are clear. Let 
	$$
	\ell = \min_{a.o \in \Omega} \{ (-v)(\chi_\alpha(a)) \} \quad \text{and} \quad \ell' = \min_{a.o \in \Omega} \{ (-v)(\chi_{-\alpha}(a)) \}.
	$$
	Then $\ell \leq -\ell'$ since $\Omega \neq \emptyset$, $U_{\alpha,\ell} = U_{\alpha,\Omega}$ and $U_{-\alpha, \ell'} = U_{-\alpha,\Omega}$. Now if $g\in \operatorname{Stab}_L(\Omega)$, then there are $u\in (U_\alpha)_\F , u' \in (U_{-\alpha})_\F , n \in N$ such that $g.p = uu'n.p = p$ for all $p\in \Omega$, by Lemma \ref{lem:BTL_l+l'=-infty}. We distinguish three cases for $\varphi_\alpha(u)$. If $\varphi_{\alpha}(u) \leq \ell$, then $u \in U_{\alpha,\ell}$ and $u'n.p = u^{-1}.p = p$ for all $p \in \Omega$, so $u' \in U_{-\alpha, \ell'}$ and $n\in N_\Omega$ by Proposition \ref{prop:UonA}, so $g \in U_{\alpha,\ell} U_{-\alpha,\ell'} N_{\Omega}$.
	
	If $\varphi_\alpha(u) \geq -\ell'$, then we write $u = u'' m(u) u''' $ for some $u'',u''' \in (U_{-\alpha})_\F $ with $\varphi_{-\alpha}(u'') = \varphi_{-\alpha}(u''') = -\varphi_\alpha(u) < \ell'$. Abbreviating $\bar{u}: = m(u) u'''u' m(u)^{-1} \in (U_{\alpha})_\F $ we obtain
	$g = uu'n = u'' \bar{u} m(u) n \in U_{-\alpha,\ell'} \cdot (U_{\alpha})_\F  \cdot N$ by Lemma \ref{lem:BTmUm}. Then $\bar{u}m(u)n.p = p$ for all $p \in \Omega$, whence $\varphi_{\alpha}(\bar{u}) \leq \ell$ by Proposition \ref{prop:UonA} and then $m(u)n.p = p$ for all $p\in \Omega$, so $m(u)n \in N_\Omega$. Then $g \in U_{\alpha,\ell}  U_{-\alpha,\ell'} N_\Omega$.
	
	We claim that the final case $\ell < \varphi_\alpha(u) < -\ell'$ cannot happen. For this we first notice that for $p\in \Omega$ we have $u.p = p$ if and only if $u'n.p = p$, if and only if $u'.p = p$. Now let $a_1,a_2 \in A_\F$ with $(-v)(\chi_\alpha(a_1)) = \ell$ and $(-v)(\chi_{\alpha}(a_2)) = -\ell'$. Now $u.a_1.o \neq a_1.o$, so $u'.a_1.o \neq a_1.o$, so $\varphi_{-\alpha}(u') > -\ell$, but $u.a_2.o = a_2.o$, so $u'.a_2.o = a_2.o$, so $\varphi_{-\alpha}(u')\leq \ell'$. But this is not possible, since then $-\ell \leq \varphi_{-\alpha}(u') \leq \ell'$ contradicting $\ell < \ell'$.
	We have shown that $\operatorname{Stab}_{L}(\Omega) \subseteq U_{\alpha,\ell} U_{-\alpha,\ell'} N_{\Omega}$.  
\end{proof}

\subsubsection{Bruhat-Tits theory in higher rank}\label{sec:BT_higher_rank}

In this section, we continue to assume that $\Sigma$ is a reduced root system. Let $>$ be an order on $\Sigma$ and $\Omega \subseteq \A$. Now for any subset $\Theta \subseteq \Sigma_{>0}$ closed under addition, let
$$
U_{\Theta,\Omega} :=\left\{ u \in \exp\left( \bigoplus_{\alpha \in \Theta} (\frakg_\alpha)_\F \right) \colon  u.p=p \text{ for all } p \in \Omega \right\} 
$$
and in particular
$$
U_{\Omega}^+ := U_{\Sigma_{>0}, \Omega} = \left\{ u \in U_{\F} \colon u.p=p \text{ for all } p\in \Omega\right\}
$$
for $\Theta = \Sigma_{>0}$. Note that $U_{\emptyset}^+ = U_\F$. Analogously we define
$$
U_{\Omega}^- := \left\{ u \in \exp \left( \bigoplus_{\alpha <0} (\frakg_\alpha)_\F \right) \colon  u.p=p \text{ for all } p\in \Omega \right\}.
$$
Let
$$
N_\Omega : = \left\{ n \in \operatorname{Nor}_{G_\F}(A_\F) \colon n.p = p \text{ for all } p \in \Omega \right\}.
$$
When $\Omega = \{p\}$ consists of a single point, we abbreviate the notation by omitting the brackets such as in 
$U_{p}^+ := U_{\{p\}}^+$, $U_{p}^- := U_{\{p\}}^-$ and $N_p := N_{\{p\}}$. The goal of this subsection is to prove in Theorem \ref{thm:BTstab_fin} that if $\Omega$ is finite, then the pointwise stabilizer of $\Omega$ satisfies
$$
\operatorname{Stab}_{G_\F}(\Omega) = U_{\Omega}^+ U_\Omega^- N_\Omega.
$$
\begin{proposition}\label{prop:BTUOmega}
	For any subset $\Omega \subseteq \A$ and subset $\Theta \subseteq \Sigma_{>0}$ closed under addition,
	$
	U_{\Theta, \Omega} = \left\langle u \in U_{\alpha,\Omega} \colon \alpha \in \Theta \right\rangle
	$
	and in particular, 
	$
	U_\Omega^+ = \langle  u \in U_{\alpha, \Omega} \colon  \alpha >0 \rangle.
	$
	More precisely, if $\Sigma_{>0} = \left\{ \alpha_1, \ldots , \alpha_r \right\}$ with $\alpha_1 < \ldots < \alpha_r$, then the product map
	\begin{align*}
		\prod_{\alpha \in \Sigma_{>0}} U_{\alpha, \Omega} & \to U_\Omega^+ \\
			(u_{\alpha_1}, \ldots , u_{\alpha_r}) & \mapsto u_{\alpha_1} \cdots u_{\alpha_r}
	\end{align*}
    is a bijection.
\end{proposition}
\begin{proof}
	For any $\Theta = \{\alpha_1, \ldots, \alpha_r\} \subseteq \Sigma_{>0}$ closed under addition with $\alpha_1 > \ldots > \alpha_k$, the image of $\prod_{\alpha \in \Theta} U_{\alpha,\Omega} $ under the product map is contained in $U_{\Theta,\Omega}$. On the other hand, if we start with $u \in U_{\Theta, \Omega}$, we can apply Lemma \ref{lem:BCH_consequence} to obtain $(u_1, \cdots , u_k) \in (U_{\alpha_1}) \times \ldots (U_{\alpha_k})_\F $ with $u = u_1 \cdot \ldots \cdot u_k$. A point $a.o \in \Omega$ is fixed by $u$ if and only if $a^{-1}ua \in U_\F(O)$, so
	$$
	a^{-1}ua \cdot \ldots \cdot a^{-1}ua \in U_\F(O).
	$$
	We can apply Lemma \ref{lem:nn'} repeatedly to obtain $a^{-1}u_1a \in U_\F(O)$, \ldots, $a^{-1}u_ka \in U_\F(O)$ for all $a.o \in \Omega$. This exactly means that $u_i \in U_{\alpha,\Omega}$. This shows that the product map is injective. Note that the order in the statement of the Proposition is the inverse of the one used in the proof, but one follows from the other by applying the inverse.
\end{proof}

We notice that the previous propositions hold for any chosen order on $\Sigma$. We may in particular invert the order to obtain
\begin{align*}
	U_\Omega^- = \langle u \in U_{\alpha,\Omega} \colon \alpha < 0 \rangle 
\end{align*}
from Proposition \ref{prop:BTUOmega}. 
We now define the groups
$$
P_\Omega := \langle  U_{\alpha,\Omega} \colon \alpha \in \Sigma \rangle \quad \text{ and } \quad 
\hat{P}_\Omega := \langle N_\Omega, U_{\alpha,\Omega}\colon \alpha \in \Sigma \rangle
$$
generated by all elements of $U_{\alpha, \Omega}$, not just those with positive $\alpha$. The following is a generalization of Lemma \ref{prop:BTL_l+l'} to higher rank.
\begin{proposition}\label{prop:BTPOmegaUUN} 
	Let $\Omega \subseteq \A$ finite. Then $P_\Omega$ and $\hat{P}_\Omega$ decompose as
	\begin{align*}
		P_\Omega &= U_\Omega^- \cdot U_\Omega^+ \cdot(\operatorname{Nor}_{G_\F}(A_\F)\cap P_{\Omega}) \\
		\hat{P}_\Omega &= U_\Omega^- \cdot U_\Omega^+ \cdot N_{\Omega}.
	\end{align*}
\end{proposition}
\begin{proof}
	The inclusions $\supseteq$ are clear. Recall that the set $\Sigma_{>0}$ is determined by a basis $\Delta \subseteq \Sigma$, or equivalently a chamber. After choosing a chamber, an ordering $>$ of $\Sigma$ can be obtained by choosing an order $\alpha_1 < \alpha_2 < \ldots < \alpha_r$ on $\Delta = \{ \alpha_1, \alpha_2, \ldots , \alpha_r \}$; the ordering on $\Sigma$ is then the lexicographical ordering \cite[VI.1.6, p. 174-175]{Bou08}. 
	
	For the $\subseteq$ direction of the description of $P_\Omega$, we will first show that 
	$$
	U_\Omega^- \cdot U_\Omega^+ \cdot(\operatorname{Nor}_{G_\F}(A_\F)\cap P_{\Omega})
	$$
	is independent of the chamber defining $\Sigma_{>0}$.
	Let $<_1, <_2$ be two orderings on $\Sigma$ whose chambers $C_1, C_2$ are related by a reflection determined by a simple root $\alpha \in \Delta_1 \subseteq \Sigma$, so $\Delta_2 = r_\alpha(\Delta_1)$. We may then assume that $<_1$ and $<_2$ are determined by lexicographical orderings such that $0 <_1 \alpha <_1 \beta$ for all $\beta \in \Delta_{1} \setminus \{\alpha\}$ and $0 <_2 -\alpha <_2 \beta$ for all $\beta \in \Delta_2\setminus \{\alpha\}$. We notice that then $\Sigma_{>_1 0} \setminus \{\alpha\} = \Sigma_{>_2 0} \setminus \{-\alpha\}$, since
	for $\beta \in \Sigma_{>_1 0}\setminus \{\alpha\}$ there are $\lambda_{\delta} \in \Z_{\geq 0}$, at least one of which is strictly positive for some $\delta \in \Delta_1 \setminus \{\alpha\}$, such that
	\begin{align*}
		\beta &= \sum_{\delta\in \Delta_1} \lambda_\delta \delta 
		= \sum_{\delta \in \Delta_1} \lambda_\delta \left(r_\alpha(\delta) - 2\frac{\langle \delta, \alpha\rangle}{\langle \alpha, \alpha\rangle} r_\alpha(\alpha)\right) \\
		&= \sum_{\delta \in \Delta_1 \setminus \{ \alpha\}} \lambda_\delta r_\alpha(\delta) + \left(\sum_{\delta \in \Delta_1 \setminus \{\alpha\}} - 2 \frac{\langle\delta,\alpha\rangle}{\langle \alpha, \alpha \rangle} \lambda_\delta - \lambda_\alpha \right) r_\alpha(\alpha)  \\
		& \in \sum_{\delta \in \Delta_1\setminus \{\alpha\}} \mathbb{Z}_{\geq 0} r_\alpha(\delta) + \mathbb{Z} r_\alpha(\alpha) = \sum_{\delta \in \Delta_2\setminus \{-\alpha \}} \mathbb{Z}_{\geq 0} \delta + \mathbb{Z} (-\alpha) ,
	\end{align*}
	where we first used that $\langle \delta, \alpha \rangle$ is nonpositive for all $\delta \in \Delta_1 \setminus \{\alpha\}$. Bases of a root system have the property that any element of the root system written in that basis has either all non-negative or all non-positive coefficients. Since there is a strictly positive coefficient for the element $\beta$ written in the basis $\Delta_2$ as above, all coefficients have to be non-negative, so $\beta \in \Sigma_{>_20} \setminus \{-\alpha\}$. 
	 It now suffices to show that 
	 $$
	 U_\Omega^- \cdot U_\Omega^+ \cdot(\operatorname{Nor}_{G_\F}(A_\F)\cap P_{\Omega})
	 $$ is the same for $<_1$ and $<_2$ to deduce that it is independent of the chamber used to define $\Sigma_{>0}$, since the reflections determined by the simple roots $\Delta$ generate the spherical Weyl group which acts transitively on the chambers.
	 
	 We abbreviate $Y := (\operatorname{Nor}_{G_\F}(A_\F)\cap P_{\Omega})$. We use Proposition \ref{prop:BTUOmega} repeatedly to obtain 
	 \begin{align*}
	 	U_{\Sigma_{<_1 0},\Omega} \cdot U_{\Sigma_{>_1 0},\Omega} \cdot Y 
	 	&= \prod_{\beta <_1 0} U_{\beta, \Omega} \cdot \prod_{\beta >_1 0} U_{\beta, \Omega} \cdot Y \\
	 	&= \prod_{\substack{\beta <_1 0\\ -\alpha \neq \beta}} U_{\beta, \Omega} \cdot  U_{-\alpha, \Omega} \cdot \prod_{\substack{\beta >_1 0\\ \alpha \neq \beta}} U_{\beta, \Omega}  \cdot U_{\alpha, \Omega} \cdot Y \\
	 	&= \prod_{\substack{\beta <_1 0\\ -\alpha \neq \beta}} U_{\beta, \Omega} \cdot U_{\Sigma_{>_2 0},\Omega}  \cdot U_{\alpha, \Omega} \cdot Y \\
	 	&= \prod_{\substack{\beta <_1 0\\ -\alpha \neq \beta}} U_{\beta, \Omega} \cdot  \prod_{\substack{\beta >_1 0\\ \alpha \neq \beta}} U_{\beta, \Omega} \cdot U_{-\alpha, \Omega} \cdot U_{\alpha, \Omega} \cdot Y 
	 \end{align*}
	at which point we invoke Lemma \ref{prop:BTL_l+l'} to continue
	\begin{align*}
		&= \prod_{\substack{\beta <_1 0\\ -\alpha \neq \beta}} U_{\beta, \Omega} \cdot  \prod_{\substack{\beta >_1 0\\ \alpha \neq \beta}} U_{\beta, \Omega} \cdot U_{\alpha, \Omega} \cdot U_{-\alpha, \Omega} \cdot Y \\
		&= \prod_{\substack{\beta <_1 0\\ -\alpha \neq \beta}} U_{\beta, \Omega} \cdot  \prod_{\beta >_1 0} U_{\beta, \Omega} \cdot U_{-\alpha, \Omega} \cdot Y \\
		&= \prod_{\substack{\beta <_1 0\\ -\alpha \neq \beta}} U_{\beta, \Omega} \cdot  U_{\alpha, \Omega}  \cdot \prod_{\substack{\beta >_1 0\\ \alpha \neq \beta}} U_{\beta, \Omega} \cdot  U_{-\alpha, \Omega} \cdot Y \\
		&=  \prod_{\beta <_2 0} U_{\beta, \Omega}  \cdot \prod_{\beta >_2 0} U_{\beta, \Omega} \cdot Y \\
		&= U_{\Sigma_{<_2 0}, \Omega} \cdot U_{\Sigma_{>_2 0},\Omega} \cdot Y.
	\end{align*}
	Now that we have shown that it is independent of the chamber defining the order on $\Sigma$, we show the direction $\subseteq$ by induction on word length. The base case is clear. Now let $u_1 \in U_{\alpha, \Omega}$ for some $\alpha \in \Sigma$ and let $u \in P_\Omega$. We may choose the order on $\Sigma$ such that $u_1 \in U_\Omega^-$. Then by the induction assumption, we have $u \in U_\Omega^- \cdot U_\Omega^+ \cdot Y$ and hence also $u_1u \in U_\Omega^- \cdot U_\Omega^+ \cdot Y$.
	
	It remains to show $\hat{P} \subseteq U_\Omega^- U_{\Omega}^+ N_\Omega$. We use the fact that for every $n \in N_\Omega$ we have $n U_{\alpha,\Omega} n^{-1} = U_{[n](\alpha),n.\Omega} U_{[n](\alpha), \Omega}$ where $[n]$ is the representative of the spherical Weyl group corresponding to $n$. This can be used to show
	$$
	\hat{P}_\Omega = \langle N_\Omega, U_{\alpha,\Omega} \colon \alpha \in \Sigma \rangle = \langle U_{\alpha, \Omega} \colon \alpha \in \Sigma \rangle \cdot N_\Omega = P_\Omega N_\Omega = U_\Omega^- U_\Omega^+ N_\Omega, 
	$$ 
	making use of the description of $P_\Omega$.
\end{proof}

Next, we will obtain a 'mixed Iwasawa'-decomposition for the group $G_\F$, upgrading the result of Proposition \ref{prop:BT_mixed_Iwasawa_rank_1} in rank 1.
\begin{theorem}\label{thm:BT_mixed_Iwasawa}
	Fix an order on the root system $\Sigma$. Then
	$G_\F = U^+ \cdot \operatorname{Nor}_{G_\F}(A_\F) \cdot \hat{P}_o$.
\end{theorem}
\begin{proof}
	We abbreviate $\tilde{G} := U^+ \cdot \operatorname{Nor}_{G_\F}(A_\F) \cdot \hat{P}_o$. It suffices to show that $g\tilde{G} \subseteq \tilde{G}$ for every $g\in G_\F$, since then $G_\F \subseteq G_\F \tilde{G} \subseteq \tilde{G} \subseteq G_\F$. By Theorem \ref{thm:BWB}, and the fact that non-trivial elements of the spherical Weyl group $W_s$ of the form $m(u)$ can be obtained from elements in $U^+$ and $U^-$, $G_\F = \langle U^+, U^-, T_\F  \rangle$, where 
	\begin{align*}
		U^+ &= \langle (U_\alpha)_\F \colon \alpha >0 \rangle \\
		U^- &= \langle (U_\alpha)_\F \colon \alpha <0 \rangle \\
		T_\F &= \operatorname{Cen}_{G_\F} (A_\F)
	\end{align*}
	as earlier. Since $U^+ \tilde{G} \subseteq \tilde{G}$ and $T_\F \tilde{G} \subseteq \tilde{G}$, it suffices to show $U^- \tilde{G} \subseteq \tilde{G}$. Since equality holds in $[(\frakg_\alpha)_\F, (\frakg_\beta)_\F] = (\frakg_{\alpha + \beta})_\F$, 
	$$
	U^- = \langle (U_{-\alpha})_\F \colon \alpha \in \Delta \rangle
	$$
	for the basis $\Delta \subseteq \Sigma$ associated to the chosen order, compare \cite[Remark 14.5(2)]{Bor}. Therefore it suffices to show $(U_{-\alpha})_\F \tilde{G} \subseteq \tilde{G}$ for all $\alpha\in \Delta$ to show the theorem.
	
	Now let $\alpha \in \Delta$ and consider the complete lexicographical order on $\Sigma$ such that $\alpha < \delta$ for all $\delta \in \Delta$. Let $U_\alpha' := \langle (U_\beta)_\F \colon \beta >0, \beta \neq \alpha\rangle \subseteq U^+$. Then by Lemma \ref{lem:BCH_normalizer}, $U^+ = (U_\alpha)_\F \cdot U_\alpha' = U_\alpha' \cdot (U_\alpha)_\F$. We note that when we instead consider the order with basis $r_\alpha(\Delta)$, $-\alpha$ is the smallest positive element as in the proof of Proposition \ref{prop:BTUOmega} and then Lemma \ref{lem:BCH_normalizer} gives $(U_{-\alpha})_\F \cdot U_\alpha' = U_\alpha' \cdot (U_{-\alpha})_\F$. We use Proposition \ref{prop:BT_mixed_Iwasawa_rank_1} to show
	\begin{align*}
		(U_{-\alpha})_\F \tilde{G} &= (U_{-\alpha})_\F U_\alpha' (U_\alpha)_\F \operatorname{Nor}_{G_\F}(A_\F) \hat{P}_o \\
		&= U_\alpha' (U_{-\alpha})_\F (U_\alpha)_\F \operatorname{Nor}_{G_\F}(A_\F) \hat{P}_o \\
		&\subseteq U_\alpha' (U_\alpha)_\F (\operatorname{Nor}_{G_\F}(A_\F) \cap (L_{\pm\alpha})_\F) L_o \cdot \operatorname{Nor}_{G_\F}(A_\F) \hat{P}_o \\
		&\subseteq  U^+ T_\F L_o \cdot \operatorname{Nor}_{G_\F}(A_\F) \hat{P}_o \ \amalg \ U^+ T_\F  m L_o \cdot \operatorname{Nor}_{G_\F}(A_\F) \hat{P}_o,
	\end{align*}
	where $L_o = \langle U_{\alpha, o}, U_{-\alpha, o} \rangle$ and $m\in \operatorname{Nor}_{G_\F}(A_\F)$ represents the reflection $r_\alpha$ in the spherical Weyl group $W_s$. By Lemma \ref{prop:BTL_l+l'}, $L_o \subseteq U_{\alpha, o} U_{-\alpha,o} \operatorname{Nor}_{G_\F}(A_\F)$ and $L_o \subseteq U_{-\alpha, o} U_{\alpha,o} \operatorname{Nor}_{G_\F}(A_\F)$. Then
	\begin{align*}
		(U_{-\alpha})_\F \tilde{G} & \subseteq U^+ T_\F (U_\alpha)_\F (U_{-\alpha})_\F \operatorname{Nor}_{G_\F}(A_\F) \hat{P}_o \ \amalg \ U^+ T_\F m (U_{-\alpha})_\F (U_{\alpha})_\F \operatorname{Nor}_{G_\F}(A_\F) \hat{P}_o \\
		&=  U^+ T_\F (U_\alpha)_\F (U_{-\alpha})_\F \operatorname{Nor}_{G_\F}(A_\F) \hat{P}_o \\
		&= U^+ (U_\alpha)_\F (U_{-\alpha})_\F T_\F \operatorname{Nor}_{G_\F}(A_\F) \hat{P}_o \\
		&= U^+ (U_{-\alpha})_\F \operatorname{Nor}_{G_\F}(A_\F) \hat{P}_o
	\end{align*}
	where we used $m (U_{\alpha})_\F m^{-1} = (U_{-\alpha})_\F$ and $m (U_{-\alpha})_\F m^{-1} = (U_{\alpha})_\F$. We claim that for every $u' \in (U_{-\alpha})_\F$ and $n \in \operatorname{Nor}_{G_\F}(A_\F)$, $u'n \in \tilde{G}$. To to see this, consider $\beta = [n^{-1}](\alpha)$ and $v := n^{-1}u'n\in (U_{[n^{-1}](-\alpha)})_\F = (U_{-\beta})_\F $. Then we can apply the rank one mixed Iwasawa decomposition, Proposition \ref{prop:BT_mixed_Iwasawa_rank_1}, for $\beta$ to obtain 
	\begin{align*}
		u'n = nv &\in n(L_{\pm \beta})_\F \subseteq n(U_{\beta})_\F \operatorname{Nor}_{G_\F}(A_\F) \hat{P}_o \\
		&= n (U_{\beta})_\F n^{-1} n \operatorname{Nor}_{G_\F}(A_\F) \hat{P}_o  \\
		&= (U_\alpha)_\F \operatorname{Nor}_{G_\F}(A_\F) \hat{P}_o \subseteq \tilde{G}.
	\end{align*}
	Assuming the claim, $(U_{-\alpha})_\F \tilde{G} \subseteq U^+ (U_{-\alpha})_\F \operatorname{Nor}_{G_\F}(A_\F) \hat{P}_o \subseteq U^+ \tilde{G} \hat{P}_o = \tilde{G}$ completes the proof. 
\end{proof}

The mixed Iwasawa decomposition can be used to show that $\hat{P}_o = \operatorname{Stab}_{G_\F}(o)$. In fact we will eventually show that $\hat{P}_\Omega$ is the whole pointwise stabilizer of $\Omega$ in Theorem \ref{thm:BTstab}.
\begin{proposition}\label{prop:BTstab_point}
	The stabilizer of a single point $p \in \A$ satisfies
	$$
	\operatorname{Stab}_{G_\F}(p) = \hat{P}_p = \langle N_{p} , U_{\alpha, p} \colon \alpha \in \Sigma \rangle = U_p^- U_{p}^+ N_p.
	$$
\end{proposition}
\begin{proof} The two expressions on the right coincide with $\hat{P}_p$, see Proposition \ref{prop:BTPOmegaUUN}.	The inclusion $\operatorname{Stab}_{G_\F}(p)  \supseteq \hat{P}_p$ is clear. For the other direction, we first consider $p=o$ and use the mixed Iwasawa decomposition, Theorem \ref{thm:BT_mixed_Iwasawa}, $G_\F = U^+ \operatorname{Nor}_{G_\F}(A_\F) \hat{P}_o$. For $g \in \operatorname{Stab}_{G_\F}(o)$, let $u \in U^+, n \in \operatorname{Nor}_{G_\F}(A_\F)$ and $p \in \hat{P}_o$ with $g=unp$. Then $o = g.o = unp.o = un.o = n.o$ by Proposition \ref{prop:UonA}. Thus $n \in N_o$ and hence also $u \in U_o^+ = \langle U_{\alpha,o} \colon \alpha>0 \rangle$.
	
	For general $p=a.o\in \A$ and elements $g\in G_\F$, $g.p=p$ if and only if $a^{-1}ga.o=o$. Use Proposition \ref{prop:BTPOmegaUUN} to write $a^{-1}ga = uu'n$ as a product of elements $u\in U^+_o, u' \in U_o^-$ and $n \in N_o$. We note that $aua^{-1} \in U_p^+, au'a^{-1} \in U_p^-$ and $ana^{-1} \in N_p$, whence $g \in \hat{P}_p$.
\end{proof}

We state a Lemma that allows us to prove that $\hat{P}_\Omega$ is the pointwise stabilizer of any finite subset $\Omega\subseteq \A$. 
\begin{lemma}\label{lem:BTstab_dir_lemma}
	For $p,q\in \A$ there is an order on $\Sigma$ such that $U_q^+ \subseteq U_p^+$.
\end{lemma}
\begin{proof}
	If $q=o$ and $p$ lies in the fundamental Weyl chamber $C_0 \subseteq \A$, then we can take the standard order $>$ associated to $C_0$. Elements $g\in U_q^+ = U_\F(O)$ stabilize all of $C_0$ pointwise, see Theorem \ref{thm:NO_fixes_s0}, so $U_q^+ \subseteq U_p^+$.
	
	In general, there is an element $n \in \operatorname{Nor}_{G_\F}(A_\F)$ with $n.q = 0$ and $n.p \in C_0$ ($n$ is a translation by $-q$ followed by a representative of a unique element in the spherical Weyl group). Then by the above $U_{n.q}^+ \subseteq U_{n.p}^+$ with respect to the standard order associated to $C_0$. If $>_2$ and $+_2$ denote the order with positive roots $\Sigma_{>_2 0} = \left\{ \alpha \in \Sigma \colon [n^{-1}](\alpha) > 0 \right\} $, then $U_{n.p}^+ =n^{-1} U_{p}^{+_2} n$, so
	$$
	U_q^{+_2} = n U_{n.q}^+ n^{-1} \subseteq n U_{n.p}^+ n^{-1}= U_p^{+_2}. 
	$$
\end{proof}

\begin{theorem}\label{thm:BTstab_fin}
	Let $\Omega \subseteq \A$ be a finite subset. Then the pointwise stabilizer of $\Omega$ satisfies
	$$
	\operatorname{Stab}_{G_\F}(\Omega) = \hat{P}_\Omega = \langle N_\Omega, U_{\alpha,\Omega} \colon \alpha \in \Sigma \rangle =  U_\Omega^+ U_\Omega^- N_\Omega.
	$$
\end{theorem}
\begin{proof}
	We claim that
	$$
	\hat{P}_{\Omega \cup \{p\}} = \hat{P}_\Omega \cap \hat{P}_p
	$$
	for all non-empty, finite $\Omega \subseteq \A$ and $p \in \A$. The inclusion $\subseteq$ is clear. Now let $g \in \hat{P}_\Omega \cap \hat{P}_p$. For some $q \in \Omega$, choose an ordering of $\Sigma$ such that $U_q^+ \subseteq U_p^+$ by Lemma \ref{lem:BTstab_dir_lemma}. 
	Now use Proposition \ref{prop:BTPOmegaUUN} to write $g = nu'u$ with $n \in N_{\Omega}$, $u' \in U_{\Omega}^-$, $u\in U_{\Omega}^+$. We have $u \in U_{\Omega}^+ \subseteq U_{q}^+ \subseteq U_{p}^+$, so $gu^{-1} = nu' \in \hat{P}_p$. Since now $u'.p = n^{-1}.p$, we have $u'.p = p$ by Proposition \ref{prop:UonA}. This means that all three elements $u, u'$ and $n$ fix $\Omega \cup \{p\}$, so $g \in N_{\Omega \cup \{p\}} U_{\Omega \cup \{p\}}^- U_{\Omega \cup \{p\}}^+ = \hat{P}_{\Omega \cup \{p\}}$. 
	
	We now show $\operatorname{Stab}(\Omega) = \hat{P}_{\Omega} $ by induction over the size of $\Omega$. If $|\Omega|=1$, this is Proposition \ref{prop:BTstab_point}. Now assume the statement holds for $\Omega$. Given any $p\in \A$, we have
	$$
	\operatorname{Stab}_{G_\F}(\Omega \cup \{p\}) = \operatorname{Stab}_{G_\F}(\Omega) \cap \operatorname{Stab}_{G_\F}(p) 
	= \hat{P}_\Omega \cap \hat{P}_p = \hat{P}_{\Omega \cup \{p\}},
	$$
	where we used the induction assumption and the claim.
\end{proof}

\subsection{Axiom (A2)} \label{sec:A2_subsubsection}

In this section we assume $\Sigma$ to be reduced. We are now equipped to prove axiom (A2). We first show a special case.

\begin{proposition}\label{prop:A2_b}
	Let $g\in G_\F$ and $\Omega = g^{-1}\A \cap \A$. Then there exists an element $w\in W_a$ such that for all $p \in \Omega$, $g.p = w(p)$.
\end{proposition}
\begin{proof}
	We may assume that $\Omega \neq \emptyset$. For subsets $Y \subseteq \Omega$ we consider
	$$
	N_{g,Y} := \left\{ n \in \operatorname{Nor}_{G_\F}(A_\F) \colon  g.p = n.p \text{ for all } p\in Y \right\},
	$$
	so that our goal is to prove that $N_{g,\Omega} \neq \emptyset$.
	If $Y=\{p\}$ is a singleton, then $N_{g,\left\{p\right\}}$ is non-empty, since if $p=a.o$ and $g.p = b.o$ for $a,b\in A_\F$, then $ba^{-1} \in N_{\{p\}}$. We will now show by induction on the size of $Y$, that $N_{g,Y} \neq \emptyset$ for all \emph{finite} sets $Y \subseteq \Omega$.
	
	Let $Y\subseteq \Omega$ finite and $p \in \Omega$. Assume that there are $n_{Y} \in N_{g,Y}$ and $n_p \in N_{g,\{p\}} $, then $g^{-1}n_Y$ and $g^{-1}n_p $ lie in the pointwise stabilizers $\operatorname{Stab}_{G_\F}(Y)$ and $ \operatorname{Stab}_{G_\F}(\{p\})$. We may choose an order on $\Sigma$ such that $U_{Y}^+ \subseteq U_{\{p\}}^+$, by making sure that the defining chamber for the	order based at some point in $Y$ contains $p$, see Lemma \ref{lem:BTstab_dir_lemma}.
	  Then, by Theorem \ref{thm:BTstab_fin},
	\begin{align*}
		n_Y^{-1} n_p &\in N_Y U_Y^- U_Y^+ U_{\{p\}}^+ U_{\{p\}}^- N_{\{p\}} \\
		&=  N_Y U_{Y}^- U_{\{p\}}^+ U_{\{p\}}^- N_{\{p\}}\\
		&=  N_Y U_{Y}^- U_{\{p\}}^- U_{\{p\}}^+ N_{\{p\}}\\
		&\subseteq N_Y U^- U^+ N_{\{p\}}
	\end{align*}
	Let $n_Y' \in N_Y, n_p' \in N_{\{p\}}$ such that $(n_Y')^{-1}n_Y^{-1}n_pn_p' \in U^{-} U^{+} \cap \operatorname{Nor}_{G_\F}(A_\F)$. 
	Since every element of $U^+$ stabilizes some affine chamber pointwise, so does every element of $U^-U^+$, by Proposition \ref{prop:UonA}. The element of the affine Weyl group represented by $(n_Y')^{-1}n_Y^{-1}n_pn_p'$ thus acts trivially on some affine chamber, hence acts trivially on all of $\A$. 
	Therefore $[n_Yn_Y'] = [n_pn_p'] \in W_a$ and $ n_Y n_Y' \in N_{g,Y} \cap N_{g,\{p\}} = N_{g,Y\cup \{p\}}$, concluding the induction.

	 Recall that the affine Weyl group $W_a = W_s \ltimes \A$ is isomorphic to the quotient $\operatorname{Nor}_{G_\F}(A_\F) / \operatorname{Cen}_{G_\F}(A_\F)$. Let $\pi \colon \operatorname{Nor}_{G_\F}(A_\F) \to W_s$ be the induced map to the finite group $W_s$. Let $Y_0 \subseteq \Omega$ be a finite subset, so that $|\pi(N_{g,Y_0})|$ is minimal. For any $p \in \Omega$,
	 $$
	 N_{g,Y_0 \cup \{p\}} = N_{g,Y_0} \cap N_{g,\{p\}} \subseteq N_{g,Y_0}
	 $$
	 and thus by minimality $\pi(N_{g,Y_0 \cup \{p\}}) = \pi(N_{g,Y_0})$. We claim that in fact $N_{g,Y_0 \cup \{p\}} = N_{g,Y_0}$ for every $p \in \Omega$.
	 Pick $n_0\in N_{g,Y_0}$. 
	 We decompose $w := [n_0] =  (t_a,\pi(n_0)) \in W_a = \A \ltimes W_s$, so that $t_a$ is translation given by multiplication of some $a\in A_\F$. For any $p \in \Omega$, there exists $n' \in N_{g,Y_0 \cup \{p\}}$ such that $\pi(n_0) = \pi(n')$ and we decompose similarly $w' := [n'] = (t_{a'},\pi(n')) \in W_a$. Acting on some $q\in Y_0$, we have
	 $$
	 a.(\pi(n_0)(q)) = n_0.q = g.q = n'.q = a'.(\pi(n')(q)) = a'.(\pi(n_0)(q)),
	 $$
	 which implies $t_a = t_{a'}$ and thus $w= (t_a,\pi(n_0)) = (t_{a}', \pi(n')) = w'$. For any $p \in \Omega$ we thus have
	 $$
	 n_0.p = w(p) = w'(p) =  n'.p = g.p,
	 $$
	 so $n_0 \in N_{g, Y_0\cup\{p\}}$ as required. Finally
	 $$
	 N_{g,\Omega}
	 = \bigcap_{p \in \Omega} N_{g, Y_0\cup \{p\}} = N_{g,Y_0} \neq \emptyset,
	 $$
	 so taking any $n \in N_{g,Y_0}$ provides the required element $w = [n] \in W_a$ with $w(p) = g.p$ for all $p \in \Omega$.
\end{proof}
From the above proof, we can also extract the following Lemma, by taking $g=\operatorname{Id}$ and noting that then $N_\Omega =  N_{g,\Omega} =  N_{g,Y_0} = N_{Y_0}$.
\begin{lemma}\label{lem:BTNY0=NOmega}
	Let $\Omega \subseteq \A$. There exists a finite subset $Y_0 \subseteq \Omega$ such that $N_{Y_{0}} = N_{\Omega}$.
\end{lemma}
To be able to prove the $W_a$-convexity part of axiom (A2), we want to describe the stabilizer of a (possibly infinite) subset $\Omega$ as in Theorem \ref{thm:BTstab_fin}. For this we first need a Lemma.

\begin{lemma}\label{lem:U+U-fix_U+fix}
	Let $u^+ \in U^+$, $u^- \in U^-$ and $p\in \A$. If $u^+.p = u^-.p$, then $u^+.p = p = u^-.p$.
\end{lemma}
\begin{proof}
	We start by noting that if $a \in \F^{n\times n}$ is an upper triangular matrix and $b \in \F^{n\times n}$ is a lower triangular matrix, both with ones on the diagonal such that $ab \in O^{n\times n}$, then $a, b \in O^{n\times n}$, as can be checked by matrix calculations. 
	
	We will now prove the Lemma in the case that $p = o \in \B$. We know that $u^+.p = u^-.p$ which is equivalent to $(u^+)^{-1}u^- \in G_\F(O)$. The adjoint map $\operatorname{Ad} \colon G_\F \to \operatorname{GL}(\frakg_\F)$ is a $\K$-morphism that sends elements of $U^+$ to upper triangular matrices and elements of $U^-$ to lower triangular matrices, see Lemma \ref{lem:kan_form}. Moreover, since $\operatorname{Ad}(g)X = gXg^{-1}$ for $X \in \frakg_\F$, we can use Lemma \ref{lem:orthogonal_valuation} on a basis to see that $\operatorname{Ad}(g)$ is defined by polynomials and if $g \in G_\F(O)$, then these polynomials have coefficients in $O$.
	Thus $\operatorname{Ad}(G_\F(O)) \subseteq \operatorname{Ad}(G_\F)(O)$ and $\operatorname{Ad}((u^+)^{-1})\operatorname{Ad}(u^-) \in O^{n\times n}$ and by the preceding remark, $\operatorname{Ad}(u^+) \in O^{n\times n}$ and $\operatorname{Ad}(u^-) \in O^{n\times n}$.
	
	 While $\operatorname{Ad} \colon G_\F \to \operatorname{GL}(\frakg_\F)$ may have a nontrivial finite kernel, its restriction $\operatorname{Ad}|_{U^+} \colon U^+ \to \operatorname{GL}(\frakg_\F) $ is an isomorphism onto its image, since the exponential map gives an isomorphism $U^+ \cong \bigoplus_{\alpha>0} \frakg_\alpha$ and since $\operatorname{ad}$ restricted to the Lie algebra of $U^+$ is an isomorphism to its image defined over $\K$. Then, the inverse map $\operatorname{Ad}(U^+) \to U^+$ is also defined by polynomials with coefficients in $\K$, so $u^+ \in O^{n\times n}$. Similarly $u^- \in O^{n\times n}$. This means that $u^+.o = o$ and $u^-.o = o$ as required.
	
	If now $p = a.o$ for some $a \in A_\F$, then $(u^+)^{-1}u^-.p = p$ if and only if $$a^{-1}(u^+)^{-1}aa^{-1}u^- a \in G_\F(O).$$
	By the above argument $a^{-1}u^+ a, a^{-1}u^-a \in G_\F(O)$ and thus $u^+.p = p$ and $u^-.p = p$.
\end{proof}

We upgrade the description of $\hat{P}_{\Omega}$ as the pointwise stabilizer of a finite subset in Theorem \ref{thm:BTstab_fin}, to arbitrary subsets $\Omega \subseteq \A$. 
\begin{theorem}\label{thm:BTstab}
	The pointwise stabilizer of any subset $\Omega \subseteq \A$ satisfies
	$$
	\operatorname{Stab}_{G_\F}(\Omega) = \hat{P}_\Omega = U_\Omega^+ U_{\Omega}^- N_{\Omega} .
	$$
\end{theorem}
\begin{proof} The inclusions $\supseteq$ are clear. We apply Lemma \ref{lem:BTNY0=NOmega} to obtain a finite subset $Y_0 \subseteq \Omega$ such that $N_{Y_0} = N_{\Omega}$. Then $\operatorname{Stab}(\Omega) \subseteq \operatorname{Stab}(Y_0) = U_{Y_0}^+ U_{Y_0}^- N_{Y_0}$ by Theorem \ref{thm:BTstab_fin}. Let $g\in \operatorname{Stab}(\Omega)$ and $u^+ \in U_{Y_0}^+, u^- \in  U_{Y_0}^-, n \in N_{Y_0} = N_\Omega$ with $g=u^+u^-n$. Then $u^+ u^- = gn^{-1} \in \operatorname{Stab}(\Omega)$. Thus, $(u^+).p = u^-.p$ for all $p\in \Omega$ and by Lemma \ref{lem:U+U-fix_U+fix}, $u^+.p = p = u^-.p$, in particular $u^+ \in U_\Omega^+$ and $u^- \in U_{\Omega}^-$. We now know
	$$
	\operatorname{Stab}_{G_\F}(\Omega) \subseteq U_\Omega^+ U_\Omega^- N_\Omega \subseteq \hat{P}_\Omega \subseteq \operatorname{Stab}_{G_\F}(\Omega)
	$$
	concluding the proof.
\end{proof}
We now prove a special case of $W_a$-convexity.

\begin{proposition}\label{prop:A2_a}
	Let $g\in G_\F$. Then $\Omega = g^{-1}\A \cap \A$ is a finite intersection of affine half-apartments. 
\end{proposition}
\begin{proof}
	By Proposition \ref{prop:A2_b}, we know that there is an $n\in \operatorname{Nor}_{G_\F}(A_\F)$ such that $g^{-1}n \in \operatorname{Stab}(\Omega)$ pointwise and by Theorem \ref{thm:BTstab} $g^{-1}n \in U_{\Omega}^+ U_{\Omega}^- N_{\Omega}$. So let $u^+ \in  U_{\Omega}^+, u^- \in  U_{\Omega}^- $ and $n' \in  N_{\Omega}$ with $g^{-1} n = u^+ u^- n'$. We note that for $\tilde{n} := n(n')^{-1}$ we have $g^{-1}\tilde{n} = u^+ u^- $.
	
	Recall that the affine half-space given by $\alpha \in \Sigma$ and $k\in \Lambda$ is
	$$
	H_{\alpha,k}^+ = \left\{ a.o \in \A \colon (-v)(\chi_\alpha(a)) \geq k  \right\}.
	$$
	In Proposition \ref{prop:Uconv}, we showed that for $u^+$ there are $k_\alpha \in \Lambda \cup \{-\infty\}$ for every $\alpha >0$ such that 
	$$
	\left\{ p\in \A \colon u^+.p \in \A \right\} = \bigcap_{\alpha >0} H_{\alpha, k_\alpha}^+.
	$$	
	where we take the convention that $H_{\alpha, -\infty}^+ = \A$. By changing the order on $\Sigma$, we similarly obtain for $u^-$ some $k_\alpha\in \Lambda \cup \{-\infty\}$ for $\alpha <0$ such that
	$$
	\left\{ p\in \A \colon u^-.p \in \A \right\} = \bigcap_{\alpha <0} H_{\alpha, k_\alpha}^+.
	$$
	We now show that
	$$
	\Omega = \bigcap_{\alpha \in \Sigma} H_{\alpha, k_\alpha}^+.
	$$
	By Lemma \ref{lem:U+U-fix_U+fix}, for any $p\in \Omega$ we have $u^+.p =p$  and $u^{-}.p=p$, hence $p \in \bigcap H_{\alpha, k_\alpha}^+$. If on the other hand $p \in \bigcap H_{\alpha, k_\alpha}^+$, then $u^+u^-.p=u^+.p=p$, in particular $p \in u^+u^- \A \cap \A = g^{-1}\tilde{n}\A \cap \A = g^{-1} \A \cap \A = \Omega$. This concludes the proof that $\Omega$ is a finite intersection of half-spaces.
\end{proof}

We now put together Propositions \ref{prop:A2_b} and \ref{prop:A2_a} to prove axiom (A2).
\begin{theorem}\label{thm:A2}
	Axiom (A2) holds: 
	\begin{enumerate}
	\item [(A2)] For every $f,f' \in \Fun$, the set $B := f^{-1}(f(\A) \cap f'(\A)) \subseteq \A$ is a finite intersection of affine half-apartments and there is $w\in W_a$ such that $f|_B = f'\circ w |_{B}$. 
	\end{enumerate}
\end{theorem}
\begin{proof}
	By definition of $\Fun$, there are $h,h' \in G_\F$ such that $f=h.f_0$ and $f' = h'.f_0$, where $f_0 \colon \A \to \B$ is the inclusion. Then for $g:= (h')^{-1}h$ we have $B = g^{-1}\A \cap \A$ since
	$$
	f(g^{-1}\A \cap \A) = h.(g^{-1}\A \cap \A) = h'.\A \cap h.\A = f'(\A) \cap f(\A) = f(B).
	$$ 
	By Proposition \ref{prop:A2_a}, $B$ is a finite intersection of affine half-apartments. By Proposition \ref{prop:A2_b}, there is $w\in W_a$ such that $g.p = w(p)$ for all $p\in B$, so that for all $p \in B$
	$$
	f(p) = h.p = h'g.p = h'.w(p) = h'.f_0(w(p)) = (f' \circ w)(p),
	$$
	concluding the proof.
\end{proof}

\subsection{Axiom (A4)}\label{sec:A4}

In this section the root system $\Sigma$ does not need to be reduced, except for being able to apply axiom (A2) in Lemma \ref{lem:subsector_nice}. Axiom (A4) is a statement about sectors. Let $s_0 = f_0(C_0) \subseteq \B$ be the fundamental sector corresponding to the fundamental Weyl chamber $C_0 \subseteq \A$. All sectors are of the form $g.s_0$ for some $g\in G_\F$. If a sector $s$ is a subset of another sector $s'$, then $s$ is called a \emph{subsector} of $s'$. 

From (A2) we get that subsectors of $s_0$ are of the form $a.s_0$ for $a\in A_\F$.
\begin{lemma}\label{lem:subsector_nice}
	For every subsector $s \subseteq s_0$ there exists $a\in A_\F$ such that $s=a.s_0$.
\end{lemma}
\begin{proof}
	A general subsector of $s_0$ is of the form $g.s_0$ for some $g\in G_\F$. By Axiom (A2) we know that $s_0 \to g.s_0$ is realized by an element $w\in W_a$ of the affine Weyl group, $g.b.o = w(b.o)$ for all $b.o \in s_0$. We decompose $w = ( t_a,w_s) \in \A \rtimes W_s$ for some $a\in A_\F$. Since $o\in s_0$ and $a.o = w(o) \in g.s_0 \subseteq s_0$, we know that $\chi_\alpha(a) \geq 1$ for all $\alpha >0$. We also know that $w_s(s_0)$ is one of the finitely many sectors based at $o$. If $w_s(s_0) \neq s_0$, then there is some $\alpha >0$ with $\chi_\alpha(b) \leq 1$ for all $b.o \in w_s(s_0)$. Since $w_s(s_0)$ is a cone with open interior, there are $b.o\in w_s(s_0)$ with arbitrary negative $\chi_\alpha(b)$, in particular there is some $b.o\in w_s(s_0)$ with $(-v)(\chi_\alpha(b)) < (-v)(\chi_\alpha(a)^{-1})$, so that $(-v)(\chi_\alpha(ab)) < 0$. But this contradicts $a.w_s(s_0) \subseteq s_0$, since $a.b.o \notin s_0$. We conclude that $w_s(s_0) = s_0$ and thus $g.s_0 = a.s_0$.
\end{proof}

 While studying the model apartment $\A$, Bennett \cite[Prop 2.9]{Ben1} proved the following lemma.

\begin{lemma}
	For all $p \in \A$ there is a $q\in \A$ such that $q+C_0 \subseteq (p+C_0) \cap C_0$.
\end{lemma}
In our setting this lemma translates in terms of $A_\F.o = \A\subseteq \B$.
\begin{lemma} \label{lem:A4:a_sub}
	For all $a\in A_\F$ there is a $b\in A_\F$ such that $b.s_0 \subseteq a.s_0 \cap s_0$.
\end{lemma}

\begin{figure}[h]
	\centering
	\includegraphics[width=0.5\linewidth]{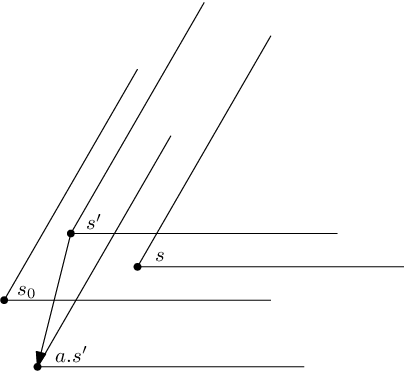}
	\caption{Lemma \ref{lem:A4:a_sub_cor} states that for every subsector $s' \subseteq s_0$ and $a \in A_\F$ there is a sector $s$ contained in both $s'$ and $a.s'$.   }
	\label{fig:subsectors_lemma}
\end{figure}

We also get a slightly more general statement, illustrated in Figure \ref{fig:subsectors_lemma}.
\begin{lemma}\label{lem:A4:a_sub_cor}
	For all subsectors $s' \subseteq s_0$ and for all $a\in A_\F$ there is a subsector $s \subseteq a.s' \cap s_0$.
\end{lemma}
\begin{proof}
	Let $s'=\tilde{a}.s_0$ for some $\tilde{a} \in A_\F$, see Lemma \ref{lem:subsector_nice}. We apply Lemma \ref{lem:A4:a_sub} with $a\tilde{a}$ to get $b\in A_\F$ with $s:= b.s_0 \subseteq a\tilde{a}.s_0 \cap s_0 = a.s' \cap s_0$.
\end{proof}

We now use Proposition \ref{prop:exists_a_NO} to show that while elements of $U_\F$ may not fix $s_0$ itself, they at least fix a subsector.
\begin{lemma}\label{lem:N_fixes_subsector}
	For every $u\in U_\F$ there is a subsector $s\subseteq s_0$ with $u.p=p$ for all $p\in s$.	
\end{lemma} 
\begin{proof}
	Let $a\in A_\F$ with $ a^{-1}ua \in U_\F(O)$ as in Proposition \ref{prop:exists_a_NO}. For all $b.o \in s_0$ we have $a^{-1}ua.b.o = b.o$ by Corollary \ref{cor:NO_fixes_s0}.
	Therefore $u$ fixes $a.s_0$ pointwise. We don't know whether $a.s_0$ is a subsector of $s_0$, but we can apply Lemma \ref{lem:A4:a_sub} to find a subsector $s$ of $s_0$, which is also a subsector of $a.s_0$ and therefore is fixed pointwise by $u$. 
\end{proof}
Now that we understand the action of $U_\F$ better, let us turn to the group $B_\F=M_\F A_\F U_\F$ introduced in Section \ref{sec:BWB}, where $M_\F=\operatorname{Cen}_{K_\F}(A_\F) ,
A_\F$ and $U_\F$ are as before.
\begin{lemma}\label{lem:A4:B_acts}
	For all $b\in B_\F = M_\F A_\F U_\F$ there is a subsector $s\subseteq s_0$ with $b.s \subseteq s_0$.
\end{lemma}
\begin{proof}
	We first recall the elements of $M_\F$ fix all of $\A$ pointwise. Let $m\in M_\F$, $u \in U_\F$ and $a\in A_\F$ such that $b=mua$. Using Lemma \ref{lem:N_fixes_subsector}, we find a subsector $s'\subseteq s_0$ which is fixed by $u$ pointwise. Applying Lemma \ref{lem:A4:a_sub_cor} to $s'$ and $a^{-1}$ we get a subsector $s\subseteq s_0 \cap a^{-1}.s'$. We now have
	$$
	b.s \subseteq b.a^{-1}.s' = m u .s' = m .s' = s' \subseteq s_0,
	$$
	as claimed.
\end{proof}
We are now ready to prove axiom (A4) using the Bruhat decomposition $G_\F=B_\F W_sB_\F$, Theorem \ref{thm:BWB}.
\begin{theorem}\label{thm:A4}
	Axiom 
	\begin{itemize}
 \item [(A4)] For any sectors $s_1,s_2 \subseteq \B$ there are subsectors $s_1' \subseteq s_1, s_2' \subseteq s_2$ such that there is an $f\in \Fun$ with $s_1', s_2' \subseteq f(\A)$.
	\end{itemize}
	holds for $(\B,\Fun)$. 
\end{theorem}

\begin{proof}
	The action of $G_\F$ on the sectors is transitive by definition and we may hence assume without loss of generality that one of the sectors in (A4) is $s_0$ and the other is given by $g.s_0$ for some $g \in G_\F$. We have to prove 
	\begin{itemize}
		\item [(A4)'] For all $g\in G_\F$, there are subsectors $s \subseteq s_0, s' \subseteq g.s_0$ such that there is an $f\in \Fun$ with $s, s' \subseteq f(\A)$.
	\end{itemize}
	\begin{figure}[h]
		\centering
		\includegraphics[width=0.8\linewidth]{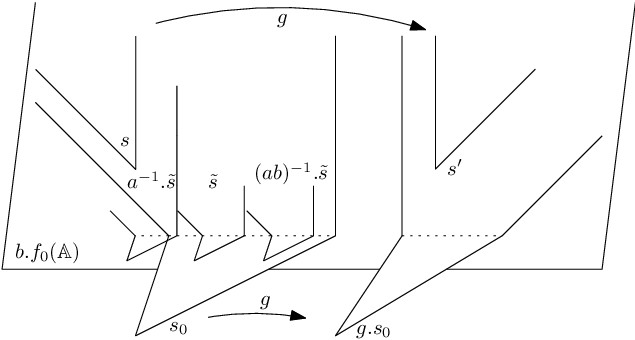}
		\caption{Axiom (A4) states that while the sectors $s_0,g.s_0$ may not lie in a common flat, they contain subsectors $s,s'$ contained in a common flat $f(\A)$. }
		\label{fig:A4proof}
	\end{figure}
	The situation is illustrated in Figure \ref{fig:A4proof}. The Bruhat-decomposition, Theorem \ref{thm:BWB}, states that $G_\F$ is a disjoint union 
	$$
	G_\F = \bigcup_{w=[k] \in W_s} B_\F k B_\F
	$$
	of double cosets. 
	Let $g \in G_\F$ and $b,b' \in B_\F, k \in N_\F := \operatorname{Nor}_{K_\F}(A_\F)$ with $g=bkb'$. We further decompose $b'=ma'ua$ for $m\in M_\F, u\in U_\F(O)$ and $a,a'\in A_\F$ using $B_\F=M_\F U_\F A_\F$ and Proposition \ref{prop:exists_a_NO}. Since $(ab)^{-1} \in B_\F$, apply Lemma \ref{lem:A4:B_acts} to obtain a subsector $\tilde{s} \subseteq s_0$ with $(ab)^{-1}.\tilde{s} \subset s_0$. We apply Lemma \ref{lem:A4:a_sub_cor} to $\tilde{s}$ and $a^{-1}\in A_\F$ to obtain a subsector $s \subset a^{-1}.\tilde{s} \cap s_0$. We claim that $s$ and $s':=g.s$ are the required subsectors in (A4)' and that $f= b.f_0$ defines the common apartment.
	
	We have by construction that $s\subseteq s_0$ and clearly $s'=g.s \subseteq g.s_0$. It remains to show that $s$ and $s'$ are in the apartment $b.f_0(\A)$. We see
	$$
	b^{-1}.s \subseteq b^{-1}.a^{-1}.\tilde{s} \subseteq s_0 \subseteq f_0(\A)
	$$
	hence $s \subseteq b.f_0(\A)$ and
	\begin{align*}
		s' &= g.s = bkma'ua.s \subseteq bkma'u.\tilde{s} = bkma'.\tilde{s} \\
		&\subseteq bkm.f_0(\A) = bkm.f_0(\A) = bk.f_0(\A) = b.f_0(\A),
	\end{align*}
	where we used Corollary \ref{cor:NO_fixes_s0} for $u.\tilde{s} = \tilde{s}$ since $\tilde{s} \subseteq s_0$ and $u\in U_\F(O)$ . 
\end{proof}

\subsection{Axiom (EC)}\label{sec:EC}
\begin{figure}[h]
 	\centering
 	\includegraphics[width=0.5\linewidth]{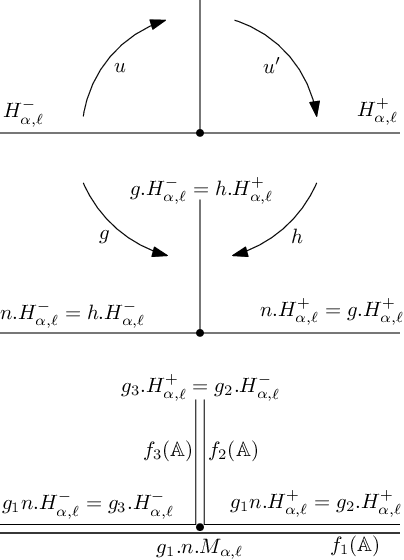}
 	\caption{Axiom (EC) states that if $f_1(\A) \cap f_2(\A)$ is a half-apartment, then there exists $f_3 \in \F$ situated as illustrated. Proposition \ref{prop:EC_alpha} deals with the case where the half-apartment is contained in $\A$, illustrated in the first two parts. Theorem \ref{thm:EC} then tackles the general case of axiom (EC). }
 	\label{fig:EC}
 \end{figure}

In this section, the Jacobson-Morozov-construction as well as axiom (A2) is used, so we require $\Sigma$ to be reduced. We recall that a \emph{half-apartment in the model apartment} is a set of the form
$$
H_{\alpha,\ell}^+ = \left\{ a.o \in \A \colon (-v)(\chi_\alpha(a)) \geq \ell \right\}
$$
for some $\alpha \in \Sigma$ and some $\ell \in \Lambda$. A \emph{half-apartment in the building} is $g.H_\alpha^+$ for any $g \in G_\F$. We first use axiom (A2) to make sure that if a half-apartment of the building is included in $\A$, then it is an affine half-apartment in the model apartment.
\begin{lemma}\label{lem:subhalfapartment_nice}
	Let $g\in G_\F, \alpha \in \Sigma, \ell \in \Lambda$ such that $g.H_{\alpha,\ell}^+ \subseteq \A$. Then there is $\alpha \in \Sigma, \ell' \in \Lambda$ such that 
	$$
	g.H_\alpha^+ = H_{\alpha', \ell'}^+ .
	$$
\end{lemma}
\begin{proof}
	By axiom (A2), Theorem \ref{thm:A2}, there is $w \in W_a$ such that $g.H_\alpha^+ = w(H_\alpha^+)$. We decompose $w = (t_a, w_s) \in \A \rtimes W_s = W_a$ and define $\alpha' := w_s(\alpha)$. Then $g.H_{\alpha,\ell}^+ = w(H_{\alpha,\ell}^+) = a.H_{\alpha',\ell}^+) = H_{\alpha', \ell'}^+$ for $\ell' := (-v)(\chi_{\alpha'}(a)) + \ell'$. 
\end{proof}
We start by proving axiom
\begin{enumerate}
	\item [(EC)] For $f_1,f_2 \in \Fun$, if $f_1(\A)\cap f_2(\A)$ is a half-apartment, then there exists $f_3 \in \Fun$ such that $f_i(\A)\cap f_3(\A)$ are half-apartments for $i\in \{1,2\}$. Moreover $f_3(\A)$ is the symmetric difference of $f_1(\A)$ and $f_2(\A)$ together with the boundary wall of $f_1(\A) \cap f_2(\A)$.
\end{enumerate}
in the special case where 
$$
f_1(\A)\cap f_2(\A) = H_{\alpha,\ell}^+ :=  \left\{ a.o \in \A \colon (-v)(\chi_\alpha(a)) \geq \ell \right\}
$$
for some $\alpha \in \Sigma$ and $\ell \in \Lambda$, before deducing the full statement in Theorem \ref{thm:EC}. The situation is illustrated in Figure \ref{fig:EC}.
 \begin{proposition}\label{prop:EC_alpha}
	Let $g \in G_\F$ such that $g^{-1}.\A \cap \A = H_{\alpha, \ell}^+$ for some $\alpha \in \Sigma, \ell \in \Lambda$. Then there exists an $h\in G_\F$ such that $h^{-1}.\A \cap \A = H_{\alpha,\ell}^-$ and $ h.H_{\alpha,\ell}^+ = g.H_{\alpha,\ell}^-$. Moreover there is some $n\in \operatorname{Nor}_{G_\F}(A_\F)$ such that pointwise $g.H_{\alpha,\ell}^+ = n.H_{\alpha,\ell}^+$ and $h.H_{\alpha,\ell}^- = n.H_{\alpha,\ell}^-$.
\end{proposition}
\begin{proof}
	Since $g.H_{\alpha,\ell}^+ = \A \cap g.\A \subseteq \A$, we can apply axiom (A2), Theorem \ref{thm:A2}, to obtain $n \in \operatorname{Nor}_{G_\F}(A_\F)$ such that $g.H_{\alpha,\ell}^+ = n.H_{\alpha,\ell}^+$ pointwise. Hence $n^{-1}g$ fixes $H_{\alpha,\ell}^+$ pointwise and by Corollary \ref{cor:NalphaO_fixes_H_affine}, $n^{-1}g \in U_{\alpha, \ell} A_\F(O) M_\F$. So let $u \in U_{\alpha, \ell}$ such that $n^{-1}g.p = u.p$ for all $p \in \A$. By Lemma \ref{lem:BTm_uuu}, the Jacobson-Morozov-homomorphism can be used to define
	\begin{align*}
		u' = \varphi_\F \begin{pmatrix}
			1 & 0 \\ - 1/t & 0
		\end{pmatrix} \in (U_{-\alpha})_\F
	\end{align*}
	with $\varphi_{-\alpha} (u') = -\varphi_\alpha(u) = -\ell$ and $m(u) = u'uu' \in \operatorname{Nor}_{G_\F}(A_\F)$ which acts as the reflection along $M_{\alpha, \ell} = \{ a.o \in \A \colon (-v)(\chi_\alpha(a)) = \ell \}$ by Proposition \ref{prop:BTmu_in_Wa}. For $h := n(u')^{-1} \in G_\F$ we then have
	$$
	h.H_{\alpha,\ell}^+ = n (u')^{-1} m(u) .H_{\alpha,\ell}^- = n u u' .H_{\alpha,\ell}^- = nu.H_{\alpha,\ell}^- = g.H_{\alpha,\ell}^-
	$$
	with $(\A \cap h.H_{\alpha,\ell}^+) = g.M_{\alpha,\ell} = n.M_{\alpha,\ell} = h.M_{\alpha,\ell} \subseteq h.H_{\alpha,\ell}^-$ and therefore 
	$$
	h.(h^{-1}.\A\cap \A) = \A \cap h.\A = (\A \cap h.H_{\alpha,\ell}^+ )\cup(\A \cap h.H_{\alpha,\ell}^-) = \A \cap h.H_{\alpha,\ell}^- = n.H_{\alpha,\ell}^-
	$$
    which implies $h^{-1}.\A \cap \A = h^{-1}n.H_{\alpha,\ell}^- = u'.H_{\alpha,\ell}^- = H_{\alpha,\ell}^- $ as required.
\end{proof}

\begin{theorem}\label{thm:EC}
	Axiom 
	\begin{enumerate}
		\item [(EC)] For $f_1,f_2 \in \Fun$, if $f_1(\A)\cap f_2(\A)$ is a half-apartment, then there exists $f_3 \in \Fun$ such that $f_i(\A)\cap f_3(\A)$ are half-apartments for $i\in \{1,2\}$. Moreover $f_3(\A) = f_1(\A)\cap f_3(\A) \cup f_2(\A)\cap f_3(\A)$ and $\partial (f_1(\A)\cap f_2(\A) ) = f_1(\A)\cap f_3(\A) \cap f_2(\A)\cap f_3(\A)$. 
	\end{enumerate}
	holds.
\end{theorem}
\begin{proof}
	Let $g_1, g_2 \in G_\F$ such that $f_1 = g_1.f_0$ and $f_2 = g_2.f_0$. For $g := g_1^{-1} g_2$, we have $f_2 = g_1 g.f_0$. Since $f_1(\A)\cap f_2(\A)$ is a half-apartment, so is
	$$
	(g_2)^{-1}(f_1(\A)\cap f_2(\A)) = g_2^{-1}(g_1.\A \cap g_2.\A) = g^{-1}.\A \cap \A \subseteq \A
	$$
	and by Lemma \ref{lem:subhalfapartment_nice}, there is some $\alpha \in \Sigma$ and $\ell \in \Lambda$ such that $g^{-1}.\A \cap \A = H_{\alpha, \ell}^+$. Now we use Proposition \ref{prop:EC_alpha} to obtain $h\in G_\F$ and $n\in\operatorname{Nor}_{G_\F}(A_\F)$ with $g.H_{\alpha,\ell}^+ = n.H_{\alpha,\ell}^+$, $h.H_{\alpha,\ell}^- = n.H_{\alpha,\ell}^-$ and $h.H_{\alpha,\ell}^+ = g.H_{\alpha,\ell}^-$. We define $g_3 := g_1h $ and $f_3 = g_3.f_0$.
	We have as in Figure \ref{fig:EC}.
	\begin{align*}
		g_3.H_{\alpha,\ell}^- &= g_1h.H_{\alpha,\ell}^- = g_1.n.H_{\alpha,\ell}^- \\
		g_2.H_{\alpha,\ell}^+ &= g_1g.H_{\alpha,\ell}^+ = g_1.n.H_{\alpha,\ell}^+ \\
		g_3.H_{\alpha,\ell}^+ &= g_1h.H_{\alpha,\ell}^+ = g_1g.H_{\alpha,\ell}^- = g_2.H_{\alpha,\ell}^- 
	\end{align*}
	so $f_1(\A)\cap f_3(\A) = f_3(H_{\alpha,\ell}^-)$ and $f_2(\A)\cap f_3(\A) = f_3(H_{\alpha,\ell}^+)$ are half-apartments and
	$$
	f_3(\A) = f_3(H_{\alpha,\ell}^+) \cup f_3(H_{\alpha,\ell}^-) = f_1(\A)\cap f_3(\A) \cup f_2(\A)\cap f_3(\A).
	$$
	The wall of the half-apartment $f_1(\A)\cap f_2(\A) = g_2.H_{\alpha,\ell}^+ $ is given by 
	\begin{align*}
		\partial (f_1(\A)\cap f_2(\A) ) &= g_2.M_{\alpha,\ell} = g_1g.M_{\alpha,\ell} = g_1.n.M_{\alpha,\ell} =
		g_1n.H_{\alpha,\ell}^+ \cap g_1n.H_{\alpha,\ell}^- \\
		&= f_1(\A)\cap f_2(\A) \cap f_1(\A)\cap f_3(\A).
	\end{align*}
\end{proof} 

\noindent
This concludes the proof of the last remaining axiom and the proof of Theorem \ref{thm:B_is_building}.

\subsection{Beyond reduced root systems}

The main theorem of this thesis relies on the assumption that the root system $\Sigma$ is reduced.
\begin{reptheorem}{thm:B_is_building}
	If the root system $\Sigma$ is reduced, then the pair $(\B,\Fun)$ is an affine $\Lambda$-building.
\end{reptheorem}
We remark that $\B$ is defined for any self-adjoint semisimple linear algebraic $\K$-group, independent of its root system. We expect Theorem \ref{thm:B_is_building} to still hold without the assumption on $\Sigma$. 
\begin{question}
	Is the pair $(\B,\Fun)$ an affine $\Lambda$-building even when $\Sigma$ is not reduced?
\end{question}
We outline here how our proof relies on the assumption, how the assumption cannot be removed using our strategy and a possible alternative proof strategy that might be of use to eliminate the assumption. 

In our proof, the assumption first comes up in the Jacobson-Morozov Lemma, both in the setting of Lie algebras, see Lemma \ref{lem:JM_basic} and in the semialgebraic setting, see Proposition \ref{prop:Jacobson_Morozov_real_closed}. Explicit calculations in $\mathfrak{su}_{1,2}$ show that Lemma \ref{lem:JM_basic} does not hold and similarly, Proposition \ref{prop:Jacobson_Morozov_real_closed} does not hold for $\operatorname{SU}(1,2)$. For given $\alpha \in \Sigma$ and $X\in \frakg_\alpha \oplus \frakg_{2\alpha}$, the task of the Jacobson-Morozov Lemma is to find $Y \in \frakg_{-\alpha}\oplus \frakg_{-2\alpha}$ and $H\in \fraka \subseteq \frakg_0$ such that $(X,Y,H)$ is an $\mathfrak{sl}_2$-triplet. While $Y$ can be found as a multiple of $\theta(X)$, there is no guarantee that $H :=[X,Y] \in \fraka$ when $\frakg_{2\alpha}\neq 0$.

In the proof of Theorem \ref{thm:B_is_building}, the Jacobson-Morozov Lemma is used to associate to each $u \in (U_\alpha)_\F$ an element $m(u) \in \operatorname{Nor}_{G_\F}(A_\F)$ representing a reflection as in Proposition \ref{prop:BTmu_in_Wa}. We expect the following construction in the rank one subgroup $L_{\pm \alpha}$ to give a definition of $m(u)$ also when $\Sigma$ is not reduced. We use the Bruhat-decomposition of $L_{\pm\alpha}$.

\begin{repcorollary}{cor:levi_Bruhat}
	Let $(B_{\alpha})_\F := (M_{\pm \alpha})_\F (A_{\pm \alpha})_\F (U_\alpha)_\F$. Then there is a representative $m \in (N_{\pm \alpha})_\F$ of the unique non-trivial element in $W_{\pm \alpha}$, so that
	$$
	(L_{\pm \alpha})_\F = (B_{\alpha})_\F \ \amalg \   (B_{\alpha})_\F \cdot m \cdot (B_{ \alpha})_\F.
	$$
\end{repcorollary}
If $u \in (U_\alpha)_\F$, let $u \in (B_{-\alpha})_\F \ \amalg \   (B_{-\alpha})_\F \cdot m \cdot (B_{ -\alpha})_\F$ and since $(B_{\alpha})_\F \cap (B_{-\alpha})_\F = \{\operatorname{Id}\}$, we can find $b,b' \in (B_{-\alpha})_\F$ such that $u = bmb'$. Writing $b = \bar{u}\bar{a}\bar{m}$ and $b' = \bar{m}'\bar{a}'\bar{u}'$ with $\bar{m}, \bar{m}' \in (M_{\pm \alpha})_\F, \bar{a},\bar{a}' \in (A_{\pm \alpha})_\F$ and $\bar{u}, \bar{u}'\in (U_{-\alpha})_\F$, we define
$$
m(u):= \bar{a} \bar{m} m \bar{m}' \bar{a}' = \bar{u}^{-1} u (\bar{u}')^{-1}  \in \operatorname{Nor}_{(L_{\pm\alpha})_\F }((A_{\pm \alpha})_\F ) \cap (U_{-\alpha})_\F (U_{\alpha})_\F (U_{-\alpha})_\F.
$$
We note that as a consequence of Lemma \ref{lem:levi_fixes_A}, $\operatorname{Nor}_{(L_{\pm\alpha})_\F }((A_{\pm \alpha})_\F ) \subseteq \operatorname{Nor}_{G_\F}(A_\F)$. In the case where $\Sigma$ is reduced, the uniqueness statement of Lemma \ref{lem:BTUUUm} then shows that $m(u)$ defined here agrees with $m(u)$ defined via the Jacobson-Morozov Lemma. This suggests that the definition of $m(u)$ as outlined here should be used when $\Sigma$ is not reduced. However, our proof also relies on the explicit description of $m(u)$ as $m(u) = u'uu''$ with $\varphi_{-\alpha}(u') = \varphi_{-\alpha}(u'') = -\varphi_\alpha(u)$, see the second part of Lemma \ref{lem:BTUUUm}, which follows from the explicit description of the root group valuation in Lemma \ref{lem:BTphiu_is_vt} relying heavily on the Jacobson-Morozov description. One way to allow a similar level of understanding of the root group valuations in the case of non-reduced root systems could be to do a case by case analysis of all rank one groups.

\newpage
	
	\begin{appendices}

\section{Appendix: The building for $\operatorname{SL}(n,\F)$}\label{sec:appendixSLn}

To obtain Theorem \ref{thm:B_is_building}, that $\B$ is an affine $\Lambda$-building, a relatively large amount of effort goes into proving axiom (A2). The development of the theory following \cite{BrTi} in Subsections \ref{sec:BT_root_groups}, \ref{sec:BT_rank_1} and \ref{sec:BT_higher_rank} is not needed when the group $G$ is well understood. In this appendix, we give an alternative proof of axiom (A2) in the case where $G_\F = \operatorname{SL}_n(\F)$. The proof still relies on the general theory developed for semialgebraic groups in Sections \ref{sec:split_tori} and \ref{sec:decompositions}, but is significantly shorter.

Let $\F$ be a non-Archimedean real closed field with order compatible valuation $(-v) \colon \F \to \Lambda \cup \{ -\infty\}$ and valuation ring $O = \{a \in \F \colon (-v)(a) \leq 0\}$. We consider the semisimple linear algebraic group $G= \operatorname{SL}_n$ with maximal $\K$-split torus
$$
S = \left\{ \begin{pmatrix}
	\star & & \\ & \ddots & \\ && \star
\end{pmatrix} \in \operatorname{SL}_n \right\}. .
$$ 
Then the groups showing up in the various decompositions of Section \ref{sec:decompositions} are given by 
\begin{align*}
	K_\F &= \operatorname{SO}_n(\F)\\
	A_{\F} &= \left\{ a = (a_{ij}) \in S_{\F} \colon a_{ii}>0   \right\} \\
	U_{\F} &= \left\{g=(g_{ij}) \in \operatorname{SL}_{n}(\F) \colon  g_{ii}= 1, \, g_{ij} = 0 \text{ for } i>j  \right\} \\
	N_{\F} &= \left\{ \text{ permutation matrices with entries in }\pm 1 \right\} \\
	M_{\F} & = \left\{ a=(a_{ij}) \in S_{\F} \colon a_{ii} \in \{ \pm 1\}    \right\} \\
	B_{\F} &= \left\{ g = (g_{ij}) \in \operatorname{SL}_n(\F) \colon g_{ij} = 0 \text{ for } i > j \right\}.
\end{align*}
Then $\fraka = \operatorname{Lie}(A_\R) = \{ H \in \R^{n\times n} \colon \operatorname{tr}(H) = 0 \text{ and $H$ is diagonal }  \}$ and the root system $\Sigma$ associated with $\operatorname{SL}(n,\R)$ is given by $\Sigma = \{ \alpha_{ij} \in \fraka^\star \colon i\neq j \in \{1, \ldots , n\} \}$ for the roots
\begin{align*}
	\alpha_{ij} \colon \fraka &\to \R \\
	H& \mapsto H_{ii} - H_{jj}.
\end{align*}
The spherical Weyl group $W_s$ is then isomorphic to the symmetric group $S_n$ on $n$ letters. If we choose the ordered basis $\Delta = \{ \alpha_{12}, \ldots , \alpha_{(n-1)n}\}$, we obtain the positive Weyl chamber
$$
A_\F^+ = \{  \operatorname{Diag}(a_1, \ldots , a_n) \in A_\F \colon a_1 \geq a_2 \geq \ldots  \geq a_n  \}.
$$
As in Section \ref{sec:building_def} we define the non-standard symmetric space
$$
P_{1}(n,\F) := \left\{ A \in \F^{n\times n} \colon A\tran = A, \ \det(A) = 1 ,\ A \text{ is positive definite } \right\}
$$
and note that $\operatorname{SL}(n,\F)$ acts transitively on $P_1(n,\F)$. The multiplicative norm \begin{align*}
	N_\F \colon \A_\F &\to \F_{\geq 1} \\
	\operatorname{Diag}(a_1, \ldots , a_n) &\mapsto \prod_{i\neq j} \max \left\{ \frac{a_i}{a_j}, \frac{a_j}{a_i}   \right\}
\end{align*}
then gives a $G_\F$-invariant $\Lambda$-pseudo distance $d  = (-v) \circ N_\F \circ \delta_\F$ which allows us to define the $\Lambda$-metric space $\B := P_1(n,\F)/\!\!\sim$. We endow $\B$ with the apartment structure $\Fun = \{ g.f_0 \colon g \in \operatorname{SL}(n,\F) \}$, where $f_0 \colon \A  \to \B$ is the inclusion of the apartment $\A := A_\F.o$ for the basepoint $o = [\operatorname{Id}] \in \B$. The following is a special case of Theorem \ref{thm:B_is_building}.
\begin{theorem}
	The $\Lambda$-metric space $\B$ is an affine $\Lambda$-building of type $\A = \A(\Sigma,\Lambda,\Lambda^n)$.
\end{theorem} 
The axioms (A1), (A3) and (TI) follow as described in Section \ref{sec:axiomsA1A3TI}. We will give a hands on proof of axiom (A2) that only relies on Theorem \ref{thm:stab} about the stabilizer of $o \in \B$. Axiom (A4) can then be proved as in Section \ref{sec:A4}, relying on (A2) and the fact that every $u \in U_\F$ stabilizes a sector, which can be proven as in Section \ref{sec:Wconvexity_for_U}, or looking at matrix entries directly. Finally, axiom (EC) relies on (A2) and the fact that Jacobson-Morozov morphisms can be found, which can be seen for $\operatorname{SL}(2,\F)$ explicitly. 

We first look at the axiom
\begin{enumerate}
	\item [(A2)] For all $f_1,f_2 \in \Fun$, if $f_1(\A)\cap f_2(\A)\neq \emptyset$, then the set $\Omega := f_2^{-1}\circ f_1(\A)$ is a $W_a$-convex set and there exists $w \in W_a$ such that $f_2|_\Omega = f_1\circ w |_\Omega$. 
\end{enumerate}
in the special case where $f_2=f_0$ and $f_1=g.f_0$ for some $g\in G_\F$ with $g.o=o$. The set $\Omega$ is defined by $f_0(\Omega)=f_0(\A)\cap g.f_0(\A)$. Let $B\subset A_\F$ be the set of all elements $a\in A_\F$ with $a.o\in f_0(\Omega)$. Note that elements in $B$ are diagonal matrices whose diagonal entries are positive elements in the real closed field $\F$, we can thus take their $n$-th roots. We would like to prove that $\Omega$ is a finite intersection of closed affine half-apartments. The reason we define $W_a$-convexity this way is that general convex combinations of points in the apartment $\A$ are not possible, since we do not have a vector space-structure on $\A$. Some convex combinations however are still possible and $\Omega$ contains them.  

\begin{lemma}
	\label{lem:A2conv}
	For all $a,a' \in B$ and all $n,m \in \N$, $
	\sqrt[n+m]{a^n a'^m} \in B$.
\end{lemma}
\begin{proof}
	Let $b,b' \in A_\F$ such that $g.a.o = b.o$ and $g.a'.o = b'.o$. By Theorem \ref{thm:stab} $b^{-1}ga, (b')^{-1}g a' \in G_\F(O)$. We can now exploit the explicit structure of $A_\F$ to write in coordinates
	$$
	g_{ij}\frac{a_j}{b_i} \in O \quad \text{and} \quad g_{ij}\frac{a'_j}{b'_i} \in O.
	$$
	Thus also 
	$$
	g_{ij}^{n+m} \frac{a_j^n {a'}_j^m}{b_i^n {b'}_i^m} \in O \quad \text{and so}\quad
	g_{ij} \frac{\sqrt[n+m]{a_j^n {a'}_j^m}}{\sqrt[n+m]{b_i^n {b'}_i^m}} \in O.
	$$
	We note that if $a = \operatorname{Diag}\left(a_1, \ldots , a_n\right) \in A_\F$, then also 
	$$
	\sqrt[n+m]{a} :=\operatorname{Diag}(\sqrt[n+m]{a_1} , \ldots , \sqrt[n+m]{a_n}  ) \in A_\F,
	$$
	since this is a first order statement that is true over the reals. Then
	$$
	g.\sqrt[n+m]{a^n a'^m}.o = \sqrt[n+m]{b^n b'^m}.o \in \A
	$$
completes the proof.
\end{proof}

We now introduce a notion of regularity. Elements $a.o \in f_0(\A)$ with $a\in A_\F$ are represented by diagonal matrices with entries $a_1, \ldots, a_n \in \F$. Consider the amount of distinct entries $(-v)(a_1), \ldots ,(-v)(a_n)$ as elements in $\Lambda = (-v)(\F_{>0})$. 
Let $\pmb{a}\in B$ be a maximal element with respect to the amount of distinct entries. The intuition is that elements with $a_i=a_j$ for $i\neq j$ lie on a wall. The more walls they lie on, the less regular they are. 
Other elements in $B$ may have as many distinct entries as $\pmb{a}$ but the next Lemma states that those are the same ones as the ones for $\pmb{a}$. 

\begin{lemma}
	\label{lem:A2regular}
	Let $b \in B$. If $(-v)(b_i)\neq(-v)(b_j)$ for some $i,j$, then $(-v)(\pmb{a}_i)\neq(-v)(\pmb{a}_j)$. Equivalently if $(-v)(\pmb{a}_i) =(-v)(\pmb{a}_j)$, then $(-v)(b_i) = (-v)(b_j)$.
\end{lemma}
\begin{proof}
	For any element $b \in B$, we define $I_b = \{ (i,j) \colon (-v)(b_i) \neq (-v)(b_j) \}$. The maximality of $\pmb{a}$ means that $|I_b|\leq |I_{\pmb{a}} |$ for all $b \in B$. We want to prove that $I_b \subseteq I_{\pmb{a}}$ for all $b\in B$. We note that there are certainly less than $n^2$ elements in $I_{\pmb{a}}$, and define for $k \in \{0 , 1 , \ldots , n^2  \}$ the element
	$$
	c(k)=\sqrt[n^2]{\pmb{a}^{k} b^{n^2-k}} \in B 
	$$
	which is in $B$ by Lemma \ref{lem:A2conv}. Let $(i,j) \in I_{\pmb{a}} \cup I_b$. Thus $(-v)(\pmb{a_i}/\pmb{a_j}) \neq 0$ and $(-v)(b_i/b_j) \neq 0$. We define the map
	$$
	k \mapsto (-v)( c(k)_i / c(k)_j ) \in \Lambda
	$$
	which is either constant, but not $=0$, or strictly monotonous ascending or descending. Either way, there is at most one value of $k$ for which the map is $0$. Let $k_{ij}$ be this value, if it exists. 
	
	Now pick a $k$ which does not appear as $k_{ij}$ for any $(i,j) \in I_{\pmb{a}} \cup I_b$. This means that $c(k)_i \neq c(k)_j$ for all $(i,j) \in I_{\pmb{a}} \cup I_b$. We conclude that $I_{\pmb{a}} \cup I_b \subseteq I_{c(k)}$. But since $\pmb{a}$ is maximal, we have $I_{\pmb{a}} = I_{c(k)}$ and $I_b \subseteq I_{\pmb{a}}$.

\end{proof}

Recall that $O\subsetneqq \F$ is a strict subring of a non-Archimedean real closed field $\F$. Denote the units of $O$ by $O^\times$. The following is a standard fact about valuation rings.
\begin{lemma}
	In any valuation ring $O$, the set of nonunits $O\setminus O^\times$ is an ideal.
\end{lemma}
\begin{proof}
	Let $x,y\in O\setminus O^\times, r \in O $. Since $x$ is not a unit, $x^{-1}=r\cdot (rx)^{-1} \not\in O$. Therefore $(rx)^{-1}\not\in O$, so $rx \in O \setminus O^\times$. If $x$ and $y$ are nonzero, then $x/y\in O$ or $y/x\in O$, since $O$ is a valuation ring. Writing $x+y=y\cdot(1+x/y)=x\cdot(1+y/x)$, we see that also $x+y\in O\setminus O^\times$.
\end{proof}
The determinant of a matrix is a polynomial of its entries. By the previous Lemma, if all entries of a matrix are in $O\setminus O^\times$, then so is its determinant. We conclude that if the determinant of a matrix with entries in $O$ is $1$, then at least one entry has to be a unit. Actually, even more holds.
\begin{lemma}
	\label{lem:A2unit}
	Let $g$ be a matrix with entries in $O$ and $\operatorname{det}(g)=1$. Then in every row and every column there exists at least one entry that is a unit. In fact, there is a permutation $\sigma$, such that $g_{i\sigma(i)}\in O^\times$.
\end{lemma}
\begin{proof}
	Consider the $i$'th row $(g_{i1}, g_{i2}, \ldots, g_{in})$ of $g$. Using Laplace's formula for the determinant
	$$
	\det (g) = \sum_{j=1}^n (-1)^{i+j} g_{ij} \cdot \det M_{ij},
	$$
	where $M_{ij}$ are the minor matrices after deleting row $i$ and column $j$, we see that since $\det (g)=1 \in O^\times$, at least one of the $g_{ij}$ for $j=1, \ldots , n$ has to be a unit. Equivalently there has to be a unit in every column. We can extract a permutation by alternating the row and column argument. 
\end{proof}

We will use Lemma \ref{lem:A2regular} and Lemma \ref{lem:A2unit} to prove a certain rigidity of the action of $G_\F$ on $f_0(\A)$, namely when an isometry $g\in G_\F$ fixes the identity, then its action on $f_0(\A)$ can be described by an element in the spherical Weyl group $W_s$. The apartment $\A$ itself could have many other symmetries, but the isometry can only be extended to all of $\B$, when it is one of the finitely many elements in $W_s$. Consider the fundamental sector
$$
s_0 = \{ a.o \in \A \colon a_1 \geq \ldots \geq a_n \}
$$
based at $o$. The sectors in $f_0(\A)$ based at $o$ are fundamental domains of the action of $W_s$, see for instance \cite[Chapter 10.3]{Hum2}. So two points that lie inside a common sector based at $o$ can only be related by an element in $W_s$ if they are the same. This implies a first restricted version of axiom (A2). 

\begin{proposition}
	\label{pro:A2id}
	Let as before $\pmb{a}\in \B$ be maximal with respect to the amount of distinct entries and $g \in G_\F$ such that $g.o=o$. If both $\pmb{a}.o$ and $g.\pmb{a}.o$ lie in the same sector based at $o$, then $\pmb{a}.o=g.\pmb{a}.o$ and moreover the action of $g$ is the identity on the set $\Omega := f_0^{-1}\left( f_0(\A)\cap g.f_0(\A) \right)$, i.e. $f_0|_\Omega = g.f_0|_\Omega$.
\end{proposition}
\begin{proof} Let $\pmb{a} = \diag{\pmb{a}}$ and $\pmb{b} = \diag{\pmb{b}} \in A_\F$ such that $g.\pmb{a}.o = \pmb{b}.o$. Since $\pmb{a}.o$ and $\pmb{b}.o$ lie in the same sector based at $o$, their entries satisfy the same ordering. Formally we can find a permutation $\rho$, such that 
	\begin{align}\label{eq:ord}
		(-v)(\pmb{a}_{\rho(1)}) & \leq (-v)(\pmb{a}_{\rho(2)}) \leq \ldots \leq (-v)(\pmb{a}_{\rho(n)}) \quad \text{and} \\
		(-v)(\pmb{b}_{\rho(1)}) & \leq (-v)(\pmb{b}_{\rho(2)}) \leq \ldots \leq (-v)(\pmb{b}_{\rho(n)}). \nonumber
	\end{align}
	From $g.o=o$ follows that $g\in G_\F(O)$ by Proposition \ref{thm:stab}. By the previous Lemma \ref{lem:A2unit} we get a permutation $\sigma$, such that $g_{i\sigma(i)}\in O^\times$. In coefficients, $\pmb{b}^{-1}g\pmb{a} \in G_\F(O)$ implies $\pmb{b}_i^{-1}g_{i\sigma(i)}\pmb{a}_{\sigma(i)} \in O$ and therefore also $\pmb{b}_i^{-1}\pmb{a}_{\sigma(i)} \in O$. This means that $\operatorname{Diag}( \pmb{a}_{\sigma(1)}, \ldots , \pmb{a}_{\sigma(n)}).o=\diag{\pmb{b}}.o$. So although the elements $\pmb{a}_{\sigma(i)}$ and $\pmb{b}_{i}$ may not be exactly the same, they satisfy $(-v)(\pmb{a}_{\sigma(i)}) = (-v)(\pmb{b}_i)$. So we can order the entries as
	\begin{alignat}{3}
		\label{eq:ineq}
		(-v)(\pmb{b}_{\rho(1)}) & \leq (-v)(\pmb{b}_{\rho(2)})&& \leq \ldots \leq (-v)(\pmb{b}_{\rho(n)}) \quad \text{and thus} \\
		\rotatebox{90}{$=$}\ \ \ \ & \quad \ \quad \ \rotatebox{90}{$=$} && \quad \quad \quad \quad \quad \ \rotatebox{90}{$=$} \nonumber
		\\
		(-v)(\pmb{a}_{\rho(\sigma(1))}) & \leq (-v)(\pmb{a}_{\rho(\sigma(2))})&& \leq \ldots \leq (-v)(\pmb{a}_{\rho(\sigma(n))}). \nonumber
	\end{alignat}
	But we already have a decreasing ordering of the diagonal entries of $\pmb{a}$ in equation (\ref{eq:ord}), so they have to be the same, i.e. for all $i$, $(-v)(\pmb{a}_{\rho(i)})=(-v)(\pmb{a}_{\rho(\sigma(i))})$ and thus $(-v)(\pmb{a}_i) = (-v)(\pmb{a}_{\sigma(i)})=(-v)(\pmb{b}_i) \in \Lambda$. We now know that $\pmb{a}.o=\pmb{b}.o$, but we still need to show that $g$ is the identity on all of $\Omega$. 
	
	If all the inequalities in (\ref{eq:ineq}) were strict, then $\sigma$ would necessarily be the identity. For indices $i,j$ where there is an equality, $\sigma$ could change these entries. However by Lemma \ref{lem:A2regular} we know that whenever $(-v)(\pmb{a}_i)=(-v)(\pmb{a}_j)$, then also $(-v)(a_i) = (-v)(a_j)$ for any $a=\diag{a} \in B$. So for all $i$, $(-v)(a_{\sigma(i)})=(-v)(a_i)$. In fact we see
	$$
	g_{i\sigma(i)}\frac{a_{\sigma(i)}}{a_i}=g_{i\sigma(i)} \in O
	$$
	and thus $g.a.o=a.o$ for all $a \in B$, so $f_0|_\Omega = g.f_0|_\Omega$.
\end{proof}
We will now allow slightly more general points $\pmb{a}.o$ and $\pmb{b}.o$.
\begin{proposition}
	\label{pro:A2W}
	Let $g\in G_\F$ with $g.o = o \in \B$, then there exists an element $w \in W_s$ such that $f_0|_\Omega = g.f_0 \circ w|_\Omega$, where $\Omega := f_0^{-1}\left( f_0(\A)\cap g.f_0(\A) \right)$.
\end{proposition}
\begin{proof}
	Use Lemma \ref{lem:A2regular} to get $\pmb{a}\in B$ maximal with respect to the amount of distinct entries. Let $\pmb{b}\in A_\F$ such that $\pmb{b}.o=g.\pmb{a}.o \in f_0(\Omega)$. Since the sectors in $f_0(\A)$ are fundamental domains for $W_s$, there is a $w \in W_s$ such that $w.\pmb{b}.o$ and $\pmb{a}.o$ are in the same sector. Equivalently $\pmb{b}.o$ and $w^{-1}.\pmb{a}.o$ are in the same sector. Note that $gw.o = o$, $\pmb{b}.o = (gw).(w^{-1}\pmb{a}).o$, $f_0(\Omega)=f_0(\A)\cap g.f_0(\A)=f_0(\A)\cap (gw).f_0(\A)$ (since $w.\A = \A$) and $w^{-1}\pmb{a}$ is maximal with respect to the amount of distinct entries, so we can apply Proposition \ref{pro:A2id} to get $f_0|_\Omega = (gw).f_0|_\Omega=g.f_0 \circ w |_\Omega$.
\end{proof}
To show that axiom (A2) holds, there is one more obstacle, namely $g$ may not preserve $o$ in general. In fact, $\Omega$ may not even contain $o$, so we have to first translate $\Omega$ to be able to use the previous propositions. 
\begin{proposition}
	\label{pro:A2aW} 
	Let $g \in G_\F$ and define $\Omega := f_0^{-1}\left( f_0(\A)\cap g.f_0(\A) \right)$. Then there exists an element $\overline{w} \in W_a$ such that $f_0|_\Omega = g.f_0 \circ \overline{w}|_\Omega$.
\end{proposition}
\begin{proof}
	If $\Omega=\emptyset$, then any element of $W_a$ will suffice. Otherwise choose and fix $\pmb{a}\in B$ (not necessarily maximal) and $\pmb{b}\in A_\F$ with $\pmb{b}.o =g.\pmb{a}.o \in f_0(\Omega)$. Now consider any $a,b \in A_\F$ with $b.o = g.a.o \in f_0(\Omega)$. We translate the problem by $\pmb{b}^{-1}$ and get
	$$
	\pmb{b}^{-1}.b.o = (\pmb{b}^{-1}.g.\pmb{a}).\pmb{a}^{-1}.a.o
	$$
	where $(\pmb{b}^{-1}.g.\pmb{a}).o = o$. We are almost in the situation of Proposition \ref{pro:A2W}, but we have a different $\Omega$. In fact since $\pmb{b}^{-1}.f_0(\A)=f_0(\A)=\pmb{a}.f_0(\A)$,
	\begin{align*}
		f_0(\A)\cap (\pmb{b}^{-1}.g.\pmb{a}).f_0(\A) &=\pmb{b}^{-1}.f_0(\A)\cap \pmb{b}^{-1}.g.f_0(\A) \\ &=\pmb{b}^{-1}.f_0(\Omega)=f_0(\pmb{b}^{-1}.\Omega),
	\end{align*} 
	we can apply Proposition \ref{pro:A2W} to get a $w \in W_s$ such that
	$$
	f_0|_{\pmb{b}^{-1}.\Omega}= (\pmb{b}^{-1}.g.\pmb{a}).f_0 \circ w |_{\pmb{b}^{-1}.\Omega} \ \ ,
	$$
	which can be rewritten as
	$$
	\pmb{b}^{-1}.f_0|_\Omega =  \pmb{b}^{-1}.g.f_0\circ (\pmb{a} w \pmb{b}^{-1}) |_\Omega.
	$$
	Renaming $\overline{w}=\pmb{a} w \pmb{b}^{-1} \in W_a$ and translating back results in the required formula
	$$
	f_0|_\Omega =  g.f_0\circ \overline{w} |_\Omega.
	$$
\end{proof}
We can now conclude the proof of axiom (A2).
\begin{proposition}\label{prop:A2_SLn}
	\label{pro:A2} The axiom
	\begin{enumerate}
		\item [(A2)] For all $f_1,f_2 \in \Fun$, if $f_1(\A)\cap f_2(\A)\neq \emptyset$, then the set $\Omega := f_2^{-1}\circ f_1(\A)$ is a $W_a$-convex set and there exists $w \in W_a$ such that $f_2|_\Omega = f_1\circ w |_\Omega$. 
	\end{enumerate}
	holds for $(\B,\Fun)$.
\end{proposition}
\begin{proof}
	Let $g,h \in G_\F$ with $f_1=g.f_0$ and $f_2=h.f_0$ such that $g.f_0(\A)\cap h.f_0(\A)\neq \emptyset$. We define $\Omega := (h.f_0)^{-1}(h.f_0(\A)\cap g.f_0(\A))$. We will first show the second part of (A2). To be able to apply Proposition \ref{pro:A2aW}, we act with $h^{-1}$ on
	$$
	h.f_0(\Omega)=h.f_0(\A)\cap g.f_0(\A),
	$$
	to get
	$$
	f_0(\Omega)=f_0(\A)\cap (h^{-1}g).f_0(\A).
	$$
	We can apply Proposition \ref{pro:A2aW} with $h^{-1}g$ to get an element $w \in W_a$ such that $
	f_0|_\Omega = (h^{-1}g).f_0 \circ w |_\Omega, 
	$
	and thus $f_2|_\Omega = h.f_0|_\Omega = g.f_0\circ w |_\Omega=f_1\circ w |_\Omega$.
	
	It remains to show that $\Omega$ is $W_a$-convex. Elements $a.o \in f_0(\Omega)$ are exactly those elements which have the special property that
	$
	h.a.o = g.w.a.o,
	$
	which is equivalent to $a^{-1}h^{-1}gwa \in G_\F(O)$ by Theorem \ref{thm:stab}. If $a = \diag{a} \in A_\F$, then
	$$
	(h^{-1}gw)_{ij}\frac{a_j}{a_i} \in O
	$$
	in matrix entries. Taking the valuation we obtain
	$$
	(-v)(\chi_{\alpha_{ij}}(a)) = (-v)(a_i/a_j) \geq (-v)((h^{-1}gw)_{ij}) \in \Lambda
	$$
	and these are exactly the inequalities that define affine half-apartments
	$$
	H_{\alpha_{ij},k}^+ := \left\{  a.o \in \A \colon (-v)(\chi_{\alpha_{ij}}(a)) \geq k  \right\}
	$$
	for $k = (-v)((h^{-1}gw)_{ij})  $. Since there are only finitely many pairs $(i,j)$, $\Omega$ is a finite intersection of half-apartments, i.e. $W_a$-convex.
\end{proof}

	\end{appendices}

 \newpage
 
    \bibliographystyle{alpha_noand} 
    \bibliography{refs} 
 
	
\end{document}